\documentclass[11pt]{amsart}

\usepackage{amssymb,epsfig,xcolor,upgreek,cancel,enumitem}

\usepackage{pgfplots}
\pgfplotsset{/pgf/number format/use comma, compat=newest}

\newcommand{\R}{{\bf R}}
\newcommand{\N}{{\bf N}}

\newfont{\bbten}{bbold12}

\newcommand{\miu}{\leqslant}
\newcommand{\mau}{\geqslant}

\newcommand{\dint}{\displaystyle \int}

\def\media#1{{\int\!\!\!\!\!\!-}_{\!\!\!{#1}}}

\newtheorem{teorema}{Theorem}[section]

\newtheorem{oss}[teorema]{\sc Remark}
\newtheorem{nota}[teorema]{\sc Nota}
\newtheorem{esempio}[teorema]{\sc Example}
\newtheorem{definition}[teorema]{Definition}
\newtheorem{prop}[teorema]{Proposition}
\newtheorem{lemma}[teorema]{Lemma}
\newtheorem{cor}[teorema]{Corollary}
\newtheorem{theorem}[teorema]{Theorem}

\newcommand{\be}{\begin{equation}}
\newcommand{\ee}{\end{equation}}
\newcommand{\ba}{\begin{array}}
\newcommand{\ea}{\end{array}}
\newcommand{\bea}{\begin{eqnarray}}
\newcommand{\eea}{\end{eqnarray}}
\newcommand{\arst}{\renewcommand{\arraystretch}{1.5}}
\def\arrst#1{\renewcommand{\arraystretch}{#1}}

\newcommand{\bd}{\begin{definition}}
\newcommand{\ed}{\end{definition}}
\newcommand{\bprop}{\begin{prop}}
\newcommand{\eprop}{\end{prop}}
\newcommand{\bthm}{\begin{teorema}}
\newcommand{\ethm}{\end{teorema}}
\newcommand{\boss}{ \begin{oss}  \rm -\ }
\newcommand{\eoss}{ \end{oss} }
\newcommand{\bex}{ \begin{esempio}  \rm -\ }
\newcommand{\eex}{ \end{esempio} }
\newcommand{\bcor}{\begin{cor}}
\newcommand{\ecor}{\end{cor}}
\newcommand{\besempio}{ \begin{esempio}  \rm -\ }
\newcommand{\eesempio}{ \end{esempio} }
\newcommand{\bnota}{ \begin{nota}  \rm -\ }
\newcommand{\enota}{ \end{nota} }

\newcommand{\dimo}{{\it Proof}\ \ -\ \ }

\newcommand{\finedimo}{\hfill $\square$ \\}
\newcommand{\fine}{ \\}

\newcommand{\e}{\epsilon}
\newcommand{\la}{\lambda}

\newcommand{\Om}{\Omega}

\newcommand{\om}{\omega}

\newcommand{\B}{\mathcal{B}}

\newcommand{\V}{\mathcal{V}}
\newcommand{\W}{\mathcal{W}}

\newcommand{\osc}{\mathop{\rm osc}\limits}

\newcommand{\esssup}{\mathop{\rm ess\hspace{0.1 cm}sup}\limits}

\def\media#1{{\dint\!\!\!\!\!\!-}_{\!\!\!{#1}}}

\newcommand{\q}{\mathfrak q}

\newcommand{\emm}{\mathpzc{M}}
\newcommand{\enn}{\mathpzc{N}}
\newcommand{\anna}{\uptau}
\newcommand{\carlo}{N}
\newcommand{\lela}{\upomega}

\newcommand{\loc}{\rm loc}

\makeglossary
\makeindex

\title{A Harnack's inequality for mixed type evolution equations}

\author{\sc Fabio Paronetto}
\thanks{F. Paronetto - Dipartimento di Matematica, Universit\`a di Padova, via 
Trieste 63, 35121, Padova, Italy.  e-mail: fabio.paronetto@unipd.it. Fax +39 049 8271333, phone +39 049 8271318}

\date{\today}

\DeclareMathAlphabet{\mathpzc}{OT1}{pzc}{m}{it}

\begin{document}

\maketitle
\bibliographystyle{fabiostile}

\begin{abstract}
We define a homogeneous parabolic De Giorgi classes of order 2 which suits a mixed type class of evolution equations whose simplest example is
$\mu (x) \frac{\partial u}{\partial t} - \Delta u = 0$ where $\mu$ can be positive, null and negative, so in particular elliptic-parabolic
and forward-backward parabolic equations are included.
For functions belonging to this class we prove local boundedness and show a Harnack inequality which, as by-products, gives
H\"older-continuity, in particular in the interface $I$ where $\mu$ change sign,  and a maximum principle.
\end{abstract}

\ \\


\section{Introduction}

The purpose of this paper is to study problems related to equations of mixed type whose simplest example may be
\begin{equation}
\label{lequazione}
\mu (x) \frac{\partial u}{\partial t} - \Delta u = 0 \qquad \qquad \text{in } \Omega \times (0,T)
\end{equation}
where $\mu$ is a function changing sign and possibly taking also the value zero in some region of positive measure, $\Omega$ an open subset of $\R^n$ and $T > 0$.
This means that the equation can be forward parabolic in a subregion $\Omega_+ \times (0,T)$, 
backward parabolic in another subregion $\Omega_- \times (0,T)$ and also a family of elliptic equations depending on the parameter $t$ in a third
subregion $\Omega_0 \times (0,T)$ of $\Omega \times (0,T)$.
For the existence of solutions to such equations we refer to \cite{fabio4} and the forthcoming paper \cite{fabio9}.
In these papers coefficient $\mu$ is considered depending also on time, but here we confine to $\mu$ depending only on the spatial variable. \\
We want to remember that equations of this kind were already been considered, but always in some simple cases. We recall, among the many papers
about mixed type equations,
\cite{baogris} and \cite{pag-tal}, were the two following equations are considered
\begin{gather*}
x \frac{\partial u}{\partial t} - \frac{\partial^{2m} u}{\partial x^{2m}} = 0 \quad (m \geqslant 1 , \, m \in \N) \, ,	\qquad
\text{sgn}(x) \frac{\partial u}{\partial t} - \frac{\partial^{2} u}{\partial x^{2}} + k u = f \, , 
\end{gather*}
and one of the many cases considered by Beals, \cite{beals}, where the following equation is considered
$$
x \frac{\partial u}{\partial t} - \frac{\partial}{\partial x} \left( (1-x^2) \frac{\partial u}{\partial x} \right) \ = 0 
$$
and, for a general situation like that we consider in \cite{fabio4} and \cite{fabio9},
but confined to $\mu \geqslant 0$, we recall \cite{show1}. \\
For the many applications we refer to \cite{show1} and \cite{fabio9} and the references therein. \\  [0.4em]
To come back to the content of the present paper,
we give a Harnack type inequality (see Theorem \ref{Harnack1} and Theorem \ref{Harnack2}) for a wide class of functions belonging to a proper De Giorgi class.
By this result, on one side we give a generalized Harnack inequality which includes the classical ones for elliptic equations and for parabolic equations,
on the other we study regularity and maximum principles
of solutions of equations like \eqref{lequazione}; in particular we get some local H\"older continuity on the interfaces where $\mu$ change sign
(see the examples at the end of the paper). \\
We recall that a result of regularity in a very general setting was already considered in \cite{fabio7} and \cite{fabio9}, but this regards only
regularity in time. \\  [0.4em]
Just to avoid to confine
to consider equations with $\mu : \Omega \to \{ -1, 0 , 1 \}$ and, on the contrary, to consider also, for instance, $\mu$ continuous,
one is forced to consider weighted spaces.
For this reason we consider a more general De Giorgi class suitably defined to contain quasi-minima
(see Section \ref{De Giorgi classes and Q-minima}) for the equation
\begin{equation}
\label{lequazione2}
\mu \frac{\partial u}{\partial t} - \text{div} (\lambda D u) = 0 \qquad \qquad \text{in } \Omega \times (0,T)
\end{equation}
with $\mu$ and $\lambda$ functions in $L^1_{\loc} (\Omega)$, $\lambda > 0$ while $\mu$ is valued in $\R$. Indeed the De Giorgi class we consider
contains also solutions of more general equations, like
\begin{equation}
\label{equazionegeneralissima}
\mu (x) \frac{\partial u}{\partial t} - \text{div} \, A (x,t,u,Du) = B(x,t,u,Du)
\end{equation}
with ($L \geqslant 1, M > 0$)
\begin{align}
\label{condizionissime}
& \big(A (x,t,u,Du) , Du \big) \geqslant \lambda (x) |Du|^2			\, , 		\nonumber	\\
& | A (x,t,u,Du) | \leqslant L \, \lambda (x) |Du|					\, , 					\\
& | B (x,t,u,Du) | \leqslant M \, \lambda (x) |Du|					\, .		\nonumber
\end{align}
To give our main result we follow \cite{dib1} and \cite{gianazza-vespri}, but we want to stress that
the De Giorgi class we consider is different from the one considered in those papers, also when $\mu \equiv 1$ and that
not only because of the more complicate nature of the equations we consider (the reason lies in Lemma \ref{stimettaDG}). \\ [0.4em]
Since our class contains parabolic quasi-minima we want recall that
quasi-minima or quasi-minimizers (briefly $Q$-minima) were introduced by Giaquinta and Giusti in \cite{giagiu}, where they prove local H\"older continuity
extending the result due to De Giorgi for the solution of elliptic equations, while Harnack inequality for quasi-minima was proved by
DiBenedetto and Trudinger in \cite{dib-tru}.
In the parabolic setting the definition of quasi-minima is due to Wieser (see \cite{wieser}), who proves H\"older continuity for a suitable
parabolic De Giorgi class. \\  [0.5em]
Going back to degenerate elliptic and parabolic equations, where by ``degenerate'' we mean where some weights are involved like in \eqref{lequazione2},
we precise that we consider $\mu$ and $\lambda$ such that
$$
|\mu|_{\lambda} :=
\left\{
\begin{array}{ll}
|\mu|		&	\text{ if } \mu \not = 0,		\\
\lambda	&	\text{ if } \mu = 0 
\end{array}
\right.
\quad \text{and} \quad \lambda \qquad \text{are Muckenhoupt weights} ,
$$
a class of weights we introduce in Section \ref{paragrafo2}. Precisely $|\mu|_{\lambda} \in A_{\infty}$ and $\lambda \in A_2$.
Moreover we assume a condition relating $|\mu|_{\lambda}$ and $\lambda$, assumption (H.2), which is the existence of two constants $q > 2$ and $K > 0$
such that ($x \in \R^n$, $\rho > r > 0$)
\begin{equation}
\label{siamoinritardo}
\left(\frac{|B_r(x)|}{|B_\rho(x)|}\right)^{1/n} 
\left( \frac{|\mu|_{\lambda} (B_r(x))}{|\mu|_{\lambda} (B_\rho(x))} \right)^{1/q} 
\left( \frac{\lambda (B_r(x))}{\lambda (B_\rho(x))} \right)^{-1/2} \leqslant K \, .
\end{equation}
We stress that we are forced to introduce the weight $|\mu|_{\lambda}$ extending $|\mu|$ because
the weight $|\mu|$ could be zero in some region with positive measure and in that case the measure associated to $|\mu|$, even if non-negative,
would not be {\em doubling}. We recall that $\omega \in L^1_{\loc}(\Omega)$, $\omega : \Omega \to [0, +\infty]$, satisfies a doubling condition if there is a positive constant $c$ such that
$$
\omega (B_{2r} (x_0)) \leqslant c \, \omega (B_{r} (x_0))
$$
for every $x_0 \in \Omega$ and $r > 0$ such that $B_{2r} (x_0) \subset \Omega$ (and where $\omega (A)$ denotes $\int_A \omega (x) dx$).
Assumption we need for the weights $|\mu|_{\lambda}$ and $\lambda$ are summarized in (H.1), (H.2), (H.3), (H.4) in Section \ref{De Giorgi classes and Q-minima}.
In particular (H.4) gives also a condition about the geometry of the interface separating the regions
$\Omega_+ = \{ \mu > 0 \}, \Omega_0 = \{ \mu = 0 \}, \Omega_- = \{ \mu < 0 \}$, condition which turns out to be sufficient to get the Harnack inequality. 
We do not know if this is sharp and are not able to give a counterexample to this condition. \\ [0.4em]
Harnack's inequality for parabolic equations was first proved separately by Hadamard and Pini in 1954 just for the heat equations, then Moser, Aronson, Serrin, Trudinger
gave some generalizations of this result. But 
among the many papers studying Harnack's inequality and regularity of partial differential equations, both parabolic and elliptic,
we confine to mention some papers regarding degenerate cases similar to the
one we are considering, referring also to the references contained in them for the more classical results. \\
First we recall \cite{fa-ke-se} where for the first time, at least for our knowledge,
a Muckenhoupt condition on $\lambda$, and precisely $\lambda \in A_2$, was consider to study regularity of the solutions of equations like
$$
\text{div} (\lambda \, D u) = 0
$$
or more generally $\text{div} (a \cdot D u) = 0$ where $a$ satisfies
$$
\lambda (x) |\xi|^2 \leqslant \big( a(x) \cdot \xi, \xi \big) \leqslant L \, \lambda(x) |\xi|^2
$$
In this regard we also recall \cite{tru1} and \cite{tru2}, where some sommability conditions, and not some local conditions, on the weight were requested. \\
As regards the parabolic case, we recall that equations like \eqref{lequazione2} are considered in \cite{chia-se1}, where $\mu \equiv 1$ is considered,
and in \cite{chia-se2}, where
$\mu = \lambda$ is considered. In both these papers a condition $A_2$ on the weight $\lambda$ is considered, when $\mu \equiv 1$ to show that
$L^{\infty}$ bounds and Harnack inequality are impossible, in the second paper where $\mu = \lambda$ to show $L^{\infty}$ bounds and a Harnack inequality.
To get the Harnack inequality with $\mu \equiv 1$ a stronger request has to be made, i.e. $\lambda$ has to belong to $A_{1+2/n}$ which is
a proper subclass of $A_2$ (see \cite{chia-se3}). \\
A more recent paper we mention about linear elliptic equation with principal part in divergence form \cite{mohammed},
where the matrix $a$ defining the principal part satisfies
\begin{equation}
\label{veronastation}
\lambda_1 (x) |\xi|^2 \leqslant \big( a(x,t) \cdot \xi, \xi \big) \leqslant \lambda_2 (x) |\xi|^2 \, ,
\end{equation}
but satisfying \eqref{siamoinritardo} with $\lambda_1$ in the place of $\lambda$ and $\lambda_2$ in the place of $|\mu|_{\lambda}$;
this implies the Sobolev-Poincar\'e inequality
$$
\Big[ \frac{1}{\nu(B_\rho)} {\int_{B_\rho}} |u(x)|^q \lambda_2 (x) dx \Big]^{1/q} \leqslant C \, \rho \, 
	\Big[ \frac{1}{\omega(B_\rho)} {\int_{B_\rho}} |Du(x)|^2 \lambda_1 (x) dx  \Big]^{1/2}
$$
for every Lipschith function with either support contained in $B_{\rho}$ or with null mean value.
About parabolic equations with some $\mu$ in front of $\partial_t$
we also mention \cite{ishige}, where an equation with $\mu = \lambda$ is considered,
\cite{fernandes}, where the author considers $\mu \, \partial_t u - \text{div} \big( a(x,t) Du \big) = 0$ with $a$ satisfying \eqref{veronastation},
and \cite{gut-wheeden2} where $\lambda_1$ and $\lambda_2$ are depending also on time.
Finally we quote the recent paper \cite{surnachev}, where the technique used is the one developed
by DiBenedetto, Gianazza and Vespri in \cite{dib1} and \cite{gianazza-vespri} 
(see also \cite{articolo} for a result regarding non-linear equations)
and the result is analogous to that in \cite{chia-se3}, but 
it concerns monotone operators
with $(p-1)$-growth and the condition about $\lambda$ is $A_{1+p/n}$. \\ [0.4em]
Coming back to our result, we want to stress that our condition (H.2) on the pair $(|\mu|_{\lambda}, \lambda)$
simply reduces to require $\lambda \in A_2$ when $\mu \equiv \lambda$, while is sharp to get,
among the Muckenhoupt weights, $\lambda \in A_{1+2/n}$ when $\mu \equiv 1$
(for this see Remark \ref{notaimportante}, point $\mathpzc{D}$, and Remark \ref{cortona}), so our result cover the result obtained in
\cite{chia-se2} and \cite{chia-se3} and, confining to $\mu > 0$ almost everywhere, generalizes them to doubly weighted parabolic equations. \\ [0.4em]
As regards some result concerning mixed type equations we want to recall a recent result contained in \cite{alk-liske},
where the authors prove H\"older continuity for the limit of a family of functions, the solutions $(u_{\epsilon})_{\epsilon > 0}$ of a class of parabolic equations like
$$
\partial_t (\omega_{\epsilon} u) - \text{div} \, ( \omega_{\epsilon} a (x,t) Du) = 0 \, ,
$$
with $a$ satisfying \eqref{veronastation} with $\lambda_1$ and $\lambda_2$ positive constants and
$\omega_{\epsilon} = 1$ in one region, $\omega_{\epsilon} = \epsilon$ in another.  \\ [0.4em]
Before concluding the introduction we want to stress some difficulties and some interesting thing regarding the main results
(Theorem \ref{Harnack1} and Theorem \ref{Harnack2}). A first comment is the following:
given a ball $B_{\rho}(x_o) \subset \Omega$ and once defined
$B_{\rho}^+(x_o) := B_{\rho}(x_o) \cap \{ \mu > 0 \}$,
$B_{\rho}^-(x_o) := B_{\rho}(x_o) \cap \{ \mu < 0 \}$, 
$B_{\rho}^0(x_o) := B_{\rho}(x_o) \cap \{ \mu = 0 \}$, 
we (in particular) show there is a positive constant $c$  such that for every $u$ in a proper class
$$
u(x_o, t_o)\leqslant c \inf_{B_{\rho} (x_o)} 
	\left\{
	\begin{array}{ll}
	u \left(x, t_o + \vartheta \, \rho^2 \, \frac{|\mu|_{\lambda} (B_{\rho}(x_o))}{\lambda (B_{\rho}(x_o))} \right)
															&	\text{ if } x \in \overline{B_{\rho}^+(x_o)}	\\	[0.5em]
	u \left(x, t_o - \vartheta \, \rho^2 \, \frac{|\mu|_{\lambda} (B_{\rho}(x_o))}{\lambda (B_{\rho}(x_o))} \right)
															&	\text{ if } x \in \overline{B_{\rho}^-(x_o)}	\\	[0.5em]
	u(x, t_o)													&	\text{ if } x \in \overline{B_{\rho}^0(x_o)} .
	\end{array}
	\right.
$$
Notice that the temporal interval where $\mu \not= 0$ is proportional to
$$
\rho^2 \, \frac{|\mu|_{\lambda} (B_{\rho}(x_o))}{\lambda (B_{\rho}(x_o))}
$$
where ($\mu_+$ and $\mu_-$ the positive and negative parts of $\mu$)
$$
|\mu|_{\lambda} (B_{\rho}(x_o)) = \mu_+ (B_{\rho}^+(x_o)) + \mu_- (B_{\rho}^-(x_o)) + \lambda (B_{\rho}^0(x_o)) \, ;
$$
what we want to stress is then that the natural temporal delay, for instance where $\mu > 0$, depends
also on the measure of the regions where $\mu < 0$ and $\mu = 0$. \\
This causes a difficulty in proving our result, in particular Theorem \ref{Harnack1}, because the natural cylinders are alike
$$
B_{\rho}(x) \times \left(t, t +\rho^2 \, \frac{|\mu|_{\lambda} (B_{\rho}(x))}{\lambda (B_{\rho}(x))} \right)
$$
and so (in general) it not true that
$$
B_{r}(x) \times \left(t, t + r^2 \, \frac{|\mu|_{\lambda} (B_{r}(x))}{\lambda (B_{r}(x))} \right)
	\subset B_{R}(x) \times \left(t, t + R^2 \, \frac{|\mu|_{\lambda} (B_{R}(x))}{\lambda (B_{R}(x))} \right) , \quad \text{with } \quad 0 < r < R \, .
$$
Other natural difficulties are due to the equation, which can change its nature around an interface, and so every result used by DiBenedetto, Gianazza and Vespri
is to be suitably modified and adapted. \\ [0.4em]
The paper is organized as follows: in Section \ref{paragrafo2} we introduce the class of Muchenhoupt weights and prove
some results needed in the following; in particular a simple, but fundamental, extension of a classical lemma will be needed (see Lemma \ref{lemmuzzofurbo-quinquies}).
Section \ref{paragrafo3} is devoted to a brief comment about mixed type equations, needed to explain a requirement we make in the De Giorgi class.
In Section \ref{De Giorgi classes and Q-minima} we introduce a degenerate mixed type evolution equation, the Q-minima for that equation, assumptions about weights involved in that equation and the De Giorgi class suited to that equations which, as already mentioned, turns out to be different
from the one introduced in \cite{gianazza-vespri} or in \cite{wieser} also when $\mu \equiv 1$; we also show that Q-minima (and then a large class of solutions)
are contained in the De Giorgi class we define. In the following three sections we prove local boundedness, the fundamental step so-called
{\em expansion of positivity} (see Section \ref{secPositivity}) and a Harnack type inequality stated in Theorem \ref{Harnack1} and Theorem \ref{Harnack2}.
Finally, we give some natural consequences of the inequality we obtain and,
due to the particular nature of the equation, some examples in the hope to help comprehension.   \\ [0.4em]
Finally we want to mention that a short version of the paper, without proofs and with a simpler situation where we consider $\mu$ bounded and $\lambda \equiv 1$,
can be found in \cite{fabio11}. \\ [0.4em]
{\sc Acknowledgments - } The author is pleased to thank R. Serapioni and V. Recupero for some nice and interesting discussions on the subject. \\

\section{Preliminaries on weights}
\label{paragrafo2}

In this section we remind and introduce some definitions and results about $A_p$ weights needed in the following. 
For most of the results we refer to \cite{gc-rdf}. \\
By $B_{\rho}(x_0)$ we will denote the open ball
$\{x \in \R^n \, | \, |x-x_0| < \rho \}$, and sometimes we will simply write $B_{\rho}$ or $B$ if it is
not there is no need to specify further.
With the word {\it weight} we will mean a function $\eta$ such that
$$
\eta \text{ weight if: } \quad \eta > 0 \text{ a.e. in } \R^n \quad \text{ and } \quad \eta \in L^1_{\loc}(\R^n).
$$
Given a weight $\eta$ and a function $u \in L^p(\Omega, \eta)$ with $\Omega$ open set of $\R^n$ we will write
$$
\eta(B) := \int_B \eta \, dx \, , \hskip20pt
{{\int\!\!\!\!\!\!-}_{\!\!{B}}} |u|^p \eta \, dx := \frac{1}{\eta(B)} \int_B |u|^p \eta \, dx \, .
$$

\bd
\label{ap}
Let $p >1$, $K > 0$ be constants, $\omega$ a weight.
We say that $\omega$ belongs to the class $A_p(K)$ if 
\begin{equation}
\label{adueduealfa}
\bigg( {{\int\!\!\!\!\!\!-}_{\!\!{B}}} \omega \, dx \bigg)^{1/p}
	\bigg( {{\int\!\!\!\!\!\!-}_{\!\!{B}}} \omega^{-1/(p-1)} dx \bigg)^{(p-1)/p} \leqslant K
\hskip10pt {\it for\ every\ ball\ }B \subset \R^n
\end{equation}
We say that $\omega$ belongs to the class $A_{\infty} (K, \varsigma)$ if
\begin{equation}
\label{Ainfinito}
\frac{\omega (S)}{\omega (B)} \leqslant K \left( \frac{|S|}{|B|} \right)^{\varsigma}
\end{equation}
for every ball $B$ and every measurable set $S \subset B$. \\
We denote by $A_p = \bigcup_{K \mau 1} A_p(K)$. It turns out $($see, e.g., \cite{gc-rdf}$)$
that $A_{\infty} = \bigcup_{p > 1} A_p$. \\
Given a positive weight $\eta$, a class $A_{p}(K;\eta)$ and all the previous classes may be defined
in a analogous way simply replacing the measure $dx$ with $\eta \, dx$. \\
More generally a pair $(\nu , \omega)$ of weights belong to $A_{p,q}^{\alpha} (B_0, K)$, $\alpha \in [0,n)$, $B_0$ ball $($possibly $\R^n)$ if
\begin{equation}
\label{adueduealfa2}
|B|^{\alpha/n}
\bigg( {{\int\!\!\!\!\!\!-}_{\!\!{B}}} \nu \, dx \bigg)^{1/q}
	\bigg( {{\int\!\!\!\!\!\!-}_{\!\!{B}}} \omega^{-1/(p-1)} dx \bigg)^{(p-1)/p} \leqslant K
\hskip10pt {\it for\ every\ ball\ }B \subset B_0 \, .
\end{equation}
We simply write $A_{p,q}^{\alpha} (K)$ if $B_0 = \R^n$.
For $\alpha = 0$ we get the classical Muckenhoupt class of pairs $($for more details we refer to \cite{gc-rdf}$)$;
for $\alpha = 0$, $q = p$, $\nu = \omega$, $B_0 = \R^n$ we get the class $A_p$.
\ed

\noindent
We remind some properties of $A_p$ weights
(the same properties hold for $A_p(\eta)$ weights), for which we refer to \cite{gc-rdf}.
$A_p$ weights verify the {\it doubling property} which is the following:
for a fixed $t > 1$ there exists a constant $c_d > 1$ which we denote by $c_d(\omega)$, such that
\be
\label{doubling_property}
\int_{tB}\omega \, dx \leqslant c_d (\omega) \int_B \omega \, dx
\ee
for every ball $B$ of $\R^n$,
where by $t B$ we mean the ball concentric to $B$ and whose radius
is $t$ times the lenght of the side of $B$. 
If $\omega \in A_p(K)$ one has that for every $t > 0$ the constant $c_d$
depends (only) on $t, n, p, K$. \\
Moreover $\omega \in A_p(K)$ satisfies the following {\it reverse H\"older's inequality}:
there is $\delta = \delta (n,p,K) > 0$ and a constant $c_{\textit{rh}} =  c_{\textit{rh}}(p,K) \geqslant 1$ such that
\be
\label{maggiore_sommabilita`}
\ba{c}
\bigg( \displaystyle\media{B} \omega^{1+\delta}dx
\bigg)^{1/{(1+\delta)}} \leqslant
		c_{\textit{rh}} \ \bigg( \displaystyle\media{B} \omega \, dx \bigg)\ , \\
\bigg(\displaystyle\media{B} \omega^{ {-{\frac{1}{p-1}}} (1+\delta)}dx
		\bigg)^{1/{(1+\delta)}} \leqslant
	c_{\textit{rh}} \ \bigg( \displaystyle\media{B} \omega ^{ {-{\frac{1}{p-1}}} }dx \bigg),
\ea
\ee
for every ball $B$. A consequence of the definition of $A_p$ weights and of \eqref{maggiore_sommabilita`} are the two
following inequalities. If $\omega \in A_p(K)$ then, called $\varsigma$ the quantity $\delta/(1+\delta)$, one has
\begin{align}
\label{ecomeservono!}
\left(\frac{|S|}{|B|}\right)^p \leqslant K \, \frac{\omega(S)}{\omega(B)} \, , \hskip20pt
\frac{\omega(S)}{\omega(B)} \leqslant c_{\textit{rh}} \left(\frac{|S|}{|B|}\right)^{\varsigma}
\end{align}
for every measurable $S \subset B$, for every $B$ ball of $\R^n$.

\boss
\label{minoreq}
Another interesting property of $A_p$ weights is the following. \\
If $\omega \in A_p(K)$ then there is $p' < p$,
$p' = p'(n,p,K)$, and $K' = K' (n,p,K)$ such that $\omega \in A_{p'}(K')$.
To prove this fact take $\omega \in A_p(K)$, $\delta, c_{\textit{rh}}$ considered in \eqref{maggiore_sommabilita`},
choose $\bar{p}$ in such a way that
$$
\frac{1}{\bar{p}-1} = \frac{1}{p-1} (1 + \delta)
$$
(precisely $\bar{p} = (p+\delta)(1+\delta)^{-1} < p$) and using \eqref{maggiore_sommabilita`} we get
\begin{align*}
\media{B} \omega \, dx \ \bigg( \media{B} \omega^{-\frac{1}{p'-1}} dx \bigg)^{p'-1} & \leqslant
\media{B} \omega \, dx \ \bigg( \media{B} \omega^{-\frac{1}{p-1}(1+\delta)} dx \bigg)^{\frac{p-1}{1+\delta}} \leqslant \\
& \leqslant c_{\textit{rh}}^{p-1} \, \media{B} \omega \, dx \ \bigg( \media{B} \omega^{-\frac{1}{p-1}} dx \bigg)^{p-1} \leqslant 
	c_{\textit{rh}}^{p-1}  K^p \, .
\end{align*}
for every $p' \in [\bar{p}, p]$.
\eoss

\boss
\label{remarcuccia}
Suppose to have $\nu, \omega \in A_{\infty}$, i.e. there are $s_1, s_2, K_1, K_2 > 1$ such that
$\omega \in A_{s_1}(K_1)$ and $\nu \in A_{s_2}(K_2)$. \\
Then the weight $\omega/\nu \in A_{\infty}(\nu)$, i.e. there is $r > 1$ such that
$\omega/\nu \in A_r(c; \nu)$ or
\be
\label{stimazza}
\frac{\dint_B \omega \, dx}{\dint_B \nu \, dx}
\left( \frac{\dint_B \Big(\frac{\omega}{\nu}\Big)^{-1/(r-1)} \nu \, dx}{\dint_B \nu \, dx} \right)^{r-1} \miu c
\hskip20pt
\textrm{for every }B \textrm{ ball in }\R^n.
\ee
Indeed multiplying and dividing by $|B|^r$ we get that the above inequality is equivalent to
$$
\frac{1}{|B|^r} \, \dint_B \omega dx \left(\int_B \omega^{-1/(r-1)} \nu^{r/(r-1)} dx \right)^{r-1} \leqslant
c \, \frac{1}{|B|^r} \left( \int_B \nu \, dx \right)^r \, .
$$
Now by H\"older's inequality, if $a^{-1} + b^{-1} = 1$, $a, b > 1$, we get that
$$
\ba{l}
{\displaystyle \frac{1}{|B|^r}}
\dint_B \omega \, dx \left( \dint_B \omega^{-1/(r-1)} \nu^{r/(r-1)} dx \right)^{r-1} \miu		\\ [0.5em]
\hskip40pt \miu \media{B} \omega \, dx \left(\media{B} \omega^{-a/(r-1)} dx \right)^{(r-1)/a}
	\left(\media{B} \nu^{rb/(r-1)} dx \right)^{(r-1)/b} \, .
\ea
$$
Since $a$ and $r$ are arbitrary we can choose $1 + (r-1)/a = s_1$, 
so that
$\omega \in A_{1+(r-1)/a}(K)$ and consequently
$$
\media{B} \omega \, dx \left(\media{B} \omega^{-a/(r-1)} dx \right)^{(r-1)/a} \miu K_1 \, .
$$
Moreover if $\nu \in A_{\infty}$ 
by the higher summability property of $A_{\infty}$ weights, there is $\delta = \delta (s_2, n, K_2) > 0$ 
such that \eqref{maggiore_sommabilita`} holds. 
Notice that it is possible to choose $a, b, r > 1$ in such a way
$$
\frac{1}{a} + \frac{1}{b} = 1\, , \hskip20pt
\frac{r-1}{a} = s_1 - 1 \, , \hskip20pt
\frac{r b}{r-1} = 1 + \delta \, .
$$
With these choices
there is $c_1 = c_1(s_2, n, K_2)$ such that
$$
\left(\media{B} \nu^{rb/(r-1)} dx \right)^{(r-1)/b} \miu c_1 \left( \media{B} \nu \, dx \right)^r \, .
$$
Then \eqref{stimazza} holds with $c = K_1 c_1$, $c = c(s_2, n, K_1, K_2)$, $r = r(s_1, s_2, n, K_1, K_2)$.
\eoss

We recall now some classical results about weighted inequalities.
The following one in particular can be found in \cite{chanillo-wheeden} and is the weighted version of the standard Sobolev-Poincar\'e inequality. 
Given two weights $\nu, \omega$ in $\R^n$ and $p,q$ with $1 < p < q$ the following condition: \\
\begin{center}
\begin{minipage}{11cm}
{there is a constant $K > 0$ such that}
\end{minipage} 
\begin{equation}
\label{sob-poin-cond}
\left(\frac{|B_r(\bar{x})|}{|B_\rho(\bar{x})|}\right)^{\alpha/n} 
\left( \frac{\nu(B_r(\bar{x}))}{\nu(B_\rho(\bar{x}))} \right)^{1/q} 
\left( \frac{\omega(B_r(\bar{x}))}{\omega(B_\rho(\bar{x}))} \right)^{-1/p} \leqslant K
\end{equation}
\begin{minipage}{11cm}
for every pair of concentric balls $B_r$ and $B_\rho$ with $0 < r < \rho$,
\end{minipage} 
\end{center}
\ \\
with $\alpha = 1$ is essentially necessary and sufficient to have the Sobolev-Poincar\'e inequality. 
Below we confine to state only the result we need.
For more details we refer to \cite{chanillo-wheeden}. 

\bd
\label{bipiqualfa}
For a pair of weights $\nu,\omega$ and $\alpha \in [0,n)$ we will write $($this is not a standard notation$)$
$$
(\nu,\omega) \in B_{p,q}^{\alpha} (K)
$$
if it satisfies \eqref{sob-poin-cond} for every pair of balls $B_r (\bar{x}), B_\rho(\bar{x})$ with $r < \rho$ and $\bar{x} \in \R^n$.
\ed

\bthm
\label{chanillo-wheeden}
Consider $p,q$ such that $1 < p < q$,
$\rho > 0$, $x_0 \in \R^n$, two weights $\nu ,\omega$ in $\R^n$
such that $\omega \in A_p(K_1)$, $(\nu,\omega) \in B_{p,q}^1(K_2)$ and $\nu$ satisfies \eqref{doubling_property}.
Then there is a constant $\gamma_1$ depending $($only$)$ on $n, p, q, K_1, K_2, c_{d}(\nu)$
$($the doubling constants of the weight $\nu)$ such that
\begin{equation}
\label{disuguaglianzadisobolev}
\Big[ \frac{1}{\nu(B_\rho)} {\int_{B_\rho}} |u(x)|^q \nu (x) dx \Big]^{1/q} \leqslant \gamma_1 \, \rho \, 
	\Big[ \frac{1}{\omega(B_\rho)} {\int_{B_\rho}} |Du(x)|^p \omega (x) dx  \Big]^{1/p}
\end{equation}
for every $u$ Lipschitz continuous function
defined in $B_{\rho} = B_{\rho}(x_0)$, with either support contained in $B_{\rho}(x_0)$ or with null mean value.
\ethm

\boss
\label{rmkipotesi}
Notice that the previous theorem holds also for every $q' \in [1, q]$ in the place of $q$ with the same constant $\gamma_1$.
Indeed condition \eqref{sob-poin-cond} holds with the same constant $K$ for every $q' \in [1, q]$. \\
Moreover, using \eqref{ecomeservono!}, one gets that in particular Theorem \ref{chanillo-wheeden} holds when $\nu = \omega \in A_p$ with
$q = np/(n-1) > p$ (and in fact also with some greater value thanks to Remark \ref{minoreq}).
\eoss

\boss
\label{notaimportante}
Here we want to stress some important facts we will need later. \\
$\mathpzc{A}$ {\em - If $(\nu,\omega) \in A_{p,q}^1(K, B_0)$ with $1 < p < q$, $\nu \in A_{\infty}$, then there are $\alpha \in (0,1)$,
$\tilde{q} \in (p,q)$, $\tilde{K} \geqslant K$ such that $(\nu,\omega) \in A_{p,\tilde{q}}^{\alpha}(\tilde{K}, B_0)$}. \\
By \eqref{ecomeservono!}, since $\nu \in A_{\infty}$, we get the existence of $\varsigma > 0$ such that for every $\delta > 0$
$$
\left(\frac{\nu(B_r)}{\nu(B_R)}\right)^{\delta} \leqslant (c_{\textit{rh}}(\nu))^{\delta} \left(\frac{|B_r|}{|B_R|}\right)^{\varsigma \delta}.
$$
Now we choose $\delta$ and consequently define $\tilde{q}$ in such a way that
$$
\frac{1}{q} + \delta < \frac{1}{p} \hskip20pt \textrm{and} \hskip20pt \frac{1}{\tilde{q}} := \frac{1}{q} + \delta \, .
$$
Now we can fix $\alpha \in (0,1)$ and we do that in such a way that $\varsigma \delta = (1-\alpha)/n$. Then we have for $r < R$
\begin{align*}
K \left( \frac{\omega(B_r(\bar{x}))}{\omega(B_R(\bar{x}))} \right)^{1/p} & \geqslant
	\left(\frac{|B_r(\bar{x})|}{|B_R(\bar{x})|}\right)^{1/n} 
	\left( \frac{\nu(B_r(\bar{x}))}{\nu(B_R(\bar{x}))} \right)^{1/q} 	=			\\
& = 	\left(\frac{|B_r(\bar{x})|}{|B_R(\bar{x})|}\right)^{\alpha/n}
	\left(\frac{|B_r(\bar{x})|}{|B_R(\bar{x})|}\right)^{(1-\alpha)/n} 
	\left( \frac{\nu(B_r(\bar{x}))}{\nu(B_R(\bar{x}))} \right)^{1/q} 	\geqslant		\\
& \geqslant \frac{1}{(c_{\textit{rh}}(\nu))^{\delta}}\left(\frac{|B_r(\bar{x})|}{|B_R(\bar{x})|}\right)^{\alpha/n}
	\left( \frac{\nu(B_r(\bar{x}))}{\nu(B_R(\bar{x}))} \right)^{1/q'} 	\, .
\end{align*}
$\mathpzc{B}$ {\em - If $(\nu,\omega) \in A_{p,q}^1(K_2, B_0)$ with $1 < p < q$, $\nu \in A_{\infty}$, $\omega \in A_p(K_1)$,
then there are $p' \in (1,p)$, $q' \in (p,q)$, $K_2' \geqslant K_2$ such that $(\nu,\omega) \in A_{p',q'}^{1}(K_2', B_0)$}. \\
Consider the values of $\alpha, q', K'$ ($K' \geqslant K_2$) of point $\mathpzc{A}$: then we know that $(\nu,\omega) \in A_{p,q'}^{\alpha}(K', B_0)$.
If we consider $p'$ in such a way that
$$
\frac{p-p'}{p'} = \frac{1-\alpha}{n}
$$
we get, using the assumptions, the fact $\omega \in A_p(K_1)$ and \eqref{ecomeservono!}, for $r < R$
\begin{align*}
K' \left( \frac{\omega(B_r(\bar{x}))}{\omega(B_R(\bar{x}))} \right)^{1/p'} & \geqslant
	\left( \frac{\omega(B_r(\bar{x}))}{\omega(B_R(\bar{x}))} \right)^{1/p' - 1/p}
	\left(\frac{|B_r(\bar{x})|}{|B_R(\bar{x})|}\right)^{\alpha/n} 
	\left( \frac{\nu(B_r(\bar{x}))}{\nu(B_R(\bar{x}))} \right)^{1/q'} 	\geqslant			\\
& \geqslant \left(\frac{1}{K_1}\right)^{\frac{p-p'}{p'p}}
	\left(\frac{|B_r(\bar{x})|}{|B_R(\bar{x})|}\right)^{\frac{p-p'}{p'}}
	\left(\frac{|B_r(\bar{x})|}{|B_R(\bar{x})|}\right)^{\alpha/n}
	\left( \frac{\nu(B_r(\bar{x}))}{\nu(B_R(\bar{x}))} \right)^{1/q'} 	=				\\
& = \left(\frac{1}{K_1}\right)^{\frac{p-p'}{p'p}}
	\left(\frac{|B_r(\bar{x})|}{|B_R(\bar{x})|}\right)^{1/n}
	\left( \frac{\nu(B_r(\bar{x}))}{\nu(B_R(\bar{x}))} \right)^{1/q'} 	\, .
\end{align*}
Taking $K_2' = K' K_1^{\frac{p-p'}{p'p}}$ (which depends on $K_1, K_2, c_{\textit{rh}}(\nu), p, q, n$) we conclude. \\
Actually one can require simply $\omega \in A_{\infty}$ and get not only $(\nu,\omega) \in A_{p',q'}^{1}(K_2', B_0)$,
but in fact $(\nu,\omega) \in A_{p',q'}^{\alpha'}(K_2', B_0)$ with $p' \in (1,p)$, $q' \in (p,q)$, $\alpha' \in (\alpha,1)$. \\ [0.3em]
$\mathpzc{C}$ {\em - If $(\nu,\omega) \in A_{2,q}^1(K, B_0)$ with $q > 2$ the function
$f(\bar{x},r) = r^{2 \alpha} \frac{\nu (B_r(\bar{x}))}{\omega (B_r(\bar{x}))}$ satisfies the following inequality: by point $\mathpzc{A}$
and Remark $\ref{rmkipotesi}$
we get that there are $\alpha \in (0,1)$ and $\tilde{K} = \tilde{K} (K, c_{\textit{rh}}(\nu))$ such that
$$
f(\bar{x},r) \leqslant {\tilde{K}}^2 f(\bar{x},R)
$$ 
for every $\bar{x} \in B_0$ and $r,R$ satisfying $0 < r < R$.} \\
Indeed by assumptions we derive
\begin{gather*}
\left(\frac{|B_r(\bar{x})|}{|B_R(\bar{x})|}\right)^{\alpha/n} 
\left( \frac{\nu(B_r(\bar{x}))}{\nu(B_R(\bar{x}))} \right)^{1/2} 
\left( \frac{\omega(B_r(\bar{x}))}{\omega(B_R(\bar{x}))} \right)^{-1/2} \leqslant \tilde{K} \, .
\end{gather*}
Taking the power $2$ we immediately get the thesis.  \\ [0.3em]
$\mathpzc{D}$ {\em -  Consider $\nu \equiv 1$. Then there are $q > p$ and $\hat{K}$ depending on $n, p, K$ such that}
$$
\omega \in
\left\{
\begin{array}{ll}
A_{1+ p/n} (K)		&		\text{ for } n \geqslant {\displaystyle \frac{p}{p-1}} \, ,	\\ [0.5em] 
A_p (K)			&		\text{ for } n \leqslant {\displaystyle \frac{p}{p-1}}  \, .
\end{array}
\right.
\quad \Longrightarrow \quad
(1, \omega) \in B_{q, p}^{1} ( \hat{K}) \, .
$$
First of all notice that for every $n$ we have indeed $\omega \in A_{1+ p/n} (K)$. 
In particular, by Remark \ref{minoreq}, there is $K'$ and $\varepsilon$ such that $\omega \in A_{1+ p/n - \varepsilon} (K')$. Using
the first one in \eqref{ecomeservono!} with $S = B_r(x)$ and $B = B_{\rho}(x)$ ($x \in \Omega$ and $\rho > r > 0$) we get
$$
\left(\frac{| B_r |}{| B_{\rho} |}\right)^{1+ \frac{p}{n} - \varepsilon} \leqslant K' \, \frac{\omega (B_r)}{\omega(B_{\rho})} .
$$
What we want to prove, since $\nu \equiv 1$, is
$$
\left(\frac{| B_r |}{| B_{\rho} |}\right)^{\frac{1}{n} + \frac{1}{q}} \leqslant \hat{K} \, \left( \frac{\omega (B_r)}{\omega(B_{\rho})} \right)^{1/p}.
$$
Taking the power $p$ we get the thesis with $\hat{K} = (K')^p$ and choosing some $q \in \left( p, \frac{p}{1 - \varepsilon}\right)$.
\eoss

\boss
\label{cortona}
In this remark we want to stress that the request $\omega \in A_{1+p/n}$ is optimal among the Muckenhoupt class to get that
$(1, \omega) \in B_{q, p}^{1}$ for some $q > p$. \\
Indeed consider
$\omega (x) = |x|^{\beta}$ which is $A_r$ if and only if $- n < \beta < (r - 1)n$. If we consider $r > 1 + p/n$ then it is possible to choose $\beta > p$ and in this case
to get $(1, \omega) \in B_{q, p}^{1}$ for some $q > p$ we should consider $p < \beta < p + (p/q - 1)n$, but this is clearly impossible.
\eoss

\noindent

\noindent
We state now a slight generalization of a result about Muckenhoupt type weights,  (see \cite{gut-wheeden3}
and \cite{gut-wheeden2}).

\bthm
\label{gut-whee}
Consider $B_{\rho} = B_{\rho}(x_0)$ a ball of $\R^n$ whose radius's lenght is $\rho$,
$\omega \in A_2(K_1)$, $(\nu,\omega) \in B_{2,q}^1(K_2)$ for some $q > 2$, $\nu \in A_{\infty}$.
Then there is $\upsigma_1 \in (1,q)$ $($see also the remark below$)$
such that for every $A \subset B_{\rho}(x_0)$, for every Lipschitz continuous function $u$
defined in $B_{\rho}(x_0)$, with either support contained in $B_{\rho}(x_0)$ or with null mean value
and for every $\kappa \in (1, \upsigma_1]$
$$
\frac{1}{\upsilon (B_{\rho})} \int_A |u|^{2\kappa} \upsilon \, dx \leqslant \gamma_1^{2} \, \rho^2 \, 
	\Big( \frac{1}{\nu (B_{\rho})} \int_A  |u|^{2} \nu \, dx \Big)^{\kappa -1} 
	\Big( \frac{1}{\omega (B_{\rho})} \int_{B_{\rho}}  |D u|^2 \omega \, dx \Big)
$$
where the inequality holds both with $\upsilon = \nu$ and $\upsilon = \omega$
$($and in fact with every weight for which Theorem $\ref{chanillo-wheeden}$ holds$)$.
\ethm
\boss
\label{dipendenza}
The assumption $\nu \in A_{\infty}$ means that there is $s > 1, K_3 \geqslant 1$ such that
$\nu \in A_s(K_3)$. Following the proof of Theorem \ref{gut-whee} (and thanks to Remark \ref{remarcuccia})
one can see that the constant $\kappa$ depends $($only$)$ on $n, q, s, K_1, K_3$.
\eoss
\noindent
\dimo
Consider $\kappa > 1$ (to be chosen) and consider $h_0, r > 1$ in such a way that
$$
(\kappa -1) + \frac{1}{h_0} + \frac{1}{r} = 1 \, . 
$$
Writing $|u|^{2\kappa} \upsilon$ as $|u|^{2(\kappa-1)} \nu^{k-1} u^2 \upsilon ^{1/h_0} \upsilon^{1-1/h_0} \nu^{1-\kappa}$ we get
$$
\int_{A} |u|^{2\kappa} \upsilon \, dx \leqslant
\left( \int_{A} u^2 \nu \, dx \right)^{\kappa-1} \left( \int_{B_{\rho}}|u|^{2h_0} \upsilon \, dx \right)^{\frac{1}{h_0}}
\left(\int_{B_{\rho}} \upsilon^{(1-1/h_0)r} \nu^{(1-\kappa)r} dx \right)^{\frac{1}{r}} \, .
$$
Now we chose $h_0 = q/2$ in such a way Theorem \ref{chanillo-wheeden} holds both with $\upsilon = \nu$ and $\upsilon = \omega$
on the left hand side of the inequality. 
For such a $h_0$ we get (we have not chosen $k$ and $r$ yet)
$$
\left( {{\int\!\!\!\!\!\!-}_{\!\!B_{\rho}}} |u|^{2h_0} \, \upsilon \, dx \right)^{\frac{1}{2h_0}} \leqslant
	\gamma_1 \, \rho \, 
	\left( {{\int\!\!\!\!\!\!-}_{\!\!B_{\rho}}} |D u|^2 \omega \, dx \right)^{1/2} \, .
$$
Now consider $\upsilon = \omega$. The previous inequality becomes
$$
\left( \int_{B_{\rho}} |u|^{2h_0} \, \upsilon \, dx \right)^{\frac{1}{h_0}} \leqslant
	\gamma_1^2 \, \rho^2 \, \frac{1}{(\omega(B_{\rho}))^{1-\frac{1}{h_0}}} 
	\left( \int_{B_{\rho}} |D u|^2 \omega \, dx \right) \, .
$$
Since $(1-h_0^{-1})r = r (\kappa-1) + 1$ we may write
$$
\int_{B_{\rho}} \omega^{(1-1/h_0)r} \nu^{(1-\kappa)r} \, dx =
\int_{B_{\rho}} \left(\frac{\omega}{\nu}\right)^{(\kappa-1)r + 1} \nu \, dx \, .
$$
Since $\omega/{\nu} \in A_{\infty}({\nu})$ (see Remark \ref{remarcuccia}) the function $\omega/{\nu}$ satisfies
a reverse H\"older inequality. Then there are
two positive constants $\delta, c_{\textit{rh}}$ such that, for every ball $B$,
$$
\frac{1}{\nu (B)} \int_B \left(\frac{\omega}{\nu}\right)^{1 + \delta} \nu \, dx 
\leqslant c_{\textit{rh}} \left[ \frac{1}{{\nu} (B)} \int_B \frac{\omega}{\nu} \, \nu \, dx \right]^{1 + \delta} =
c_{\textit{rh}} \left[ \frac{\omega (B)}{{\nu} (B)} \right]^{1 + \delta}
$$
(the constants $c_{\textit{rh}}, \delta$ depend on $n, s, K_1, K_3$ if $\nu \in A_s(K_3)$).
Then we will choose $\kappa , r$ in such a way that $(\kappa -1)r = \delta$ and consequently, by what remarked above, we get
\begin{align*}
\int_{B_{\rho}} \omega^{(1-1/h_0)r} \nu^{(1-\kappa)r} \, dx \leqslant
c_{\textit{rh}} \, {\nu}(B_{\rho}) \left[ \frac{\omega (B_{\rho})}{{\nu} (B_{\rho})} \right]^{1 + (\kappa-1)r}
\miu c_{\textit{rh}} \, \frac{\big(\omega (B_{\rho})\big)^{1 + (\kappa-1)r}}{\big({\nu} (B_{\rho})\big)^{(\kappa-1)r}} \, .
\end{align*}
Then we get the thesis when $\upsilon = \omega$. If $\upsilon = \nu$ the proof is easier since
the quantity $\upsilon^{(1-1/h_0)r} \nu^{(1-\kappa)r}$ reduces to $\nu$.
\finedimo

\noindent
We briefly recall the definition (a possible definition, in our case equivalent to the other possible ones)
of weighted Sobolev spaces for $\nu \in A_{\infty}$ and $\omega \in A_p$.
Given an open and bounded set $\Omega \subset \R^n$ by $L^p(\Omega, \nu)$
we denote the set of measurable functions $u : \Omega \to \R$ such that $\int_{\Omega} |u|^p \nu \, dx$ is finite.
By $W^{1,p}(\Omega, \nu, \omega)$ we denote the space
$\{ u \in L^p(\Omega, \nu) \cap W^{1,1}_{\text{loc}}(\Omega) \, | \, D_i u \in  L^p(\Omega, \omega) \}$
endowed with the obvious norm;
by $W^{1,p}_0(\Omega, \nu, \omega)$ we denote the closure of $C^1_c(\Omega)$ in $W^{1,p}(\Omega, \nu, \omega)$. 
Indeed we will write $H^1(\Omega, \nu, \omega)$ for $W^{1,2}(\Omega, \nu, \omega)$. \\

\noindent
Coming back to the result stated in Theorem \ref{gut-whee}, integrating in time one immediately gets what follows.

\bcor
\label{cor-gut-whee}
With the same assumptions of Theorem $\ref{gut-whee}$, consider moreover $s_1, s_2 \in (0,T)$.
Consider a family of open sets $A(t)$, $t \in (s_1,s_2)$ in such a way $E = \cup_{t \in (s_1,s_2)} A(t)$
is a an open subset of $B_{\rho} \times (s_1, s_2)$.
For every $v \in C^0([s_1,s_2]; L^2(B_{\rho},\nu)) \cap L^2(s_1,s_2; W^{1,2}_0(B_{\rho}, \nu, \omega))$ it holds
\begin{align*}
\frac{1}{\upsilon (B_{\rho})} 
\int \!\!\! & \int_{E} |u|^{2\kappa} (x,t) \upsilon (x) \, dx dt \leqslant
	\gamma_1^{2} \, \rho^2 \, \left( \frac{1}{\nu (B_{\rho})} \right)^{\kappa-1}	\cdot					\\
& \cdot \Big( \sup_{s_1 < t < s_2} \int_{A(t)} |u|^{2}(x,t) \nu (x) \, dx \Big)^{\kappa-1}
	\frac{1}{\omega (B_{\rho})} \int_{s_1}^{s_2}\!\! \int_{B_{\rho}} |D u|^2 (x,t) \, \omega (x) \, dx dt
\end{align*}
where the inequality holds both with $\upsilon = \nu$ and $\upsilon = \omega$.
\ecor

\begin{lemma}
\label{lemma2.2}
Consider $B_{\rho} = B_{\rho}(x_0)$ a ball, $p \in (1,+\infty)$, $q \in [1,+\infty)$, $\nu$, $\omega$
and $v \in W^{1,p}(B_{\rho}, \nu, \omega)$ for which assumptions of Theorem $\ref{chanillo-wheeden}$ hold,
$k, l \in \R$ with $k < l$. Consider also a subset $Z$ of $B_{\rho}$ and denote by $\bar\nu$ the function
taking value $0$ in $Z$ and $\bar\nu \equiv \nu$ in $B_{\rho} \setminus Z$. Then
$$
(l-k)^q \, \bar\nu( \{ v < k \} ) \, \bar\nu( \{ v > l \}) \leqslant
	\, 2^q \, \gamma_1^q \, \rho^q \, \bar\nu(B_{\rho}) \, \nu(B_{\rho}) \, \omega(B_{\rho})^{-\frac{q}{p}} \,  
	\left( \int_{B_{\rho} \cap \{ k < v < l \}} |Dv|^p \, \omega \, dx \right)^{q/p} \, .
$$
\end{lemma}
\boss
The previous result holds in every open set $\Omega$, provided that
Theorem \ref{chanillo-wheeden} holds with $\Omega$ in the place of $B_{\rho}$.
\eoss
\noindent
\dimo
Denoted by $A$ the set $\{x \in B_{\rho} \setminus Z \, | \, v(x) < k \}$ and
suppose $\bar{\nu}(A) > 0$, otherwise there is nothing to prove.
Following the proof of Theorem 3.16 in \cite{giusti} we have that for every $u$ which takes the value zero in $A$
\begin{equation}
\label{misuradiA}
\int_{B_{\rho}} |u - u_{B_{\rho}}|^q \bar\nu \, dx = 
\int_{B_{\rho}\setminus A} |u - u_{B_{\rho}}|^q \bar\nu \, dx + \int_{A} |u_{B_{\rho}}|^q \bar\nu \, dx
\mau |u_{B_{\rho}}|^q \int_{A} \bar\nu \, dx \, ,
\end{equation}
where $u_{B_{\rho}} = |B_{\rho}|^{-1} \int_{B_{\rho}} u(x) \, dx$.
Consider the function
$$
u:= 	\left\{
	\ba{ll}
		\min \{v , l \} - k	&	\textrm{if } v > k	\\
		0			&	\textrm{if } v \miu k \, .
	\ea
	\right.
$$
and estimate, first from below
$$
\int_{B_{\rho}} |u|^q \bar\nu \, dx = 
	\int_{\{ v > l \}} (l-k)^q \bar\nu \, dx + \int_{\{ k < v < l \}} (v-k)^q \bar\nu \, dx
	\mau  (l-k)^q \int_{\{ v > l \}} \bar\nu\, dx \, ,
$$
and then, using \eqref{misuradiA}, from above
\begin{align*}
\left( \int_{B_{\rho}} |u|^q \bar\nu \, dx \right)^{\frac{1}{q}} & \miu
	\left( \int_{B_{\rho}} \big[ |u - u_{B_{\rho}}| + |u_{B_{\rho}}| \big]^q \bar\nu \, dx	\right)^{\frac{1}{q}}\\
& \miu \left( \int_{B_{\rho}} |u - u_{B_{\rho}}|^q \bar\nu \, dx \right)^{\frac{1}{q}} + 
	\left( |u_{B_{\rho}}|^q \int_{B_{\rho}} \bar\nu \, dx \right)^{\frac{1}{q}}	\\
& \miu 2 \, \left( \frac{\bar\nu(B_{\rho})}{\bar\nu(A)} \int_{B_{\rho}} |u - u_{B_{\rho}}|^q \bar\nu \, dx \right)^{\frac{1}{q}} \, .
\end{align*}
Now, if $q > p$ we can apply Theorem \ref{chanillo-wheeden}; if $q \leqslant p$, notice that
$(\nu ((B_{\rho}))^{-1} \int_{B_{\rho}} |u - u_{B_{\rho}}|^q \bar\nu \, dx )^{1/q} \leqslant
(\nu ((B_{\rho}))^{-1} \int_{B_{\rho}} |u - u_{B_{\rho}}|^{q'} \bar\nu \, dx )^{1/q'}$ for $q' > q$.
Then, by Theorem \ref{chanillo-wheeden} used if necessary with $q' > p$, we finally get
$$
\displaylines{
\hfill
(l-k)^q \int_{\{ v > l \}} \bar\nu \, dx \leqslant
	2^q \, \gamma_1^q \, \rho^q \, \frac{\bar\nu(B_{\rho})}{\bar\nu(A)} \, \frac{\nu(B_{\rho})}{\omega(B_{\rho})^{q/p}}
	\, \left( \int_{\{k < v < l\}} |Dv|^p \omega \, dx \right)^{q/p} \, .
\hfill\llap{$\square$}}
$$
\fine

\begin{lemma}
\label{lemmaMisVar}
Consider $x_0 \in \Omega$ and $\rho > 0$ such that $B_{2\rho}(x_0) \subset \Omega$, $\sigma \in (0,\rho)$, $\omega \in A_2(K_1)$,
$(\nu,\omega) \in B_{2,q}^1(K_2)$, $\nu \in A_{\infty}$, $q > 2$, $\alpha, \beta > 0$. 
Consider $\B$ an open and non-empty subset of $B_{\rho}(x_0)$ such that $\B^{\sigma} = \{ x \in \Omega \, | \, \text{\rm dist}(x, \B) < \sigma \}$
is a subset of $B_{\rho}(x_0)$.
Then, for every $\varepsilon , \delta\in (0,1)$ there exists $\eta \in (0,1)$ such that for
every $u\in W^{1,2}_{\rm loc}(\Omega, \nu, \omega)$ satisfying
$$
\int_{\B^{\sigma}} |Du|^2 \, \omega \, dx \leqslant \beta \, \frac{\omega(B_{\rho}(x_0))}{\rho^2},
$$
and
$$
\nu(\{ u > 1\} \cap \B )\geqslant \alpha \, \nu (B_\rho(x_0)), 
$$
there exists $x^*\in \B$ with $B_{\eta\rho}(x^*) \subset \B$ such that
$$
\nu(\{ u > \varepsilon \}\cap B_{\eta\rho}(x^*)) > (1-\delta) \, \nu (B_{\eta\rho}(x^*)).
$$
\end{lemma}
\noindent
\dimo
For any positive $\eta$ satisfying $\eta \rho < \sigma/2$, we can consider a finite disjoint 
family of balls $(B_{\eta\rho}(x_i))_{i\in I}$ with the property that 
$$
\bigcup _{i\in I} B_{\eta\rho}(x_i) \subset \B \subset \bigcup_{i\in I}  B_{2\eta\rho}(x_i) \subset \B^{\sigma} \, .
$$
Again for simplicity, we denote by $B_i$ the ball $B_{\eta\rho}(x_i)$ and by $B_{ii}$ the ball $B_{2\eta\rho}(x_i)$.
We denote by $I^+$ and $I^-$ the sets
$$
I^+=\{i\in I: \nu(\{u>1\} \cap B_{ii}) > \frac{\alpha}{2 \, c_d(\nu)} \, \nu(B_{ii}) \},
$$
$$
I^-=\{i\in I: \nu(\{u>1\} \cap B_{ii}) \leqslant \frac{\alpha}{2 \, c_d(\nu)} \, \nu (B_{ii}) \}
$$
where $c_d(\nu)$ is the doubling constant of the weight $\nu$, which, from now on, we will simply denote by $c_d$.
By assumption we then get 
\begin{align*}
\alpha \, \nu(B_{\rho}(x_0)) & \leqslant \nu (\{u>1\}\cap \B)	
\leqslant \sum_{i \in I^+} \nu (\{u >1\}\cap B_{ii}) 
	+\frac{\alpha}{2\, c_d} \sum_{i\in I^-}\nu (B_{ii})	\leqslant						\\
&\leqslant \sum_{i \in I^+}\nu (\{u>1\}\cap B_{ii}) 
	+\frac{\alpha}{2} \sum_{i\in I^-}\nu (B_{i}) 		\leqslant						\\
&\leqslant \sum_{i \in I^+} \nu(\{u>1\}\cap B_{ii}) + \frac{\alpha}{2} \, \nu (\B) \leqslant	\\
&\leqslant \sum_{i \in I^+} \nu(\{u>1\}\cap B_{ii}) + \frac{\alpha}{2} \, \nu (B_{\rho}(x_0)) \, . 
\end{align*}
By this we get that
\begin{equation}
\label{esti1}
\frac{\alpha}{2} \, \nu (B_{\rho}(x_0)) \leqslant  \sum_{i\in I^+}\nu (\{u>1\}\cap B_{ii}) \, .
\end{equation}
Now fix $\varepsilon , \delta\in (0,1)$ and assume by contradiction that
\begin{equation}
\label{estDelta}
\nu(\{u > \varepsilon \} \cap B_i) \leqslant (1-\delta) \, \nu (B_i),\qquad
	\textrm{ for every } i \in I := I^+ \cup I^- \, .
\end{equation}
This clearly would imply in particular that
$$
\frac{\nu(\{u\leqslant \varepsilon\}\cap B_{ii})}{\nu (B_{ii})}\geqslant \frac{\delta}{c_d} =: \delta'
	\qquad \textrm{ for every } i \in I^+ \, .
$$
By this last inequality,
Lemma \ref{lemma2.2} with $p = q= 2$, $k = \varepsilon$ and $l = 1$,
$\bar{\nu} = \nu$ we would obtain that
\begin{align}
\label{estOmega}
\nonumber
\delta' \, \nu(\{u>1\}\cap B_{ii}) \leqslant & \, 
	\frac{\nu(\{u\leqslant \varepsilon\}\cap B_{ii})}{\nu(B_{ii})} \, \nu(\{u>1\}\cap B_{ii}) \leqslant			\\
\leqslant & \,  \frac{4 \gamma_1^2}{(1-\varepsilon)^2} \, (\eta \rho)^2 \, \frac{\nu(B_{ii})}{\omega(B_{ii})} \, 
	\int_{\{\varepsilon < u < 1\}\cap B_{ii}} |Du|^2 \, \omega \, dx  \, .
\end{align}
By Remark \ref{notaimportante}, point $\mathpzc{A}$, we get the existence of $a \in (0,1)$, $K_2'$, such that
(see also Remark \ref{rmkipotesi})
$$
\left(\frac{|B_{2\eta\rho}(x_i)|}{|B_{2\rho}(x_i)|}\right)^{a/n} 
\left( \frac{\nu (B_{2\eta\rho}(x_i))}{\nu (B_{2\rho}(x_i))} \right)^{1/2} 
\left( \frac{\omega (B_{2\eta\rho}(x_i))}{\omega (B_{2\rho}(x_i))} \right)^{-1/2}  \leqslant K_2' \, ,
$$
i.e.
$$
\eta^{2a} \, \frac{\nu (B_{2\eta\rho}(x_i))}{\omega (B_{2 \eta \rho}(x_i))} \leqslant 
	(K_2')^2 \, \frac{\nu (B_{2\rho}(x_i))}{\omega (B_{2\rho}(x_i))} \, .
$$
Notice that
\begin{gather*}
\nu (B_{2\rho}(x_i)) \leqslant c_d (\nu) \, \nu (B_{\rho}(x_i)) \leqslant c_d (\nu) \, \nu (B_{2\rho}(x_0)) \leqslant (c_d (\nu))^2 \, \nu (B_{\rho}(x_0))	\\
\omega (B_{2\rho}(x_0)) \leqslant \omega (B_{4\rho}(x_i)) \leqslant (c_d (\omega))^2 \omega (B_{\rho}(x_i)) \leqslant (c_d (\omega))^2 \omega (B_{2\rho}(x_i)) 
\end{gather*}
by which we get
$$
\frac{\nu (B_{2\rho}(x_i))}{\omega (B_{2\rho}(x_i))} \leqslant \frac{(c_d (\nu))^2}{(c_d (\omega))^2} \, \frac{\nu (B_{\rho}(x_0))}{\omega (B_{\rho}(x_0))} \, .
$$
Summing up on $I^+$, from \eqref{esti1} and \eqref{estOmega} we get
\begin{align*}
\frac{\alpha}{2} \, \delta' \, \nu (B_{\rho}(x_0)) \leqslant & \, \sum_{i\in I^+}
\frac{4 \gamma_1^2}{(1-\varepsilon)^2} \, (\eta \rho)^2 \, \frac{\nu(B_{ii})}{\omega(B_{ii})} \, 
	\int_{\{\varepsilon < u < 1\}\cap B_{ii}} |Du|^2 \, \omega \, dx		\leqslant														\\
\leqslant & \, 
\frac{4 \gamma_1^2}{(1-\varepsilon)^2} \, \eta^{2(1-a)} \rho^2 (K_2')^2 \frac{(c_d (\nu))^2}{(c_d (\omega))^2} \, 
\frac{\nu (B_{\rho}(x_0))}{\omega (B_{\rho}(x_0))} \, 
	\sum_{i\in I^+} \int_{\{\varepsilon < u < 1\}\cap B_{ii}} |Du|^2 \, \omega \, dx		\leqslant											\\
\leqslant & \, 
\frac{4 \gamma_1^2}{(1-\varepsilon)^2} \, \eta^{2(1-a)}  (K_2')^2 \frac{(c_d (\nu))^2}{(c_d (\omega))^2} \, \beta \, \nu (B_{\rho}(x_0))		\, .
\end{align*}
The conclusion follows by taking the limit $\eta\to 0$.
\finedimo

\noindent
Here we state three results, which are corollaries rispectively of Theorem \ref{chanillo-wheeden}, Lemma \ref{lemma2.2}, Lemma \ref{lemmaMisVar}.

\bcor
\label{corollario1}
In the same assumptions of Theorem $\ref{chanillo-wheeden}$ suppose moreover $a, b \in \R$, $a < b$. Then
\begin{equation}
\label{disuguaglianzadisobolevcontempo}
\Big[ \frac{1}{\nu(B_\rho)} {\int_a^b \!\!\! \int_{B_\rho}} |u(x,t)|^p \nu (x) dx dt \Big]^{1/p} \leqslant \gamma_1 \, \rho \, 
	\Big[ \frac{1}{\omega(B_\rho)} {\int_a^b \!\!\!  \int_{B_\rho}} |Du(x,t)|^p \omega (x) dx dt \Big]^{1/p}
\end{equation}
for every $u$ Lipschitz continuous function in $B_{\rho}(x_0) \times (a,b)$ such that for every $t \in (a,b)$ $u(\cdot, t)$ has
either support contained in $B_{\rho}(x_0)$ or null mean value $($with respect to the variable $x)$.
\ecor
\noindent
\dimo
It is sufficient first to observe that
$(\nu ((B_{\rho}))^{-1} \int_{B_{\rho}} |u - u_{B_{\rho}}|^p \nu \, dx )^{1/p} \leqslant
(\nu ((B_{\rho}))^{-1} \int_{B_{\rho}} |u - u_{B_{\rho}}|^{q'} \nu \, dx )^{1/q'}$ for $q > p$,
then to take the power $p$ and integrate in time.
\finedimo

\bcor
\label{corollario2}
Consider $B_{\rho} = B_{\rho}(x_0)$ a ball, $a, b \in \R$, $a < b$, $p \in (1,+\infty)$, $\nu$, $\omega$
and $v \in L^p (a, b; W^{1,p}(B_{\rho}, \nu, \omega))$ for which assumptions of Theorem $\ref{chanillo-wheeden}$ hold,
$k, l \in \R$ with $k < l$. Consider also a subset $Z$ of $B_{\rho}$ and denote by $\bar\nu$ the function
taking value $0$ in $Z$ and $\bar\nu \equiv \nu$ in $B_{\rho} \setminus Z$. Then
\begin{align*}
(l - k)^p \, & \bar\nu \otimes \mathcal{L}^1 \, ( \{ v < k \} ) \, \bar\nu \otimes \mathcal{L}^1 \, ( \{ v > l \}) \leqslant				\\
&	\leqslant \, 2^p \, \gamma_1^p \, \rho^p \, \bar\nu \otimes \mathcal{L}^1 \, \big(B_{\rho} \times (a,b) \big) \, 
	\frac{\nu (B_{\rho})}{\omega (B_{\rho})} \,  
	\iint_{( B_{\rho} \times (a,b)) \cap \{ k < v < l \}} |Dv|^p \, \omega \, dx dt \, .
\end{align*}
\ecor
\noindent
\dimo
One can follow the proof of Lemma \ref{lemma2.2} integrating in space and time and finally applying Corollary \ref{corollario1}.
\finedimo

\bcor
\label{corollario3}
Consider $x_0 \in \Omega$ and $\rho > 0$ such that $B_{2\rho}(x_0) \subset \Omega$,  $a, b \in \R$, $a < b$,
$\sigma \in (0,\rho)$, $\omega \in A_2(K_1)$,
$(\nu,\omega) \in B_{2,q}^1(K_2)$, $\nu \in A_{\infty}$, $q > 2$, $\alpha, \beta > 0$. 
Consider $\B$ an open and non-empty subset of $B_{\rho}(x_0)$ such that also $\B^{\sigma} = \{ x \in \Omega \, | \, \text{\rm dist}(x, \B) < \sigma \}$
is a subset of $B_{\rho}(x_0)$, $a, b \in \R$, $a < b$.
Then, for every $\varepsilon , \delta\in (0,1)$ there exists $\eta \in (0,1)$ such that for
every $u\in L^2(a,b; W^{1,2}_{\rm loc}(\Omega, \nu, \omega))$ satisfying
$$
\int_a^b \!\!\! \int_{\B^{\sigma}} |Du|^2 \, \omega \, dx dt \leqslant \beta \, (b - a) \, \frac{\omega(B_{\rho}(x_0))}{\rho^2}
$$
and
$$
\nu \otimes {\mathcal L}^1 \big(\{ u > 1\} \cap (\B \times (a,b) ) \big) \geqslant \alpha \, (b - a) \, \nu (B_\rho(x_0)), 
$$
there exists $x^*\in \B$ with $B_{\eta\rho}(x^*) \subset \B$ such that
$$
\nu \otimes {\mathcal L}^1 \big( \{ u > \varepsilon \}\cap (B_{\eta\rho}(x^*) \times (a,b) ) \big) > (1-\delta) \, (b - a) \, \nu (B_{\eta\rho}(x^*)).
$$
\ecor
\noindent
\dimo
One can repeat the proof of Lemma \ref{lemmaMisVar} using a family of disjoint cylinders
$(B_{\eta\rho}(x_i) \times (a,b))_{i\in I}$ with the property that 
$$
\bigcup _{i\in I} B_{\eta\rho}(x_i) \subset \B \subset \bigcup_{i\in I}  B_{2\eta\rho}(x_i) \subset \B^{\sigma} \, ,
$$
taking the measure $\nu \otimes {\mathcal L}^1$ instead of $\nu$ and finally using Corollary \ref{corollario2} to conclude.
\finedimo

\noindent
We conclude stating a standard lemma (see, for instance, Lemma 7.1 in \cite{giusti}) and one of its possible generalizations which will be needed later.

\begin{lemma}
\label{giusti}
Let $(y_h)_h$ be a sequence of positive real numbers such that
$$
y_{h+1} \leqslant c \, b^h \, y_h^{1+\alpha}
$$
with $c, \alpha > 0$, $b > 1$.
If $y_0 \miu c^{-1/\alpha} b^{-1/\alpha^2}$ then
$$
\lim_{h\to +\infty} y_h = 0 \, .
$$
\end{lemma}

\begin{lemma}
\label{lemmuzzofurbo-quinquies}
Let $(y_h)_h$ and $(\epsilon_h)_h$ two sequences of non-negative real numbers such that
\begin{equation}
\label{ipotesi}
y_{h+1} \leqslant c \, b^h \, (y_h + \epsilon_h) \, y_h^{\alpha} \, , \hskip20pt y_{h+1} \leqslant y_h \, , \hskip20pt  \lim_{h \to + \infty} \epsilon_h = 0 \, ,
\end{equation}
$c, \alpha > 0$, $b > 1$.
If $y_0 < c^{-1/\alpha} b^{-1/\alpha^2}$ then
$$
\lim_{h\to +\infty} y_h = 0 \, .
$$
\end{lemma}
\noindent
\dimo
If $\epsilon_h = 0$ for every $h$ we reduce to Lemma \ref{giusti}.
Otherwise, say $\bar{y}$ the limit $\lim_h y_h$ which exists by the monotonicity of $(y_h)_h$ and suppose that
$$
y_0 < c^{-1/\alpha} b^{-1/\alpha^2} \, .
$$
Now, by contradiction, assume that
$$
\bar{y} > 0  \, .
$$
By assumptions we have that for each $\varepsilon > 0$ there is $\bar{h} = \bar{h}(\varepsilon)$ such that
\begin{equation}
\label{assurdo}
\epsilon_h \leqslant \varepsilon \hskip30pt \text{for every } h \geqslant \bar{h}  \, .
\end{equation}
Now for each $\delta > 0$ we choose $\varepsilon$ such that
$\varepsilon < \delta \, \bar{y}$ so that we get $\delta \, y_h \geqslant \varepsilon$ for every $h$.
In particular for $h \geqslant \bar{h}$ we get
\begin{align*}
y_{h+1}	& \leqslant \, c \, b^h \, (y_h + \epsilon_h) \, y_h^{\alpha} \leqslant						\\
		& \leqslant \, c \, b^h \, (y_h + \varepsilon) \, y_h^{\alpha} \leqslant				\\
		& \leqslant \, c \, b^h \, (y_h + \delta \, y_h) \, y_h^{\alpha}		=				\\
		& = (1+\delta) \, c \, b^h \, y_h^{1 + \alpha} \, .
\end{align*}
Using the lemma above we have that if $y_{\bar{h}} \leqslant (1+\delta)^{-1/\alpha} \, c^{-1/\alpha} b^{-1/\alpha^2}$
than $\lim_{h} y_h = \bar{y} = 0$, where $\bar{h}$ depends on $\varepsilon$ which depends on the choice of $\delta$.
By the monotonicity of $(y_h)_h$
if $y_0 \leqslant (1+\delta)^{-1/\alpha} \, c^{-1/\alpha} b^{-1/\alpha^2}$ the condition on $y_{\bar{h}}$ is garanteed
whatever the value of $\bar{h}$. Since $y_0 < c^{-1/\alpha} b^{-1/\alpha^2}$ there is $\delta > 0$
such that $y_0 \leqslant (1+\delta)^{-1/\alpha} \, c^{-1/\alpha} b^{-1/\alpha^2}$
and so we would derive that $\bar{y} = 0$, which contradicts the assumption $\bar{y} > 0$.
\finedimo

\section{Preliminaries about mixed type equations}
\label{paragrafo3}

This brief section is devoted to a remark 
about equations of mixed type, like for example
\begin{equation}
\label{equazionegenerale}
\mu (x) \frac{\partial u}{\partial t} - \textrm{div} (a (x,t,Du)) = 0 ,
\end{equation}
where 
$a$ is a Caratheodory function such that
\begin{align}
& a(x,t,0) = 0 \, ,																		\nonumber	\\
\label{proprieta`}
& (a(x,t,\xi) - a(x,t,\eta), \xi - \eta) \geqslant \lambda (x) |\xi - \eta|^2	\, ,										\\
& |a(x,t, \xi) - a(x,t,\eta) | \leqslant L \, \lambda (x) |\xi - \eta|	\, ,									\nonumber
\end{align}
for every $\xi, \eta \in \R^n$, where $L$ is a positive constant and $\mu = \mu(x), \lambda = \lambda(x)$ are functions,
$\lambda$ positive, while $\mu$ may change sign (and also be zero in some positive measure regions). \\
\noindent
Before talking about mixed type equations we want to recall that a weighted Sobolev space $H^1 (\Omega, |\mu|, \lambda)$ endowed with the norm
$$
\| u \|^2 := \int_{\Omega} u^2 |\mu| dx + \int_{\Omega} |D u|^2 \lambda dx
$$
can be defined even if the function $|\mu|$ takes the value zero in a subset whose measure is positive
(we refer to \cite{fabio3} for the definition and the completeness of this space).
If we denote the space $L^2(0,T; H^1_0(\Omega, |\mu|, \lambda)) $ by $\V$ and the space
$\left\{ u \in \V \, | \, \mu u' \in \V' \right\}$ by $\W$
($u'$ denotes the derivative of $u$, $\V'$ the dual space of $\V$) one has that a solution
of \eqref{equazionegenerale} belong to $\W$ and (see  \cite{fabio4})
$$
u \in \W \hskip10pt \Longrightarrow \hskip10pt t \mapsto  \int_{\Omega} u^2(x,t) \mu(x) dx
\hskip10pt \text{is continuous in } [0,T] \, ,
$$
and
\begin{equation}
\label{finitezzamu}
\int_{\Omega} u^2(x,t) |\mu|(x) dx \hskip10pt \text{is finite for every } t \in [0,T] \, .
\end{equation}
On the other hand, for $u$ solution of \eqref{equazionegenerale}, the function 
\begin{equation}
\label{finitezzalambda}
t \mapsto \int_{\Omega} u^2(x,t) \lambda (x) dx \hskip20pt \text{is not necessarily  } L^{\infty}_{\loc}(0,T)
\end{equation}
even if it is finite for almost every $t$
since ${\displaystyle \int_0^T \!\! \int_{\Omega} u^2(x,t) \lambda (x) dx}$ is finite.
In the next section we will define a De Giorgi type class of functions requiring
\begin{equation}
\label{ipotesipoconaturale}
t \mapsto \int_{A} u^2(x,t) |\mu|_{\lambda} (x) dx \hskip10pt \text{belongs to } L^{\infty}_{\loc}(0,T) \hskip10pt \text{for every } A \subset\subset \Omega \, .
\end{equation}
This is something more of the natural requirement \eqref{finitezzamu} and
%
%
%
this a priori is not guaranteed by the equation in a general situation, but in many cases it is true, as we mention below.
This condition will be needed only if there is a region in which the equation reduces to a family of elliptic equations, i.e. if there is an open set in which $\mu = 0$. \\
More in general, using a corollary of Theorem 2.1 in \cite{fabio7} one can prove that, if $u$ is the solution of the problem
\begin{equation}
\label{aaa}
\left\{
\arst
\ba{l}
{\displaystyle \mu \frac{\partial u}{\partial t}}
				 - \textrm{div} (a(x,t)\cdot Du) = 0 	\hskip30pt	\hfill		\text{in } \Omega \times (0,T)					\\
u = \phi													\hfill		\text{in } \partial\Omega \times (0,T)			\\
u (x,0) = \varphi (x)											\hfill		\text{in } \{ x \in \Omega \, | \, \mu(x) > 0 \}		\\
u (x,0) = \psi (x)											\hfill		\text{in } \{ x \in \Omega \, | \, \mu(x) < 0 \}
\ea
\right.
\end{equation}
for some $\phi \in \W$, $\varphi, \psi \in L^2(\Omega)$, if
$$
\phi_t \in \W \hskip15pt \text{and} \hskip15pt a \hskip5pt \text{is regular in time }
$$
(we refer to \cite{fabio7} for the precise requirement about regularity of $a$) we derive that the function
$w = \eta (u-\phi) \in H^1(0,T; H^1(\Omega, |\mu|, \lambda))$, and then in  particular
\begin{equation}
\label{cont}
u \in C^0((0,T); H^1(\Omega, |\mu|, \lambda))
\end{equation}
and as a by-product one gets that $u$ satisfies \eqref{ipotesipoconaturale} since $H^1(\Omega, |\mu|, \lambda) \subset L^2(\Omega, \lambda)$. \\ [0.3em]
Analogous considerations hold for Neumann boundary conditions. \\ [0.3em]
We observe that in general a solution of a family of elliptic equation will be not regular in time
(if, e.g., $a$ is not regular in time) as we will show with an example in the last section. \\

\section{De Giorgi classes and Q-minima}
\label{De Giorgi classes and Q-minima}

From this section on we will focus our attention on a class of functions which contains the solutions of some forward-backward evolution equations,
also possibly a family of elliptic equations, whose simplest example is the following ($\lambda$ is positive, but $\mu$ is valued in $\R$, both may be unbounded)
\begin{equation}
\label{equazione}
\mu \frac{\partial u}{\partial t} - \text{div} (\lambda D u) = 0 \hskip30pt
\textrm{in } \Omega \times (0,T) \, ,
\end{equation}
but one can think to \eqref{equazionegenerale} or to \eqref{equazionegeneralissima}.
The connection of the class we are going to define and this equation will be clarified below. We will show that solutions of such
a homogeneous equation, of equation \eqref{equazionegenerale}
and also of a wider class of homogeneous equations are {\it quasi-minimizers} (from now on we will call them more simply,
and according to the original definition, $Q$-minima, see Definition \ref{quasi-min})
for equation \eqref{equazione}, and $Q$-minima are contained in the De Giorgi class we are going to define. \\ 
\ \\
{\bf Assumptions about $\mu$ and $\lambda$ - } Given $\mu$ and $\lambda$ defined in $\R^n$, $\lambda$ positive almost everywhere,
while $\mu$ may be positive, null and negative, we define
$$
\mu_{\lambda} :=
\left\{
\ba{ll}
\mu		&	\text{ if } \mu \not = 0,		\\	[0.5em]
\lambda	&	\text{ if } \mu = 0 .
\ea
\right.
$$
Once considered $\Omega$ on open subset of $\R^n$ and $T > 0$ we require $\mu$ and $\lambda$ to satisfy what follows:
there is $q > 2$ such that
\begin{align*}
\text{(H.1) -  } & \lambda \in A_2(K_1) \, ,													\\
\text{(H.2) -  } & (|\mu|_{\lambda},\lambda) \in B_{2,q}^1(K_2)  \, ,									\\
\text{(H.3) -  } & |\mu|_{\lambda} \in A_{\infty} (K_3, \varsigma) \, .
\end{align*}
This conditions (see Theorem \ref{chanillo-wheeden}) garantees the validity of the Sobolev-Poincar\'e type inequality
$$
\Big[ \frac{1}{|\mu|_{\lambda} (B_\rho)} {\int_{B_\rho}} |u(x)|^q |\mu|_{\lambda} (x) dx \Big]^{1/q} \leqslant \gamma_1 \, \rho \, 
	\Big[ \frac{1}{\lambda (B_\rho)} {\int_{B_\rho}} |Du(x)|^2 \lambda (x) dx  \Big]^{1/2}
$$
and of all the results which follows (in particular Theorem \ref{gut-whee} and Corollary \ref{cor-gut-whee}). \ \\
The condition (H.2) (see Remark \ref{notaimportante}, point $\mathpzc{A}$) garantees the existence of $\alpha \in (0,1)$, $\tilde{K}_2 > K_2$
depending on $K_2$ and $c_{\textit{rh}}(|\mu|_{\lambda})$ and
$\tilde{q} \in (2, q)$ such that, thanks also to Remark \ref{rmkipotesi},
\begin{align*}
\text{(H.2)}' \text{ - }   & (|\mu|_{\lambda},\lambda) \in B_{2,\tilde{q}}^{\alpha}(\tilde{K}_2) \subset B_{2,2}^{\alpha}(\tilde{K}_2)  \, .
\end{align*}
We will suppose that the sets
$$
\Omega_+ := \{ x \in \Omega \, | \, \mu(x) > 0 \} , \hskip10pt 
\Omega_- := \{ x \in \Omega \, | \, \mu(x) < 0 \} \hskip10pt \textrm{and} \hskip10pt
\Omega_0 := \Omega \setminus \big( \Omega_+ \cup \Omega_-  \big)
$$
are the union of a finite number of open and connected subsets of $\Omega$.
This means, for instance, that $\mu$ cannot change sign in a Cantor type set with positive measure. \\
Beyond to $\mu_+$ and $\mu_-$, which will denote respectively the positive and negative part of $\mu$,
we define
\begin{equation}
\label{lambda}
\lambda_+ := \left\{
\ba{ll}
\lambda		&	\text{ in } \Omega_+						\\	[0.5em]
0			&	\text{ in } \Omega \setminus \Omega_+ 
\ea
\right. 
,	\quad 
\lambda_- := \left\{
\ba{ll}
\lambda		&	\text{ in } \Omega_-						\\	[0.5em]
0			&	\text{ in } \Omega \setminus \Omega_- 
\ea
\right.
, \quad
\lambda_0 := \left\{
\ba{ll}
\lambda		&	\text{ in } \Omega_0						\\	[0.5em]
0			&	\text{ in } \Omega \setminus \Omega_0
\ea
\right.
\, .
\end{equation}
In this way notice that $$|\mu|_{\lambda} = |\mu| + \lambda_0 = \mu_+ + \mu_- + \lambda_0 .$$
Notice that hypotheses (H.1) and (H.3) (see \eqref{doubling_property})
implies that $\lambda$ and $|\mu|_{\lambda}$ are doubling, i.e. there is a constant $\q$ such that
\begin{equation}
\label{doublingmula}
\begin{array}{c}
|\mu|_{\lambda} \big( B_{2\rho}(x) \big) \leqslant \q \, |\mu|_{\lambda} \big( B_{\rho}(x) \big) ,			\\	[0.5em]
{\lambda} \big( B_{2\rho}(x) \big) \leqslant \q \, {\lambda} \big( B_{\rho}(x) \big)
\end{array}
\end{equation}
for every $x \in \Omega$ and $\rho > 0$ for which $B_{2\rho}(x) \subset \Omega$. \\
Moreover by \eqref{ecomeservono!}, once denoted by $c_{\textit rh}(\lambda)$ the constant satisfying \eqref{maggiore_sommabilita`} with the weight $\lambda$
and $\varsigma (\lambda)$ the constant appearing in \eqref{ecomeservono!} with $\omega = \lambda$ and
$c_{\textit rh}(|\mu|_{\lambda})$ and $\varsigma (|\mu|_{\lambda})$ the analogous with $\omega = |\mu|_{\lambda}$, we get that
\begin{equation}
\label{carlettomio}
\frac{\lambda(S)}{\lambda(Q)} \leqslant \upkappa \, \left(\frac{|\mu|_{\lambda}(S)}{|\mu|_{\lambda}(Q)}\right)^{\uptau} \, ,
\qquad
\frac{|\mu|_{\lambda}(S)}{|\mu|_{\lambda}(Q)} \leqslant \upkappa \, \left( \frac{\lambda(S)}{\lambda(Q)} \right)^{\uptau}
\end{equation}
where $\uptau = \min \{ {\varsigma (\lambda)/r}, {\varsigma (|\mu|_{\lambda})/2} \}$ and
$\upkappa = \max \{ c_{\textit{rh}} (\lambda) \, K_3^{\varsigma (\lambda)/r} , c_{\textit{rh}} (|\mu|_{\lambda}) \, K_1^{\varsigma (|\mu|_{\lambda})/2} \}$. \\
Once defined $I$, the set of ``interfaces'' as follows:
$$
I_+ = \partial\Omega_+ \cap \Omega \, , \hskip10pt 
I_- = \partial\Omega_- \cap \Omega \, , \hskip10pt 
I_0 = \partial\Omega_0 \cap \Omega \, , \hskip20pt 
I := I_+ \cup I_- \cup I_0 \, ,
$$
we moreover will assume the following additional assumptions where, for simplicity, we assume the first holds with the the same constant $\q$ as before:
\begin{align*}
\text{(H.4) -  }	&
\left|
\begin{array}{ll}
\mu_+ \big( B_{2\rho}(x) \big) \leqslant \q \, \mu_+ \big( B_{\rho}(x) \big)
										&	\quad \qquad \text{for every } x \in {\Omega}_+ \cup I_+	,			\\	[0.5em]
\mu_- \big( B_{2\rho}(y) \big) \leqslant \q \, \mu_- \big( B_{\rho}(y) \big)
										&	\quad \qquad \text{for every } y \in {\Omega}_- \cup I_-	,			\\	[0.5em]
\lambda_0 \big( B_{2\rho}(z) \big) \leqslant \q \, \lambda_0 \big( B_{\rho}(z) \big)
										&	\quad \qquad \text{for every } z \in {\Omega}_0 \cup I_0 ,
\end{array}
\right.																										\\
\ 			&		\ 																						\\
\text{(H.5) -  }	&
I \text{ is a such that } \lim_{\varepsilon \to 0^+} |I^{\varepsilon}| = 0 ,
\end{align*}
where (H.4) holds for every $\rho > 0$ for which $B_{2\rho}(x) \subset \Omega$ and
$I^{\varepsilon}$ is the open $\varepsilon$-neighbourhood of $I$ and is defined in \eqref{ingr.dimagr}. \\ [0.3em]
Some comments about (H.4) and (H.5) are in order.
First notice that since $|\mu|_{\lambda}$ satisfies \eqref{doublingmula}, at least one of the three requirements in (H.4) holds for every $x \in \Omega$. \\
Notice moreover that assumption (H.4) is deeply connected to a geometric requirement about the set $I$ of interfaces,
indeed (H.4) has to hold in particular for points belonging to $I$. Finally, about the set $I$, notice that (H.5) is weaker than
the requirement that $I$ is a $\mathcal{H}^{n-1}$-rectifiable set because $I$ could be also not rectifiable.
For all these comments we refer to the last section, in which some examples are shown. \\
\ \\
{\bf Some notations - } 
By $u_+(y)$ we define the function $\max \{ u(y), 0 \}$ and by $u_-(y)$ $ \max \{ -u(y), 0 \}$.
We will write $u_+^2$ or $u_-^2$ to denote
$$
u_+^2(y) := ( u_+(y))^2 \, , \hskip20pt u_-^2(y) := ( u_-(y))^2 \, .
$$
Given $A \subset \Omega$ we will denote,
for a given $\varepsilon > 0$,
\begin{equation}
\label{ingr.dimagr}
\begin{array}{c}
A^{\varepsilon} := \big\{ x \in \Omega \, \big| \, \textrm{dist} (x,A) < \varepsilon \big\} \, , \hskip20pt
	A_{\varepsilon} := \big\{ x \in \Omega \, \big| \, \textrm{dist} (x,A^c) < \varepsilon \big\} \, ,							\\ [0.5em]
\text{while for } \varepsilon = 0 \hskip10pt A^{\varepsilon} = A_{\varepsilon} := A \, .
\end{array}
\end{equation}
Fix, beyond $x_0$, $t_0 \in (0,T)$.
For a given $\varepsilon > 0$ and a ball $B_{\rho}(x_0)$ we define the sets
\begin{align*}
I_{\rho,\varepsilon}(x_0) := (I \cap B_{\rho}(x_0))^{\varepsilon} \, ,	 & \hskip 20pt
	B_{\rho}^0(x_0) := B_{\rho}(x_0) \cap \Omega_0																\\
B_{\rho}^+(x_0) := B_{\rho}(x_0) \cap \Omega_+ \, , & \hskip 20pt
	B_{\rho}^-(x_0) := B_{\rho}(x_0) \cap \Omega_- \, , 															\\
I_{\rho}^+ (x_0) := I \, \cap \, \overline{B_{\rho}^+}(x_0) \, ,	\hskip 20pt
	I_{\rho}^- (x_0) := I \, \cap & \, \overline{B_{\rho}^-}(x_0) \, , \hskip 20pt
		I_{\rho}^0 (x_0) := I \, \cap \, \overline{B_{\rho}^0}(x_0) \, ,													\\
I_{\rho,\varepsilon}^+(x_0) :=  (I_{\rho}^+(x_0))^{\varepsilon} \cap B_{\rho}^+(x_0)\, ,
	& \hskip 20pt	I^{\rho,\varepsilon}_+(x_0) :=	(I_{\rho}^+(x_0))^{\varepsilon} \setminus I_{\rho,\varepsilon}^+(x_0) 	\, ,		\\
I_{\rho,\varepsilon}^-(x_0) :=  (I_{\rho}^-(x_0))^{\varepsilon} \cap B_{\rho}^-(x_0)\, ,
	& \hskip 20pt	I^{\rho,\varepsilon}_-(x_0) :=	(I_{\rho}^-(x_0))^{\varepsilon} \setminus I_{\rho,\varepsilon}^-(x_0) 	\, ,		\\
I_{\rho,\varepsilon}^0(x_0) :=  (I_{\rho}^0(x_0))^{\varepsilon} \cap B_{\rho}^0(x_0)\, ,
	& \hskip 20pt	I^{\rho,\varepsilon}_0(x_0) :=	(I_{\rho}^0(x_0))^{\varepsilon} \setminus I_{\rho,\varepsilon}^0(x_0) 	\, .
\end{align*}
We define the following functions
\begin{align}
\label{funzioneacca}
& h(x_0, \rho) := \frac{|\mu|_{\lambda} \left(B_{\rho}(x_0)\right)}{\lambda \left(B_{\rho}(x_0)\right)} \, , 
						\hskip20pt  f (x_0, \rho) := h(x_0,\rho) \rho^2 \, .
\end{align}
These functions depend a priori on $x_0$, but just for simplicity we will not specify this dependence writing only
$h(\rho)$ and $f(\rho)$ if not strictly necessary. \\
Notice that the function $h$ satisfies, if $\mu \not= 0$ almost everywhere, the following inequalities
\begin{equation}
\label{stimeacca}
h(x_0, \rho) \leqslant \q \, h(x_0, 2\rho) \, , \hskip20pt h(x_0, 2\rho) \leqslant \q \, h(x_0, \rho) \, .
\end{equation}
Other sets we define are the following: fix $x_0 \in \Omega$ and $t_0 \in (0,T)$, $R > 0$, $\upbeta > 0$
and $s_1, s_2 \in (0,T)$ with
$s_1 < t_0 < s_2$ and satisfying
\begin{equation}
\label{esseunoeessedue}
\left.
\begin{array}{lll}
i \, ) 	& s_2 - t_0 = t_0 - s_1 = \upbeta \, h(x_0, R) R^2 & \qquad \text{ when we consider } B_R^+(x_0) \text{ or } B_R^-(x_0) \, , \\ 	[0.5em]
ii \, ) & s_1, s_2 \quad \text{arbitrary} & \qquad \text{ when we consider } B_R^0(x_0) \, .
\end{array}
\right.
\end{equation}
Inside the cylinder $B_R(x_0) \times (s_1, s_2)$ for 
$$
\theta \in [0,1)
$$
we define
\begin{equation}
\label{sigmateta}
\sigma_\theta := \theta \, \upbeta \, h(x_0,R) \, R^2 \, .
\end{equation}
in such a way that
$\sigma_\theta \in [0, \upbeta \, h(x_0,R) R^2)$; then for $\rho \in (0, R)$ and $\varepsilon > 0$ 
and taking $s_1, s_2$ as in \eqref{esseunoeessedue}, point $i \, )$, we define the sets
\begin{equation}
\label{notazione1}
\begin{array}{l}
Q_R^{\upbeta,\texttt{\,>}}(x_0,t_0) := B_R (x_0) \times (t_0, s_2)		\, ,		\qquad
						Q_R^{\upbeta,\texttt{\,<}}(x_0,t_0) := B_R (x_0) \times (s_1,t_0)		\, , 						\\ [0.5em]
Q_R^{\upbeta,+} (x_0,t_0) := B_R^+ (x_0) \times (t_0, s_2)		\, ,	\qquad
	Q_R^{\upbeta,-} (x_0,t_0) := B_R^- (x_0) \times (s_1,t_0)		\, , 												\\ [0.5em]
Q_{R;\rho, \theta}^{\upbeta,+} (x_0,t_0) := B_{\rho}^+ (x_0) \times (t_0 + \sigma_\theta, s_2)	\, ,							\\ [0.5em]
Q_{R;\rho, \theta}^{\upbeta,-} (x_0,t_0) := B_{\rho}^- (x_0) \times (s_1, t_0 - \sigma_{\theta}) \, ,							\\ [0.5em]
Q_{R;\rho, \theta}^{\upbeta,+,\varepsilon} (x_0,t_0) := 
	\left\{
	\begin{array}{ll}
	B_{\rho + \varepsilon} (x_0) \times (t_0 + \sigma_\theta, s_2) 		
							&		\hspace{-2cm} \text{ if } B_{\rho + \varepsilon}^+ (x_0) = B_{\rho + \varepsilon} (x_0) ,	\\	[0.2em]
	\left( (B_{\rho}^+ (x_0))^{\varepsilon} \times (t_0 + \sigma_\theta, s_2) \right) \cup
	\big( (I_{\rho}^+(x_0))^{\varepsilon} \times (t_0, s_2) \big)
							& \text{ otherwise} ,
	\end{array}
	\right.																								\\	[1em]
Q_{R;\rho, \theta}^{\upbeta,-,\varepsilon} (x_0,t_0) := 
	\left\{
	\begin{array}{ll}
	B_{\rho + \varepsilon} (x_0) \times (s_1, t_0 - \sigma_\theta) 		
							&	\hspace{-2cm} \text{ if } B_{\rho + \varepsilon}^- (x_0) = B_{\rho + \varepsilon} (x_0) ,		\\	[0.2em]
	\left( (B_{\rho}^- (x_0))^{\varepsilon} \times (s_1, t_0 - \sigma_\theta) \right) \cup
	\big( (I_{\rho}^-(x_0))^{\varepsilon} \times (s_1, t_0) \big)
							& \text{ otherwise} 	,
	\end{array}
	\right.																								\\	[1em]
\end{array}
\end{equation}
and with $s_1, s_2$ arbitrary (see \eqref{esseunoeessedue}) we define
\begin{equation}
\label{notazione2}
\begin{array}{l}
Q_{R;\rho; s_1, s_2}^{0} (x_0) := B_{\rho}^0 (x_0) \times (s_1, s_2)	\quad \text{for } \rho \leqslant R \, ,					\\	[0.5em]
Q_{R;\rho; s_1, s_2}^{0,\varepsilon} (x_0) := ( B_{\rho}^0(x_0))^{\varepsilon} \times (s_1, s_2) 
\end{array}
\end{equation}
The first subscript $R$ below $Q$ denotes that
$s_2 - t_0$ and $t_0 - s_1$ are proportional to $R^2$ and that we consider subsets of $B_R \times (0,T)$. \\
\ \\
\noindent
We now introduce the De Giorgi class for equation \eqref{equazionegenerale}. \\
In the following definition we will use the measures $\mu_+$ and $\mu_-$ rescaled by the factor $h(x_0,R)$.
We will make the implicit assumption that the support of these measures (or functions) is the same of $\mu_+$ and $\mu_-$, i.e.
$$
\frac{\mu_+}{h(x_0, {R})} (x) := 
	\left\{
	\begin{array}{ll}
	{\displaystyle \frac{\mu_+}{h(x_0, {R})}	}			&	\text{if } \mu_+ (x) > 0	\, ,		\\ 	[1em]
	0											&	\text{if } \mu_+ (x) = 0	\, ,
	\end{array}
	\right.
\quad \frac{\mu_-}{h(x_0, {R})} := 
	\left\{
	\begin{array}{ll}
	{\displaystyle \frac{\mu_-}{h(x_0, {R})} }				&	\text{if } \mu_- (x) > 0	\, ,		\\	[1em]
	0											&	\text{if } \mu_- (x) = 0	\, .
	\end{array}
	\right.
$$
Moreover in the definition which follows we require that $u \in L^{\infty}_{\loc} ((0,T); L^2_{\rm{loc}} (\Omega, |\mu|_{\lambda}))$
even if only the terms
$$
\int_{B_{\rho}} u^2 (x,t) \mu_+ (x) dx	\quad \text{and} \quad \int_{B_{\rho}} u^2 (x,t) \mu_- (x) dx
$$
are, a priori, bounded (see Section \ref{paragrafo3}). 
The fact that also ${\displaystyle \int_{B_{\rho}^0} u^2 (x,t) \lambda (x) dx}$ is to be finite will be needed,
for instance, to prove point $iii \,)$ of Theorem \ref{Linfinity}.

\bd[De Giorgi classes]
\label{classiDG}
Consider $\Omega$ an open subset of $\R^n$ and $T > 0$ and a point $(x_0, t_0) \in \Omega \times (0,T)$.
Consider $R, r, \tilde{r} > 0$, $r < \tilde{r} \leqslant R$, $\upbeta > 0$,
$\theta, \tilde\theta$ such that $0 \leqslant \tilde\theta < \theta < 1$,
$s_1, s_2, t_0 \in (0,T)$, $s_1 < t_0 < s_2$ satisfying \eqref{esseunoeessedue}.
We say that a function
$$
u \in L^2_{\rm{loc}}(0,T; H^1_{\rm{loc}} (\Omega, |\mu|, \lambda)) \cap L^{\infty}_{\loc} ((0,T); L^2_{\rm{loc}} (\Omega, |\mu|_{\lambda}))
$$
belongs to the De Giorgi class $DG_+(\Omega, T, \mu, \la, \gamma)$, being $\gamma$ a positive constant, if for every
$\varepsilon \in [ 0, R - \tilde{r}]$ and $\theta - \tilde\theta = (\tilde{r} - r)^2/R^2$ and every $k \in \R$
the following inequalities hold
$(\sigma_\theta$ is defined in \eqref{sigmateta}$)$: \\ [0.5em]
$i \, )$ for $s_2 = t_0 + \upbeta \, h(x_0, {R}) R^2$ and $B_R(x_0) \times [t_0, s_2] \subset \Omega \times (0,T)$
\begin{align}
\label{DGgamma+}
\sup_{t \in (t_0 + \sigma_\theta, s_2)} & \int_{B_{r + \varepsilon}^+} (u-k)_+^2 (x,t) \mu_+ (x) dx +	
\sup_{t \in (t_0, t_0 + \sigma_{\tilde\theta})} \int_{I^{r, \varepsilon}_+} (u-k)_+^2 (x,t) \mu_-(x) \, dx								\nonumber	\\
& \hskip130pt + \iint_{Q_{R;r, \theta}^{\upbeta, +,\varepsilon}}  |D(u-k)_+|^2\, \la  \, dx ds \leqslant								\nonumber	\\
\leqslant & \, \gamma \Bigg[
		\sup_{t \in (t_0,t_0 + \sigma_{\tilde\theta})} \int_{I_{r, \tilde{r}-r + \varepsilon}^+} (u-k)_+^2 (x,t) \mu_+ (x) \, dx +						\\
& \hskip30pt + \sup_{t \in (t_0 + \sigma_\theta, s_2)} \int_{I^{r, \tilde{r}-r + \varepsilon}_+} (u-k)_+^2 (x,t) \mu_- (x) \, dx +			\nonumber	\\
& \hskip30pt + \frac{1}{(\tilde{r} - r)^2}
		\iint_{Q_{R;r , {\tilde\theta}}^{\upbeta, +,\tilde{r}-r + \varepsilon}} 
		(u-k)_+^2\, \left( \frac{\mu_+}{\upbeta \,  h(x_0, R)} + \la \right) \, dx dt  \Bigg] 	;									\nonumber
\end{align}
$ii \, )$ for $s_1 = t_0 - \upbeta \, h(x_0, R) R^2$ and $B_R(x_0) \times [s_1, t_0] \subset \Omega \times (0,T)$
\begin{align}
\label{DGgamma-}
\sup_{t \in (s_1, t_0 - \sigma_\theta)} & \int_{B_{r + \varepsilon}^-} (u-k)_+^2 (x,t) \mu_- (x) dx +	
\sup_{t \in (t_0 - \sigma_{\tilde\theta}, t_0)} \int_{I^{r, \varepsilon}_-} (u-k)_+^2 (x,t) \mu_+(x) \, dx								\nonumber	\\
& \hskip130pt + \iint_{Q_{R;r, \theta}^{\upbeta,-,\varepsilon}}  |D(u-k)_+|^2\, \la  \, dx ds \leqslant								\nonumber	\\
\leqslant & \, \gamma \Bigg[
	\sup_{t \in (t_0 - \sigma_{\tilde\theta}, t_0)} \int_{I_{r, \tilde{r}-r + \varepsilon}^-} (u-k)_+^2 (x,t) \mu_- (x) \, dx +							\\
& \hskip30pt + \sup_{t \in (s_1, t_0 - \sigma_\theta)} \int_{I^{r, \tilde{r}-r + \varepsilon}_-} (u-k)_+^2 (x,t) \mu_+ (x) \, dx + 			\nonumber	\\
& \hskip30pt + \frac{1}{(\tilde{r} - r)^2}
		\iint_{Q_{R;r, {\tilde\theta}}^{\upbeta, -,\tilde{r}-r + \varepsilon}} 
		(u-k)_+^2\, \left( \frac{\mu_-}{\upbeta \, h(x_0, R)} + \la \right) \, dx dt  \Bigg] ;										\nonumber
\end{align}
$iii \, )$ for $s_1$ and $s_2$ arbitrary and $B_R(x_0) \times [s_1, s_2] \subset \Omega \times (0,T)$
\begin{align}
\label{DGgamma0}
\iint_{Q_{R;r; s_1, s_2}^{0,\varepsilon} (x_0)} & |D(u-k)_+|^2 \la \, dx dt \leqslant												\nonumber	\\
& \leqslant \gamma \Bigg[ \sup_{t \in (s_1, s_2)} \int_{ I_0^{r, \tilde{r} - r + \varepsilon}} (u-k)_+^2(x,t) \mu_- (x) \, dx +				\nonumber	\\
& \hskip50pt + \sup_{t \in (s_1, s_2)}\int_{ I_0^{r, \tilde{r} - r + \varepsilon}} (u-k)_+^2(x,t) \mu_+ (x) \, dx \, +									\\
& \hskip50pt + \frac{1}{(\tilde{r} - r)^2} 
					\iint_{Q_{R; r; s_1, s_2}^{0,\tilde{r}-r + \varepsilon} } (u-k)_+^2 \, \lambda \, dx dt	\Bigg] \, ;					\nonumber	
\end{align}
$iv \, )$ for every $s_2 > t_0$ such that $B_R(x_0) \times [t_0, s_2] \subset \Omega \times (0,T)$
\begin{align}
\label{DGgamma+_1}
\sup_{t \in (t_0, s_2)} \int_{B_r^+} & (u-k)_+^2 (x,t) \mu_+ (x) dx	 \leqslant \int_{B_{\tilde{r}}^+} (u-k)_+^2 (x,t_0) \mu_+ (x) dx	\, +	\nonumber	\\
+ &  \sup_{t \in (t_0, s_2)} \int_{I^{r, \tilde{r}-r}_+} (u-k)_+^2 (x,t) \mu_-(x) \, dx	 +														\\
& \hskip50pt + \, \gamma \, \frac{1}{(\tilde{r} - r)^2}
		\int_{t_0}^{s_2} \!\!\!\! \int_{B_{\tilde{r}}^+ \cup I^{r, \tilde{r}-r}_+} (u-k)_+^2\, \la \, dx dt  	;						\nonumber
\end{align}
$v \, )$ for every $s_1 < t_0$ such that $B_R(x_0) \times [s_1, t_0] \subset \Omega \times (0,T)$
\begin{align}
\label{DGgamma+_2}
\sup_{t \in (s_1, t_0)} \int_{B_r^-} & (u-k)_+^2 (x,t) \mu_- (x) dx	\leqslant 	\int_{B_{\tilde{r}}^-} (u-k)_+^2 (x,t_0) \mu_- (x) dx	\, +	\nonumber	\\
+ & \sup_{t \in (t_0, s_2)} \int_{I^{r, \tilde{r}-r}_-} (u-k)_+^2 (x,t) \mu_+(x) \, dx	 +														\\
& \hskip50pt + \, \gamma \, \frac{1}{(\tilde{r} - r)^2}
		\int_{s_1}^{t_0} \!\!\!\! \int_{B_{\tilde{r}}^- \cup I^{r, \tilde{r}-r}_-} 	(u-k)_+^2\, \la \, dx dt \, .						\nonumber
\end{align}
We will say that $u$ belongs to $DG_-(\Omega, T, \mu, \la, \gamma)$ if the estimates above 
holds for $(u-k)_-$ in the place of  $(u-k)_+$.
We will say that $u$ belongs to $DG(\Omega, T, \mu, \la, \gamma)$ if
$u \in DG_+(\Omega, T, \mu, \la, \gamma) \cap DG_-(\Omega, T, \mu, \la, \gamma)$.
\ed


\boss
\label{notachesegueladefinizione}
Notice that if $|\mu|(B_R(x_0)) = 0$, that is $B_R (x_0) \subset \Omega_0$, 
\eqref{DGgamma+}, \eqref{DGgamma-} and \eqref{DGgamma0} coincide and reduce to
\begin{align*}
\int_{s_1}^{s_2} \!\!\! \int_{B_{r}}  |D(u-k)_+|^2\, \la  \, dx dt \leqslant	 \, \gamma \, \frac{1}{(\tilde{r} - r)^2}
		\int_{s_1}^{s_2} \!\!\! \int_{B_{r}} (u-k)_+^2\,\lambda \, dx dt 
\end{align*}
by which we can derive
\begin{align}
\label{tempofissato}
\int_{B_{r}(x_0)} |D(u-k)_+|^2 (x,t)\, \la(x)  \, dx	 	\leqslant 
	\gamma \, \frac{1}{(\tilde{r} - r)^2} 	\int_{B_{\tilde{r}}(x_0)} (u-k)_+^2 (x,t)\, \lambda (x) \, dx
\end{align}
for {\em almost} every $t \in [s_1, s_2]$.
Since by assumption $u \in L^{\infty}_{\loc} ((0,T); L^2_{\rm{loc}} (\Omega, |\mu|_{\lambda}))$
we get as a by-product that $u \in L^{\infty}_{\loc} ((0,T); H^1_{\rm{loc}} (\Omega_0, \lambda, \lambda))$. \\
In some cases we can derive that \eqref{tempofissato} can hold for {\em every} $t \in [s_1, s_2]$ (see the previous section).
\eoss

\noindent
The estimates given in Definition \ref{classiDG} are also known as
\emph{energy estimates} or \emph{Caccioppoli's estimates} and we will often 
refer to them in this way. \\ \ \\
\noindent
Now denote by $\mathcal{K}(\Omega \times (0,T))$ the set $\{ K \subset \Omega 
\times (0,T) \, | \, K \textrm{ compact} \}$ and consider the functional
$$
E: L^2(0,T; H^{1}(\Omega)) \times \mathcal{K}(\Omega \times (0,T)) \to \R \, , \hskip30pt
E(w, K) = \frac{1}{2} \, \int\!\!\!\int_K |Dw|^2 \la \, dx dt \, .
$$
We are going to define a $Q$-minimum following the definition given in \cite{wieser} (see also \cite{giagiu} for the elliptic case).

\bd
\label{quasi-min}
We will call a function $u:\Omega \times (0,T) \to \R$ a $Q$-minimum for the equation \eqref{equazione}
if $u \in L^2_{\rm loc}(0,T; H^{1}_{\rm loc}(\Omega,$ $|\mu|, \la)) \cap L^{\infty}_{\loc} ((0,T); L^2_{\rm{loc}} (\Omega, |\mu|_{\lambda}))$
and there is a constant $Q \geqslant 1$ such that
\be
\label{Qminpar}
- \int \!\!\! \int_{{\rm supp} (\phi)} u \frac{\partial \phi}{\partial t} \mu \, dx dt + 
	E(u, {\rm supp}(\phi)) \leqslant Q \, E (u-\phi,{\rm supp}(\phi))
\ee
for every $\phi \in C^1_c (\Omega \times (0,T))$.
\ed

\boss
\label{notasuiQminimi}
It is easy to verify that if $u \in L^2(0,T; H^{1}(\Omega, |\mu|, \la))$ is a $Q$-minimum for 
equation \eqref{equazione}
than the map $L \phi := - \int \!\! \int_{{\rm supp} (\phi)} u \frac{\partial \phi}{\partial t} \mu \, dx dt$
with $\phi \in C^1_c (\Omega \times (0,T))$ turns out to be a linear and continuous form in
$L^2(0,T; H^{1}_0(\Omega,|\mu|, \la))$, i.e. $L$ belongs to the dual space $L^2(0,T; (H^{1}(\Omega,|\mu|, \la))')$
(the proof can be obtained following the analogous one in \cite{wieser}).
\eoss

\noindent
{\bf Solutions are $Q$-minima} -
Following the analogous proof in \cite{wieser} one can verify
that $u$ is a solution of \eqref{equazione} if and only if $u$ is a
$1$-minimum for \eqref{equazione}. \\
A second interesting fact is that a solution of \eqref{equazionegenerale} is a $Q$-minimum for the equation \eqref{equazione}.
Indeed using \eqref{proprieta`} it is easy to see that a solution of \eqref{equazionegenerale}
satisfies \eqref{Qminpar} with $Q = 2 L M$. \\

\noindent
{\bf $Q$-minima belong to the class $DG$} -
We now want to show that the De Giorgi class defined above contains $Q$-minima
and in particular solutions of \eqref{equazione}.
In Section \ref{secHarnack} we will show a Harnack type inequality, 
and then H\"older continuity, for functions in the De Giorgi classes,
and consequently for $Q$-minima and solutions of \eqref{equazione}.
To show this, first of all 
notice that if $u$ satisfies \eqref{Qminpar} for every 
$\phi \in C^1_c (\Omega \times (0,T))$ then, by density of $C^1_c (\Omega \times (0,T))$ in $\W$, $u$ satisfies \eqref{Qminpar} 
also for $\phi \in \W$; then in particular we could choose 
$\phi = (u-k)_+ \zeta^2$ with $\zeta$ 
a Lipschitz continuous and non-negative function such that
$\zeta(\cdot, t) \in \text{Lip}_0(B_{R}(x_0))$,
$|\nabla \zeta|, \zeta_t \in L^{\infty}$, $\zeta_t \mu \geqslant 0$. \\
To show this fact it is sufficient to
consider a point $(x_0, t_0) \in \Omega \times (0,T)$, a function $(u-k)_+ \zeta^2$ with $\zeta$ defined in $[s_1, s_2] \times B_{R}(x_0)$ with
$0 < s_1 < t_0 < s_2 < T$ and
$s_2 - t_0 = \upbeta \, h(x_0, R) R^2$ if $\mu_+ (B_R(x_0)) > 0$,
$t_0 - s_1 = \upbeta \, h(x_0, R) R^2$ if $\mu_- (B_R(x_0)) > 0$, 
while if $B_R(x_0) \subset \Omega_0$ $s_1$ and $s_2$ arbitrary;
then for arbitrary $\sigma_1, \sigma_2$ satisfying $s_1 \leqslant \sigma_1 < \sigma_2 \leqslant s_2$
choose $\phi_\e = (u-k)_+ \zeta^2 \tau_\epsilon$ where
$$
\tau_\epsilon(t)	= \left\{	\begin{array}{ll}
1	&	t \in [\sigma_1, \sigma_2]						\\
\epsilon^{-1}(t - \sigma_1 + \epsilon)
&	t \in [\sigma_1 - \epsilon, \sigma_1]						\\
- \epsilon^{-1}(t - \sigma_2 - \epsilon)
&	t \in [\sigma_2, \sigma_2 + \epsilon]					\\
0	&	t \not\in [\sigma_1 - \epsilon, \sigma_2 + \epsilon]
\end{array}
\right.
$$
for a suitable $\epsilon > 0$. Taking such a $\phi_\e$ in 
\eqref{Qminpar} and letting $\epsilon$ go to zero one gets that
\begin{equation}
\label{prelim_DG}
\begin{array}{l}
{\displaystyle
\frac{1}{2} \int_{B_{R}} (u-k)_+^2(x,\sigma_2) \zeta^2(x,\sigma_2) \, \mu(x) \, dx + E(u, K)  \leqslant Q \, E (u-\phi,K)	+ }	\\	[1.5em]
{\displaystyle
\hskip20pt + \frac{1}{2} \int_{B_{R}} (u-k)_+^2(x,\sigma_1) \zeta^2(x,\sigma_1) \, \mu (x) \, dx +
\int_{\sigma_1}^{\sigma_2} \!\!\! \int_{B_{R}} (u-k)_+^2 \zeta \zeta_t \, \mu \, dx dt }
\end{array}
\end{equation}
where we simply denote $B_{R}$ instead of $B_{R}(x_0)$ and $K$ denotes the part of the support of $\zeta$ contained in $B_R \times [\sigma_1, \sigma_2]$. \\
\ \\
\noindent
$1^{\circ}$ - First suppose $\mu_+(B_{R}(x_0)) > 0$ and show \eqref{DGgamma+} and \eqref{DGgamma+_1}.
We proceed as follows: consider 
$\phi = (u-k)_+ \zeta^2$ with $\zeta$ a Lipschitz continuous function to be choosen later.
Since we have that
$$
u - \phi =	\left\{
\begin{array}{ll}
u			&	u \leqslant k	\\
(u-k)(1-\zeta^2) + k	&	u > k \, .
\end{array}
\right.
$$
and ${\rm supp}(\phi) \subset \{u > k\}$ we have that
\begin{equation}
\label{prelim_DG_2}
\begin{array}{l}
{\displaystyle
E ( u-\phi, {\rm supp}(\phi)) = 
\frac{1}{2}
\int \!\!\! \int_{{\rm supp}(\phi)} \left| D \left[ (u-k)_ + (1-\zeta^2) \right] \right|^2 \la \, dx dt \leqslant }		\\	[1.5em]
\hskip20pt	{\displaystyle \leqslant\int \!\!\! \int_{{\rm supp}(\phi)} \left[ (1-\zeta^2)^2 |D(u-k)_+|^2  
	+ 4 (u-k)_+^2 \zeta^2 |D \zeta|^2 \right] \la \, dx dt \, . }
\end{array}
\end{equation}
We first prove \eqref{DGgamma+}.
We consider $r, \tilde{r} > 0$ with $r < \tilde{r} < R$, $t_0, s_2 \in (0,T)$ with $s_2 - t_0 = \upbeta \, h(x_0,R) R^2$, $\theta, \tilde\theta$ such that
$0 \leqslant \tilde\theta < \theta < 1$.
By assuming in addition that for $\varepsilon \geqslant 0$ (and sufficiently small, say $\varepsilon < R - \tilde{r}$)
$$
K := \text{supp}(\zeta) \cap \big( B_R (x_0) \times [s_1, s_2] \big) \subset Q_{R;{r}, \tilde\theta}^{\upbeta, +,\tilde{r} - r + \varepsilon} (x_0,t_0) 
$$
and that $|\zeta| \leqslant 1$, on the right hand
side we estimate $(1-\zeta^2)^2$ by $1-\zeta^2$ and the second term by $4 (u-k)_+^2 |D \zeta|^2$.
Moreover using the assumption that $u$ is a $Q$-minimum and since
$E ( u, K) = E ( (u-k)_+, K)$ we get that for every $\tau_1, \tau_2 \in [t_0, s_2]$ with $\tau_1 < \tau_2$
\begin{align}
\int_{B_{\tilde{r} + \varepsilon}} & (u-k)_+^2(x,\tau_2) \zeta^2(x,\tau_2) \mu (x) \, dx -
\int_{B_{\tilde{r} + \varepsilon}} (u-k)_+^2(x,\tau_1) \zeta^2(x,\tau_1) \mu (x) \, dx +											\nonumber	\\
& \hskip180pt		+ 2Q \int_{\tau_1}^{\tau_2} \!\!\! \int_{B_{\tilde{r} + \varepsilon}} |D(u-k)_+|^2 \zeta^2 \la \, dx dt \leqslant		\nonumber	\\
& \leqslant 
2 \int_{\tau_1}^{\tau_2} \!\!\! \int_{B_{\tilde{r} + \varepsilon}} (u-k)_+^2 \, \zeta  \zeta_t \, \mu \, dx dt +
				8Q \int_{\tau_1}^{\tau_2} \!\!\! \int_{B_{\tilde{r} + \varepsilon}} (u-k)_+^2 |D \zeta|^2 \la \, dx dt +				\nonumber	\\
& \hskip20pt + (2Q-1) \iint_{Q_{R; r, \tilde\theta}^{\upbeta,+,\tilde{r} - r + \varepsilon} \cap (B_{R} \times [\tau_1, \tau_2])}
																		 |D(u-k)_+|^2 \la \, dx dt \, .				\nonumber
\end{align}
We then choose a Lipschitz continuous function $\zeta$ (see also Figure A below where we show an example where $\mu > 0$ and $\mu < 0$) satisfying also
\begin{equation}
\begin{array}{c}
\zeta = 1 \hskip10pt \text{ in } Q_{R; r, \theta}^{\upbeta,+,\varepsilon} (x_0,t_0) \, ,		\hskip20pt
\zeta = 0 \hskip10pt \text{ in } Q_R^{\upbeta,\texttt{\,>}}(x_0,t_0) \setminus 
												Q_{R; r, \tilde\theta}^{\upbeta,+,\tilde{r} - r + \varepsilon} (x_0,t_0) \, ,			\\	[1em]
{\displaystyle \hskip8pt |D \zeta| \leqslant \frac{1}{\tilde{r} - r} } \, ,	
	\hskip20pt \theta - \tilde\theta = {\displaystyle \frac{(\tilde{r} - r)^2}{R^2} }	\, ,												\\	[1em]
\label{puredifave}
	{\displaystyle |\zeta_t| \leqslant \, \frac{1}{\sigma_{\theta} - \sigma_{\tilde\theta}} =
	\frac{1}{\upbeta \, h(x_0,R) (\tilde{r} - r)^2 } }\, , {\displaystyle \hskip15pt \zeta_t \mu \geqslant 0 }	\, ,	\hskip12pt \zeta_t \mu_- = 0 \, .
\end{array}
\end{equation}
\ \\
\ \\
\ \\
\begin{picture}(150,200)(-180,0)
\put (-105,200){\linethickness{1pt}\line(1,0){210}}
\put (-105,50){\linethickness{1pt}\line(1,0){210}}
\put (-105,50){\linethickness{1pt}\line(0,1){150}}
\put (105,50){\linethickness{1pt}\line(0,1){150}}

\put (30,40){\line(0,1){170}}

\put (-180,125){\line(1,0){320}}
\put (-170,40){\line(0,1){170}}



\put (-95,200){\line(0,-1){4}}
\put (-95,195){\line(0,-1){4}}
\put (-95,190){\line(0,-1){4}}
\put (-95,185){\line(0,-1){4}}
\put (-95,180){\line(0,-1){4}}

\put (-95,175){\line(1,0){4}}
\put (-90,175){\line(1,0){4}}
\put (-85,175){\line(1,0){4}}
\put (-80,175){\line(1,0){4}}
\put (-75,175){\line(1,0){4}}
\put (-70,175){\line(1,0){4}}
\put (-65,175){\line(1,0){4}}
\put (-60,175){\line(1,0){4}}
\put (-55,175){\line(1,0){4}}
\put (-50,175){\line(1,0){4}}
\put (-45,175){\line(1,0){4}}
\put (-40,175){\line(1,0){4}}
\put (-35,175){\line(1,0){4}}
\put (-30,175){\line(1,0){4}}
\put (-25,175){\line(1,0){4}}
\put (-20,175){\line(1,0){4}}
\put (-15,175){\line(1,0){4}}
\put (-10,175){\line(1,0){4}}
\put (-5,175){\line(1,0){4}}
\put (0,175){\line(1,0){4}}
\put (5,175){\line(1,0){4}}
\put (10,175){\line(1,0){4}}
\put (15,175){\line(1,0){4}}

\put (20,175){\line(0,-1){4}}
\put (20,170){\line(0,-1){4}}
\put (20,165){\line(0,-1){4}}
\put (20,160){\line(0,-1){4}}
\put (20,155){\line(0,-1){4}}
\put (20,150){\line(0,-1){4}}
\put (20,145){\line(0,-1){4}}
\put (20,140){\line(0,-1){4}}
\put (20,135){\line(0,-1){4}}
\put (20,130){\line(0,-1){4}}


\put (40,175){\line(0,1){4}}
\put (40,180){\line(0,1){4}}
\put (40,185){\line(0,1){4}}
\put (40,190){\line(0,1){4}}
\put (40,195){\line(0,1){4}}
\put (40,170){\line(0,1){4}}
\put (40,165){\line(0,1){4}}
\put (40,160){\line(0,1){4}}
\put (40,155){\line(0,1){4}}
\put (40,150){\line(0,1){4}}
\put (40,145){\line(0,1){4}}
\put (40,140){\line(0,1){4}}
\put (40,135){\line(0,1){4}}
\put (40,130){\line(0,1){4}}
\put (40,125){\line(0,1){4}}

\put (50,191){\line(0,1){1}}
\put (50,194){\line(0,1){1}}
\put (50,197){\line(0,1){1}}
\put (50,188){\line(0,1){1}}
\put (50,185){\line(0,1){1}}
\put (50,182){\line(0,1){1}}
\put (50,179){\line(0,1){1}}
\put (50,176){\line(0,1){1}}
\put (50,173){\line(0,1){1}}
\put (50,170){\line(0,1){1}}
\put (50,167){\line(0,1){1}}
\put (50,164){\line(0,1){1}}
\put (50,161){\line(0,1){1}}
\put (50,158){\line(0,1){1}}
\put (50,155){\line(0,1){1}}
\put (50,152){\line(0,1){1}}
\put (50,149){\line(0,1){1}}
\put (50,146){\line(0,1){1}}
\put (50,143){\line(0,1){1}}
\put (50,140){\line(0,1){1}}
\put (50,137){\line(0,1){1}}
\put (50,134){\line(0,1){1}}
\put (50,131){\line(0,1){1}}
\put (50,128){\line(0,1){1}}
\put (50,125){\line(0,1){1}}

\put (10,200){\linethickness{2pt}\line(1,0){40}}
\put (35,210){$I_{R,\tilde{r} - r + \varepsilon}(x_0) \times \{ s_2 \}$}
\put (20,50){\linethickness{2pt}\line(1,0){20}}
\put (20,30){$I_{R,\varepsilon}(x_0) \times \{ s_1 \}$}

\put (-40,140){\tiny$\mu > 0$}
\put (40,110){\tiny$\mu < 0$ \text{ or } $\mu = 0$}

\put (0,125){\linethickness{2pt}\line(1,0){1}}
\put (-10,130){\tiny$(x_0,t_0)$}

\put (-185,198){$s_2$}
\put (-169,200){\line(1,0){1}}
\put (-166,200){\line(1,0){1}}
\put (-163,200){\line(1,0){1}}
\put (-160,200){\line(1,0){1}}
\put (-157,200){\line(1,0){1}}
\put (-154,200){\line(1,0){1}}
\put (-151,200){\line(1,0){1}}
\put (-148,200){\line(1,0){1}}
\put (-145,200){\line(1,0){1}}
\put (-141,200){\line(1,0){1}}
\put (-138,200){\line(1,0){1}}
\put (-135,200){\line(1,0){1}}
\put (-132,200){\line(1,0){1}}
\put (-129,200){\line(1,0){1}}
\put (-126,200){\line(1,0){1}}
\put (-123,200){\line(1,0){1}}
\put (-120,200){\line(1,0){1}}
\put (-117,200){\line(1,0){1}}
\put (-117,200){\line(1,0){1}}
\put (-114,200){\line(1,0){1}}
\put (-111,200){\line(1,0){1}}

\put (10,158){\line(0,-1){1}}
\put (10,155){\line(0,-1){1}}
\put (10,152){\line(0,-1){1}}
\put (10,149){\line(0,-1){1}}
\put (10,146){\line(0,-1){1}}
\put (10,143){\line(0,-1){1}}
\put (10,140){\line(0,-1){1}}
\put (10,137){\line(0,-1){1}}

\put (-185,48){$s_1$}
\put (-169,50){\line(1,0){1}}
\put (-166,50){\line(1,0){1}}
\put (-163,50){\line(1,0){1}}
\put (-160,50){\line(1,0){1}}
\put (-157,50){\line(1,0){1}}
\put (-154,50){\line(1,0){1}}
\put (-151,50){\line(1,0){1}}
\put (-148,50){\line(1,0){1}}
\put (-145,50){\line(1,0){1}}
\put (-141,50){\line(1,0){1}}
\put (-138,50){\line(1,0){1}}
\put (-135,50){\line(1,0){1}}
\put (-132,50){\line(1,0){1}}
\put (-129,50){\line(1,0){1}}
\put (-126,50){\line(1,0){1}}
\put (-123,50){\line(1,0){1}}
\put (-120,50){\line(1,0){1}}
\put (-117,50){\line(1,0){1}}
\put (-114,50){\line(1,0){1}}
\put (-111,50){\line(1,0){1}}

\put (-105,160){\line(1,0){1}}
\put (-102,160){\line(1,0){1}}
\put (-99,160){\line(1,0){1}}
\put (-96,160){\line(1,0){1}}
\put (-93,160){\line(1,0){1}}
\put (-90,160){\line(1,0){1}}
\put (-87,160){\line(1,0){1}}
\put (-84,160){\line(1,0){1}}
\put (-81,160){\line(1,0){1}}
\put (-78,160){\line(1,0){1}}
\put (-75,160){\line(1,0){1}}
\put (-72,160){\line(1,0){1}}
\put (-69,160){\line(1,0){1}}
\put (-66,160){\line(1,0){1}}
\put (-63,160){\line(1,0){1}}
\put (-60,160){\line(1,0){1}}
\put (-57,160){\line(1,0){1}}
\put (-54,160){\line(1,0){1}}
\put (-51,160){\line(1,0){1}}
\put (-48,160){\line(1,0){1}}
\put (-45,160){\line(1,0){1}}
\put (-42,160){\line(1,0){1}}
\put (-39,160){\line(1,0){1}}
\put (-36,160){\line(1,0){1}}
\put (-33,160){\line(1,0){1}}
\put (-30,160){\line(1,0){1}}
\put (-27,160){\line(1,0){1}}
\put (-24,160){\line(1,0){1}}
\put (-21,160){\line(1,0){1}}
\put (-18,160){\line(1,0){1}}
\put (-15,160){\line(1,0){1}}
\put (-12,160){\line(1,0){1}}
\put (-9,160){\line(1,0){1}}
\put (-6,160){\line(1,0){1}}
\put (-3,160){\line(1,0){1}}
\put (0,160){\line(1,0){1}}
\put (3,160){\line(1,0){1}}
\put (6,160){\line(1,0){1}}
\put (9,160){\line(1,0){1}}

\put (-30,10){Figure A}

\end{picture}

\ \\
\ \\
\ \\
\noindent
Plugging such a $\zeta$ into the last inequality and dividing by $2Q$ we get that
\begin{align}
\label{stimaprimadellemma}
\frac{1}{2Q} & \int_{B_{r + \varepsilon}^+ \cup I^{r, \tilde{r} - r + \varepsilon}_+} (u-k)_+^2(x,\tau_2) \mu (x) \, dx -
\frac{1}{2Q} \int_{ I^+_{r, \tilde{r} - r + \varepsilon} \cup  I^{r, \varepsilon}_+} (u-k)_+^2(x,\tau_1) \mu (x) \, dx +						\nonumber		\\
& \hskip80pt
	+ \iint_{Q_{R; r, \theta}^{\upbeta, +,\varepsilon} \cap (B_{R} \times [\tau_1, \tau_2])} |D(u-k)_+|^2 \la \, dx dt \leqslant			\nonumber		\\
& \hskip20pt \leqslant \ \frac{1}{2Q} \, 
	\frac{1}{(\tilde{r} - r)^2} \iint_{Q_{R; r, \tilde\theta}^{\upbeta, +,\tilde{r} - r + \varepsilon} \cap (B_{R} \times [\tau_1, \tau_2])} 
			(u-k)_+^2 \,	\left( 8Q \la + \frac{2}{\upbeta \, h(x_0, R)} \mu_+ \right) dx dt +													\\
& \hskip40pt + \frac{2Q-1}{2Q} \iint_{Q_{R; r, \tilde\theta}^{\upbeta, +,\tilde{r} - r + \varepsilon} \cap (B_{R} \times [\tau_1, \tau_2]) } 
																				|D(u-k)_+|^2 \la \, dx dt			\nonumber
\end{align}
with
$$
\tau_1 \in [t_0, t_0 + \sigma_{\tilde\theta}(R)]  \quad \text{and} \quad 
\tau_2 \in [t_0 + \sigma_{\theta}(R), s_2] \, .
$$

\noindent
Before going on with the proof we state two lemmas, the first result is a slight generalization of Lemma 5.1 in \cite{giaquinta}
(see also Section 4 in \cite{wieser}).

\begin{lemma}
\label{giaq}
Consider some non-negative functions 
$f, g_1, g_2 : [t_0, s_2] \times (0,R] \times [0,R] \to [0,M]$,
$F , G : [t_0, s_2] \times (0,R] \times [0, 1) \times [0,R] \to (0, M]$, $M$ positive constant, satisfying
\begin{align}
\label{iterLemma}
f(\tau_2, \rho, \varepsilon) \, + & \, g_2(\tau_1, \rho, \varepsilon) + F(\tau_1, \tau_2; \rho, \vartheta, \varepsilon)
\leqslant	g_1(\tau_1, \rho, \tilde\varepsilon) + g_2(\tau_2, \rho, \tilde\varepsilon) \, + 									\nonumber	\\	
& + \ {\displaystyle \frac{1}{(\tilde\varepsilon - \varepsilon)^2}} \, G(\tau_1, \tau_2; \rho, \tilde\vartheta, \tilde\varepsilon) +
	\delta \,F (\tau_1, \tau_2; \rho, \tilde\vartheta, \tilde\varepsilon) 
\end{align}
and
\begin{align*}
g_1(\tau_1, \rho, \varepsilon) & \leqslant g_1(\tau_1, \tilde\rho, \tilde\varepsilon) \, , \hskip10pt 
	g_2(\tau_2, \rho, \varepsilon) \leqslant g_2(\tau_2, \tilde\rho, \tilde\varepsilon) \, ,					\\
& F (\tau_1, \tau_2; \rho, \vartheta, \varepsilon) \leqslant F (\tau_1, \tau_2; \tilde\rho, \tilde\vartheta, \tilde\varepsilon)
\end{align*}
for every $\tau_1, \tau_2 \in [t_0, s_2]$, $\tau_1 < \tau_2$,
for every $\rho \leqslant \tilde\rho, \tilde\vartheta \leqslant \vartheta , \varepsilon \leqslant \tilde\varepsilon$ and $\delta \in (0,1)$.
Then there is a constant $c > 1$ depending only on $\delta$ such that
\begin{align*}
f(\tau_2, \rho, \varepsilon) & + g_2(\tau_1, \rho, \varepsilon) + F(\tau_1, \tau_2; \rho, \vartheta, \varepsilon)	\leqslant		\\
& \ \leqslant \frac{1}{1 - \delta} \big[ g_1(\tau_1, \rho, \tilde\varepsilon) + g_2(\tau_2, \rho, \tilde\varepsilon) \big] +
{\displaystyle \frac{c}{(\tilde\varepsilon - \varepsilon)^2}} \, G(\tau_1, \tau_2; \rho, \tilde\vartheta, \tilde\varepsilon) \, .
\end{align*}
\end{lemma}
\noindent
\dimo
We take the sequences $\vartheta_n$ and $\varepsilon_n$ defined by ($\eta$ to be chosen)
\begin{align*}
& \vartheta_0 = \vartheta \, , \hskip10pt \vartheta_{n+1} = \vartheta_n + (1-\eta) (\tilde\vartheta - \vartheta) \eta^n ,\qquad \eta \in (0,1) \, ,			\\
& \varepsilon_0 = \varepsilon \, , \hskip10pt \varepsilon_{n+1} = \varepsilon_n + (1-\eta) (\tilde\varepsilon - \varepsilon) \eta^n ,\qquad \eta \in (0,1) \, .
\end{align*}
Notice that
\begin{align*}
\varepsilon_{n+1} - \varepsilon_0 & = \varepsilon_{n+1} - \varepsilon = (\tilde\varepsilon - \varepsilon) (1 - \eta^{n+1}) \, , 						\\
\varepsilon_0 + \sum_{n=0}^{\infty} (\varepsilon_{n+1} & - \varepsilon_n) = \tilde\varepsilon \, , \quad \quad
\vartheta_0 + \sum_{n=0}^{\infty} (\vartheta_{n+1} - \vartheta_n) = \tilde\vartheta \, .
\end{align*}
By \eqref{iterLemma} we have
\begin{align*}
f(\tau_2, \rho, \varepsilon_0) & + g_2(\tau_1, \rho, \varepsilon_0) + F (\tau_1, \tau_2; \rho, \vartheta_0, \varepsilon_0) \leqslant 			\\
\leqslant & \ 	g_1(\tau_1, \rho, \varepsilon_1) + g_2(\tau_2, \rho, \varepsilon_1) \ +												\\
& + \ {\displaystyle \frac{1}{(\varepsilon_1 - \varepsilon_0)^2}} \, G(\tau_1, \tau_2; \rho, \vartheta_{1}, \varepsilon_1) +
	\delta \,F (\tau_1, \tau_2; \rho, \vartheta_{1}, \varepsilon_1) \leqslant														\\
\leqslant & \ g_1(\tau_1, \rho, \varepsilon_1) + g_2(\tau_2, \rho, \varepsilon_1) \ +
	{\displaystyle \frac{1}{(\varepsilon_1 - \varepsilon_0)^2}} \, G(\tau_1, \tau_2; \rho, \vartheta_{1}, \varepsilon_1) +					\\
& + \ \delta \Bigg[ g_1(\tau_1, \rho, \varepsilon_2) + g_2(\tau_2, \rho, \varepsilon_2) +												\\
& + {\displaystyle \frac{1}{(\varepsilon_2 - \varepsilon_1)^2}} \, G(\tau_1, \tau_2; \rho, \vartheta_{2}, \varepsilon_2) +
	\delta F (\tau_1, \tau_2; \rho, \vartheta_{2}, \varepsilon_2) \Bigg] \, .
\end{align*}
By the monotonicity property of the functions we have in fact
\begin{align*}
f(\tau_2, \rho, \varepsilon_0) & + g_2(\tau_1, \rho, \varepsilon_0) + F (\tau_1, \tau_2; \rho, \vartheta_0, \varepsilon_0) \leqslant 		\\
\leqslant & \ (1 + \delta) \Big[ g_1(\tau_1, \rho, \varepsilon_2 ) + g_2(\tau_2, \rho, \varepsilon_2) \Big]\ +							\\
& \ +  \left(\frac{1}{(\varepsilon_1 - \varepsilon_0)^2} + \frac{\delta}{(\varepsilon_2 - \varepsilon_1)^2} \right) 
									G(\tau_1, \tau_2; \rho, \vartheta_{2}, \varepsilon_2) + 								\\
& \ + \delta^2 F (\tau_1, \tau_2; \rho, \vartheta_{2}, \varepsilon_2) \, .
\end{align*}
Iterating $N$ times these inequalities we first get
\begin{align*}
f(\tau_2, \rho, \varepsilon_0) & + g_2(\tau_1, \rho, \varepsilon_0) + F (\tau_1, \tau_2; \rho, \vartheta_0, \varepsilon_0) \leqslant 		
																									\dots \, \leqslant		\\
\leqslant & \ \Big[ g_1(\tau_1, \rho, \varepsilon_{N+1}) + 
		g_2(\tau_2, \rho, \varepsilon_{N+1}) \Big] \sum_{n=0}^N \delta^n +															\\
& \ +  G(\tau_1, \tau_2; \rho, \vartheta_{N+1}, \varepsilon_{N+1}) \sum_{n=0}^N \frac{\delta^n}{(\varepsilon_{n+1} - \varepsilon_n)^2}+ 				\\
& \ + \delta^{N+1} F (\tau_1, \tau_2; \rho, \vartheta_{N+1}, \varepsilon_{N+1}) \, ;
\end{align*}
then taking the limit as $N \to +\infty$ we finally obtain
\begin{align*}
f(\tau_2, \rho, \varepsilon) & + g_2(\tau_1, \rho, \varepsilon) + F (\tau_1, \tau_2; \rho, \vartheta, \varepsilon) \leqslant					\\
& \leqslant \ \frac{1}{1 - \delta} \ \Big[ g_1(\tau_1, \rho, \tilde\varepsilon) + 
								g_2(\tau_1, \rho, \tilde\varepsilon) \Big]	+												\\
&	\quad \quad \quad + G(\tau_1, \tau_2; \rho, \tilde\vartheta, \tilde\varepsilon) 
	 \frac{1}{(\tilde\varepsilon - \varepsilon)^2}\, \frac{1}{(1-\eta)^2} \, \sum_{n=0}^{\infty} \left( \frac{\delta}{\eta^2}\right)^n \, .
\end{align*}
Taking $\eta \in (\sqrt{\delta}, 1)$ we are done. Taking for instance $\eta = \sqrt{(1+\delta)/2}$ one could have $c = (1 + \delta)/(1 - \delta)$.
\finedimo

\noindent
Call
\begin{align}
\label{funzioncine}
f (\tau, \rho, \varepsilon) 		& := \frac{1}{2Q} \int_{B_{\rho + \varepsilon}^+} (u-k)_+^2(x,\tau) \mu_+ (x) \, dx	\, ,					\nonumber	\\
g_2 (\tau, \rho, \varepsilon) 	& := \frac{1}{2Q} \int_{I^{\rho, \varepsilon}_+} (u-k)_+^2(x,\tau) \mu_- (x) \, dx \, , 					\nonumber	\\
g_1 (\tau, \rho, \varepsilon) 	& := \frac{1}{2Q} \int_{I^+_{\rho, \varepsilon}} (u-k)_+^2(x,\tau) \mu_+ (x) \, dx \, ,								\\
F (\tau_1, \tau_2; \rho, \vartheta, \varepsilon) 
		& := \iint_{Q_{R; \rho, \vartheta}^{\upbeta, +,\varepsilon} \cap (B_{R} \times [\tau_1, \tau_2])} |D(u-k)_+|^2 \la \, dx dt \, ,		\nonumber	\\
G (\tau_1, \tau_2; \rho, \vartheta, \varepsilon)
		& := \frac{1}{2Q} \, \iint_{Q_{R;\rho, \vartheta}^{\upbeta, +,\varepsilon} \cap (B_{R} \times [\tau_1, \tau_2])}
								(u-k)_+^2 \, \left( 8Q \la + \frac{2}{\upbeta \, h(x_0, R)} \mu_+ \right) dx dt \, ,				\nonumber
\end{align}
for $ \rho, \vartheta, \varepsilon \geqslant 0$; now we apply the previous lemma in \eqref{stimaprimadellemma}
with $\delta = \frac{2Q - 1}{2Q}$, $\rho = r$, $\tilde{\varepsilon} = \tilde{r} - r + \varepsilon$
and since $(1 - \delta)^{-1} = 2Q$ we derive the existence
of a positive constant $c_Q$ depending only on $Q$ (for instance, as shown at the end of the proof, one could consider
$c_Q = 4Q - 1$) such that
\begin{align}
\label{devoproseguire}
\frac{1}{2Q} \int_{B_{r + \varepsilon}^+} & (u-k)_+^2(x,\tau_2) \mu_+ (x) \, dx +
\frac{1}{2Q} \int_{I^{r, \varepsilon}_+} (u-k)_+^2(x,\tau_1) \mu_- (x) \, dx +													\nonumber	\\
& \hskip80pt
	+ \iint_{Q_{R; r, \theta}^{\upbeta, +,\varepsilon} \cap (B_{R} \times [\tau_1, \tau_2])} |D(u-k)_+|^2 \la \, dx dt \leqslant			\nonumber	\\
\leqslant & \ \int_{I^{r, \tilde{r} - r + \varepsilon}_+} (u-k)_+^2(x,\tau_2) \mu_- (x) \, dx +
	\int_{I^+_{r, \tilde{r} - r + \varepsilon}} (u-k)_+^2(x,\tau_1) \mu_+ (x) \, dx +															\\
& + \ \frac{c_Q}{2Q} \, \frac{1}{(\tilde{r} - r)^2} 
	\iint_{Q_{R;r, \tilde\theta}^{\upbeta, +,\tilde{r} - r + \varepsilon} \cap (B_R \times [\tau_1, \tau_2])} 
	(u-k)_+^2 \,	\left( 8Q \la + \frac{2}{\upbeta \, h(x_0, R)} \, \mu_+ \right) dx dt		\, . 									\nonumber
\end{align}

\noindent
Here is the second lemma, a simple but important lemma.

\begin{lemma}
\label{stimettaDG}
Consider some non-negative functions $f, g_1, g_2, g_3: [t_0, s_2] \to [0,M]$,
$F , G : [s_1, s_2] \to (0, M]$, $M$ positive constant, satisfying
$$
f(\tau_2) + g_3(\tau_1) + \int_{\tau_1}^{\tau_2} F(t) dt  \leqslant g_2(\tau_2) + g_1(\tau_1) + \int_{\tau_1}^{\tau_2} G(t) dt
$$
for every $\tau_1 <\tau_2$. Let $\theta$ and $\tilde\theta$ be the values considered in \eqref{puredifave},
$\sigma_{\theta} = \theta \, \upbeta \, h(x_0, R) R^2$, $\sigma_{\tilde\theta} = \tilde\theta \, \upbeta \, h(x_0, R) R^2$ for some positive $\upbeta$. Then
\begin{align*}
\sup_{t \in (t_0 + \sigma_{\theta}, s_2)} & f(t) + \sup_{t \in (t_0, t_0 + \sigma_{\tilde\theta})} g_3(t) + \int_{t_0}^{s_2} F(t) dt \leqslant			\\
& \leqslant \ 2 \left[ \sup_{t \in (t_0 + \sigma_{\theta}, s_2)} g_2 (t) + \sup_{t \in (t_0, t_0 + \sigma_{\tilde\theta})} g_1 (t) +  \int_{s_1}^{s_2} G(t) dt \right] \, .
\end{align*}
\end{lemma}
\noindent
\dimo
By the assumptions in particular we have
\begin{align*}
f(\tau_2) + g_3(\tau_1) \leqslant g_2(\tau_2) + g_1(\tau_1) + \int_{\tau_1}^{\tau_2} G(t) dt \, ,		\\
\int_{\tau_1}^{\tau_2} F(t) dt \leqslant g_2(\tau_2) + g_1(\tau_1) + \int_{\tau_1}^{\tau_2} G(t) dt \, .
\end{align*}
Taking the supremum in both the inequalities we get
\begin{align*}
\sup_{\tiny \begin{array}{c} 	\tau_1 \in (t_0, t_0 + \sigma_{\tilde\theta}) \\ 
						\tau_2 \in (t_0 + \sigma_{\theta}, s_2) 
		\end{array}} & 
	\big[ f(\tau_2) + g_3(\tau_1) \big] =
	\sup_{\tau_2 \in (t_0 + \sigma_{\theta}, s_2)} f(\tau_2) + \sup_{\tau_1 \in (t_0, t_0 + \sigma_{\tilde\theta})} g_3(\tau_1) \leqslant					\\
& \leqslant \ \sup_{\tiny \begin{array}{c} \tau_1 \in (t_0, t_0 + \sigma_{\tilde\theta}) \\ 
								\tau_2 \in (t_0 + \sigma_{\theta}, s_2) 
		\end{array}}
	\left[ g_2(\tau_2) + g_1(\tau_1) + \int_{\tau_1}^{\tau_2} G(t) dt \right] 	\leqslant														\\
& \leqslant \ \sup_{\tau_2 \in (t_0 + \sigma_{\theta}, s_2)} g_2(\tau_2) + 
	\sup_{\tau_1 \in (t_0, t_0 + \sigma_{\tilde\theta})} g_1(\tau_1) + \int_{t_0}^{s_2} G(t) dt
\end{align*}
and
\begin{align*}
\int_{t_0}^{s_2} F(t) dt
\leqslant \ \sup_{\tau_2 \in (t_0 + \sigma_{\theta}, s_2)} g_2(\tau_2) + 
	\sup_{\tau_1 \in (t_0, t_0 + \sigma_{\tilde\theta})} g_1(\tau_1) + \int_{t_0}^{s_2} G(t) dt  \, .
\end{align*}
Summing the two inequalities we get the thesis.
\finedimo

\noindent
Now we multiply by $2Q$ the inequality \eqref{devoproseguire} and apply the previous lemma. We get
\begin{align*}
\sup_{t \in (t_0 + \sigma_{\theta}, s_2)} & \int_{B_{r + \varepsilon}^+} (u-k)_+^2(x,t) \mu_+ (x) \, dx +  
							\sup_{t \in (t_0, t_0 + \sigma_{\tilde\theta})} \int_{I^{r, \varepsilon}_+} (u-k)_+^2(x,t) \mu_- (x) \, dx +			\\
& \hskip80pt + 2 Q \iint_{Q_{R; r, \theta}^{\upbeta, +,\varepsilon}} |D(u-k)_+|^2 \la \, dx dt \leqslant											\\
\leqslant & \ 4Q \, \sup_{t \in (t_0 + \sigma_{\theta}, s_2)} \int_{I^{r, \tilde{r} - r + \varepsilon}_+} (u-k)_+^2(x,t) \mu_- (x) \, dx +					\\
& \hskip10pt + 4Q \,  \sup_{t \in (t_0, t_0 +  \sigma_{\tilde\theta})} \int_{I^+_{r, \tilde{r} - r + \varepsilon}} (u-k)_+^2(x,t) \mu_+ (x) \, dx +			\\
& \hskip10pt + \, \frac{2 \, c_Q}{(\tilde{r} - r)^2} \iint_{Q_{R;r, \tilde\theta}^{\upbeta, +,\tilde{r} - r + \varepsilon}} 
	(u-k)_+^2 \,	\left( 8Q \la + \frac{2}{\upbeta \, h(x_0, R)} \mu_+ \right) dx dt \, .
\end{align*}
Finally, calling $\gamma$ the quantity $16 \, c_Q \, Q$ (which turns out to be greater than $1$) we get \eqref{DGgamma+}
\begin{align*}
\sup_{t \in (t_0 + \sigma_{\theta}, s_2)} & \int_{B_{r + \varepsilon}^+} (u-k)_+^2(x,t) \mu_+ (x) \, dx + 
								\sup_{t \in (t_0, t_0 +  \sigma_{\tilde\theta})} \int_{I^{r, \varepsilon}_+} (u-k)_+^2(x,t) \mu_- (x) \, dx +		\\
& \hskip150pt + \, \iint_{Q_{R; r, \theta}^{\upbeta, +,\varepsilon}} |D(u-k)_+|^2 \la \, dx dt \leqslant											\\
\leqslant & \ \gamma \Bigg[ \sup_{t \in (t_0,t_0 + \sigma_{\tilde\theta})} \int_{I^+_{r, \tilde{r} - r + \varepsilon}} (u-k)_+^2(x,t) \mu_+ (x) \, dx +		\\
& \hskip25pt + \sup_{t \in (t_0 + \sigma_{\theta}, s_2)} \int_{I^{r, \tilde{r} - r + \varepsilon}_+} (u-k)_+^2(x,t) \mu_- (x) \, dx +						\\
& \hskip25pt +  \, \frac{1}{(\tilde{r} - r)^2} \iint_{Q_{R;r, \tilde\theta}^{\upbeta, +,\tilde{r} - r + \varepsilon}} 
	(u-k)_+^2 \, \left( \la + \frac{1}{\upbeta \, h(x_0, R)} \mu_+ \right) dx dt \Bigg] \, .
\end{align*}
\ \\
\noindent
Now we prove \eqref{DGgamma+_1}. We integrate in $B_R (x_0) \times [\tau_1, \tau_2]$ with
$[\tau_1, \tau_2] \subset [t_0, s_2]$ for an arbitrary $s_2$ (we mean that it is not necessary
to consider $s_2 = t_0 + \upbeta \, h(x_0, R) R^2$) and, as done before to obtain \eqref{prelim_DG},
we get for every $[\tau_1, \tau_2] \subset [t_0, s_2]$
\begin{align*}
\frac{1}{2} \int_{B_{R}} & (u-k)_+^2(x,\tau_2) \zeta^2(x,\tau_2) \, \mu(x) \, dx + E (u,K)
	\leqslant Q \, E (u-\phi,K)	+	\\
& + \frac{1}{2} \int_{B_{R}} (u-k)_+^2(x,\tau_1) \zeta^2(x,\tau_1) \, \mu (x) \, dx +
\int_{\tau_1}^{\tau_2} \!\!\! \int_{B_{R}} (u-k)_+^2 \zeta \zeta_t \, \mu \, dx dt .
\end{align*}
Now choosing $\zeta$ (whose support depends on $\tau$) such that 
\begin{align*}
\zeta = 1 \hskip10pt \text{ in } B_{r}^+(x_0) \times [t_0, \tau]	\, ,	&	\hskip20pt
\zeta = \, 0 \hskip10pt \text{ in } B_R (x_0) \setminus \big( B_{\tilde{r}}^+ (x_0) \cup I^{r, \tilde{r}-r}_+ \big) \times [t_0, \tau] \, ,		\nonumber	\\
\zeta_t \equiv 0 , & \hskip50pt  |D \zeta| \leqslant \frac{1}{\tilde{r} - r}  	\, ,
\end{align*}
using the estimate \eqref{prelim_DG_2} and the inequality which follows it
and taking $\tau_1 = t_0$, we get that for every $\tau \in [t_0, s_2]$
\begin{align*}
\frac{1}{2Q} \int_{B_{r}^+} & (u-k)_+^2(x,\tau) \mu_+ (x) \, dx + \int_{t_0}^{\tau} \!\! \int_{B_{r}^+} |D(u-k)_+|^2 \lambda \, dx dt	\leqslant			\\
\leqslant & \ \frac{1}{2Q} \int_{B_{\tilde{r}}^+} (u-k)_+^2(x,t_0) \mu_+ (x) \, dx +
		\frac{1}{2Q} \int_{I^{r, \tilde{r}-r}_+} (u-k)_+^2(x,\tau) \mu_- (x) \, dx	 +														\\
&	+ \frac{4}{(\tilde{r} - r)^2} \int_{t_0}^{\tau} \!\! \int_{B_{\tilde{r}}^+ \cup I^{r, \tilde{r}-r}_+} (u-k)_+^2 \lambda \, dx dt
	+ \frac{2Q - 1}{2Q} \int_{t_0}^{\tau} \!\! \int_{B_{\tilde{r}}^+ \cup I^{r, \tilde{r}-r}_+} |D(u-k)_+|^2 \la \, dx dt \, .
\end{align*}
As done to obtain \eqref{DGgamma+}, we first use Lemma \ref{giaq} with the analogous functions considered in \eqref{funzioncine}
(notice that with $\varepsilon = 0$ we get $g_2(t_0,r,0) = 0$), then we use Lemma \ref{stimettaDG} to conclude and get \eqref{DGgamma+_1}. \\
\ \\
In an analogous way one can prove \eqref{DGgamma-} and \eqref{DGgamma+_2}, provided that $\mu_- (B_{R}(x_0)) > 0$. \\
\ \\
\noindent
$2^{\circ}$ - We now drop the assumptions $\mu_+ (B_{R}(x_0)) > 0$ and $\mu_- (B_{R}(x_0)) > 0$ and prove \eqref{DGgamma0}.
We recall that in this case we consider $K =  B_R (x_0) \times [s_1, s_2]$ with $s_1$ and $s_2$ arbitrary (but belonging to $[0,T]$).
Now proceeding similarly as before, taking $\phi = (u-k)_+ \zeta^2$ with $\zeta$ independent of $t$ and satisfying
$$
\begin{array}{c}
\zeta \equiv 1 \hskip10pt \text{in } (B_{r}^0 (x_0))^{\varepsilon} \, ,	\hskip20pt
	\zeta \equiv 0 \hskip10pt \text{in } B_{R} (x_0) \setminus (B_r^0 (x_0))^{\tilde{r} - r + \varepsilon} \, ,	\\ [1em]
0 \leqslant \zeta \leqslant 1	\, , 		\hskip20pt
0 \leqslant {\displaystyle |D \zeta| \leqslant \frac{1}{\tilde{r}-r} } \, ,
\end{array}
$$
from \eqref{prelim_DG}, integrating over $(B_r^0)^{\tilde{r}-r + \varepsilon} \times (\tau_1, \tau_2)$,
we derive for every $\tau_1, \tau_2 \in [s_1, s_2]$, $\tau_1 < \tau_2$,
\begin{align}
\frac{1}{2Q} & \int_{I_0^{r, \varepsilon}} (u-k)_+^2(x,\tau_2) \mu_+ (x) \, dx +
	\frac{1}{2Q} \int_{I_0^{r, \varepsilon}} (u-k)_+^2(x,\tau_1) \mu_- (x) \, dx +												\nonumber		\\
& \hskip80pt
	+ \iint_{Q_{R;r; \tau_1, \tau_2}^{0,\varepsilon} } |D(u-k)_+|^2 \la \, dx dt \leqslant											\nonumber		\\
& \hskip20pt \leqslant \ \frac{1}{2Q} \int_{ I_0^{r, \tilde{r} - r + \varepsilon}} (u-k)_+^2(x,\tau_2) \mu_- (x) \, dx +
	\frac{1}{2Q} \int_{ I_0^{r, \tilde{r} - r + \varepsilon}} (u-k)_+^2(x,\tau_1) \mu_+ (x) \, dx \, +								\nonumber		\\
& \hskip40pt + \frac{4}{(\tilde{r} - r)^2} \iint_{Q_{R; r; \tau_1, \tau_2}^{0,\tilde{r}-r + \varepsilon} } (u-k)_+^2 \, \lambda \, dx dt +
	\frac{2Q-1}{2Q} \iint_{Q_{R; r ; \tau_1, \tau_2}^{0,\tilde{r}-r + \varepsilon} } |D(u-k)_+|^2 \la \, dx dt	\, .						\nonumber	
\end{align}
We can apply Lemma \ref{giaq} with $\vartheta = \tilde\vartheta = 0$, $\rho = r$, $\tilde{\rho} = \tilde{r}$,
$\varepsilon \geqslant 0$,  $\tilde\varepsilon = \tilde{r} - r$,
$\delta = (2Q-1)/2Q$ and
\begin{align*}
g_2 (\tau, \rho, \epsilon) 	& := \frac{1}{2Q} \int_{I_0^{\rho, \epsilon}} (u-k)_+^2(x,\tau) \mu_- (x) \, dx \, , 						\\
f (\tau, \rho, \epsilon) =
g_1 (\tau, \rho, \epsilon) 	& := \frac{1}{2Q} \int_{I_0^{\rho, \epsilon}} (u-k)_+^2(x,\tau) \mu_+ (x) \, dx \, ,						\\
F (\tau_1, \tau_2; \rho, \vartheta, \epsilon) 
						& := \iint_{Q_{R;\rho; \tau_1, \tau_2}^{0,\epsilon} } |D(u-k)_+|^2 \la \, dx dt \, , 					\\
G (\tau_1, \tau_2; \rho, \vartheta, \epsilon)
						& := 4 \iint_{Q_{R;\rho; \tau_1, \tau_2}^{0,\epsilon} }   (u-k)_+^2 \, \lambda \,  dx dt \, ,
\end{align*}
and get the existence of $c_Q$ such that
\begin{align}
& \frac{1}{2Q} \int_{I_0^{r, \varepsilon}} (u-k)_+^2(x,\tau_2) \mu_+ (x) \, dx + 
	\frac{1}{2Q} \int_{I_0^{r, \varepsilon}} (u-k)_+^2(x,\tau_1) \mu_- (x) \, dx +												\nonumber	\\
& \hskip150pt
	+ \iint_{Q_{R;r; \tau_1, \tau_2}^{0,\varepsilon} } |D(u-k)_+|^2 \la \, dx dt \leqslant											\nonumber	\\
& \hskip20pt \leqslant \ \int_{ I_0^{r, \tilde{r} - r + \varepsilon}} (u-k)_+^2(x,\tau_2) \mu_- (x) \, dx +
			\int_{ I_0^{r, \tilde{r} - r + \varepsilon}} (u-k)_+^2(x,\tau_1) \mu_+ (x) \, dx \, +									\nonumber	\\
& \hskip100pt + \frac{4 c_Q}{(\tilde{r} - r)^2} 	
						\iint_{Q_{R; r ; \tau_1, \tau_2}^{0,\tilde{r}-r + \varepsilon} } (u-k)_+^2 \, \lambda \, dx dt \, .				\nonumber	
\end{align}
Taking the supremum for $\tau_1, \tau_2 \in (s_1, s_2)$ we get that $u$ satisfies \eqref{DGgamma0} with $\gamma = 4 c_Q$. \\

\section{Local boundedness for functions in $DG$ }
\label{paragrafo5}

In this section we prove that functions belonging to the De Giorgi class are locally bounded in $\Omega \times (0,T)$. \\
We start proving that a generic function $u \in DG (\Omega, T, \mu, \lambda, \gamma)$ is bounded in 
$(B_{\rho} \times (a,b) )\cap (\Omega_+ \times (0,T))$ for some set $B_{\rho} \times (a,b) \subset \subset \Omega \times (0,T)$. \\
Fix $x_0 \in \Omega$, $t_0 \in (0,T)$, $R > 0$ and in what follows assume
$$
\mu_+ (B_{R}(x_0)) > 0 \, .
$$
Then consider $\upbeta > 0$ and $s_2 \in (0,T)$ with
$$
s_2 - t_0 = \upbeta \, h(x_0, R) R^2 \, , \hskip20pt B_R(x_0) \times (t_0, s_2) \subset \Omega \times (0,T) \, . 
$$
Consider now $r , \tilde{r} , \hat{r} \in (0,R]$ such that
$$
\frac{R}{2} \leqslant r < \tilde{r}  < \hat{r} \leqslant  R \hskip10pt \text{and} \hskip10pt  \tilde{r} - r = \frac{\hat{r} - \tilde{r}}{2}
$$
and $\theta, \tilde\theta, \hat\theta$ such that
$$
0 \leqslant \hat\theta < \tilde\theta < \theta < 1\hskip12pt \text{and} \hskip12pt
	 \tilde\theta - \hat\theta = \frac{(\hat{r} - \tilde{r})^2}{R^2} \, , \hskip5pt 
	 \theta - \tilde\theta = \frac{(\tilde{r} - r)^2}{R^2}
$$
and define analogously as done in \eqref{sigmateta} (but here we simplify the notation)
$$
\sigma := \theta \, \upbeta \, h(x_0,R) \, R^2 \, , \hskip20pt
\tilde\sigma := \tilde\theta \, \upbeta \, h(x_0,R) \, R^2 \, ,
\hskip20pt
\hat\sigma := \hat\theta \, \upbeta \, h(x_0,R) \, R^2 \, ,
$$
in such a way that
$$
0 \leqslant \hat\sigma < \tilde\sigma < \sigma < s_2 - t_0 \, .
$$
Since $t_0, x_0$ will remain fixed we will often use the following simplified notations:
we will write
$$
h(\rho), \, B_{\rho}, \, Q_R^{\upbeta,+}, \, Q_R^{\upbeta,\texttt{\,>}}, Q_{R; \rho, \theta}^{\upbeta,+,\delta}, \, Q_{R; \rho, \theta}^{\upbeta,+}
$$
instead of respectively
$$
h(x_0, \rho), \,  B_{\rho}(x_0), \, Q_R^{\upbeta,+}(x_0, t_0), \, Q_R^{\upbeta,\texttt{\,>}}(x_0, t_0), \, Q_{R; \rho, \theta}^{\upbeta,+,\delta}(x_0, t_0), \, 
Q_{R; \rho, \theta}^{\upbeta,+} (x_0,t_0).
$$
In fact, to further simplify the notations, we will suppose that (it is always possible, up to a translation)
$$
t_0 = 0 \, .
$$
Finally, from now on, we will use this short notations for the following measures
$$
\begin{array}{c}
M := \mu \otimes {\mathcal L}^1 \, , \hskip20pt \Lambda := \lambda \otimes {\mathcal L}^1 \, ,
	\hskip20pt |M|_{\Lambda} :=  |\mu|_{\lambda} \otimes {\mathcal L}^1 \, ,						\\	[0.5em]
M_+ := \mu_+ \otimes {\mathcal L}^1 \, , \hskip20pt M_- := \mu_- \otimes {\mathcal L}^1 \, ,			\\	[0.5em]
\Lambda_+ := \lambda_+ \otimes {\mathcal L}^1 \, , \hskip20pt
	\Lambda_- := \lambda_- \otimes {\mathcal L}^1 \, , \hskip20pt
			\Lambda_0 := \lambda_0 \otimes {\mathcal L}^1
\end{array}
$$
where we recall that $\lambda_+, \lambda_-, \lambda_0$ have been defined in \eqref{lambda}. \\
\ \\
Now fix a function $u \in DG (\Omega, T, \mu, \lambda, \gamma)$ and define (since $\upbeta$ will remain fix we omit it in the definition of the following set)
$$
\begin{array}{c}
A_R^{+,\delta}(k; \rho, \theta) = \{ (x,t) \in Q_{R; \rho, \theta}^{\upbeta,+,\delta} \, | \, u(x,t) > k \} \, .
\end{array}
$$
Consider a function $\zeta \in \text{Lip} (B_{\tilde{r}} (x_0) \times [t_0, s_2])$ such that
$\zeta (\cdot, t) \in \text{Lip}_c (B_{\tilde{r}} (x_0))$ for every $t$ such that
(notice that $\tilde{r} - \frac{R}{2} = r - \frac{R}{2} + (\tilde{r} - r)$ and
$\hat{r} - \frac{R}{2} = \tilde{r} - \frac{R}{2} + (\hat{r} - \tilde{r})$)
$$
\begin{array}{c}
\zeta \equiv 1 \hskip10pt \text{in } Q_{R; \frac{R}{2}, \theta}^{\upbeta,+, r - \frac{R}{2}} (x_0,t_0) \, ,	\hskip20pt
	\zeta \equiv 0 \hskip10pt \text{in } Q_R^{\upbeta,\texttt{\,>}} (x_0,t_0) \setminus 
	Q_{R; \frac{R}{2}, \tilde\theta}^{\upbeta,+,\tilde{r} - \frac{R}{2}} (x_0,t_0) \, ,		\, , 								\\ [1em]
0 \leqslant \zeta \leqslant 1 \, , \hskip15pt {\displaystyle |D \zeta| \leqslant \frac{1}{\tilde{r}-r} \, , 
	\hskip15pt 0 \leqslant \zeta_t \mu \, , \hskip15pt \zeta_t \mu_- = 0 \, , \hskip15pt
	|\zeta_t| \leqslant \frac{1}{\upbeta \, h(x_0,R)} \frac{1}{(\tilde{r}-r)^2} } \, .
\end{array}
$$
In what follows we will denote by $Q_{R; R/2, \tilde\theta}^{\upbeta,+,\tilde{r} - R/2}(s)$ the set
$\{ (x,t) \in Q_{R; R/2, \tilde\theta}^{\upbeta,+,\tilde{r} - R/2} \, | \, t = s \}$. \\
First using H\"older's inequality, then
applying Corollary \ref{cor-gut-whee} to the function $(u-k)_+\zeta$
with $\upsilon = \nu = |\mu|_{\lambda}$ and $\omega = \lambda$, $E = Q_{R; R/2, \tilde\theta}^{\upbeta,+,\tilde{r} - R/2} \cap \Omega_+$
(we integrate first in $Q_{R; R/2, \theta}^{\upbeta, +,r - R/2}$, then in $Q_{R; R/2, \tilde\theta}^{\upbeta,+,\tilde{r} - R/2}$, 
with respect to the measure $\mu_+ dx dt$ which is supported in $E$), we estimate
\begin{align*}
\frac{1}{|\mu|_{\lambda} (B_R)} &
\int\!\!\!\int_{Q_{R; R/2, \theta}^{\upbeta,+,r - R/2}}  (u-k)_+^2 \mu_+ \, dx dt
	\leqslant \frac{1}{|\mu|_{\lambda} (B_R)}
	\iint_{Q_{R; R/2, \tilde\theta}^{\upbeta,+,\tilde{r} - R/2}} (u-k)_+^2 \zeta^2 \mu_+ \, dx dt 		\leqslant										\\
\leqslant & \, \frac{\big( M_+ (A_R^{+,\tilde{r} - R/2}(k; R/2, \tilde\theta)) \big)^{\frac{\kappa - 1}{\kappa}}}
			{(|\mu|_{\lambda} (B_R))^{\frac{\kappa - 1}{\kappa}}} 
		\left[\frac{1}{|\mu|_{\lambda} (B_R)} \iint_{Q_{R; R/2, \tilde\theta}^{\upbeta, +,\tilde{r} - R/2}} 
		(u-k)_+^{2\kappa} \zeta^{2\kappa} \mu_+ \, dx dt \right]^{\frac{1}{\kappa}} \leqslant													\\
\leqslant & \, \frac{\big( M_+ (A_R^{+,\tilde{r} - R/2}(k; R/2, \tilde\theta)) \big)^{\frac{\kappa - 1}{\kappa}}}
			{(|\mu|_{\lambda} (B_R))^{\frac{\kappa - 1}{\kappa}}} 
	\, \gamma_1^{2/\kappa} R^{2/\kappa}
	\left(\frac{1}{|\mu|_{\lambda} (B_R)}\right)^{\frac{\kappa -1}{\kappa}}  \frac{1}{(\lambda(B_R))^{1/\kappa}}  \cdot 								\\
& \hskip40pt \cdot \Bigg( \sup_{0 < t < s_2} \int_{Q_{R; R/2, \tilde\theta}^{\upbeta, +,\tilde{r} - R/2}(t)} (u-k)_+^{2}(x,t) \zeta^2(x,t) \mu_+ (x) \, dx 
		\Bigg)^{\frac{\kappa-1}{\kappa}} \!\!\!\!\!\!\!  \cdot																				\\
& \hskip60pt \cdot \Bigg( \iint_{Q_{R; R/2, \tilde\theta}^{\upbeta, +,\tilde{r} - R/2}} |D ((u-k)_+\zeta)|^2 (x,t) \, \lambda (x) \, dx dt
		\Bigg)^{\frac{1}{\kappa}} \leqslant																							\\
\leqslant & \, \frac{\big( M_+ (A_R^{+,\tilde{r} - R/2}(k; R/2, \tilde\theta)) \big)^{\frac{\kappa - 1}{\kappa}}}
			{(|\mu|_{\lambda} (B_R))^{\frac{\kappa - 1}{\kappa}}} 
\hskip5pt \gamma_1^{2/\kappa} \frac{R^{2/\kappa}}{(\lambda(B_R))^{1/\kappa}}
	 	\left(\frac{1}{|\mu|_{\lambda} (B_R)}\right)^{\frac{\kappa -1}{\kappa}} \cdot															\\
& \hskip40pt \cdot \Bigg( \sup_{0 < t < s_2} \int_{Q_{R; R/2, \tilde\theta}^{\upbeta, +,\tilde{r} - R/2}(t)} (u-k)_+^{2}(x,t) \zeta^2(x,t) \mu_+ (x) \, dx +		\\
& \hskip80pt + \iint_{Q_{R; R/2, \tilde\theta}^{\upbeta, +,\tilde{r} - R/2}} |D ((u-k)_+\zeta)|^2 (x,t) \, \lambda (x) \, dx dt \Bigg)	\leqslant					\\
\leqslant & \, \frac{\big( M_+ (A_R^{+,\tilde{r} - R/2}(k; R/2, \tilde\theta)) \big)^{\frac{\kappa - 1}{\kappa}}}
			{(|\mu|_{\lambda} (B_R))^{\frac{\kappa - 1}{\kappa}}} 
\hskip5pt \gamma_1^{2/\kappa} \frac{R^{2/\kappa}}{(\lambda(B_R))^{1/\kappa}}
	 	\left(\frac{1}{|\mu|_{\lambda} (B_R)}\right)^{\frac{\kappa -1}{\kappa}}  \cdot 															\\
& \cdot \Bigg( \sup_{0 < t < s_2} \int_{Q_{R; R/2, \tilde\theta}^{\upbeta, +,\tilde{r} - R/2}(t)} (u-k)_+^{2}(x,t) \mu_+ (x) \, dx +							\\
& \quad + 2 \iint_{Q_{R; R/2, \tilde\theta}^{\upbeta, +,\tilde{r} - R/2}} |D (u-k)_+|^2 (x,t) \, \lambda (x) \, dx dt +
		\frac{2}{(\tilde{r} - r)^2} \iint_{Q_{R; R/2, \tilde\theta}^{\upbeta, +,\tilde{r} - R/2}} (u-k)_+^2 (x,t) \, \lambda (x) \, dx dt \Bigg)					\\
\leqslant & \, \frac{\big( M_+ (A_R^{+,\tilde{r} - R/2}(k; R/2, \tilde\theta)) \big)^{\frac{\kappa - 1}{\kappa}}}
			{(|\mu|_{\lambda} (B_R))^{\frac{\kappa - 1}{\kappa}}} 
\hskip5pt \gamma_1^{2/\kappa} \frac{R^{2/\kappa}}{(\lambda(B_R))^{1/\kappa}}
	 	\left(\frac{1}{|\mu|_{\lambda} (B_R)}\right)^{\frac{\kappa -1}{\kappa}}  \cdot 															\\
& \cdot \Bigg( \sup_{t \in (\tilde\sigma, s_2)} \int_{B_{\tilde{r}}^+} (u-k)_+^2 (x,t) \mu_+ (x) dx +	
	\sup_{t \in (0, \tilde\sigma)} \int_{I_{R/2, \tilde{r} - R/2}^+} (u-k)_+^2 (x,t) \mu_+ (x) dx \, +												\\
& \hskip10pt + 2 \iint_{Q_{R; R/2, \tilde\theta}^{\upbeta, +,\tilde{r} - R/2}} |D (u-k)_+|^2 (x,t) \, \lambda (x) \, dx dt +
	\frac{8}{(\hat{r} - \tilde{r})^2} \iint_{Q_{R; R/2, \tilde\theta}^{\upbeta, +,\tilde{r} - R/2}} (u-k)_+^2 (x,t) \, \lambda (x) \, dx dt \Bigg)
\end{align*}
where in the last inequality we have used the fact that $2 (\tilde{r} - r) = \hat{r} - \tilde{r}$. \\
Now we can continue using the energy estimates \eqref{DGgamma+} (with $\varepsilon = \tilde{r} - R/2$)
\begin{align*}
& \frac{1}{|\mu|_{\lambda} (B_R)}
\int\!\!\!\int_{Q_{R;R/2, \theta}^{\upbeta, +,r - R/2}}  (u-k)_+^2 \mu_+ \, dx dt \leqslant
		\frac{\big( M_+ (A_R^{+,\tilde{r} - R/2}(k; R/2, \tilde\theta)) \big)^{\frac{\kappa - 1}{\kappa}}}{(|\mu|_{\lambda} (B_R))^{\frac{\kappa - 1}{\kappa}}}
		\hskip5pt \gamma_1^{2/\kappa} \frac{R^{2/\kappa}}{(\lambda(B_R))^{1/\kappa}}
	 	\left(\frac{1}{|\mu|_{\lambda} (B_R)}\right)^{\frac{\kappa -1}{\kappa}}  \cdot 														\\
& \hskip25pt \cdot \Bigg[ 2\gamma \sup_{t \in (0, \hat\sigma)} \int_{I_{R/2, \hat{r} - R/2}^+} (u-k)_+^2 (x,t) \mu_+(x) \, dx +
2 \gamma \sup_{t \in (\tilde\sigma, s_2)} \int_{I^{R/2, \hat{r} - R/2}_+} (u-k)_+^2 (x,t) \mu_-(x) \, dx +										\\
& \hskip35pt + \frac{2 \gamma}{(\hat{r} - \tilde{r})^2} \iint_{Q_{R; R/2, {\hat\theta}}^{\upbeta, +,\hat{r} - R/2}} 
	(u-k)_+^2\, \left( \frac{\mu_+}{\upbeta \, h(R)} + \la \right) \, dx ds + 
	\sup_{t \in (\hat\sigma, \tilde\sigma)} \int_{I_{R/2,\tilde{r} - R/2}^+} (u-k)_+^2 (x,t) \mu_+ (x) dx +										\\
& \hskip35pt + \frac{8}{(\hat{r} - \tilde{r})^2} 
		\iint_{Q_{R; \tilde{r}, \tilde\theta}^{\upbeta, +,\tilde{r} - R/2}} (u-k)_+^2 (x,t) \, \lambda (x) \, dx dt \Bigg] \leqslant						\\
& \hskip10pt \leqslant
		\frac{\big( M_+ (A_R^{+,\tilde{r} - R/2}(k; R/2, \tilde\theta)) \big)^{\frac{\kappa - 1}{\kappa}}}{(|\mu|_{\lambda} (B_R))^{\frac{\kappa - 1}{\kappa}}}
		\hskip5pt \gamma_1^{2/\kappa} \frac{R^{2/\kappa}}{(\lambda(B_R))^{1/\kappa}}
	 	\left(\frac{1}{|\mu|_{\lambda} (B_R)}\right)^{\frac{\kappa -1}{\kappa}}  \cdot 														\\
& \hskip20pt \cdot \left[
\frac{2 \gamma + 8}{(\hat{r} - \tilde{r})^2} \iint_{Q^{\upbeta, +,\hat{r} - R/2}_{R; R/2, \hat\theta}}
		(u-k)_+^2 \left( \frac{\mu_+}{\upbeta \, h(R)} + \lambda \right) dx dt +
(2 \gamma+1)  \sup_{t \in (0, s_2)} \int_{(I_{R/2}^+)^{\hat{r} - R/2}} (u-k)_+^2 (x,t) |\mu| (x) dx \right] =										\\
& \hskip10pt = \frac{\big( M_+ (A_R^{+,\tilde{r} - R/2}(k; R/2, \tilde\theta)) \big)^{\frac{\kappa - 1}{\kappa}}}{(|\mu|_{\lambda} (B_R))^{\frac{\kappa - 1}{\kappa}}}
		\hskip5pt \gamma_1^{2/\kappa} \frac{R^{2/\kappa}}{(\lambda(B_R))^{1/\kappa}}
	 	\left(\frac{1}{|\mu|_{\lambda} (B_R)}\right)^{\frac{\kappa -1}{\kappa}} \frac{2 \gamma + 8}{(\hat{r} - \tilde{r})^2} \, \cdot					\\
& \hskip20pt \cdot {\lambda (B_{R}) }
	\Bigg[ \frac{1}{\upbeta \, |\mu|_{\lambda} (B_{R})} \iint_{Q^{\upbeta, +,\hat{r} - R/2}_{R; R/2, \hat\theta}} (u-k)_+^2 \mu_+ \, dx dt
		+ \frac{1}{\lambda (B_{R})}	\iint_{Q^{\upbeta, +,\hat{r} - R/2}_{R; R/2, \hat\theta}} (u-k)_+^2 \lambda \, dx dt \, +						\\
& \hskip60pt	+ \frac{2\gamma + 1}{2\gamma +8} \, (\hat{r} - \tilde{r})^2 \, \frac{1}{\lambda (B_{R})}
	\sup_{t \in (0, s_2)} \int_{(I_{R/2}^+)^{\hat{r} - R/2}} (u-k)_+^2 (x,t) |\mu| (x) dx \Bigg] \, .
\end{align*}
Now we divide by $s_2 - t_0 = \upbeta \, h({R})R^2$,
estimate $\frac{2\gamma + 1}{2\gamma +8}$ by $1$
and finally multiply and divide in the right hand side by $(\upbeta \, h({R}) R^2)^{\frac{\kappa - 1}{\kappa}}$. We get
\begin{align}
\frac{1}{|M|_{\Lambda} (Q_R^{\upbeta, \texttt{\,>}})} &
\int\!\!\!\int_{Q_{R; R/2, \theta}^{\upbeta, +,r - R/2}}  (u-k)_+^2 \mu_+ \, dx dt \leqslant										\nonumber	\\
\leqslant & \hskip5pt \gamma_1^{2/\kappa} \, 
	R^2 \, \upbeta^{\frac{\kappa - 1}{\kappa}} \, \, \frac{2 \gamma + 8}{(\hat{r} - \tilde{r})^2} \, 
	\frac{\big( M_+ (A_R^{+,\tilde{r} - R/2}(k; R/2, \tilde\theta)) \big)^{\frac{\kappa - 1}{\kappa}}}
		{(|M|_{\Lambda} (Q_R^{\upbeta, \texttt{\,>}}))^{\frac{\kappa - 1}{\kappa}} } \,  \left( \frac{1}{\upbeta} + 1 \right) \cdot		\nonumber	\\
\label{mitoccanum1}
& \cdot \Bigg[ \frac{1}{|M|_{\Lambda} (Q_R^{\upbeta, \texttt{\,>}})} \iint_{Q^{\upbeta, +,\hat{r} - R/2}_{R; R/2, \hat\theta}} (u-k)_+^2 \mu_+ \, dx dt
	+ \frac{1}{\Lambda (Q_R^{\upbeta, \texttt{\,>}})} \iint_{Q^{\upbeta, +,\hat{r} - R/2}_{R; R/2, \hat\theta}} (u-k)_+^2 \lambda_+ \, dx dt \, + 		\\
& \hskip50pt + \frac{1}{\Lambda (Q_R^{\upbeta, \texttt{\,>}})} 
			\iint_{Q^{\upbeta, +,\hat{r} - R/2}_{R; R/2, \hat\theta}} (u-k)_+^2 (\lambda_0 + \lambda_-) \, dx dt \, +					\nonumber	\\
& \hskip50pt	+  (\hat{r} - \tilde{r})^2 \, \frac{1}{\Lambda (Q_R^{\upbeta, \texttt{\,>}})}
	\sup_{t \in (0, s_2)} \int_{(I_{R/2}^+)^{\hat{r} - R/2}} (u-k)_+^2 (x,t) |\mu| (x) dx \Bigg] \, .									\nonumber
\end{align}

\noindent
Notice that
$$
\iint_{Q^{\upbeta, +,\hat{r} - R/2}_{R; R/2, \hat\theta}} (u-k)_+^2 (\lambda_0 + \lambda_-) \, dx dt
\quad \text{is in fact} \quad
\int_{0}^{s_2} \!\!\!\! \int_{(I_{R/2}^+)^{\hat{r} - R/2}} (u-k)_+^2 (\lambda_0 + \lambda_-) \, dx dt .
$$
\noindent
In a similar way one can estimate $\int\!\!\!\int_{Q_{R; R/2, \theta}^{\upbeta, +,r - R/2}}  (u-k)_+^2 \lambda_+ \, dx dt$.
The main difference is that we use Corollary \ref{cor-gut-whee} with $\nu = |\mu|_{\lambda} $ and $\upsilon = \omega = \lambda$. We get
\begin{align}
\frac{1}{\Lambda (Q_R^{\upbeta, \texttt{\,>}})} &
\int\!\!\!\int_{Q_{R; R/2, \theta}^{\upbeta, +,r - R/2}}  (u-k)_+^2 \lambda_+ \, dx dt \leqslant										\nonumber	\\
\leqslant & \hskip5pt \gamma_1^{2/\kappa} \, 
	R^2 \, \frac{1 + \upbeta}{\upbeta^{\frac{1}{\kappa}}} \, \frac{2 \gamma + 8}{(\hat{r} - \tilde{r})^2} \, 
	\frac{\big( \Lambda_+ (A_R^{+,\tilde{r} - R/2}(k; R/2, \tilde\theta)) \big)^{\frac{\kappa - 1}{\kappa}}}
										{(\Lambda (Q_R^{\upbeta, \texttt{\,>}}))^{\frac{\kappa - 1}{\kappa}}}\, \cdot 		\nonumber	\\
\label{mitoccanum2}
& \cdot \Bigg[ \frac{1}{|M|_{\Lambda} (Q_R^{\upbeta, \texttt{\,>}})} \iint_{Q^{\upbeta, +,\hat{r} - R/2}_{R; R/2, \hat\theta}} (u-k)_+^2 \mu_+ \, dx dt
	+ \frac{1}{\Lambda (Q_R^{\upbeta, \texttt{\,>}})} \iint_{Q^{\upbeta, +,\hat{r} - R/2}_{R; R/2, \hat\theta}} (u-k)_+^2 \lambda_+ \, dx dt \, + 		\\
& \hskip50pt + \frac{1}{\Lambda (Q_R^{\upbeta, \texttt{\,>}})} 
	\iint_{Q^{\upbeta, +,\hat{r} - R/2}_{R; R/2, \hat\theta}} (u-k)_+^2 (\lambda_0 + \lambda_-) \, dx dt \, +							\nonumber	\\
& \hskip50pt	+  (\hat{r} - \tilde{r})^2 \, \frac{1}{\Lambda (Q_R^{\upbeta, \texttt{\,>}})}
	\sup_{t \in (0, s_2)} \int_{(I_{R/2}^+)^{\hat{r} - R/2}} (u-k)_+^2 (x,t) |\mu| (x) dx \Bigg] \, .									\nonumber
\end{align}
Once defined (for $\rho \in [R/2, R]$)
\begin{align*}
\tilde{u}_{\mu_+} (l;\rho, \vartheta; \varepsilon) & := \left( \frac{1}{|M|_{\Lambda} (Q_{R}^{\upbeta, \texttt{\,>}})} 
				\int\!\!\!\int_{Q_{R; \rho, \vartheta}^{\upbeta, +,\varepsilon}} (u-l)_+^2 \mu_+ \, dx dt \right)^{1/2} \, ,				\\
\tilde{u}_{\lambda_+} (l;\rho, \vartheta; \varepsilon) & := \left( \frac{1}{\Lambda (Q_{R}^{\upbeta, \texttt{\,>}})} 
				\int\!\!\!\int_{Q_{R; \rho, \vartheta}^{\upbeta, +,\varepsilon}} (u-l)_+^2 \lambda_+ \, dx dt \right)^{1/2} \, ,				\\
\big( \tilde{u}_+ (l;\rho, \vartheta, \varepsilon) \big)^2 & := \big( \tilde{u}_{\mu_+}  (l;\rho, \vartheta, \varepsilon) \big)^2 + 
												\big( \tilde{u}_{\lambda_+}  (l;\rho, \vartheta, \varepsilon) \big)^2		\, ,
\end{align*}
we sum the two inequalities and get
\begin{align*}
\big( \tilde{u}_+ (k; & \, {\textstyle \frac{R}{2}}, \theta; \, r - {\textstyle \frac{R}{2}}) \big)^2 \leqslant \frac{C_1}{(\hat{r} - \tilde{r})^2} \,  
\left[ \frac{\big( M_+ (A_R^{+,\tilde{r} - R/2}(k; \frac{R}{2}, \tilde\theta)) \big)^{\frac{\kappa - 1}{\kappa}}}
		{|M|_{\Lambda} (Q_R^{\upbeta, \texttt{\,>}}))^{\frac{\kappa - 1}{\kappa}}} \, + \right.										\\
& \left. + \, \frac{\big( \Lambda_+ (A_R^{+,\tilde{r} - R/2}(k; \frac{R}{2}, \tilde\theta)) \big)^{\frac{\kappa - 1}{\kappa}}}
		{(\Lambda (Q_R^{\upbeta, \texttt{\,>}}))^{\frac{\kappa - 1}{\kappa}}} \right] \cdot
		\left[ \big( \tilde{u}_+ (k; {\textstyle \frac{R}{2}}, \hat\theta; \hat{r} - {\textstyle \frac{R}{2}}) \big)^2 + 
		\big( \omega^{\hat{r} - \tilde{r}} (u; k; \hat{r} ; \hat\theta) \big)^2 \right]
\end{align*}
where $C_1 = \gamma_1^{2/\kappa} \, R^2 \, {\displaystyle \frac{1 + \upbeta}{\upbeta^{\frac{1}{\kappa}}}} \, (2 \gamma + 8)$ and
\begin{align*}
\big( \omega^{\hat{r} - \tilde{r}} (u; k; \hat{r} ; \hat\theta) \big)^2:= & 
	\, \frac{1}{\Lambda (Q_R^{\upbeta, \texttt{\,>}})} 
			\iint_{Q^{\upbeta, +,\hat{r} - R/2}_{R; R/2, \hat\theta}} (u-k)_+^2 (\lambda_0 + \lambda_-) \, dx dt \, +						\\
	& + (\hat{r} - \tilde{r})^2 \, \frac{1}{\Lambda (Q_R^{\upbeta, \texttt{\,>}})}
	\sup_{t \in (0, s_2)} \int_{(I_{R/2}^+)^{\hat{r} - R/2}} (u-k)_+^2 (x,t) |\mu| (x) dx \, .
\end{align*}
Notice that for $h < k$ we have
\begin{align*}
(k-h)^2 & M_+ (A_R^{+,\tilde{r} - R/2}(k; {\textstyle \frac{R}{2}}, \tilde\theta)) \leqslant  
	\iint_{A_R^{+,\tilde{r} - R/2}(k; R/2, \tilde\theta)}  (u-h)_+^2 \mu_+ \, dx dt \leqslant							\\
&	\leqslant \,  \iint_{A_R^{+,\tilde{r} - R/2}(h; R/2 ,\tilde\theta)} (u-h)_+^2 \mu_+ \, dx dt \, ,
\end{align*}
that is
$$
M_+ (A_R^{+,\tilde{r} - R/2}(k; {\textstyle \frac{R}{2}}, \tilde\theta)) \leqslant
	\frac{M_+ (Q_{R}^{\upbeta, \texttt{\,>}})}{(k-h)^2} \, \, 
	\big( \tilde{u}_{\mu_+} (h; {\textstyle \frac{R}{2}},\tilde\theta ; \tilde{r} - {\textstyle \frac{R}{2}})  \big)^2 \, .
$$
From that (and the analogous estimate for $\Lambda_+ (A_R^{+,\tilde{r} - R/2}(k; R/2, \tilde\theta))$) we derive
\begin{align*}
\frac{M_+ (A_R^{+,\tilde{r} - R/2}(k; R/2, \tilde\theta))}{|M|_{\Lambda} (Q_{R}^{\upbeta, \texttt{\,>}})} \leqslant
\frac{M_+ (A_R^{+,\tilde{r} - R/2}(k; R/2, \tilde\theta))}{M_+ (Q_{R}^{\upbeta, \texttt{\,>}})} \leqslant 
	\frac{1}{(k-h)^2} \, \, \big( \tilde{u}_{\mu_+} (h; {\textstyle \frac{R}{2}}, \tilde\theta; \tilde{r} - {\textstyle \frac{R}{2}} )  \big)^2 \, ,				\\
\frac{\Lambda_+ (A_R^{+,\tilde{r} - R/2}(k; R/2, \tilde\theta))}{\Lambda (Q_{R}^{\upbeta, \texttt{\,>}})} \leqslant
\frac{\Lambda_+ (A_R^{+,\tilde{r} - R/2}(k; R/2, \tilde\theta))}{\Lambda_+ (Q_{R}^{\upbeta, \texttt{\,>}})} \leqslant 
	\frac{1}{(k-h)^2} \, \, \big( \tilde{u}_{\lambda_+} (h; {\textstyle \frac{R}{2}}, \tilde\theta; \tilde{r} - {\textstyle \frac{R}{2}})  \big)^2 \, .
\end{align*}
Then, applying these inequalities 
we get
\begin{align}
\label{daiterare}
\tilde{u}_+ (k; {\textstyle \frac{R}{2}}, \theta; r - {\textstyle \frac{R}{2}}) 
	\leqslant & \, \frac{C_1^{1/2}}{\hat{r} - \tilde{r}} \,  \frac{1}{(k-h)^{\frac{\kappa - 1}{\kappa}}} \, 
	\tilde{u}_+ (h; {\textstyle \frac{R}{2}}, \tilde\theta ; \tilde{r} - {\textstyle \frac{R}{2}})^{\frac{\kappa - 1}{\kappa}}\, 
		\left[ \tilde{u}_+ (k; {\textstyle \frac{R}{2}}, \hat\theta; \hat{r} - {\textstyle \frac{R}{2}}) + 
		\omega^{\hat{r} - \tilde{r}} (u; k; \hat{r} ; \hat\theta) \right]		\leqslant											\nonumber	\\
\leqslant & \, \frac{C_1^{1/2}}{\hat{r} - \tilde{r}} \,  \frac{1}{(k-h)^{\frac{\kappa - 1}{\kappa}}} \, 
	\tilde{u}_+ (h; {\textstyle \frac{R}{2}}, \tilde\theta ; \tilde{r} - {\textstyle \frac{R}{2}})^{\frac{\kappa - 1}{\kappa}}\, 
		\left[ \tilde{u}_+ (h; {\textstyle \frac{R}{2}}, \hat\theta; \hat{r} - {\textstyle \frac{R}{2}}) + 
		\omega^{\hat{r} - \tilde{r}} (u; h; \hat{r} ; \hat\theta) \right]		\leqslant											\nonumber	\\
\leqslant & \, \frac{C_1^{1/2}}{\hat{r} - \tilde{r}} \,  \frac{1}{(k-h)^{\frac{\kappa - 1}{\kappa}}} \, 
	\tilde{u}_+ (h; {\textstyle \frac{R}{2}}, \hat\theta; \hat{r} - {\textstyle \frac{R}{2}})^{\frac{\kappa - 1}{\kappa}} \, 
		\left[ \tilde{u}_+ (h; {\textstyle \frac{R}{2}}, \hat\theta; \hat{r} - {\textstyle \frac{R}{2}}) + 
		\omega^{\hat{r} - \tilde{r}} (u; h; \hat{r} ; \hat\theta) \right] .
\end{align}
Consider the following choices: for $n \in \N$, $k_0\in \R$ and a fixed $d$ we define
\begin{align*}
& k_n := k_0 + d \left( 1 -  \frac{1}{2^n} \right) \nearrow k_0 + d \, ,									\\
& r_n := \frac{R}{2} + \frac{R}{2^{n+1}}		\searrow \frac{R}{2}\, , 									\\
& \theta_n := \frac{1}{2} \left( 1 - \frac{1}{4^{n}} \right)\nearrow \frac{1}{2} \, ,							\\
& \sigma_n := \theta_n \, \upbeta \, h(x_0, R) \, R^2 \nearrow \frac{1}{2} \, \upbeta \, h(x_0, R) \, R^2 \, .
\end{align*}
Notice that (for these choices)
$$
2 \, (r_n - r_{n+1} ) = r_{n-1} - r_n \, . 
$$
With this choice of $\theta_n$ (and since $\upbeta \, h(x_0, R) R^2 = s_2 - t_0 = s_2$ since we are supposing $t_0 = 0$) we have that
$$
\sigma_n = \theta_n \, \upbeta \, h(x_0, R) \, R^2 = \theta_n \, s_2 \nearrow \frac{s_2}{2} \, .
$$
With this choices
we define the sequences
$$
u_n^+ := \tilde{u}_+ (k_n; {\textstyle \frac{R}{2}}, \theta_n; r_n - {\textstyle \frac{R}{2}}) 
	\, , \hskip20pt \omega_n^+ := \omega^{r_n - r_{n+1}} (u; k_n; r_n; \theta_n)
$$
and show that with the particular choices just made above the sequence $(u_n)_n$ is infinitesimal.
To get that it is sufficient to observe that from \eqref{daiterare} and using
\begin{align*}
& r_{n+1}  & \text{in the place of}  &  \hskip10pt	r	\, ,		 & \theta_{n+1} &  \hskip10pt	\text{in the place of}  & \theta	\, ,					\\
& r_{n}      & \text{in the place of}  &  \hskip10pt	\tilde{r}	\, , 	 & \theta_{n}     &  \hskip10pt	\text{in the place of}  &\tilde\theta \, ,				\\
& r_{n-1}   & \text{in the place of}  &  \hskip10pt	\hat{r}	\, ,	 & \theta_{n-1}  & \hskip10pt	\text{in the place of}   & \hat\theta \, ,				\\
& k_{n+1}  & \text{in the place of}  &  \hskip10pt	k		\, ,	 & k_{n-1}		  & \hskip10pt	\text{in the place of}   & h \, ,
\end{align*}
we derive
\begin{equation}
\label{enumeriamola}
u_{n+1}^+ \leqslant \, C_+ \, 2^{n+1} \,  \frac{2^{(n+1) \frac{\kappa - 1}{\kappa}}}{(3d)^{\frac{\kappa - 1}{\kappa}}} \, 
	\left( u_{n-1}^+ + \omega_{n-1}^+ \right)	(u_{n-1}^+)^{\frac{\kappa - 1}{\kappa}}	\, , \qquad n \geqslant 1 \, ,
\end{equation}
where $C_+ = \sqrt{C_1}/R = \gamma_1^{1/\kappa} \, (1 + \upbeta)^{1/2} \upbeta^{-\frac{1}{2\kappa}} \,  \, (2 \gamma + 8)^{1/2}$. Setting
$$
\alpha = \frac{\kappa - 1}{\kappa} \, , \hskip10pt
c = C_+ \, \frac{4^{1 + \alpha}}{3^{\alpha}d^{\alpha}} \,  , \hskip10pt 
b = 2^{1+\alpha} \,  , \hskip10pt 
y_n = u_n^+  \,  , \hskip10pt 
\epsilon_n = \omega_n^+ \,  .
$$
\eqref{enumeriamola} becomes
$$
u_{n+1}^+ \leqslant c \, b^{n-1} \left( u_{n-1}^+ + \omega_{n-1}^+ \right)	(u_{n-1}^+)^{\alpha} \, , \qquad n \geqslant 1 .
$$
In particular we get 
$$
u_{2(n+1)}^+ \leqslant c \, b^{2n} \left( u_{2n}^+ + \omega_{2n}^+ \right) (u_{2n}^+)^{\alpha} \, , \qquad n \geqslant 0 .
$$
Now notice that $(u_n^+)_n$ is decreasing. Then, using Lemma \ref{lemmuzzofurbo-quinquies}, provided that
\begin{equation}
\label{costanted}
u_0^+ < 
	\left( C_+ \, \frac{4^{1+\alpha}}{3^{\alpha}d^{\alpha}} \right)^{-1/\alpha} \, 2^{ -\frac{2}{\alpha} - \frac{2}{\alpha^2}} =
	3d\left( C_+ \right)^{-\frac{1}{\alpha}} \, 4^{ -\frac{2}{\alpha} - \frac{1}{\alpha^2} - 1} \, ,
\end{equation}
that is
\begin{align*}
& \left( \frac{1}{|M|_{\Lambda} (Q_{R}^{\upbeta, \texttt{\,>}})} 
				\int\!\!\!\int_{Q_{R; R/2, 0}^{\upbeta, +,R/2}} (u-k_0)_+^2 \mu_+ \, dx dt +
\frac{1}{\Lambda (Q_{R}^{\upbeta, \texttt{\,>}})} 
				\int\!\!\!\int_{Q_{R; R/2, 0}^{\upbeta, +,R/2}} (u-k_0)_+^2 \lambda_+ \, dx dt \right)^{1/2} <					\\
& \hskip200pt < \, 3d\left( C_+ \right)^{-\frac{1}{\alpha}} \, 4^{ -\frac{2}{\alpha} - \frac{1}{\alpha^2} - 1} 		 \, ,
\end{align*}
we get that the subsequence $(u_{2n})_n$ is infinitesimal and since $(u_n)_n$ is decreasing we finally derive
\begin{equation}
\label{limitezero!!!}
\lim_{n \to +\infty} u_n^+ = \tilde{u}_+ \left(k_0 + d; \frac{R}{2}, \frac{1}{2} \right) = 0
\end{equation}
where
\begin{align*}
\big( \tilde{u}_+ (l; \varrho , \vartheta) \big)^2 := & \, \big( \tilde{u}_+ (l; \varrho ,\vartheta ; 0) \big)^2 =							\\
= & \, \frac{1}{|M|_{\Lambda} (Q_R^{\upbeta, \texttt{\,>}})} 
		\int\!\!\!\int_{Q_{R; \varrho, \vartheta}^{\upbeta, +}} (u-l)_+^2 \mu_+ \, dx dt
		+ \frac{1}{\Lambda (Q_R^{\upbeta, \texttt{\,>}})} 
				\int\!\!\!\int_{Q_{R; \varrho, \vartheta}^{\upbeta, +}} (u-l)_+^2 \lambda_+ \, dx dt \, .
\end{align*}
In a complete analogous way, if $\mu_- (B_R) > 0$ and taking $s_1 = t_0 - \upbeta \, h(x_0, R) R^2$, one can prove that
\begin{equation}
\label{limitezero!!!-}
\begin{array}{l}
{\displaystyle \int\!\!\!\int_{Q_{R; R/2, 1/2}^{\upbeta, -} (x_0,t_0)} (u-k_0-d)_+^2 \mu_- \, dx dt = 0 \, ,		}				\\	[1em]
{\displaystyle \int\!\!\!\int_{Q_{R; R/2, 1/2}^{\upbeta, -}(x_0,t_0)} (u-k_0-d)_+^2 \lambda_- \, dx dt = 0 \, ,	}
\end{array}
\end{equation}
where $Q_{R; R/2, 1/2}^{\upbeta, -} (x_0,t_0) = B_R^- (x_0) \times (t_0 - \upbeta \, h(x_0, R) R^2, t_0 - \frac{1}{2} \upbeta \, h(x_0, R) R^2)$, provided that
\begin{align*}
& \left( \frac{1}{|M|_{\Lambda} (Q_{R}^{\upbeta, \texttt{\,<}})} 
				\int\!\!\!\int_{Q_{R; R, 0}^{\upbeta, -,R/2}} (u-k_0)_+^2 \mu_- \, dx dt +
\frac{1}{\Lambda (Q_{R}^{\upbeta, \texttt{\,<}})} 
				\int\!\!\!\int_{Q_{R; R, 0}^{\upbeta, -,R/2}} (u-k_0)_+^2 \lambda_- \, dx dt \right)^{1/2} <					\\
& \hskip200pt < \, 3d\left( C_- \right)^{-\frac{1}{\alpha}} \, 4^{ -\frac{2}{\alpha} - \frac{1}{\alpha^2} - 1} 		 \, ,
\end{align*}
where $C_- = C_+ = \gamma_1^{1/\kappa} \, \upbeta^{1/2} \, (2 \gamma + 8)^{1/2}$. \\ [0.5em]
The proof regarding the part in which $\mu \equiv 0$ is slightly different and we show it.
We define
$$
\sigma_1 := t_0 - \frac{R^2}{2}, \quad \sigma_2 := t_0 + \frac{R^2}{2}
\qquad \text{ so that } \ \sigma_2 - \sigma_1 =  R^2 .
$$
Moreover we suppose that
$$
\lambda_0 (B_{R}) > 0 ,
$$
otherwise there is nothing to prove. We consider $r , \tilde{r} , \hat{r} \in (R/2,R)$ as before.
Consider a function $\zeta \in \text{Lip}_c (B_{\tilde{r}} (x_0))$ (independent of $t$!) such that
$$
\begin{array}{c}
\zeta \equiv 1 \hskip10pt \text{in } Q_{R; \frac{R}{2}; \sigma_1, \sigma_2}^{0, r - \frac{R}{2}} (x_0) \, ,	\hskip20pt
	\zeta \equiv 0 \hskip10pt \text{in } 
	\big( B_R (x_0) \times (\sigma_1, \sigma_2) \big) \setminus Q_{R; \frac{R}{2}; \sigma_1, \sigma_2}^{0,\tilde{r} - \frac{R}{2}} (x_0) \, , 		\\ [1em]
0 \leqslant \zeta \leqslant 1 \, , \hskip15pt {\displaystyle |D \zeta| \leqslant \frac{1}{\tilde{r}-r} } \, .
\end{array}
$$
We moreover define
\begin{gather*}
A^{0,\delta}_{R} (k;\rho;\sigma_1, \sigma_2) := \{ (x,t) \in Q_{R;\rho; \sigma_1, \sigma_2}^{0, \delta} (x_0) \, | \, u(x,t) > k \} .
\end{gather*}
Then we proceed in a way similar to that above and estimate 
$(\lambda (B_R))^{-1} \int\!\!\!\int_{Q_{R; R/2; \sigma_1, \sigma_2}^{0,r - R/2}}  (u-k)_+^2 \lambda \, dx dt$
using first Corollary \ref{cor-gut-whee} with $\nu = \upsilon = \omega = \lambda$. One has 
(we write $Q_{R;\rho; s_1, s_2}^{0,\varepsilon}$ to mean $Q_{R;\rho; s_1, s_2}^{0,\varepsilon} (x_0)$)
\begin{align*}
& \frac{1}{\lambda (B_R)} \iint_{Q_{R; R/2; \sigma_1, \sigma_2}^{0,r-R/2}}  (u-k)_+^2 \lambda_0 \, dx dt \leqslant	
	\frac{1}{\lambda (B_R)} \iint_{Q_{R; R/2 ; \sigma_1, \sigma_2}^{0, \tilde{r} - R/2}}  (u-k)_+^2 \zeta^2 \lambda_0 \, dx dt \leqslant				\\
& \hskip30pt \leqslant  \frac{\big( \Lambda_0 (A^{0,\tilde{r} - R/2}_{R} (k; R/2 ;\sigma_1, \sigma_2) ) \big)^{\frac{\kappa - 1}{\kappa}}}
	{(\lambda (B_R))^{\frac{\kappa - 1}{\kappa}}} 
		\left[\frac{1}{\lambda (B_{R})} \iint_{Q_{R; R/2 ; \sigma_1, \sigma_2}^{0, \tilde{r} - R/2}}
	(u-k)_+^{2\kappa} \zeta^{2\kappa} \lambda_0 \, dx dt \right]^{\frac{1}{\kappa}} \leqslant												\\
& \hskip30pt \leqslant  \frac{\big( \Lambda_0 (A^{0,\tilde{r} - R/2}_{R} (k; R/2;\sigma_1, \sigma_2) ) \big)^{\frac{\kappa - 1}{\kappa}}}
	{(\lambda (B_R))^{\frac{\kappa - 1}{\kappa}}} \, \gamma_1^{2/\kappa} \, R^{2/\kappa} \frac{1}{\lambda (B_R)} \cdot							\\
& \hskip40pt \cdot \left[ \sup_{t \in (\sigma_1, \sigma_2)} 
		\int_{(B_{R/2}^0)^{\tilde{r} - R/2}} (u-k)_+^2 (x,t) \lambda_0 (x) dx \right]^{\frac{\kappa - 1}{\kappa}}
		\left[\iint_{Q_{R; R/2 ; \sigma_1, \sigma_2}^{0, \tilde{r} - R/2}} | D \big( (u-k)_+ \zeta \big) |^2 \lambda \, dx dt
											\right]^{\frac{1}{\kappa}} .
\end{align*}
Then using the energy estimate \eqref{DGgamma0} we get
\begin{align*}
\iint_{Q_{R; R/2 ; \sigma_1, \sigma_2}^{0, \tilde{r} - R/2}} & | D \big( (u-k)_+ \zeta \big) |^2 \lambda \, dx dt	\leqslant							\\
& \leqslant 2 \iint_{Q_{R; R/2 ; \sigma_1, \sigma_2}^{0, \tilde{r} - R/2}}
	\Big[ | D (u-k)_+ |^2 \zeta^2 + | D \zeta |^2 (u-k)_+^2 \Big] \lambda \, dx dt	\leqslant												\\
& \leqslant 2 \iint_{Q_{R; R/2 ; \sigma_1, \sigma_2}^{0, \tilde{r} - R/2}}
	\Big[ | D (u-k)_+ |^2 + \frac{1}{(\tilde{r} - r)^2} (u-k)_+^2 \Big] \lambda \, dx dt	\leqslant											\\
& \leqslant 2 \gamma \Bigg[ \sup_{t \in (\sigma_1, \sigma_2)} \int_{ I_0^{R/2, \hat{r} - R/2}} (u-k)_+^2(x,t) \mu_- (x) \, dx +					\\
& \hskip50pt + \sup_{t \in (\sigma_1, \sigma_2)}\int_{ I_0^{R/2, \hat{r} - R/2}} (u-k)_+^2(x,t) \mu_+ (x) \, dx \, +							\\
& \hskip50pt + \frac{1}{(\hat{r} - \tilde{r})^2} 
				\iint_{Q_{R; R/2 ; \sigma_1, \sigma_2}^{0, \hat{r} - R/2}} (u-k)_+^2 \, \lambda \, dx dt	\Bigg] +							\\
& \hskip80pt + \frac{2}{(\tilde{r} - r)^2} \iint_{Q_{R; R/2 ; \sigma_1, \sigma_2}^{0, \tilde{r} - R/2}} (u-k)_+^2 \lambda \, dx dt \, .
\end{align*}
Then we have, dividing by $\sigma_2 - \sigma_1$ in both sides,
\begin{align}
\label{mitoccanum3}
\frac{1}{(\sigma_2 - \sigma_1) \lambda (B_R)} & 
					\iint_{Q_{R; R/2 ; \sigma_1, \sigma_2}^{0, {r} - R/2}}  (u-k)_+^2 \lambda_0 \, dx dt	\leqslant				\nonumber	\\
& \hskip30pt \leqslant  \frac{\big( \Lambda_0 (A^{0,\tilde{r} - R/2}_{R} (k; R/2 ;\sigma_1, \sigma_2) ) \big)^{\frac{\kappa - 1}{\kappa}}}
	{(\sigma_2 - \sigma_1)^{\frac{\kappa - 1}{\kappa}}(\lambda (B_R))^{\frac{\kappa - 1}{\kappa}}} \, 
	\gamma_1^{2/\kappa} \, \frac{R^{2/\kappa}}{(\sigma_2 - \sigma_1)^{\frac{1}{\kappa}}}
					\frac{(\sigma_2 - \sigma_1)}{(\sigma_2 - \sigma_1)\lambda (B_R)} 					\cdot				\nonumber	\\
& \hskip50pt \cdot \Bigg[ \sup_{t \in (\sigma_1, \sigma_2)} 
		\int_{(B_{R/2}^0)^{\tilde{r} - R/2}} (u-k)_+^2 (x,t) \lambda_0 (x) dx +												\nonumber	\\
& \hskip70pt + 2 \gamma \sup_{t \in (\sigma_1, \sigma_2)} \int_{ I_0^{R/2, \hat{r} - R/2}} (u-k)_+^2(x,t) \mu_- (x) \, dx +							\\
& \hskip85pt + 2 \gamma \sup_{t \in (\sigma_1, \sigma_2)}\int_{ I_0^{R/2, \hat{r} - R/2}} (u-k)_+^2(x,t) \mu_+ (x) \, dx \, +			\nonumber	\\
& \hskip100pt + \frac{2 \gamma + 8}{(\hat{r} - \tilde{r})^2}
	\iint_{Q_{R; R/2 ; \sigma_1, \sigma_2}^{0, {\hat r} - R/2}} (u-k)_+^2 (\lambda_+ + \lambda_-) \, dx dt \, +						\nonumber	\\
& \hskip115pt + \frac{2 \gamma + 8}{(\hat{r} - \tilde{r})^2}
	\iint_{Q_{R; R/2 ; \sigma_1, \sigma_2}^{0, {\hat r} - R/2}} (u-k)_+^2 \lambda_0 \, dx dt \Bigg]	\, .							\nonumber
\end{align}
Now defining
$$
\big( \tilde{u}_0(l; \rho ; \varepsilon ; \sigma_1, \sigma_2) \big)^2 =
	\frac{1}{(\sigma_2 - \sigma_1) \lambda (B_R)}\iint_{Q_{R; \rho ; \sigma_1, \sigma_2}^{0, \varepsilon}}  (u - l)_+^2 \lambda_0 \, dx dt
$$
for $\varepsilon \in [0, R/2)$,
\begin{align*}
\big( \omega^{\hat{r} - \tilde{r}} (u; k; \hat{r}) \big)^2:= & 
	\, \frac{(\hat{r} - \tilde{r})^2}{(\sigma_2 - \sigma_1)\lambda (B_R)} \cdot \Bigg[ \sup_{t \in (\sigma_1, \sigma_2)} 
		\int_{(B_{R/2}^0)^{\hat{r} - R/2}} (u-k)_+^2 (x,t) \lambda_0 (x) dx +											\\
& \hskip65pt + \sup_{t \in (\sigma_1, \sigma_2)} \int_{ I_0^{R/2, \hat{r} - R/2}} (u-k)_+^2(x,t) \mu_- (x) \, dx +					\\
& \hskip85pt + \sup_{t \in (\sigma_1, \sigma_2)}\int_{ I_0^{R/2, \hat{r} - R/2}} (u-k)_+^2(x,t) \mu_+ (x) \, dx \, \Bigg] +			\\
& \hskip50pt +	\frac{1}{(\sigma_2 - \sigma_1)\lambda (B_R)}
	\iint_{Q_{R; R/2 ; \sigma_1, \sigma_2}^{0, {\hat r} - R/2}} (u-k)_+^2 (\lambda_+ + \lambda_-) \, dx dt
\end{align*}
and for $k > h$
$$
\frac{\Lambda_0 (A^{0,\tilde{r} - R/2}_{R} (k; R/2 ;\sigma_1, \sigma_2) )}{(\sigma_2 - \sigma_1)\lambda (B_R)} \leqslant
	\frac{1}{(k-h)^2} \big( \tilde{u}_0(h; {\textstyle \frac{R}{2}}; \tilde{r} - {\textstyle \frac{R}{2}} ;\sigma_1, \sigma_2) \big)^2
$$
and since $\sigma_2 - \sigma_1 =  \upbeta \, R^2$ we reach
\begin{align}
\tilde{u}_0(k; {\textstyle \frac{R}{2}}; r - {\textstyle \frac{R}{2}} ;\sigma_1, \sigma_2) & \leqslant
\frac{\gamma_1^{1/\kappa} \, \upbeta^{\frac{\kappa - 1}{2\kappa}} R}{(k - h)^{\frac{\kappa - 1}{\kappa}}} 
	\frac{(2 \gamma + 8)^{1/2}}{\hat{r} - \tilde{r}} \cdot																	\nonumber	\\
& \cdot \Big[ \omega^{\hat{r} - \tilde{r}} (u; k; \hat{r}) + 
	\tilde{u}_0 (k; {\textstyle \frac{R}{2}}; \hat{r} - {\textstyle \frac{R}{2}};\sigma_1, \sigma_2) \Big]
	\big( \tilde{u}_0 (h; {\textstyle \frac{R}{2}}; \tilde{r} - {\textstyle \frac{R}{2}} ; \sigma_1, \sigma_2) \big)^{\frac{\kappa - 1}{\kappa}} \leqslant \nonumber	\\
\label{noncasca}
& \leqslant
\frac{\gamma_1^{1/\kappa} \, \upbeta^{\frac{\kappa - 1}{2\kappa}} R}{(k - h)^{\frac{\kappa - 1}{\kappa}}} \frac{(2 \gamma + 8)^{1/2}}{\hat{r} - \tilde{r}} \cdot	\\
& \cdot \Big[ \omega^{\hat{r} - \tilde{r}} (u; k; \hat{r}) + \tilde{u}_0 (h; {\textstyle \frac{R}{2}}; \hat{r} - {\textstyle \frac{R}{2}}; \sigma_1, \sigma_2) \Big]
		\big( \tilde{u}_0 (h; {\textstyle \frac{R}{2}}; \hat{r} - {\textstyle \frac{R}{2}};\sigma_1, \sigma_2) \big)^{\frac{\kappa - 1}{\kappa}} .	\nonumber
\end{align}
As done before, consider the following choices for $n \in \N$, $k_0\in \R$ and a fixed $d$:
\begin{align*}
& k_n := k_0 + d \left( 1 -  \frac{1}{2^n} \right) \nearrow k_0 + d \, , \qquad 	
r_n := \frac{R}{2} + \frac{R}{2^{n+1}}		\searrow \frac{R}{2} \, ,
\end{align*}
and define the sequences
$$
u_n^0 := \tilde{u}_0 (k_n; {\textstyle \frac{R}{2}}; r_n - {\textstyle \frac{R}{2}}; \sigma_1, \sigma_2) \, , 
	\hskip20pt \omega_n^0 := \omega^{r_n - r_{n+1}} (u; k_n; r_n) \, .
$$
Making the following choices in \eqref{noncasca}
\begin{align*}
& r_{n+1}  & \text{in the place of}  &  \hskip10pt	r	\, ,		 & r_{n}      	& \hskip10pt \text{in the place of }  \ \ 	\tilde{r}	\, , 	 	\\
& r_{n-1}   & \text{in the place of}  &  \hskip10pt	\hat{r}	\, ,	 & \, 			&			\, 					 			\\
& k_{n+1}  & \text{in the place of}  &  \hskip10pt	k		\, ,	 & k_{n-1}		& \hskip10pt \text{in the place of } \ \ h \, ,
\end{align*}
we get
%
%
\begin{equation}
\label{enumeriamola_bis}
u_{n+1}^0 \leqslant \, 
\frac{\gamma_1^{1/\kappa}  \upbeta^{\frac{\kappa - 1}{2\kappa}} (2 \gamma + 8)^{1/2}}{(3d)^{\frac{\kappa - 1}{\kappa}}} \, 
(2^{\frac{2\kappa - 1}{\kappa}})^{n-1}
\left( u_{n-1}^0 + \omega_{n-1}^0 \right)	(u_{n-1}^0)^{\frac{\kappa - 1}{\kappa}}	\, , \qquad n \geqslant 1
\end{equation}
and then, similarly as before, we derive that 
\begin{align*}
\lim_{n \to +\infty} u_n^0 = & \, \tilde{u}_0 \left(k_0 + d; \frac{R}{2}; \sigma_1, \sigma_2 \right) :=
	\tilde{u}_0 \left(k_0 + d; \frac{R}{2}; 0; \sigma_1, \sigma_2 \right) = 													\\
= & \, \left( \frac{1}{(\sigma_2 - \sigma_1)\lambda (B_R)} 
	\int_{\sigma_1}^{\sigma_2} \!\!\!\! \int_{B_{R/2}^0}  (u - k_0 - d)_+^2 \lambda_0 \, dx dt \right)^{1/2} = 0
\end{align*}
provided that
$$
u_0^0 < 3 \, d \, \gamma_1^{- \frac{1}{\kappa - 1}} \, \upbeta^{-1/2} (2 \gamma + 8)^{- \frac{\kappa}{2(\kappa - 1)}}
																\, 2^{- \frac{6 \kappa^2 - 7 \kappa + 2}{(\kappa - 1)^2}} \, .
$$
\ \\
\noindent
Now we continue and conclude this section showing that $u$ is locally bounded in $\Omega$. 
In Figure B we show, supposing $\mu > 0$ and $\mu < 0$,  the sets involved in the stimates of points $i \,)$ and $ii \, )$.
%

\begin{theorem}
\label{Linfinity}
Suppose $u \in DG(\Omega,T, \mu, \lambda, \gamma)$ and consider $(x_0, t_0) \in \Omega \times (0,T)$, $\upbeta > 0$. 
Then there is a constant $c_\infty$ depending only on $\gamma, \gamma_1, \kappa , \upbeta$ 
such that: \\ [0.5em]
$i \, )$ for every $B_{R} (x_0) \times (t_0, t_0 + \upbeta \, h(x_0,R) R^2) \subset \Omega \times (0,T)$ if $\mu_+ (B_R (x_0)) > 0$ we have
\begin{align*}
\esssup_{Q_{R; R/2, 1/2}^{\upbeta,+}} & |u| \leqslant 
	c_{\infty}  \Bigg[ \frac{1}{|M|_{\Lambda} (Q_{R}^{\upbeta, \texttt{\,>}})} \iint_{Q_{R; R, 0}^{\upbeta,+,R/2}} u^2 \mu_+ \, dx dt +
	\frac{1}{\Lambda (Q_{R}^{\upbeta, \texttt{\,>}})} \iint_{Q_{R; R, 0}^{\upbeta,+,R/2}} u^2 \lambda_+ \, dx dt \Bigg]^{1/2} ;
\end{align*}
$ii \, )$ for every $B_{R} (x_0) \times (t_0 - \upbeta \, h(x_0,R) R^2, t_0) \subset \Omega \times (0,T)$ if $\mu_- (B_R (x_0)) > 0$ we have
\begin{align*}
\esssup_{Q_{R; R/2, 1/2}^{\upbeta,-}} & |u| \leqslant 
	c_{\infty} \Bigg[ \frac{1}{|M|_{\Lambda} (Q_{R}^{\upbeta, \texttt{\,<}})} \iint_{Q_{R; R, 0}^{\upbeta,-,R/2}} u^2 \mu_- \, dx dt +
	\frac{1}{\Lambda (Q_{R}^{\upbeta, \texttt{\,<}})} \iint_{Q_{R; R, 0}^{\upbeta,-,R/2}} u^2 \lambda_- \, dx dt \Bigg]^{1/2} ;
\end{align*}
$iii \, )$ for every $B_{R} (x_0) \times (\sigma_1, \sigma_2) \subset \Omega \times (0,T)$, $\sigma_2 - \sigma_1 = R^2$,
if $\lambda_0 (B_R (x_0)) > 0$
\begin{align*}
\esssup_{B_{R/2}^0 \times (\sigma_1, \sigma_2)} & |u| \leqslant c_{\infty}
	\left( \frac{1}{\Lambda (B_{R} \times (\sigma_1, \sigma_2))} \iint_{Q^{0,R/2}_{R;R;\sigma_1,\sigma_2}}   u^2 \lambda_0 \, dx dt \right)^{1/2} .
\end{align*}
\end{theorem}
\noindent
\dimo
We prove the first point, being the others very similar. By \eqref{limitezero!!!} we derive that
$$
\esssup_{Q_{R; R/2, 1/2}^+} u \leqslant k_0 + d
$$
and $d$ has to satisfy \eqref{costanted}. For example we can choose
$$
d = 2 \left( C_+ \right)^{\frac{1}{\alpha}} \, 3^{-1} 4^{ \frac{2}{\alpha} + \frac{1}{\alpha^2} + 1} \, u_0^+ .
$$
By definition of $u_0^+$, defining the quantity
$$
c_{\infty} := \frac{d}{u_0^+} = \frac{2}{3} \left( C_+ \right)^{\frac{1}{\alpha}} \, 4^{ \frac{2}{\alpha} + \frac{1}{\alpha^2} + 1} =
\frac{2}{3} \gamma_1^{\frac{1}{\kappa - 1}} \frac{(1+\upbeta)^{\frac{\kappa}{2(\kappa - 1)}}}{\upbeta^{\frac{1}{2(\kappa - 1)}}} \, 
(2 \gamma + 8)^{\frac{\kappa}{2(\kappa - 1)}} 4^{\frac{3 \kappa^2 - 3 \kappa + 1}{(\kappa - 1)^2}} \, ,
$$
choosing $k_0 = 0$ and estimating $u_+^2$ by $u^2$ we finally get
\begin{align*}
\esssup_{Q_{R; R/2, 1/2}^{\upbeta,+}} & u \leqslant 
	c_{\infty}  \left( \frac{1}{|M|_{\Lambda} (Q_{R}^{\upbeta, \texttt{\,>}})} \iint_{Q_{R; R, 0}^{\upbeta,+,R/2}} u^2 \mu_+ \, dx dt +
	\frac{1}{\Lambda (Q_{R}^{\upbeta, \texttt{\,>}})} \iint_{Q_{R; R, 0}^{\upbeta,+,R/2}} u^2 \lambda_+ \, dx dt \right)^{1/2} \, .
\end{align*}
Since the analogous argument can be applied to $-u$ we have the first claim. The points $ii \, )$ and $iii \, )$ are completely analogous:
the only difference is that the constant $c_{\infty}$ in point $ii \, )$ is the same as in point $i \, )$,
in point $iii \, )$ is
$3^{-1}  \, \gamma_1^{\frac{1}{\kappa - 1}} \, \upbeta^{1/2} (8 \gamma + 2)^{\frac{\kappa}{2(\kappa - 1)}} \, 2^{\frac{6 \kappa^2 - 7 \kappa + 2}{(\kappa - 1)^2}-1}$.
\finedimo

\boss
Notice that from points $i \, )$ and $ii \, )$ it is not possible to derive a pointwise (in time) estimate: indeed letting $\upbeta$ go to zero
the constant $c_{\infty}$ goes to $+\infty$. \\
Also in point $iii \, )$ we cannot obtain a pointwise estimate because $\sigma_2 - \sigma_1 = \upbeta R^2$ and the constant
$c_{\infty}$ depends on $\upbeta$. \\
Nevertheless one could obtain a pointwise estimate if $B_R \subset \Omega_0$ using \eqref{tempofissato}
and Theorem \ref{gut-whee}.
\eoss
\ \\
\noindent
The local boundedness for a function in the class $DG$ is immediatly needed in the following section.

\ \\
\ \\
\begin{picture}(150,200)(-180,0)
\put (-105,200){\linethickness{1pt}\line(1,0){200}}
\put (-105,50){\linethickness{1pt}\line(1,0){200}}
\put (-105,50){\linethickness{1pt}\line(0,1){150}}
\put (95,50){\linethickness{1pt}\line(0,1){150}}

\put (-5,40){\line(0,1){170}}


\put (-104,125){\line(1,0){1}}
\put (-101,125){\line(1,0){1}}
\put (-98,125){\line(1,0){1}}
\put (-95,125){\line(1,0){1}}
\put (-92,125){\line(1,0){1}}
\put (-89,125){\line(1,0){1}}
\put (-86,125){\line(1,0){1}}
\put (-83,125){\line(1,0){1}}
\put (-80,125){\line(1,0){1}}
\put (-77,125){\line(1,0){1}}
\put (-74,125){\line(1,0){1}}
\put (-71,125){\line(1,0){1}}
\put (-68,125){\line(1,0){1}}
\put (-65,125){\line(1,0){1}}
\put (-62,125){\line(1,0){1}}
\put (-59,125){\line(1,0){1}}
\put (-56,125){\line(1,0){1}}

\put (-55,125){\line(0, -1){1}}
\put (-55,122){\line(0, -1){1}}
\put (-55,119){\line(0, -1){1}}
\put (-55,116){\line(0, -1){1}}
\put (-55,113){\line(0, -1){1}}
\put (-55,110){\line(0, -1){1}}
\put (-55,107){\line(0, -1){1}}
\put (-55,104){\line(0, -1){1}}
\put (-55,101){\line(0, -1){1}}
\put (-55,98){\line(0, -1){1}}
\put (-55,95){\line(0, -1){1}}
\put (-55,92){\line(0, -1){1}}
\put (-55,89){\line(0, -1){1}}
\put (-55,86){\line(0, -1){1}}
\put (-55,83){\line(0, -1){1}}
\put (-55,80){\line(0, -1){1}}
\put (-55,77){\line(0, -1){1}}
\put (-55,74){\line(0, -1){1}}
\put (-55,71){\line(0, -1){1}}
\put (-55,68){\line(0, -1){1}}
\put (-55,65){\line(0, -1){1}}
\put (-55,62){\line(0, -1){1}}
\put (-55,57){\line(0, -1){1}}
\put (-55,54){\line(0, -1){1}}
\put (-55,51){\line(0, -1){1}}

\put (-55,200){\line(0,-1){4}}
\put (-55,195){\line(0,-1){4}}
\put (-55,190){\line(0,-1){4}}
\put (-55,185){\line(0,-1){4}}
\put (-55,180){\line(0,-1){4}}
\put (-55,175){\line(0,-1){4}}
\put (-55,170){\line(0,-1){4}}
\put (-55,165){\line(0,-1){2}}
\put (-55,163){\line(1,0){4}}
\put (-50,163){\line(1,0){4}}
\put (-45,163){\line(1,0){4}}
\put (-40,163){\line(1,0){4}}
\put (-35,163){\line(1,0){4}}
\put (-30,163){\line(1,0){4}}
\put (-25,163){\line(1,0){4}}
\put (-20,163){\line(1,0){4}}
\put (-15,163){\line(1,0){4}}
\put (-10,163){\line(1,0){4}}

\put (-55,187){\line(1,1){12}}
\put (-55,177){\line(1,1){22}}
\put (-55,167){\line(1,1){32}}
\put (-49,163){\line(1,1){37}}
\put (-39,163){\line(1,1){33}}
\put (-29,163){\line(1,1){23}}
\put (-19,163){\line(1,1){13}}


\put (45,199){\line(0,-1){1}}
\put (45,196){\line(0,-1){1}}
\put (45,193){\line(0,-1){1}}
\put (45,190){\line(0,-1){1}}
\put (45,187){\line(0,-1){1}}
\put (45,184){\line(0,-1){1}}
\put (45,181){\line(0,-1){1}}
\put (45,178){\line(0,-1){1}}
\put (45,175){\line(0,-1){1}}
\put (45,172){\line(0,-1){1}}
\put (45,169){\line(0,-1){1}}
\put (45,166){\line(0,-1){1}}
\put (45,163){\line(0,-1){1}}
\put (45,160){\line(0,-1){1}}
\put (45,157){\line(0,-1){1}}
\put (45,154){\line(0,-1){1}}
\put (45,151){\line(0,-1){1}}
\put (45,148){\line(0,-1){1}}
\put (45,145){\line(0,-1){1}}
\put (45,142){\line(0,-1){1}}
\put (45,139){\line(0,-1){1}}
\put (45,136){\line(0,-1){1}}
\put (45,133){\line(0,-1){1}}
\put (45,130){\line(0,-1){1}}
\put (45,127){\line(0,-1){1}}

\put (45,125){\line(1,0){1}}
\put (48,125){\line(1,0){1}}
\put (51,125){\line(1,0){1}}
\put (54,125){\line(1,0){1}}
\put (57,125){\line(1,0){1}}
\put (60,125){\line(1,0){1}}
\put (63,125){\line(1,0){1}}
\put (66,125){\line(1,0){1}}
\put (69,125){\line(1,0){1}}
\put (72,125){\line(1,0){1}}
\put (75,125){\line(1,0){1}}
\put (78,125){\line(1,0){1}}
\put (81,125){\line(1,0){1}}
\put (84,125){\line(1,0){1}}
\put (87,125){\line(1,0){1}}
\put (90,125){\line(1,0){1}}
\put (93,125){\line(1,0){1}}

\put (45,50){\line(0,1){4}}
\put (45,55){\line(0,1){4}}
\put (45,60){\line(0,1){4}}
\put (45,65){\line(0,1){4}}
\put (45,70){\line(0,1){4}}
\put (45,75){\line(0,1){4}}
\put (45,80){\line(0,1){4}}
\put (45,85){\line(0,1){2}}
\put (45,87){\line(-1,0){4}}
\put (40,87){\line(-1,0){4}}
\put (35,87){\line(-1,0){4}}
\put (30,87){\line(-1,0){4}}
\put (25,87){\line(-1,0){4}}
\put (20,87){\line(-1,0){4}}
\put (15,87){\line(-1,0){4}}
\put (10,87){\line(-1,0){4}}
\put (5,87){\line(-1,0){4}}
\put (0,87){\line(-1,0){4}}

\put (-5,77){\line(1,1){10}}
\put (-5,67){\line(1,1){20}}
\put (-5,57){\line(1,1){30}}
\put (-2,50){\line(1,1){37}}
\put (8,50){\line(1,1){37}}
\put (18,50){\line(1,1){27}}
\put (28,50){\line(1,1){17}}

\put (-90,140){\tiny$\mu > 0$}
\put (60,110){\tiny$\mu < 0$}

\put (-5,125){\linethickness{2pt}\line(1,0){1}}
\put (-10,130){\tiny$(x_0,t_0)$}

\put (-30,10){Figure B}

\end{picture}

\section{Expansion of positivity}
\label{secPositivity}

In this section we will see many preliminary results needed to prove Harnack's inequality. \\
In what follows we fix the following points and sets: given 
three points $(x^{\diamond}\!, t^{\diamond}), (x^{\circ}\!, t^{\circ}), (x^{\star}\!, t^{\star}) \in \Omega \times (0,T)$ in such a way that
\begin{align*}
& Q_{R}^{\upbeta^{\diamond}, \texttt{\,>}} (x^{\diamond}\!, t^{\diamond}) = 
	B_R(x^{\diamond}) \times (t^{\diamond}, s_2^{\diamond}) \subset \Omega \times (0,T)
	\quad & \text{where } s_2^{\diamond} = t^{\diamond} + \upbeta^{\diamond} \, h \! \left(x^{\diamond} \!, R \right) R^2 \, ,						\\
& Q_{R}^{\upbeta^{\circ}, \texttt{\,<}} (x^{\circ}\!, t^{\circ}) = 
	B_R(x^{\ast}) \times (s_1^{\circ}, t^{\circ}) \subset \Omega \times (0,T)
	\quad & \text{where } s_1^{\circ} = t^{\ast} - \upbeta^{\circ} \, h \! \left(x^{\ast} \!, R \right) R^2 \, ,										\\
& Q_{R}^{s_1^{\star} , s_2^{\star}} (x^{\star}\!, t^{\star}) :=
	B_R(x^{\star}) \times (s_1^{\star} , s_2^{\star}) \subset \Omega \times (0,T)
\quad & \text{where } s_1^{\star} = t^{\star} - \frac{\upbeta^{\star}}{2} R^2 , \ s_2^{\star} = t^{\star} + \frac{\upbeta^{\star}}{2} R^2 \, ,
\end{align*}
with $\upbeta^{\diamond}, \upbeta^{\circ} , \upbeta^{\star} > 0$. \\ [0.3em]
\noindent
We recall that, thanks to the results of the previous section, a function belonging to the De Giorgi class $DG$
is locally bounded.

\bprop
\label{prop-DeGiorgi1}
Consider three points $(x^{\diamond}\!, t^{\diamond}), (x^{\circ}\!, t^{\circ}), (x^{\star}\!, t^{\star}) \in \Omega \times (0,T)$ and
$\rho \in (0,R)$. Suppose
$Q_{R}^{\upbeta^{\diamond}, \texttt{\,>}} (x^{\diamond}\!, t^{\diamond})$,
$Q_{R}^{\upbeta^{\circ}, \texttt{\,<}} (x^{\circ}\!, t^{\circ})$,
$Q_{R}^{s_1^{\star} , s_2^{\star}} (x^{\star}\!, t^{\star})$
are contained in $\Omega \times (0,T)$. Then for every choice of $\theta^{\diamond} , \theta^{\circ} \in (0,1)$ and $a, \sigma \in (0,1)$ there are \\
$\overline{\nu}^{\diamond} \in (0,1)$, depending only on $\kappa, \gamma_1, \gamma, a$, $\theta^{\diamond}$, $\upbeta^{\diamond}$,\\
$\overline{\nu}^{\circ} \in (0,1)$, depending only on $\kappa, \gamma_1, \gamma, a$, $\theta^{\circ}$, $\upbeta^{\circ}$, \\
$\overline{\nu}^{\star} \in (0,1)$, depending only on $\kappa, \gamma_1, \gamma, a, (R - \rho)/R$, $\max\{ 1, 1/\upbeta^{\star} \}$, \\
$\overline{\nu} \in (0,1)$, depending only on $\kappa, \gamma_1, \gamma, a, (R - \rho)/R$, \\
such that for every $u \in DG_+(\Omega, T, \mu, \lambda, \gamma)$ and fixed $\overline{m}, \omega$ satisfying \\ [0.5em]
$i \, )$ $\overline{m} \geqslant \sup_{Q_{R; R, 0}^{\upbeta^{\diamond},+} (x^{\diamond} \!, t^{\diamond})} u, \hskip10pt \omega \geqslant 
\osc_{Q_{R; R, 0}^{\upbeta^{\diamond}+} (x^{\diamond} \!, t^{\diamond})} u$
we have that if $\mu_+ (B_{\rho}) > 0$ and
\begin{align*}
\frac{ M_+ (A_{0}^+)}{|M|_{\Lambda} (Q_{R}^{\upbeta^{\diamond}, \texttt{\,>}} (x^{\diamond}\!, t^{\diamond}))} +
\frac{ \Lambda_+ (A_{0}^+)}{\Lambda (Q_{R}^{\upbeta^{\diamond}, \texttt{\,>}} (x^{\diamond}\!, t^{\diamond}))} \leqslant \overline{\nu}^{\diamond} ,
\end{align*}
$\hskip8pt$ where $A_0^+ = \{ (x,t) \in Q_{R; R, 0}^{\upbeta^{\diamond}, +} (x^{\diamond}\!, t^{\diamond}) \, | 
\, u(x,t) > \overline{m} - \sigma \omega \}$, then
$$
u(x,t) \leqslant \overline{m} - a \, \sigma \, \omega  \hskip30pt 
	\text{for a.e. }	(x,t) \in Q_{R; \rho, \theta^{\diamond}}^{\upbeta^{\diamond}, +} (x^{\diamond}\!, t^{\diamond}) \, ;
$$
$ii \, )$ $\overline{m} \geqslant \sup_{Q_{R; R, 0}^{\circ, -} (x^{\circ} \!, t^{\circ})} u, 
	\hskip10pt \omega \geqslant \osc_{Q_{R; R, 0}^{\circ, -} (x^{\circ} \!, t^{\circ})} u$
we have that if $\mu_- (B_{\rho}) > 0$ and
\begin{align*}
\frac{ M_- (A_{0}^-)}{|M|_{\Lambda} (Q_{R}^{\circ, \texttt{\,<}} (x^{\circ}\!, t^{\circ}))} +
\frac{ \Lambda_- (A_{0}^-)}{\Lambda (Q_{R}^{\circ, \texttt{\,<}} (x^{\circ}\!, t^{\circ}))} \leqslant \overline{\nu}^{\circ} ,
\end{align*}
$\hskip8pt$ where $A_0^- = \{ (x,t) \in Q_{R; R, 0}^{\circ, -} (x^{\circ} \!, t^{\circ}) \, | \, u(x,t) > \overline{m} - \sigma \omega \}$, then
$$
u(x,t) \leqslant \overline{m} - a \, \sigma \, \omega  \hskip30pt \text{for a.e. }	(x,t) \in Q_{R; \rho, \theta^{\circ}}^{\circ, -}  (x^{\circ}\!, t^{\circ}) \, ;
$$
$iii \, )$ $\overline{m} \geqslant \sup_{Q_{R}^{s_1^{\star} , s_2^{\star}} (x^{\star}\!, t^{\star})} u, \hskip10pt 
	\om \geqslant \osc_{B_R(x^{\star}) \times (s_1^{\star} , s_2^{\star})} u$
we have that if $\lambda_0 (B_{\rho}) > 0$ and
\begin{align*}
\Lambda_0 (A_{0}^0) \leqslant \overline{\nu}^{\star} \, \Lambda (Q_{R}^{s_1^{\star} , s_2^{\star}} (x^{\star}\!, t^{\star}))
\end{align*}
$\hskip8pt$ where $A_0^0 = \{ (x,t) \in Q_{R; R, s_1^{\star}, s_2^{\star}}^{0} (x^{\star} \!, t^{\star}) \, | \, u(x,t) > \overline{m} - \sigma \omega \}$,
then
$$
u(x,t) \leqslant \overline{m} - a \, \sigma \, \omega  \hskip30pt \text{for a.e. }	(x,t) \in Q_{R; \rho, s_1^{\star}, s_2^{\star}}^{0} (x^{\star} \!, t^{\star}) \, ;
$$
$iv \, )$ $\overline{m} \geqslant \sup_{B_R(x^{\star})} u (\cdot, t), \hskip10pt \om \geqslant \osc_{B_R(x^{\star})} u (\cdot, t)$
we have that if 
$B_R(x^{\star}) \subset \Omega_0$ and
\begin{align*}
\lambda \big(\{ x \in B_{R} (x^{\star}) \, | \, u(x,t) > \overline{m} - \sigma \omega \} \big) \leqslant 
	\overline{\nu} \ \lambda (B_R (x^{\star}) )
\end{align*}
then
$$
u(x,t) \leqslant \overline{m} - a \, \sigma \, \omega  \hskip30pt \text{for a.e. }	x \in B_{\rho} (x^{\star})
$$
for a.e. $t \in (s_1^{\star} , s_2^{\star})$.
\eprop

\boss
\label{pluto}
The requirement $\mu_+(B_{\rho}) > 0$ in point {\em i }\!)
(and analogously $\mu_-(B_{\rho}) > 0$ in point {\em ii }\!) and $\lambda_0 (B_{\rho}) > 0$ in point {\em iii }\!) is not technically needed,
for the proof it would be sufficient to have $\mu_+(B_{R}) > 0$. We require it just to give a meaning to the thesis of the theorem.
\eoss

\noindent
\dimo
We prove only the first claim, being the other similar. Often we will not write the point $(x^{\diamond} \!, t^{\diamond})$, just to simplify the notation.
First of all fix $a , \sigma \in (0,1)$ which will remain fixed for all the proof.
Choose $\theta^{\diamond} \in (0,1)$ and $\rho \in (0,R)$, assume that $\mu_+(B_{\rho}) > 0$ and
consider the following sequences ($h \in \N$)
$$
\displaystyle \rho_h = \rho + \varepsilon^{h} (R - \rho), \quad \quad \theta_h = \theta^{\diamond} - \varepsilon^{2h} \theta^{\diamond}  \, ,
$$
where $\varepsilon \in (0,1)$.
We require that $(\theta_{h+1} - \theta_h) R^2$ is to be equal to $(\rho_{h} - \rho_{h+1})^2$ (as required in 
Definition \ref{classiDG} and in the proof that
a $Q$-minimum belongs to the De Giorgi class, see \eqref{puredifave}): we derive that $\theta^{\diamond}$ has to satisfy
\begin{equation}
\label{teta}
\theta^{\diamond} = \frac{1 - \varepsilon}{1 + \varepsilon} \, \frac{(R-\rho)^2}{R^2} \, .
\end{equation}
Referring to definitions \eqref{notazione1} we will consider
$$
x_0 = x^{\diamond}\!, \qquad t_0 = t^{\diamond}\!, \qquad s_2 = s_2^{\diamond} := t^{\diamond} + \upbeta^{\diamond} \, h(x^{\diamond}\!, R) R^2 \, ,
$$
but we will often omit to write them just to simplify the notation. Now we moreover define, for $h \in \N$ and $a, \sigma \in (0,1)$,
\begin{equation}
\label{tutticonacca}
\begin{array}{c}
\displaystyle B_h = B_{\rho_h}(x^{\diamond}) \, , 																	\\
\displaystyle \delta_h : = \sum_{j = h}^{\infty} (\rho_j - \rho_{j+1}) = \rho_h - \rho = \varepsilon^{h} (R - \rho) \searrow 0 \, ,		\\ [4mm]
Q_h^+ :=	Q_{R; \rho, \theta_h}^{\upbeta^{\diamond}, +, \rho_h-\rho} (x^{\diamond} \!, t^{\diamond})							\\ [4mm]
I_h^+ := (I_{\rho}^+(x^{\diamond}))^{\delta_h} \, ,																	\\ [2mm]
\displaystyle \sigma_h = a \, \sigma + \varepsilon^h (1 - a)\, \sigma \searrow a\sigma\, ,	
		\quad\quad  k_h = \overline{m} - \sigma_h \om \nearrow  \overline{m} - a \sigma \omega 	\, ,							\\ [4mm]
A_h^+ = \{ (x,t) \in Q_h^+ \, | \, u(x,t) > k_h \}		\, .
\end{array}
\end{equation}
Notice that
$$
\begin{array}{c}
Q_{h+1}^+ \subset Q_h^+	\, , \hskip20pt				A_{h+1}^+ \subset A_h^+	\, ,				\\ [0.5em]
\rho_h - \rho_{h+1} = (1 - \varepsilon)\varepsilon^{h} (R - \rho) 	\, ,									\\ [0.5em]
\theta_{h+1} \, \upbeta^{\diamond} \, h(x^{\diamond}\!,R) \, R^2 - \theta_h \, \upbeta^{\diamond} \, h(x^{\diamond}\!,R) \, R^2 =
\theta^{\diamond} (1 - \varepsilon^2) \, \varepsilon^{2h} \, \upbeta^{\diamond} \, h(x^{\diamond}\!,R) \, R^2 .
\end{array}
$$ 
In the next picture we show some possible $Q_h^+$ marked by dashed lines, while the one marked by longer lines is the limit set (for $h \to +\infty$).
\ \\
\ \\
\ \\
\begin{picture}(150,200)(-180,0)
\put (-105,200){\linethickness{1pt}\line(1,0){210}}
\put (-105,50){\linethickness{1pt}\line(1,0){210}}
\put (-105,50){\linethickness{1pt}\line(0,1){150}}
\put (105,50){\linethickness{1pt}\line(0,1){150}}

\put (-40,220){\tiny$\mu > 0$}
\put (40,220){\tiny$\mu < 0 \text{ or }\mu = 0$}

\put (30,40){\line(0,1){170}}

\put (-180,125){\line(1,0){320}}


\put (-83,198){\line(0,-1){1}}
\put (-83,195){\line(0,-1){1}}
\put (-83,192){\line(0,-1){1}}
\put (-83,189){\line(0,-1){1}}
\put (-83,186){\line(0,-1){1}}
\put (-83,183){\line(0,-1){1}}
\put (-83,180){\line(0,-1){1}}
\put (-83,177){\line(0,-1){1}}
\put (-83,174){\line(0,-1){1}}
\put (-83,171){\line(0,-1){1}}
\put (-83,168){\line(0,-1){1}}
\put (-83,165){\line(0,-1){1}}
\put (-83,162){\line(0,-1){1}}
\put (-83,159){\line(0,-1){1}}
\put (-83,156){\line(0,-1){1}}
\put (-83,153){\line(0,-1){1}}
\put (-83,150){\line(0,-1){1}}
\put (-83,147){\line(0,-1){1}}
\put (-83,144){\line(0,-1){1}}
\put (-83,141){\line(0,-1){1}}
\put (-83,138){\line(0,-1){1}}
\put (-83,135){\line(0,-1){1}}
\put (-83,132){\line(0,-1){1}}

\put (-83,129){\line(1,0){1}}
\put (-80,129){\line(1,0){1}}
\put (-77,129){\line(1,0){1}}
\put (-74,129){\line(1,0){1}}
\put (-71,129){\line(1,0){1}}
\put (-68,129){\line(1,0){1}}
\put (-65,129){\line(1,0){1}}
\put (-62,129){\line(1,0){1}}
\put (-59,129){\line(1,0){1}}
\put (-56,129){\line(1,0){1}}
\put (-53,129){\line(1,0){1}}
\put (-50,129){\line(1,0){1}}
\put (-47,129){\line(1,0){1}}
\put (-44,129){\line(1,0){1}}
\put (-41,129){\line(1,0){1}}
\put (-38,129){\line(1,0){1}}
\put (-35,129){\line(1,0){1}}
\put (-32,129){\line(1,0){1}}
\put (-29,129){\line(1,0){1}}
\put (-26,129){\line(1,0){1}}
\put (-23,129){\line(1,0){1}}
\put (-20,129){\line(1,0){1}}
\put (-17,129){\line(1,0){1}}
\put (-14,129){\line(1,0){1}}
\put (-11,129){\line(1,0){1}}
\put (-8,129){\line(1,0){1}}
\put (-5,129){\line(1,0){1}}
\put (-2,129){\line(1,0){1}}
\put (1,129){\line(1,0){1}}
\put (4,129){\line(1,0){1}}
\put (7,129){\line(1,0){1}}
\put (10,129){\line(1,0){1}}
\put (13,129){\line(1,0){1}}
\put (16,129){\line(1,0){1}}

\put (17,129){\line(0,-1){1}}
\put (17,126){\line(0,-1){1}}

\put (-76,199){\line(0,-1){1}}
\put (-76,196){\line(0,-1){1}}
\put (-76,193){\line(0,-1){1}}
\put (-76,190){\line(0,-1){1}}
\put (-76,187){\line(0,-1){1}}
\put (-76,184){\line(0,-1){1}}
\put (-76,181){\line(0,-1){1}}
\put (-76,178){\line(0,-1){1}}
\put (-76,175){\line(0,-1){1}}
\put (-76,172){\line(0,-1){1}}
\put (-76,169){\line(0,-1){1}}
\put (-76,166){\line(0,-1){1}}
\put (-76,163){\line(0,-1){1}}
\put (-76,160){\line(0,-1){1}}
\put (-76,157){\line(0,-1){1}}
\put (-76,154){\line(0,-1){1}}
\put (-76,151){\line(0,-1){1}}
\put (-76,148){\line(0,-1){1}}
\put (-76,145){\line(0,-1){1}}
\put (-76,142){\line(0,-1){1}}

\put (-75,141){\line(1,0){1}}
\put (-72,141){\line(1,0){1}}
\put (-69,141){\line(1,0){1}}
\put (-66,141){\line(1,0){1}}
\put (-63,141){\line(1,0){1}}
\put (-60,141){\line(1,0){1}}
\put (-57,141){\line(1,0){1}}
\put (-54,141){\line(1,0){1}}
\put (-51,141){\line(1,0){1}}
\put (-48,141){\line(1,0){1}}
\put (-45,141){\line(1,0){1}}
\put (-42,141){\line(1,0){1}}
\put (-39,141){\line(1,0){1}}
\put (-36,141){\line(1,0){1}}
\put (-33,141){\line(1,0){1}}
\put (-30,141){\line(1,0){1}}
\put (-27,141){\line(1,0){1}}
\put (-24,141){\line(1,0){1}}
\put (-21,141){\line(1,0){1}}
\put (-18,141){\line(1,0){1}}
\put (-15,141){\line(1,0){1}}
\put (-12,141){\line(1,0){1}}
\put (-9,141){\line(1,0){1}}
\put (-6,141){\line(1,0){1}}
\put (-3,141){\line(1,0){1}}
\put (0,141){\line(1,0){1}}
\put (3,141){\line(1,0){1}}
\put (6,141){\line(1,0){1}}
\put (9,141){\line(1,0){1}}
\put (12,141){\line(1,0){1}}
\put (15,141){\line(1,0){1}}
\put (18,141){\line(1,0){1}}
\put (21,141){\line(1,0){1}}

\put (24,141){\line(0,-1){1}}
\put (24,138){\line(0,-1){1}}
\put (24,135){\line(0,-1){1}}
\put (24,132){\line(0,-1){1}}
\put (24,129){\line(0,-1){1}}
\put (24,126){\line(0,-1){1}}

\put (-75,210){$- \rho$}
\put (-70,200){\linethickness{2pt}\line(1,0){1}}

\put (65,210){$\rho$}
\put (70,200){\linethickness{2pt}\line(1,0){1}}


\put (-70,200){\line(0,-1){4}}
\put (-70,195){\line(0,-1){4}}
\put (-70,190){\line(0,-1){4}}
\put (-70,185){\line(0,-1){4}}
\put (-70,180){\line(0,-1){4}}
\put (-70,175){\line(0,-1){4}}
\put (-70,170){\line(0,-1){4}}
\put (-70,165){\line(0,-1){4}}
\put (-70,160){\line(0,-1){4}}
\put (-70,155){\line(0,-1){4}}
\put (-70,150){\line(1,0){4}}
\put (-65,150){\line(1,0){4}}
\put (-60,150){\line(1,0){4}}
\put (-55,150){\line(1,0){4}}
\put (-50,150){\line(1,0){4}}
\put (-45,150){\line(1,0){4}}
\put (-40,150){\line(1,0){4}}
\put (-35,150){\line(1,0){4}}
\put (-30,150){\line(1,0){4}}
\put (-25,150){\line(1,0){4}}
\put (-20,150){\line(1,0){4}}
\put (-15,150){\line(1,0){4}}
\put (-10,150){\line(1,0){4}}
\put (-5,150){\line(1,0){4}}
\put (0,150){\line(1,0){4}}
\put (5,150){\line(1,0){4}}
\put (10,150){\line(1,0){4}}
\put (15,150){\line(1,0){4}}
\put (20,150){\line(1,0){4}}
\put (25,150){\line(1,0){4}}

\put (-20,30){Figure C}
\end{picture}

\noindent
First of all notice that since
\begin{align*}
(k_{h+1} - k_h)^2 & M_+ (A_{h+1}^+) \leqslant  \iint_{A_{h+1}^+}  (u - k_h)_+^2 \mu_+ \, dx dt \leqslant 	\iint_{Q_{h+1}^+}  (u - k_h)_+^2 \mu_+ \, dx dt
\end{align*}
and $k_{h+1} - k_h = (1-a) \, \sigma \, \omega \, \varepsilon^{h+1}$ we can estimate
\begin{align}
\label{first}
\varepsilon^{2h+2} \, (1-a)^2 \sigma^2 \omega^2 \, \frac{ M_+ (A_{h+1}^+)}{|M|_{\Lambda} (Q_R^{\upbeta^{\diamond}, \texttt{\,>}})} \leqslant
	\frac{1}{|M|_{\Lambda} (Q_R^{\upbeta^{\diamond}, \texttt{\,>}})} \iint_{Q_{h+1}^+}  (u - k_h)_+^2 \mu_+ \, dx dt \, .
\end{align}
Similarly
\begin{equation}
\label{second}
{\displaystyle
\varepsilon^{2h+2} \, (1-a)^2 \sigma^2 \om^2 \, \frac{\Lambda_+ (A_{h+1}^+)}{\Lambda (Q_R^{\upbeta^{\diamond}, \texttt{\,>}})}}
\leqslant \frac{1}{\Lambda (Q_R^{\upbeta^{\diamond}, \texttt{\,>}})} \iint_{Q_{h+1}^+} (u-k_h)_+^2 \lambda_+ \, dx dt \, .
\end{equation}

\noindent
Then we can argue in a completely analogous way as done to obtain \eqref{mitoccanum1} and \eqref{mitoccanum2}. Taking in \eqref{mitoccanum1}
\begin{equation}
\label{valoriacchesimi}
\begin{array}{l}
\rho_{h+1} = r \, , \hskip20pt \rho_h = \tilde{r} \, , \hskip20pt \rho_{h-1} = \hat{r}	\, ,	\hskip20pt \rho \text{ in place of } R/2 \, ,		\\	[2mm]
\theta_{h+1} = \theta \, , \hskip20pt \theta_h = \tilde\theta \, , \hskip20pt \theta_{h-1} = \hat\theta \, , \hskip20pt k_h = k \, ,
\end{array}
\end{equation}
we get (the only difference with \eqref{mitoccanum1} is that $2(\rho_h - \rho_{h+1}) \not= \rho_{h-1} - \rho_{h}$ unless $\varepsilon = 1/2$)
\begin{align*}
\frac{1}{|M|_{\Lambda} (Q_R^{\upbeta^{\diamond}, \texttt{\,>}})} &
\iint_{Q_{h+1}^+}  (u-k_h)_+^2 \mu_+ \, dx dt \leqslant																	\nonumber	\\
\leqslant & \hskip5pt \gamma_1^{2/\kappa} \, R^2 \, \frac{1 + \upbeta^{\diamond}}{({\upbeta^{\diamond}})^{\frac{1}{\kappa}}} \, 
												\frac{2 \gamma + 2}{(\rho_{h} - \rho_{h+1})^2} \, \cdot					\nonumber 	\\
& \cdot \Bigg[ \frac{1}{|M|_{\Lambda} (Q_R^{\upbeta^{\diamond}, \texttt{\,>}})} \iint_{Q_{h-1}^+} (u - k_h)_+^2 \mu_+ \, dx dt
				+ \frac{1}{\Lambda (Q_R^{\upbeta^{\diamond}, \texttt{\,>}})} \iint_{Q_{h-1}^+} (u - k_h)_+^2 \lambda_+ \, dx dt \, +	\nonumber	\\
& \hskip50pt + \frac{1}{\Lambda (Q_R^{\upbeta^{\diamond}, \texttt{\,>}})}  \iint_{I_{h-1}^+ \times (t^{\diamond} \!, s_2^{\diamond})} 
															(u - k_h)_+^2 (\lambda_0 + \lambda_-) \, dx dt \, +			\nonumber	\\
& \hskip50pt	+  (\rho_{h} - \rho_{h+1})^2 \, \frac{1}{\Lambda (Q_R^{\upbeta^{\diamond}, \texttt{\,>}})}
	\sup_{t \in (t^{\diamond}\!, s_2^{\diamond})} \int_{I_{h-1}^+} (u - k_h)_+^2 (x,t) |\mu| (x) dx \Bigg] \, .
\end{align*}
Now since (here we use $\sigma_h \leqslant \sigma$)
\begin{align*}
\iint_{Q_{h-1}^+}  (u-k_h)_+^2 \mu_+ \, dx dt \leqslant M_+ (A_{h-1}^+) \ \sup_{Q_{h-1}^+} (u - k_h)^2 \leqslant M_+ (A_{h-1}^+) (\sigma \omega)^2 \, ,	  \\
\iint_{Q_{h-1}^+}  (u-k_h)_+^2 \lambda_+ \, dx dt \leqslant 
							\Lambda_+ (A_{h-1}^+) \ \sup_{Q_{h-1}^+} (u - k_h)^2 \leqslant \Lambda_+ (A_{h-1}^+) (\sigma \omega)^2 \, ,	
\end{align*}
by the above inequality and by \eqref{first} we get
\begin{align*}
\frac{ M_+ (A_{h+1}^+)}{|M|_{\Lambda} (Q_R^{\upbeta^{\diamond}, \texttt{\,>}})} \leqslant & \, 
		\frac{\gamma_1^{2/\kappa} \, R^2}{\varepsilon^{2h+2} (1-a)^2 \sigma^2 \omega^2} \, 
		\frac{1 + \upbeta^{\diamond}}{(\upbeta^{\diamond})^{\frac{1}{\kappa}}}
		\frac{2 \gamma + 2}{(1-\varepsilon)^2 \varepsilon^{2h}(R - \rho)^2} \, 												\\
& \cdot \left( \frac{ M_+ (A_{h}^+)} {|M|_{\Lambda} (Q_R^{\upbeta^{\diamond}, \texttt{\,>}})} \right)^{\frac{\kappa - 1}{\kappa}}\, \cdot
	\Bigg[ \frac{M_+ (A_{h-1}^+)}{|M|_{\Lambda} (Q_R^{\upbeta^{\diamond}, \texttt{\,>}})} (\sigma \omega)^2 + 
		\frac{\Lambda_+ (A_{h-1}^+)}{\Lambda (Q_R^{\upbeta^{\diamond}, \texttt{\,>}})} (\sigma \omega)^2 + 					\\
& \hskip30pt + \frac{1}{\Lambda (Q_R^{\upbeta^{\diamond}, \texttt{\,>}})} \iint_{I_{h-1}^+ \times (t^{\diamond} \!, s_2^{\diamond})}
														(u - k_h)_+^2 (\lambda_0 + \lambda_-) \, dx dt \, +				\\
& \hskip30pt + \, \frac{(R - \rho)^2 (1-\varepsilon)^2 \varepsilon^{2h}}{\Lambda (Q_R^{\upbeta^{\diamond}, \texttt{\,>}})}
	\sup_{t \in (t^{\diamond}\!, s_2^{\diamond})} \int_{I_{h-1}^+} (u - k_h)_+^2 (x,t) |\mu| (x) dx \Bigg]	\, .					
\end{align*}
Now defining first
\begin{align*}
y_h := & \ y^{M}_h + y^{\Lambda}_h \, , \hskip15pt \text{ where }
	y^{M}_h := \frac{M_+ (A_{h}^+)}{|M|_{\Lambda} (Q_R^{\upbeta^{\diamond}, \texttt{\,>}})}
			 \hskip10pt \text{and} \hskip10pt y^{\Lambda}_h := \frac{\Lambda_+ (A_{h}^+)}{\Lambda (Q_R^{\upbeta^{\diamond}, \texttt{\,>}})}  \, ,
\end{align*}
and, since $(u - k_h)_+^2$ is bounded by $(\sigma \omega)^2$ and estimating
\begin{align*}
\frac{1}{\sigma^2 \omega^2 \Lambda (Q_R^{\upbeta^{\diamond}, \texttt{\,>}})} & \iint_{I_{h}^+ \times (t^{\diamond} \!, s_2^{\diamond})}
														(u - k_h)_+^2 (\lambda_0 + \lambda_-) \, dx dt \, +						\\
& \hskip10pt + \, \frac{(R - \rho)^2 (1-\varepsilon)^2 \varepsilon^{2(h-1)}}{\sigma^2 \omega^2\Lambda (Q_R^{\upbeta^{\diamond}, \texttt{\,>}})}
	\sup_{t \in (t^{\diamond}\!, s_2^{\diamond})} \int_{I_{h}^+} (u - k_h)_+^2 (x,t) |\mu| (x) dx	\leqslant									\\
\leqslant & \ \frac{(\Lambda_0 + \Lambda_-)(I_{h}^+ \times (t^{\diamond} \!, s_2^{\diamond}))}{\Lambda (Q_R^{\upbeta^{\diamond}, \texttt{\,>}})}
	+ \, \frac{(R - \rho)^2 (1-\varepsilon)^2 \varepsilon^{2(h-1)} |\mu| (I_{h}^+)}{\Lambda (Q_R^{\upbeta^{\diamond}, \texttt{\,>}})}
\end{align*}
defining also
\begin{align*}
\epsilon_h := \frac{(\Lambda_0 + \Lambda_-)(I_{h}^+ \times (t^{\diamond} \!, s_2^{\diamond}))}{\Lambda (Q_R^{\upbeta^{\diamond}, \texttt{\,>}})}
	+ \, \frac{(R - \rho)^2 (1-\varepsilon)^2 \varepsilon^{2(h-1)} |\mu| (I_{h}^+)}{\Lambda (Q_R^{\upbeta^{\diamond}, \texttt{\,>}})}				
\end{align*}
we first get
\begin{align*}
y^M_{h+1} \leqslant
\frac{\gamma_1^{2/\kappa} R^2 \, (2\gamma + 2)}{(1-a)^2 (1-\varepsilon)^2 \varepsilon^2 (R - \rho)^2}
	\, \frac{1 + \upbeta^{\diamond}}{(\upbeta^{\diamond})^{\frac{1}{\kappa}}} \, \frac{1}{\varepsilon^{4h}} \ 
	(y^M_{h})^{\frac{\kappa - 1}{\kappa}} \left[ y^{M}_{h-1} + y^{\Lambda}_{h-1} + \epsilon_{h-1} \right] \, .
\end{align*}
Taking \eqref{valoriacchesimi} in \eqref{mitoccanum2} we can argue in a similar way to estimate $y^{\Lambda}_{h+1}$ and get
\begin{align*}
y^{\Lambda}_{h+1} \leqslant
\frac{\gamma_1^{2/\kappa} R^2\, (2\gamma + 2)}{(1-a)^2 (1 - \varepsilon)^2 \varepsilon^2 (R - \rho)^2} 
		\, \frac{1 + \upbeta^{\diamond}}{(\upbeta^{\diamond})^{\frac{1}{\kappa}}} \, 
		\frac{1}{\varepsilon^{4h}} \ 
		 (y^{\Lambda}_{h})^{\frac{\kappa - 1}{\kappa}} \left[ y^{M}_{h-1} + y^{\Lambda}_{h-1} + \epsilon_{h-1} \right] \, .
\end{align*}
Summing the two inequalities
and since the sequences $(y^{M}_h)_h$, $(y^{\Lambda}_h)_h$ are decreasing we finally get
\begin{align*}
y_{h+1} \leqslant \ 
\frac{\gamma_1^{2/\kappa} R^2 \, (2\gamma + 2)}{(1-a)^2 (1-\varepsilon)^2 \varepsilon^2 (R - \rho)^2} \, 
		\frac{1 + \upbeta^{\diamond}}{(\upbeta^{\diamond})^{\frac{1}{\kappa}}} \, 
		\frac{1}{\varepsilon^{4h}} \ y_{h-1}^{\frac{\kappa - 1}{\kappa}} \left( y_{h-1} + \epsilon_{h-1} \right)
\end{align*}
for every $h \geqslant 1$; then, for instance, 
\begin{align*}
y_{2(h+1)} \leqslant \ 
\frac{\gamma_1^{2/\kappa} R^2 \, (2\gamma + 2)}{(1-a)^2 (1-\varepsilon)^2 (R - \rho)^2 \varepsilon^{6}} \, 
		\frac{1 + \upbeta^{\diamond}}{(\upbeta^{\diamond})^{\frac{1}{\kappa}}} \, 
		\frac{1}{\varepsilon^{8h}} \ y_{2h}^{\frac{\kappa - 1}{\kappa}} \left( y_{2h} + \epsilon_{2h} \right) \, .
\end{align*}
Using \eqref{teta} to write $R^2/(R - \rho)^2$ and Lemma \ref{lemmuzzofurbo-quinquies} with
\begin{gather*}
c = \frac{\gamma_1^{2/\kappa} \, (2\gamma + 2)}{(1 - a)^2 (1 - \varepsilon^2) \, \varepsilon^{6} \, \theta^{\ast}} \, 
	\frac{1 + \upbeta^{\diamond}}{(\upbeta^{\diamond})^{\frac{1}{\kappa}}} \, , 
	\hskip10pt \alpha = \frac{\kappa - 1}{\kappa} \, , \hskip10pt b = \frac{1}{\varepsilon^8} \, ,
\end{gather*}
we derive that the subsequence $(y_{2h})_h$ of even indexes, and in fact the whole sequence $(y_h)_h$ since $(y_h)_h$ is decreasing, 
is converging to zero provided that
\begin{align*}
\frac{ M_+ (A_{0}^+)}{|M|_{\Lambda} (Q_R^{\upbeta^{\diamond}, \texttt{\,>}})} + 
	\frac{ \Lambda_+ (A_{0}^+)}{\Lambda (Q_R^{\upbeta^{\diamond}, \texttt{\,>}})} \leqslant 
\left(\frac{(1-a)^2 (1-\varepsilon^2) \, \varepsilon^{6} \, \theta^{\diamond} (\upbeta^{\diamond})^{\frac{1}{\kappa}}}
	{\gamma_1^{2/\kappa} \, (1 + \upbeta^{\diamond}) \, (2\gamma + 2)} \right)^{\frac{\kappa}{\kappa-1}} \, 
	\, \varepsilon^{\frac{8\kappa^2}{(\kappa - 1)^2}} \, .
\end{align*}
By the definition of $A_h$ we have that
\begin{gather*}
Q_0^+ = Q_{R; \rho, 0}^{\upbeta^{\diamond}, +,R - \rho} (x^{\diamond} \!, t^{\diamond}) \qquad \text{and} \qquad 
			A_0^+ = \big\{ (x,t) \in Q_0^+ \, \big| \, u(x,t) > \overline{m} - \sigma \omega \big\}
\end{gather*}
but we can consider
\begin{align*}
A_0^+ & = \big\{ (x,t) \in Q_{R; R, 0}^{\upbeta^{\diamond}, +} (x^{\diamond} \!, t^{\diamond}) \, \big| \, u(x,t) > \overline{m} - \sigma \omega \big\} =		\\	
& = \big\{ (x,t) \in B_R^+ (x^{\diamond}) \times (t^{\diamond} \!, t^{\diamond} + h(x^{\diamond} \!, R) R^2) \, \big| \, u(x,t) > \overline{m} - \sigma \omega \big\}
\end{align*}
since we will consider the measures $M_+$ and $\Lambda_+$ of this set.
Then we have derived that
$$
u(x,t) \leqslant \overline{m} - a \, \sigma \, \omega  \hskip30pt \text{for a.e. } 
	(x,t) \in Q_{R; \rho, \theta^{\diamond} }^{\upbeta^{\diamond},+} (x^{\diamond} \!, t^{\diamond})
$$
provided that
$$
\frac{ M_+ (A_{0}^+)}{|M|_{\Lambda} (Q_R^{\upbeta^{\diamond}, \texttt{\,>}})} + 
	\frac{ \Lambda_+ (A_{0}^+)}{\Lambda (Q_R^{\upbeta^{\diamond}, \texttt{\,>}})} 
		\leqslant \overline{\nu}^{\diamond}
$$
where
$$
\overline{\nu}^{\diamond} =
\left(\frac{(1-a)^2 (1-\varepsilon^2) \, \varepsilon^{6} \, \theta^{\diamond} \, (\upbeta^{\diamond})^{\frac{1}{\kappa}}}
			{\gamma_1^{2/\kappa} \, (1 + \upbeta^{\diamond}) \, (2\gamma + 2)} \right)^{\frac{\kappa}{\kappa-1}} \, 
			\, \varepsilon^{\frac{8\kappa^2}{(\kappa - 1)^2}} \, .
$$
In a complete analogous way: fix a point $(x^{\circ} \!, t^{\circ})$ such that $\mu_- (B_R(x^{\circ})) > 0$.
One gets that taking the same values as before for $\rho, a , \sigma$ and $\theta^{\circ} \in (0,1)$ there is $\overline{\nu}^{\circ} > 0$
such that if
$$
\frac{ M_- (A_{0}^-)}{|M|_{\Lambda} (Q_R^{\upbeta^{\circ}, \texttt{\,<}})} +
\frac{ \Lambda_- (A_{0}^-)}{\Lambda (Q_R^{\upbeta^{\circ}, \texttt{\,<}})} \leqslant \overline{\nu}^{\circ} \, , 
$$
where the ball $B_R$ is centred in $x^{\circ}$ and
\begin{gather*}
A_0^- = \left\{ (x,t) \in Q_{R; \rho, 0}^{\upbeta^{\circ}, -,R-\rho} (x^{\circ} \!, t^{\circ}) \, | \, u(x,t) > \overline{m} - \sigma \omega \right\} \, , 		\\
\overline{\nu}^{\circ} =
\left(\frac{(1-a)^2 (1-\varepsilon^2) \, \varepsilon^{6} \, \theta^{\circ} \, (\upbeta^{\circ})^{\frac{1}{\kappa}}}
	{\gamma_1^{2/\kappa} \, (1 + \upbeta^{\circ}) \, (2\gamma + 2)} \right)^{\frac{\kappa}{\kappa-1}} \, 
	\, \varepsilon^{\frac{8\kappa^2}{(\kappa - 1)^2}} \, ,
\end{gather*}
then
$$
u(x,t) \leqslant \overline{m} - a \, \sigma \, \omega  \hskip30pt \text{for a.e. } 
	(x,t) \in Q_{R; \rho, \theta^{\circ} }^{\upbeta^{\circ}, -} (x^{\circ} \!, t^{\circ}) \, .
$$
Finally we analyse the part in which $\mu \equiv 0$, which is slightly different.
Fix a point $(x^{\star} \!, t^{\star})$ such that $\lambda_0 (B_R(x^{\star})) > 0$, consider $k_h$ and $\sigma_h$ as in
\eqref{tutticonacca}.
Arguing as done to obtain \eqref{mitoccanum3} and taking in \eqref{mitoccanum3} for $k, r, \tilde{r}, \hat{r}$ the same values as in \eqref{valoriacchesimi}
and for $\sigma_1, \sigma_2$ respectively $s_1^{\star}$ and $s_2^{\star}$ we get
\begin{align*}
& \frac{1}{\Lambda (B_R \times (s_1^{\star}, s_2^{\star}))}
\int\!\!\!\int_{Q_{h+1}^0}  (u - k_h)_+^2 \lambda_0 \, dx dt \leqslant														\\
& \hskip40pt \leqslant \gamma_1^{2/\kappa} \, (\upbeta^{\star})^{\frac{\kappa - 1}{\kappa}}
	R^2 \, \frac{2 \gamma + 2}{(1 - \varepsilon)^2 \varepsilon^{2h}(R - \rho)^2} \, 
	\frac{\big( \Lambda_0 (A_h^0) \big)^{\frac{\kappa - 1}{\kappa}}}
									{(\Lambda (B_R \times (s_1^{\star}, s_2^{\star})))^{\frac{\kappa - 1}{\kappa}}}\, \cdot 		\\
& \hskip50pt \cdot \Bigg[ 
	\frac{1}{\Lambda (B_R \times (s_1^{\star}, s_2^{\star}))}
	\iint_{Q_{h-1}^0} (u - k_h)_+^2 \lambda_0 \, dx dt +																\\
& \hskip60pt + \frac{1}{\Lambda (B_R \times (s_1^{\star}, s_2^{\star}))}
	\iint_{I_{h-1}^0 \times (s_1^{\star} \!, s_2^{\star})} (u - k_h)_+^2 (\lambda_+ + \lambda_-) \, dx dt +							\\
& \hskip60pt + \frac{(1 - \varepsilon)^2 \varepsilon^{2h}(R - \rho)^2}{\Lambda (B_R \times (s_1^{\star}, s_2^{\star}))}
	\sup_{t \in (s_1^{\star}, s_2^{\star})}  \int_{(B_{\rho_h}^0)^{\rho_h - \rho}} (u - k_h)_+^2 (x,t) \lambda_0 (x) dx				\\
& \hskip60pt	+  \frac{(1 - \varepsilon)^2 \varepsilon^{2h}(R - \rho)^2}{\Lambda (B_R \times (s_1^{\star}, s_2^{\star}))}	
	\sup_{t \in (s_1^{\star}, s_2^{\star})} \int_{I_{h-1}^0} (u - k_h)_+^2 (x,t) \mu_+ (x) dx +									\\
& \hskip60pt	+  \frac{(1 - \varepsilon)^2 \varepsilon^{2h}(R - \rho)^2}{\Lambda (B_R \times (s_1^{\star}, s_2^{\star}))}	
	\sup_{t \in (s_1^{\star}, s_2^{\star})} \int_{I_{h-1}^0} (u - k_h)_+^2 (x,t) \mu_- (x) dx \Bigg]
\end{align*}
where
\begin{gather*}
I_h^0 := (I_{\rho}^0(x^{\star}))^{\rho_h - \rho} \setminus I_{\rho,\rho_h - \rho}^0(x^{\star}) 					\\
Q_h^0 :=	Q_{R; \rho, s_1^{\star}, s_2^{\star}}^{0, \rho_h - \rho} (x^{\star} \!, t^{\star})	\, ,					\\
A_h^0 = \{ (x,t) \in Q_h^0 \, | \, u(x,t) > k_h \}		\, .
\end{gather*}
Since, as for \eqref{first}, we have
\begin{gather*}
\varepsilon^{2h+2} \, (1-a)^2 \sigma^2 {\omega}^2 \, \frac{\Lambda_0 (A_{h+1}^0)}{\Lambda (B_R \times (s_1^{\star}, s_2^{\star}))} \leqslant
	\frac{1}{\Lambda (B_R \times (s_1^{\star}, s_2^{\star}))} \int\!\!\!\int_{Q_{h+1}^0}  (u - k_h)_+^2 \lambda_0 \, dx dt \, ,				\\
\iint_{Q_{h-1}^0}  (u-k_h)_+^2 \lambda_0 \, dx dt \leqslant 
						\Lambda_0 (A_{h-1}^0) \ \sup_{Q_{h-1}^0} (u - k_h)^2 \leqslant \Lambda_0 (A_{h-1}^0) (\sigma {\omega})^2 \, ,	
\end{gather*}
we derive
\begin{align*}
y_{h+1} \leqslant
\frac{\gamma_1^{2/\kappa} \, (\upbeta^{\star})^{\frac{\kappa - 1}{\kappa}} \, R^2 \, (2\gamma + 2)}{(1-a)^2 (1-\varepsilon)^2
	\varepsilon^2 (R - \rho)^2} \, \frac{1}{\varepsilon^{4h}} \ y_{h-1}^{\frac{\kappa - 1}{\kappa}} \left( y_{h-1} + \epsilon_{h-1} \right)
\end{align*}
where here we have defined
\begin{align*}
y_h := & \ \frac{\Lambda_0 (A_{h}^0)}{\Lambda (B_R \times (s_1^{\star}, s_2^{\star}))}									\\
\epsilon_h := & \frac{1}{\Lambda (B_R \times (s_1^{\star}, s_2^{\star}))}
	\Bigg[ (\Lambda_+ + \Lambda_-) (I_{h-1}^0 \times (s_1^{\star}, s_2^{\star})) + 										\\
& \hskip40pt	+  (R - \rho)^2 (1-\varepsilon)^2 \varepsilon^{2h} \Big( \Lambda_0 ((B_{\rho_h}^0)^{\rho_h - \rho}) +
	|\mu|	 (I_{h-1}^0) \Big) \Bigg] .
\end{align*}
Arguing similarly as before we get that $y_h$ tends to zero, that is
$$
u(x,t) \leqslant \overline{m} - a \, \sigma \, \omega  \hskip30pt \text{for a.e. } 
	(x,t) \in Q_{R; \rho, s_1^{\star}, s_2^{\star}}^{0} (x^{\star} \!, t^{\star}) \, ,
$$
provided that
$$
\frac{ \Lambda_0 (A_{0}^0)}{\Lambda (B_R \times (s_1^{\star}, s_2^{\star})} \leqslant \overline{\nu}^{\star}
$$
where
$$
\overline{\nu}^{\star} =
\left[
\frac{(1-a)^2 \, (1 - \varepsilon)^2 \, \varepsilon^6 \, (R - \rho)^2}
	{\gamma_1^{2/\kappa}  R^2 \, (2\gamma + 2)} \right]^{\frac{\kappa}{\kappa-1}}\, \frac{1}{\upbeta^{\star}} \, 
	\varepsilon^{\frac{8 \kappa^2}{(\kappa - 1)^2}} \, .
$$
Notice that ($\gamma_1 > 1$)
$$
\left[
\frac{(1-a)^2 \, (1 - \varepsilon)^2 \, \varepsilon^6 \, (R - \rho)^2}
	{\gamma_1^{2/\kappa}  R^2 \, (2\gamma + 2)} \right]^{\frac{\kappa}{\kappa-1}} \, 
	\, \varepsilon^{\frac{8 \kappa^2}{(\kappa - 1)^2}} \leqslant 1
$$
and to garantee $\overline{\nu}^{\star} \leqslant 1$ for every choice of $\upbeta^{\star}$ (say less than $1$)
we can choose $\varepsilon$ in a suitable way. For example taking $\varepsilon$ in such a way that
$\varepsilon^{\frac{8 \kappa^2}{(\kappa - 1)^2}} / \upbeta^{\star} = 1/2$, i.e.
$$
\varepsilon = \left( \frac{\upbeta^{\star}}{2} \right)^{\frac{(\kappa - 1)^2}{8 \kappa^2}}
$$
we have $\overline{\nu}^{\star} < 1$ and we get rid of the dependence of $1/\upbeta^{\star}$ for $\upbeta^{\star}$ small. \\ [0.3em]
%
%
%
%
For the last point we can proceed as follows: first notice that $B_R := B_R (x^{\star}) \subset \Omega_0$.
With the same $k_h$ and $\rho_h$ as before we consider $B_h := B_{\rho_h} (x^{\star})$, define the sequence of test functions
$$
\zeta_h : B_R \to [0,1] \, , \qquad
\zeta_h (x) =	\left\{
			\begin{array}{ll}
			1	&	\text{ in } B_{h+1}						\\
			0	&	\text{ in } B_R \setminus B_{h}
			\end{array}
			\right.
\qquad | D \zeta_h | \leqslant \frac{1}{\rho_{h} - \rho_{h+1}}
$$
and for almost every $t \in (0,T)$ we define $A_h = \{ x \in B_{\rho_h} (x^{\star}) \, | \, u(x,t) > k_h \}$.
Using Theorem \ref{chanillo-wheeden} with $2 \kappa$ in the place of $q$ (see also Remark \ref{rmkipotesi}) we have
\begin{align*}
& (1 - a)^2 \sigma^2 \omega^2 \varepsilon^{2(h+1)} \frac{\lambda (A_{h+1})}{\lambda (B_R)} \leqslant
	\frac{1}{\lambda (B_R)} \int_{B_{{h+1}}}  (u - k_h)_+^2 (x,t) \lambda (x)  \, dx \leqslant													\\
& \hskip30pt \leqslant \frac{1}{\lambda (B_R)} \int_{B_{{h}}}  (u - k_h)_+^2 (x,t) \zeta_h^2 (x) \lambda (x)  \, dx \leqslant							\\
& \hskip30pt \leqslant  \left(\frac{\lambda (A_h)}{\lambda (B_R)}\right)^{\frac{\kappa - 1}{\kappa}} 
		\left[\frac{1}{\lambda (B_{R})} \int_{B_{{h}}} 
		(u - k_h)_+^{2\kappa} (x,t) \zeta_h^{2\kappa} (x) \lambda(x) \, dx \right]^{\frac{1}{\kappa}} \leqslant									\\
& \hskip30pt \leqslant \left(\frac{\lambda (A_h)}{\lambda (B_R)}\right)^{\frac{\kappa - 1}{\kappa}}
		\frac{\gamma_1^{2} \, R^{2}}{\lambda (B_R)} \int_{B_{{h}}} | D \big( (u - k_h)_+ \zeta_h \big) |^2 \lambda \, dx \leqslant						\\
& \hskip30pt \leqslant \left(\frac{\lambda (A_h)}{\lambda (B_R)}\right)^{\frac{\kappa - 1}{\kappa}}
		\frac{2 \, \gamma_1^{2} \, R^{2}}{\lambda (B_R)} \int_{B_{{h}}} 
			\left[ | D (u - k_h)_+ |^2 + \frac{1}{(\rho_h - \rho_{h+1})^2} (u - k_h)_+^2 \right] \lambda \, dx \leqslant							\\
& \hskip30pt \leqslant \left(\frac{\lambda (A_h)}{\lambda (B_R)}\right)^{\frac{\kappa - 1}{\kappa}}
		\frac{2 \, \gamma_1^{2} \, R^{2}}{\lambda (B_R)} 
		\Bigg[ \frac{\gamma}{(\rho_{h-1} - \rho_h)^2} \int_{B_{{h-1}}} (u - k_h)_+^2 (x,t)\, \lambda (x) \, dx +									\\
& \hskip200pt		+	 \frac{1}{(\rho_h - \rho_{h+1})^2} \int_{B_{h}} (u - k_h)_+^2 \lambda \, dx \Bigg]  \leqslant								\\
& \hskip30pt \leqslant \left(\frac{\lambda (A_h)}{\lambda (B_R)}\right)^{\frac{\kappa - 1}{\kappa}}
		\frac{2 \, \gamma_1^{2} \, R^{2}}{\lambda (B_R)} \, 
		\frac{\gamma + 1}{\varepsilon^{2h} (R - \rho)^2 (1 - \varepsilon)^2}  \int_{B_{{h-1}}} (u - k_h)_+^2 (x,t)\, \lambda (x) \, dx \leqslant			\\
& \hskip30pt \leqslant 
		{2 \, \gamma_1^{2} \, R^{2}} \, 
		\frac{\gamma + 1}{\varepsilon^{2h} (R - \rho)^2 (1 - \varepsilon)^2} \, \sigma^2 \omega^2 \, 
		\left(\frac{\lambda (A_{h-1})}{\lambda (B_R)}\right)^{1 + \frac{\kappa - 1}{\kappa}}\, .
\end{align*}
We can conclude similarly as before using Lemma \ref{giusti} and provided that
$$
\displaylines{
\hfill
\frac{\lambda (A_0)}{\lambda (B_R)} \leqslant \overline{\nu}^{} =
\left[
\frac{(1-a)^2 \, (1 - \varepsilon)^2 \, \varepsilon^6 \, (R - \rho)^2}
	{\gamma_1^2  R^2 \, (2 \gamma + 2)} \right]^{\frac{\kappa}{\kappa-1}} \, 
	\, \varepsilon^{\frac{8 \kappa^2}{(\kappa - 1)^2}} \, .
\hfill\llap{$\square$}}
$$
\fine

\bprop
\label{prop-DeGiorgi2}
Consider three points $(x^{\diamond}\!, t^{\diamond}), (x^{\circ}\!, t^{\circ}), (x^{\star}\!, t^{\star}) \in \Omega \times (0,T)$ and
$r \in (0,R)$. Suppose
$Q_{R}^{\upbeta^{\diamond}, \texttt{\,>}} (x^{\diamond}\!, t^{\diamond})$
$Q_{R}^{\upbeta^{\circ}, \texttt{\,<}} (x^{\circ}\!, t^{\circ})$
$Q_{R}^{s_1^{\star} , s_2^{\star}} (x^{\star}\!, t^{\star})$
are contained in $\Omega \times (0,T)$.
Then for every choice of $\theta^{\diamond}, \theta^{\circ} \in (0,1)$ and $a, \sigma \in (0,1)$ there are \\
$\underline{\nu}^{\diamond} \in (0,1)$, depending only on $\kappa, \gamma_1, \gamma, a$, $\theta^{\diamond}$, $\upbeta^{\diamond}$, \\
$\underline{\nu}^{\circ} \in (0,1)$, depending only on $\kappa, \gamma_1, \gamma, a$, $\theta^{\circ}$, $\upbeta^{\circ}$,\\
$\underline{\nu}^{\star} \in (0,1)$, depending only on $\kappa, \gamma_1, \gamma, a, (R - r)/R$, $\max\{ 1, 1/\upbeta^{\star} \}$, \\
$\underline{\nu} \in (0,1)$, depending only on $\kappa, \gamma_1, \gamma, a, (R - r)/R$, \\
such that for every $u \in DG_-(\Omega, T, \mu, \lambda, \gamma)$ and fixed $\underline{m}, \omega$ satisfying \\ [0.5em]
$i \, )$ $\underline{m} \leqslant \inf_{Q_{R; R, 0}^{\upbeta^{\diamond},+} (x^{\diamond} \!, t^{\diamond})} u, \hskip10pt \omega \geqslant 
\osc_{Q_{R; R, 0}^{\upbeta^{\diamond}+} (x^{\diamond} \!, t^{\diamond})} u$
we have that if $\mu_+ (B_{r}) > 0$ and
\begin{align*}
\frac{ M_+ (A_{0}^+)}{|M|_{\Lambda} (Q_{R}^{\upbeta^{\diamond}, \texttt{\,>}} (x^{\diamond}\!, t^{\diamond}))} +
\frac{ \Lambda_+ (A_{0}^+)}{\Lambda (Q_{R}^{\upbeta^{\diamond}, \texttt{\,>}} (x^{\diamond}\!, t^{\diamond}))} \leqslant \underline{\nu}^{\diamond} ,
\end{align*}
$\hskip8pt$ where $A_0^+ = \{ (x,t) \in Q_{R; R, 0}^{\upbeta^{\diamond}, +} (x^{\diamond}\!, t^{\diamond}) \, | 
\, u(x,t) < \underline{m} + \sigma \omega \}$, then
$$
u(x,t) \geqslant \underline{m} + a \, \sigma \, \omega  \hskip30pt 
	\text{for a.e. }	(x,t) \in Q_{R; r, \theta^{\diamond}}^{\upbeta^{\diamond}, +} (x^{\diamond}\!, t^{\diamond}) \, ;
$$
$ii \, )$ $\underline{m} \leqslant \inf_{Q_{R; R, 0}^{\circ, -} (x^{\circ} \!, t^{\circ})} u, 
	\hskip10pt \omega \geqslant \osc_{Q_{R; R, 0}^{\circ, -} (x^{\circ} \!, t^{\circ})} u$
we have that if $\mu_- (B_{r}) > 0$ and
\begin{align*}
\frac{ M_- (A_{0}^-)}{|M|_{\Lambda} (Q_{R}^{\circ, \texttt{\,<}} (x^{\circ}\!, t^{\circ}))} +
\frac{ \Lambda_- (A_{0}^-)}{\Lambda (Q_{R}^{\circ, \texttt{\,<}} (x^{\circ}\!, t^{\circ}))} \leqslant \underline{\nu}^{\circ} ,
\end{align*}
$\hskip8pt$ where $A_0^- = \{ (x,t) \in Q_{R; R, 0}^{\circ, -} (x^{\circ} \!, t^{\circ}) \, | \, u(x,t) < \underline{m} + \sigma \omega \}$, then
$$
u(x,t) \geqslant \underline{m} + a \, \sigma \, \omega  \hskip30pt \text{for a.e. }	(x,t) \in Q_{R; r, \theta^{\circ}}^{\circ, -}  (x^{\circ}\!, t^{\circ}) \, ;
$$
$iii \, )$ $\underline{m} \leqslant \inf_{Q_{R}^{s_1^{\star} , s_2^{\star}} (x^{\star}\!, t^{\star})} u, \hskip10pt 
	\om \geqslant \osc_{Q_{R}^{s_1^{\star} , s_2^{\star}} (x^{\star}\!, t^{\star})} u$
we have that if $\lambda_0 (B_{r}) > 0$ and
\begin{align*}
\Lambda_0 (A_{0}^0) \leqslant \underline{\nu}^{\star} \, \Lambda (Q_{R}^{s_1^{\star} , s_2^{\star}} (x^{\star}\!, t^{\star}))
\end{align*}
$\hskip8pt$ where $A_0^0 = \{ (x,t) \in Q_{R; R, s_1^{\star}, s_2^{\star}}^{0} (x^{\star} \!, t^{\star}) \, | \, u(x,t) < \underline{m} + \sigma \omega \}$,
then
$$
u(x,t) \geqslant \underline{m} + a \, \sigma \, \omega  \hskip30pt \text{for a.e. }	(x,t) \in Q_{R; r, s_1^{\star}, s_2^{\star}}^{0} (x^{\star} \!, t^{\star}) \, ;
$$
$iv \, )$ $\underline{m} \leqslant \inf_{B_R(x^{\star})} u (\cdot, t), \hskip10pt \om \geqslant \osc_{B_R(x^{\star})} u (\cdot, t)$
we have that if $B_R(x^{\star}) \subset \Omega_0$ and
\begin{align*}
\lambda \big(\{ x \in B_{R} (x^{\star}) \, | \, u(x,t) < \underline{m} + \sigma \omega \} \big) \leqslant \underline{\nu} \ \lambda (B_R (x^{\star}) )
\end{align*}
then
$$
u(x,t) \geqslant \underline{m} + a \, \sigma \, \omega  \hskip30pt \text{for a.e. }	x \in B_{r} (x^{\star})
$$
for a.e. $t \in (0,T)$.
\eprop

%

\ \\
\noindent
We now need some results which are preparatory for one fundamental step in view of proving the Harnack's inequality, 
Lemma \ref{esp_positivita}, which is usually referred to as {\em expansion of positivity}. \\
\ \\
We define, for a fixed point $(\bar{y}, \bar{s}) \in \Om \times (0,T)$ and a fixed $h > 0$, the sets
\begin{align}
\label{Aacca_ro2}
& A_{h,\rho}^+ (\bar{y}, \bar{s})  = \{ x \in B_{\rho}^+(\bar{y}) \, | \, u(x,\bar{s}) < h \} \, ,	\nonumber	\\
& A_{h,\rho}^- (\bar{y}, \bar{s})  = \{ x \in B_{\rho}^-(\bar{y}) \, | \, u(x,\bar{s}) < h \} \, ,				\\
& A_{h,\rho}^0 (\bar{y}, \bar{s})  = \{ x \in B_{\rho}^0(\bar{y}) \, | \, u(x,\bar{s}) < h \} \, .	\nonumber
\end{align}

\boss
\label{Aacca_ro}
Observe that the condition $u(x,\bar{s}) \geqslant h$ for every
$x \in B_{\rho}(\bar{y})$ implies that
$A_{h,4 \rho}(\bar{y}, \bar{s}) \subset B_{4\rho}(\bar{y}) \setminus B_{\rho}(\bar{y})$
and then in particular, if $\omega$ is a doubling weight ($c_{\omega}$ denotes the doubling constant of $\omega$), 
one has
$$
\omega \big(A_{h,4\rho}(x^{\ast}, t^{\ast})\big) \leqslant \, \omega \big(B_{4\rho}(x^{\ast}) \setminus B_{\rho}(x^{\ast}) \big) 
	\leqslant \left(1 - c_{\omega}^{-2} \right) \, \omega \big(B_{4\rho}(x^{\ast}) \big)\, .
$$
In our situation this holds for $|\mu|_{\lambda}$, thanks to \eqref{doublingmula}, but also for $\mu_+$, $\mu_-$, $\lambda_0$ thanks to the
assumtpion (H.4).
\eoss

\begin{lemma}
\label{lemma1}
Given  $(x^{\ast} \!, t^{\ast})$ such that $B_{4\rho}(x^{\ast}) \subset \Omega$ then \\ [0.3em]
$i \, )$
if $\lambda_0 \big(B_{4\rho}(x^{\ast}) \big) > \lambda_0 \big(B_{\rho}(x^{\ast}) \big) > 0$ there exists $\eta \in (0,1)$, depending only on $\q$,
such that for every $\bar{t} \in (0,T)$
we have that, given $h > 0$ and $u \geqslant 0$ belonging to $DG(\Omega, T, \mu, \la, \gamma)$ for which the following holds
$$
u(x,\bar{t}) \geqslant h \hskip20pt \text{a.e. in } B_{\rho}^0(x^{\ast}) ,
$$
then
\begin{align*}
\lambda_0 (A_{\eta h,4\rho}^0 (x^{\ast} \!, \bar{t})) & \leqslant
	\left(1 - \frac{1}{2} \frac{1}{\q^2} \right) \, \lambda_0 \big(B_{4\rho}^0(x^{\ast}) \big) .
\end{align*}
If $B_{4\rho}(x^{\ast}) \times [t^{\ast} - \upbeta \, h(x^{\ast} \!, 4\rho) \, \rho^2, t^{\ast} + \upbeta \, h(x^{\ast} \!, 4\rho) \, \rho^2]
\subset \Omega \times (0,T)$ with $\upbeta \in (0,16]$ then: \\ [0.3em]
$ii \, )$
if $\mu_+ \big(B_{4\rho}(x^{\ast}) \big) > \mu_+ \big(B_{\rho}(x^{\ast}) \big) > 0$
there exists $\eta \in (0,1)$, depending only on $\gamma, \q$, and there exists
$\tilde\upbeta \in (0, \upbeta]$, depending only on $\gamma$ and $\upbeta$, such that, given
$h > 0$ and $u \geqslant 0$ belonging to $DG(\Omega, T, \mu, \la, \gamma)$ for which the following holds
$$
u(x,t^{\ast}) \geqslant h \hskip20pt \text{a.e. in } B_{\rho}^+(x^{\ast}) ,
$$
then for every $t \in [t^{\ast}, t^{\ast} + \tilde\upbeta \, h(x^{\ast} \!, 4\rho) \, \rho^2]$
\begin{align*}
\mu_+ (A_{\eta h,4\rho}^+ (x^{\ast} \!, t)) & \leqslant  
	\left(1 - \frac{1}{2} \frac{1}{\q^2} \right) \, \mu_+ \big(B_{4\rho}^+(x^{\ast}) \big) ;
\end{align*}
$iii \, )$
if $\mu_- \big(B_{4\rho}(x^{\ast}) \big) > \mu_- \big(B_{\rho}(x^{\ast}) \big) > 0$
there exists $\eta \in (0,1)$, depending only on $\gamma, \q$, and there exists
$\tilde\upbeta \in (0, \upbeta]$, depending only on $\gamma$ and $\upbeta$, such that, given
$h > 0$ and $u \geqslant 0$ belonging to $DG(\Omega, T, \mu, \la, \gamma)$ for which the following holds
$$
u(x,t^{\ast}) \geqslant h \hskip20pt \text{a.e. in } B_{\rho}^-(x^{\ast}) ,
$$
then for every $t \in [t^{\ast} - \tilde\upbeta \, h(x^{\ast} \!, 4\rho) \, \rho^2, t^{\ast}]$
\begin{align*}
\mu_- (A_{\eta h,4\rho}^- (x^{\ast} \!, s)) & \leqslant
	\left(1 - \frac{1}{2} \frac{1}{\q^2} \right) \, \mu_- \big(B_{4\rho}^-(x^{\ast}) \big) ;
\end{align*}
$iv \, )$
there exist $\eta \in (0,1)$, depending only on $\gamma$ and $\q$, and there exists
$\tilde\upbeta \in (0, \upbeta]$, depending only on $\gamma$ and $\upbeta$, such that, given
$h > 0$ and $u \geqslant 0$ belonging to $DG(\Omega, T, \mu, \la, \gamma)$ for which the following holds
$$
u(x,t^{\ast}) \geqslant h \hskip20pt \text{a.e. in } B_{\rho}(x^{\ast}) ,
$$
then 
\begin{align*}
|\mu|_{\lambda} \big( A_{\eta h,4\rho}^+ (x^{\ast} \!, t) \cup A_{\eta h,4\rho}^- (x^{\ast} \!, s) \cup A_{\eta h,4\rho}^0 (x^{\ast} \!, t^{\ast}) \big)
	\leqslant  \left(1 - \frac{1}{2} \frac{1}{\q^2} \right) \, |\mu|_{\lambda} \big(B_{4\rho}(x^{\ast}) \big)
\end{align*}
for every $t \in [t^{\ast}, t^{\ast} + \tilde\upbeta \, h(x^{\ast} \!, 4\rho) \, \rho^2]$ and
$s \in [t^{\ast} - \tilde\upbeta \, h(x^{\ast} \!, 4\rho)\, \rho^2, t^{\ast}]$.
\end{lemma}

\noindent
\dimo
First we prove point $ii \, )$. Consider 
$s_1 = t^{\ast} - \upbeta h(x^{\ast}\!, 4 \rho) \rho^2$, $s_2 = t^{\ast} + \upbeta h(x^{\ast}\!, 4 \rho) \rho^2$.
Apply the energy estimate \eqref{DGgamma+_1} to the function $(u - h)_-$ with $x_0 = x^{\ast}$, $t_0 = t^{\ast}$, $r = 4 \rho (1-\sigma)$
for an arbitrary $\sigma \in (0,1)$, $R = \tilde{r} = 4 \rho$, $\varepsilon = 0$. With this choice we have $\tilde{r} - r = 4 \rho \sigma$.
Then we get
\begin{align*}
\sup_{t \in (t^{\ast} \!, s_2)} & \int_{B_{4 \rho (1-\sigma)}^+(x^{\ast})} (u - h)_-^2 (x,t) \mu_+ (x) dx	 \leqslant 				\\
\leqslant & \int_{B_{4\rho}^+(x^{\ast})} (u - h)_-^2 (x,t^{\ast}) \mu_+ (x) dx	 \ +
		\sup_{t \in (t^{\ast}\!, s_2)} \int_{I^{4\rho, 4\rho\sigma}_+} (u-h)_-^2 (x,t) \mu_-(x) \, dx +						\\
& + \, \frac{\gamma}{(4 \rho \sigma)^2} \int_{t^{\ast}}^{s_2} \!\!\!\! \int_{B_{4\rho}^+(x^{\ast}) \cup I^{4\rho, 4\rho\sigma}_+}  (u - h)_-^2\, \lambda \, dx ds .
\end{align*}
Now, in addition to this inequality, we use the two following inequalities: first that in a set $A_{\eta h, r}$ we have that $(u-h)_- \geqslant (1-\eta) h$;
moreover, since $u \mau 0$, $(u - h)_- \miu h$.
Then, using also Remark \ref{Aacca_ro}, we get for every $t \in [t^{\ast}, s_2]$
\begin{align*}
(1-\eta)^2 h^2 & \mu_+\big( A_{\eta h,4\rho (1-\sigma)}^+ (x^{\ast} \!, t) \big) 				\leqslant 				\\
\leqslant & \int_{A_{\eta h,4\rho (1-\sigma)}^+ (x^{\ast} \!, t)} (u - h)_-^2 (x,t) \mu_+ (x) dx 	\leqslant				\\
\leqslant & \int_{B_{4 \rho (1-\sigma)}^+(x^{\ast})} (u - h)_-^2 (x,t) \mu_+ (x) dx +			\leqslant				\\
\leqslant & \ h^2 \, \mu_+ \big(B_{4\rho}(x^{\ast}) \setminus B_{\rho}(x^{\ast}) \big) +
		 h^2 \mu_-({I^{4\rho, 4\rho\sigma}_+}) +												
		\frac{\gamma h^2}{(4\rho\sigma)^2} \, \Lambda \big( (B_{4\rho}^+(x^{\ast}) \cup I^{4\rho, 4\rho\sigma}_+) \times (t^{\ast}, s_2) \big) .
\end{align*}
Using the following decomposition
\begin{align*}
A_{\eta h, 4\rho}^+(x^{\ast} \!, t) & =
A_{\eta h,4\rho (1-\sigma)}^+ (x^{\ast} \!, t) \cup \big\{x \in B_{4\rho}^+ (x^{\ast}) \setminus B_{4\rho(1-\sigma)}^+(x^{\ast}) \, \big| \, u(x,t) < \eta h \big \}	,
\end{align*}
and then the last estimate we get
\begin{align}
(1-\eta)^2 & \mu_+ \big( A_{\eta h,4\rho}^+ (x^{\ast} \!, t) \big)	\leqslant												\nonumber	\\
\leqslant & \ (1-\eta)^2 \Big[ \mu_+\big( A_{\eta h,4\rho (1-\sigma)}^+ (x^{\ast} \!, t) \big)
		+ \mu_+ \big(  B_{4\rho} (x^{\ast}) \setminus B_{4\rho(1-\sigma)}(x^{\ast}) \big) \Big] \leqslant						\nonumber	\\
\label{oralareplichiamo}
\leqslant & \ \mu_+ \big(B_{4\rho}(x^{\ast}) \setminus B_{\rho}(x^{\ast}) \big) +
		\mu_-({I^{4\rho, 4\rho\sigma}_+}) +												
		\frac{\gamma}{(4\rho\sigma)^2} \, \Lambda \big( (B_{4\rho}^+(x^{\ast}) \cup I^{4\rho, 4\rho\sigma}_+) \times (t^{\ast}, s_2) \big) +		\\
	      & + \, (1-\eta)^2 \mu_+ \big(  B_{4\rho} (x^{\ast}) \setminus B_{4\rho(1-\sigma)}(x^{\ast}) \big) .						\nonumber
\end{align}
If the thesis were false we would have that
for every $\tilde{\upbeta} \in (0, \upbeta]$ and $\eta \in (0,1)$ there would be
$\bar{t} \in [t^{\ast}, t^{\ast} + \tilde\upbeta \, h(x^{\ast} \!, 4\rho) \, \rho^2]$  such that
$$
\left(1 - \frac{1}{2} \frac{1}{\q^2} \right) \, \mu_+ \big(B_{4\rho}^+(x^{\ast}) \big)
	< \mu_+ (A_{\eta h,4\rho}^+ (x^{\ast} \!, \bar{t}))
$$
and then
\begin{align*}
(1-\eta)^2 & \left(1 - \frac{1}{2} \frac{1}{\q^2} \right) \, \mu_+ \big(B_{4\rho}^+(x^{\ast}) \big) <	\\
	<	& \ \mu_+ \big(B_{4\rho}(x^{\ast}) \setminus B_{\rho}(x^{\ast}) \big) +
		  \mu_- ({I^{4\rho, 4\rho\sigma}_+}) +												
		\frac{\gamma}{(4\rho\sigma)^2} \, \Lambda \big( (B_{4\rho}^+(x^{\ast}) \cup I^{4\rho, 4\rho\sigma}_+) \times (t^{\ast}, s_2) \big) +		\\
      & + \, (1-\eta)^2 \mu_+ \big(  B_{4\rho} (x^{\ast}) \setminus B_{4\rho(1-\sigma)}(x^{\ast}) \big) .			
\end{align*}
Then taking, for instance, $\upbeta = \sigma^3$, letting $\sigma$ and $\eta$ go to zero we would find the contradiction
(and here is needed $\mu_+ \big(B_{4\rho}(x^{\ast}) \big) > \mu_+ \big(B_{\rho}(x^{\ast}) \big) > 0$)
\begin{gather*}
\left(1 - \frac{1}{2} \frac{1}{\q^2} \right) \, \mu_+ \big(B_{4\rho}^+(x^{\ast}) \big)
	\leqslant \mu_+ \big(B_{4\rho}(x^{\ast}) \setminus B_{\rho}(x^{\ast}) \big)				\\
\Downarrow																	\\
2 \mu_+ \big( B_{\rho}(x^{\ast}) \big) \leqslant \mu_+ \big( B_{\rho}(x^{\ast}) \big) \, .
\end{gather*}
In a way analogous to \eqref{oralareplichiamo} we can derive  for every $s \in [s_1, t^{\ast}]$
\begin{align}
(1-\eta)^2 & \mu_- \big( A_{\eta h,4\rho}^- (x^{\ast} \!, s) \big)	\leqslant												\nonumber	\\
\leqslant & \ (1-\eta)^2 \Big[ \mu_- \big( A_{\eta h,4\rho (1-\sigma)}^- (x^{\ast} \!, s) \big)
		+ \mu_- \big(  B_{4\rho} (x^{\ast}) \setminus B_{4\rho(1-\sigma)}(x^{\ast}) \big) \Big] \leqslant						\nonumber	\\
\label{replica1}
\leqslant & \ \mu_- \big(B_{4\rho}(x^{\ast}) \setminus B_{\rho}(x^{\ast}) \big) +
		\mu_+ ({I^{4\rho, 4\rho\sigma}_-}) +												
		\frac{\gamma}{(4\rho\sigma)^2} \, \Lambda \big( (B_{4\rho}^-(x^{\ast}) \cup I^{4\rho, 4\rho\sigma}_-) \times (s_1, t^{\ast}) \big) +		\\
	      & + \, (1-\eta)^2 \mu_- \big(  B_{4\rho} (x^{\ast}) \setminus B_{4\rho(1-\sigma)}(x^{\ast}) \big)						\nonumber
\end{align}
by which, again by contradiction, we prove point $iii \, )$. \\ [0.3em]
Point $i\, )$ is quite immediate.
Since $(u - h)_-(x,\bar{t}) \geqslant (1 - \eta) h$ in $A_{\eta h,4\rho}^0(x^{\ast} \!, \bar{t})$ we immediately get
\begin{align*}
(1-\eta)^2 h^2 & \lambda_0\big( A_{\eta h,4\rho (1-\sigma)}^0 (x^{\ast} \!, \bar{t}) \big) 			\leqslant 				\\
\leqslant & \int_{A_{\eta h,4\rho}^0 (x^{\ast} \!, \bar{t})} (u - h)_-^2 (x,\bar{t}) \lambda_0 (x) dx 	\leqslant				\\
\leqslant & \int_{B_{4 \rho}^0(x^{\ast})} (u - h)_-^2 (x, \bar{t}) \lambda_0 (x) dx			\leqslant
	h^2 \, \lambda_0 \big(B_{4\rho}(x^{\ast}) \setminus B_{\rho}(x^{\ast}) \big)										
\end{align*}
that is
$$
(1-\eta)^2 \lambda_0\big( A_{\eta h,4\rho (1-\sigma)}^0 (x^{\ast} \!, \bar{t}) \big)	\leqslant
	\lambda_0 \big(B_{4\rho}(x^{\ast}) \setminus B_{\rho}(x^{\ast}) \big) \leqslant 
	\left(1 - \frac{1}{\q^2} \right) \lambda_0 \big(B_{4\rho}(x^{\ast}) \big)
$$
and then $\eta$ is easily found. \\ [0.3em]
Point $iv \, )$ is obtained simply summing and rearranging the previous inequalities.
\finedimo

\begin{lemma}
\label{lemma2}
Consider $\upbeta \in (0,16]$ and $(x^{\ast} \!, t^{\ast})$ such that
$B_{5\rho}(x^{\ast}) \times [t^{\ast} - \upbeta \, h(x^{\ast} \!, 4\rho) \, \rho^2, t^{\ast} + \upbeta \, h(x^{\ast} \!, 4\rho) \, \rho^2]
\subset \Omega \times (0,T)$, consider $\eta$ and $\tilde\upbeta$ to be the values determined in Lemma $\ref{lemma1}$.
Consider $\upkappa$ and $\uptau$ the constants appearing in \eqref{carlettomio}.
Consider $u \geqslant 0$ in $DG(\Omega, T, \mu, \la, \gamma)$, $h > 0$. \\ [0.3em]
$i \, )$ If $\mu_+ \big(B_{4\rho}(x^{\ast}) \big) > \mu_+ \big(B_{\rho}(x^{\ast}) \big) > 0$ and
$u(\cdot,t^{\ast}) \geqslant h$ a.e. in $B_{\rho}^+(x^{\ast})$ \\
then for every $\e > 0$ there exists $\eta_1 \in (0,\eta)$,
$\eta_1$ depending only on $\gamma_1 , \gamma , \q , \epsilon  , \eta, \tilde\upbeta$
such that
\begin{align*}
M_+ \Big( \{u < \eta_1 h \} \, \cap \, &
	\Big[B_{4\rho}^+ (x^{\ast}) \times (t^{\ast}, t^{\ast} + \tilde\upbeta \, \rho^2 h(x^{\ast}\!, 4\rho) ) \Big] \Big)  \leqslant						\\
&	\leqslant \epsilon \ |M|_{\Lambda} \Big( B_{4\rho} (x^{\ast}) \times (t^{\ast}, t^{\ast} + \tilde\upbeta \, \rho^2 h(x^{\ast} \!, 4\rho) ) \Big) ,		\\
\Lambda_+ \Big( \{u < \eta_1 h \} \, \cap \, &
	\Big[B_{4\rho}^+ (x^{\ast}) \times (t^{\ast}, t^{\ast} + \tilde\upbeta \, \rho^2 h(x^{\ast}\!, 4\rho) ) \Big] \Big)  \leqslant						\\
&	\leqslant \upkappa \, \epsilon^{\uptau} \ \Lambda \Big( B_{4\rho} (x^{\ast}) \times (t^{\ast}, t^{\ast} + \tilde\upbeta \, \rho^2 h(x^{\ast} \!, 4\rho) ) \Big) ;
\end{align*}
$ii \, )$ if $\mu_- \big(B_{4\rho}(x^{\ast}) \big) > \mu_- \big(B_{\rho}(x^{\ast}) \big) > 0$
and $u(\cdot,t^{\ast}) \geqslant h$ a.e. in $B_{\rho}^-(x^{\ast})$ \\
then for every $\e > 0$  
there exists $\eta_1 \in (0,\eta)$, $\eta_1$ depending only on $\gamma_1 , \gamma , \q , \epsilon , \eta, \tilde\upbeta$ such that
\begin{align*}
M_- \Big( \{u < \eta_1 h \} \, \cap \, & 
	\Big[B_{4\rho}^- (x^{\ast}) \times (t^{\ast} - \tilde\upbeta \, \rho^2 h(x^{\ast}\!, 4\rho), t^{\ast} ) \Big] \Big)  \leqslant						\\
&	\leqslant \epsilon \ |M|_{\Lambda} \Big( B_{4\rho} (x^{\ast}) \times (t^{\ast} - \tilde\upbeta \, \rho^2 h(x^{\ast} \!, 4\rho) , t^{\ast} ) \Big) ,		\\
\Lambda_- \Big( \{u < \eta_1 h \} \, \cap \, & 
	\Big[B_{4\rho}^- (x^{\ast}) \times (t^{\ast} - \tilde\upbeta \, \rho^2 h(x^{\ast}\!, 4\rho), t^{\ast} ) \Big] \Big)  \leqslant						\\
&	\leqslant \upkappa \, \epsilon^{\uptau} \ \Lambda \Big( B_{4\rho} (x^{\ast}) \times (t^{\ast} - \tilde\upbeta \, \rho^2 h(x^{\ast} \!, 4\rho) , t^{\ast} ) \Big) ;
\end{align*}
$iii \, )$ consider $\upbeta > 0$ such that
$B_{5\rho}(x^{\ast}) \times [t^{\ast} - \upbeta \, h(x^{\ast} \!, 4\rho) \, \rho^2, t^{\ast} + \upbeta \, h(x^{\ast} \!, 4\rho) \, \rho^2]
\subset \Omega \times (0,T)$:
if $\lambda_0 \big(B_{4 \rho}(x^{\ast}) \big) > \lambda_0 \big(B_{\rho}(x^{\ast}) \big) > 0$
and $u \geqslant h$  a.e. in
$\big( B_{\rho}^0(x^{\ast}) \times (t^{\ast} - \upbeta \, \rho^2 h(x^{\ast}\!, 4\rho), t^{\ast} + \upbeta \, \rho^2 h(x^{\ast}\!, 4\rho) ) \big)$
then for every $\e > 0$  there exists $\eta_1 \in (0,\eta)$,
$\eta_1$ depending only on $\gamma_1 , \gamma , \q , \epsilon , \eta, \upbeta$
such that
\begin{align*}
\Lambda_0 \Big( \{u < \eta_1 h \} \cap
	\Big[B_{4 \rho}^0 (x^{\ast}) \times & (t^{\ast} - \upbeta \, \rho^2 h(x^{\ast}\!, 4\rho), t^{\ast} + \upbeta \, \rho^2 h(x^{\ast}\!, 4\rho) ) \Big] \Big) 
	\leqslant																													\\
& \leqslant \epsilon \ \Lambda 
	\Big( B_{4\rho} (x^{\ast}) \times (t^{\ast} - \upbeta \, \rho^2 h(x^{\ast}\!, 4\rho), t^{\ast} + \upbeta \, \rho^2 h(x^{\ast}\!, 4\rho ) \Big) ;
\end{align*}
$iv \, )$ if $B_{5 \rho} (x^{\ast}) \subset \Omega_0$ and
$u (\cdot, t) \geqslant h$  a.e. in $B_{\rho}(x^{\ast})$
then for every $\e > 0$  there exists $\eta_1 \in (0,\eta)$,
$\eta_1$ depending only on $\gamma_1 , \gamma , \q , \epsilon , \eta$
such that for almost every $t \in [t^{\ast} - \upbeta \, h(x^{\ast} \!, 4\rho) \, \rho^2, t^{\ast} + \upbeta \, h(x^{\ast} \!, 4\rho) \, \rho^2]$
\begin{align*}
\lambda \Big( \{u < \eta_1 h \} \cap \big( B_{4\rho} (x^{\ast}) \times \{ t \} \big) \Big) \leqslant	\epsilon \ \lambda \Big( B_{4\rho} (x^{\ast}) \Big) .
\end{align*}
\end{lemma}
\noindent
\dimo
We first show point $i\, $).
Consider $\tilde\upbeta$ and $\eta$ to be the values determined in Lemma \ref{lemma1}, point $i\, $). For simplicity, by $f$ we will denote the quantity
$$
f (x^{\ast}\!, 4\rho) = h(x^{\ast}\!, 4\rho) \, \rho^2 \, .
$$
Now we consider $m \in \N$, $\tau \in [t^{\ast}\!, t^{\ast} + \tilde\upbeta \, f(x^{\ast}\!,4\rho) ]$
and $\sigma \in [t^{\ast} - \tilde\upbeta \, f(x^{\ast}\!,4\rho), t^{\ast} ]$. 
First of all notice that for every $t^{\ast}$, for every $\tau \in [t^{\ast}\!, t^{\ast} + \tilde\upbeta \, f(x^{\ast}\!,4\rho) ]$ and
$\sigma \in [t^{\ast} - \tilde\upbeta \, f(x^{\ast}\!,4\rho), t^{\ast} ]$
and every $t \in (t^{\ast} - \upalpha \, \rho^2 h(x^{\ast}\!, 4\rho), t^{\ast} + \upalpha \, \rho^2 h(x^{\ast}\!, 4\rho) )$
we derive, using Lemma \ref{lemma1} and since for $m \in \N$ it holds
$A_{\eta h 2^{-m},4\rho}^+ (x^{\ast}\!,\tau) \subset A_{\eta h,4\rho}^+ (x^{\ast}\!,\tau)$,
$A_{\eta h 2^{-m},4\rho}^- (x^{\ast}\!,\sigma) \subset A_{\eta h,4\rho}^- (x^{\ast}\!,\sigma)$,
$A_{\eta h 2^{-m},4\rho}^0 (x^{\ast}\!,t^{\ast}) \subset A_{\eta h,4\rho}^0 (x^{\ast}\!,t^{\ast})$,
that if $\mu_+ \big( B_{\rho} (x^{\ast}) \big) > 0$, $\mu_- \big( B_{\rho} (x^{\ast}) \big) > 0$, $\lambda_0 \big( B_{\rho} (x^{\ast}) \big) > 0$
\begin{equation}
\label{nuovoarrivo}
\begin{array}{l}
{\displaystyle \frac{1}{2 \q^2} } \, \mu_+ \big(B_{4\rho}^+(x^{\ast}) \big) \leqslant 
	\mu_+ \big(B_{4\rho}^+(x^{\ast}) \setminus A_{\eta h,4\rho}^+ (x^{\ast} \!, \tau)\big)
	\leqslant 
	\mu_+ \big(B_{4\rho}^+(x^{\ast}) \setminus A_{\eta h 2^{-m},4\rho}^+ (x^{\ast} \!, \tau)\big) ,		\\	[1em]
{\displaystyle \frac{1}{2 \q^2} } \, \mu_- \big(B_{4\rho}^-(x^{\ast}) \big) \leqslant 
	\mu_- \big(B_{4\rho}^-(x^{\ast}) \setminus A_{\eta h,4\rho}^- (x^{\ast} \!, \sigma)\big) 
	\leqslant 
	\mu_- \big(B_{4\rho}^-(x^{\ast}) \setminus A_{\eta h 2^{-m},4\rho}^- (x^{\ast} \!, \sigma)\big),		\\	[1em]
{\displaystyle \frac{1}{2 \q^2} } \, \lambda_0 \big(B_{4\rho}^0(x^{\ast}) \big) \leqslant 
	\lambda_0 \big(B_{4\rho}^0(x^{\ast}) \setminus A_{\eta h,4\rho}^0 (x^{\ast} \!, t)\big) 
	\leqslant 
	\lambda_0 \big(B_{4\rho}^0(x^{\ast}) \setminus A_{\eta h 2^{-m},4\rho}^0 (x^{\ast} \!, t)\big).
\end{array}
\end{equation}
Again for simplicity, we define (since $x^{\ast}$ is fixed we omit it)
\begin{align*}
A_m^+(\tau) := & A_{\eta h 2^{-m},4\rho}^+ (x^{\ast} \!, \tau) ,		\hskip30pt	
				a_m^{+} :=  \int_{t^{\ast}}^{t^{\ast} + \tilde\upbeta \, f(x^{\ast}\!,4\rho)}	\mu_+(A_m^+(\tau)) \, d \tau ,				\\
A_m^-(\sigma) := & A_{\eta h 2^{-m},4\rho}^- (x^{\ast} \!, \sigma) ,	\hskip30pt	
				a_m^- := 	\int_{t^{\ast} - \tilde\upbeta \, f(x^{\ast}\!,4\rho)}^{t^{\ast}}	\mu_-(A_m^-(\sigma)) \, d \sigma ,				\\
A_m^0(t) := & A_{\eta h 2^{-m},4 \rho}^0 (x^{\ast} \!, t) , \hskip28pt
				a_m^0 :=	\int_{t^{\ast} - \upalpha \, f(x^{\ast}\!,4\rho)}^{t^{\ast} + \upalpha \, f(x^{\ast}\!,4\rho)}
						\lambda_0 (A_m^0(t)) dt	,												 						\\
A_m (t) := & A_{\eta h 2^{-m},4\rho} (x^{\ast} \!, t) ,				\hskip32pt
				B_{4\rho} := B_{4\rho} (x^{\ast})  ,																			\\
d_m^{\texttt{\,>}} := & \int_{t^{\ast}}^{t^{\ast} + \tilde\upbeta \, f(x^{\ast}\!,4\rho)} \lambda (A_m(t)) d t ,
														\hskip10pt
				d_m^{\texttt{\,<}} := \int_{t^{\ast} - \tilde\upbeta \, f(x^{\ast}\!,4\rho)}^{t^{\ast}} \lambda (A_m(t)) d t ,					\\
& \hskip20pt d_m := \int_{t^{\ast} - \upalpha \, f(x^{\ast}\!,4\rho)}^{t^{\ast} + \upalpha \, f(x^{\ast}\!,4\rho)} \lambda (A_m(t)) d t .
\end{align*}
First we prove point $i \,)$.
Now we estimate from above and from below the quantity
\begin{align*}
\mu_+ \big( B_{4\rho}^+ \setminus A_{m-1}^+ (\tau) \big) \int_{B_{4\rho}} \Big(u - \frac{\eta h}{ 2^m}\Big)_- (x, \tau) \mu_+ (x) dx
\end{align*}
Using also \eqref{nuovoarrivo} we get that for every $\tau \in [t^{\ast}\!, t^{\ast} + \tilde\upbeta \, f(x^{\ast}\!,4\rho) ]$
\begin{align*}
& \frac{1}{2 \q^2} \, \mu_+ \big(B_{4\rho}^+ \big) \frac{\eta h}{2^{m+1}} \mu_+ \big( A_{m+1}^+ (\tau) \big) \leqslant							\\
& \ \leqslant \, \mu_+ \big( B_{4\rho}^+ \setminus A_{m-1}^+ (\tau) \big) \frac{\eta h}{2^{m+1}} \mu_+ \big( A_{m+1}^+ (\tau) \big) \leqslant 			\\
& \ \leqslant \, \mu_+ \big( B_{4\rho}^+ \setminus A_{m-1}^+ (\tau) \big) \int_{B_{4\rho}} \Big(u - \frac{\eta h}{ 2^m}\Big)_- (x, \tau) \mu_+ (x) dx =	\\
& \ \leqslant \, \mu_+ \big( B_{4\rho}^+ \setminus A_{m-1}^+ (\tau) \big) \frac{\eta h}{2^{m}} \mu_+ \big( A_{m}^+ (\tau) \big)
\end{align*}
that is we get that for every $\tau \in [t^{\ast}\!, t^{\ast} + \tilde\upbeta \, f(x^{\ast}\!,4\rho) ]$
\begin{equation}
\label{saraclara}
\begin{array}{l}
{\displaystyle \frac{1}{2 \q^2} \, \mu_+ \big(B_{4\rho}^+ \big) \frac{\eta h}{2^{m+1}} \mu_+ \big( A_{m+1}^+ (\tau) \big) \leqslant }
{\displaystyle \ \mu_+ \big( B_{4\rho}^+ \setminus A_{m-1}^+ (\tau) \big) \frac{\eta h}{2^{m}} \mu_+ \big( A_{m}^+ (\tau) \big) } \, .
\end{array}
\end{equation}
Now to estimate the right hand side of \eqref{saraclara}
we use Lemma \ref{lemma2.2} in the ball $B_{4\rho}(x^{\ast})$ with $k = \eta h /2^{m}$, $l = \eta h /2^{m-1}$, $q=1$, $p \in (1,2)$ arbitrary,
$\omega = \lambda$, $\nu = |\mu|_{\lambda}$ ($\bar\nu = \mu_+$) we get for every $\tau \in (t^{\ast}\!, t^{\ast} + \tilde\upbeta f(x^{\ast}\!,4\rho))$
\begin{align*}
& {\displaystyle \mu_+ \big( B_{4\rho}^+ \setminus A_{m-1}^+ (\tau) \big) \frac{\eta h}{2^{m}} \mu_+ \big( A_{m}^+ (\tau) \big)	} \leqslant 				\\
& {\displaystyle \qquad \leqslant  8 \, \gamma_1 \, \rho 
		\ \frac{\mu_+ (B_{4\rho}^+) \, |\mu|_{\lambda} (B_{4\rho})}{(\lambda (B_{4\rho}))^{1/p}} \cdot \left(\int_{{A_{m-1}(\tau) \setminus A_m(\tau)}} 
								\!\!\!\!\!\!\!\!\!\!\!\!\!\!\!\!\! |D u|^p (x,\tau) \, \la  \, dx \right)^{1/p} }
\end{align*}
By this last inequality and \eqref{saraclara} and integrating in time between $t^{\ast}$ and $t^{\ast} + \tilde\upbeta \, f(x^{\ast}\!,4\rho)$ we get
\begin{align}
\label{saraclara2}
\frac{1}{2 \q^2} \, & \frac{\eta h}{2^{m+1}} a_{m+1}^+ \leqslant														\nonumber	\\
& {\displaystyle \leqslant  8  \, \gamma_1 \, \rho
		\ \frac{|\mu|_{\lambda} (B_{4\rho})}{(\lambda (B_{4\rho}))^{1/p}} \cdot 
		\int_{t^{\ast}}^{t^{\ast} + \tilde\upbeta \, f(x^{\ast}\!,4\rho)} \left(\int_{{A_{m-1}(\tau) \setminus A_m(\tau)}} 
								\!\!\!\!\!\!\!\!\!\!\!\!\!\!\!\!\! |D u|^p (x,\tau) \, \la  \, dx \right)^{1/p} } d\tau \leqslant		\nonumber	\\
& {\displaystyle \leqslant  8 \, \gamma_1 \, \rho
		\ \frac{|\mu|_{\lambda} (B_{4\rho})}{(\lambda (B_{4\rho}))^{1/p}} \cdot 
		\left( \int_{t^{\ast}}^{t^{\ast} + \tilde\upbeta \, f(x^{\ast}\!,4\rho)} \int_{{A_{m-1}(\tau) \setminus A_m(\tau)}} 
								\!\!\!\!\!\!\!\!\!\!\!\!\!\!\!\!\! |D u|^p (x,\tau) \, \la  \, dx d\tau \right)^{1/p} }
								 \big( \tilde\upbeta \, f(x^{\ast}\!,4\rho) \big)^{\frac{p-1}{p}} \leqslant					\nonumber	\\
& {\displaystyle \leqslant  8 \, \gamma_1 \, \rho
		\ \frac{|\mu|_{\lambda} (B_{4\rho})}{(\lambda (B_{4\rho}))^{1/p}} \big( \tilde\upbeta \, f(x^{\ast}\!,4\rho) \big)^{\frac{p-1}{p}} \cdot
		\left( \int_{t^{\ast}}^{t^{\ast} + \tilde\upbeta \, f(x^{\ast}\!,4\rho)} 
		\big[ \lambda (A_{m-1} (\tau)) - \lambda (A_m (\tau)) \big]\, d \tau	\right)^{\frac{2-p}{2p}} }	\cdot							\\
& {\displaystyle \hskip100pt \cdot \left( \int_{t^{\ast}}^{t^{\ast} + \tilde\upbeta \, f(x^{\ast}\!,4\rho)} \!\!\! \int_{B_{4\rho}} 
			\Big|D \Big(u - \frac{\eta h}{ 2^{m-1}}\Big)_-\Big|^2 (x,\tau) \, \la  \, dx d\tau \right)^{1/2} }.					\nonumber
\end{align}
Now we want to estimate the term in the right hand side involving the gradient of $u - \frac{\eta h}{ 2^m}$ and to do this
we apply the energy estimates \eqref{DGgamma+}, \eqref{DGgamma-}, \eqref{DGgamma0} in some suitable subsets of 
$$
B_{5\rho} (x^{\ast}) \times (t^{\ast}\! - \tilde\upbeta \, f(x^{\ast}\!, 4\rho), t^{\ast} + \tilde\upbeta \, f(x^{\ast}\!, 4\rho) ) .
$$
to estimate the quantity
$
{\displaystyle
\int_{t^{\ast}}^{t^{\ast} + \tilde\upbeta \, f(x^{\ast}\!,4\rho)} \!\!\! \int_{B_{4\rho}} \Big|D \Big(u - \frac{\eta h}{ 2^{m-1}}\Big)_-\Big|^2 (x,\tau) \, \la  \, dx d\tau  .
}
$
First we estimate, taking in \eqref{DGgamma+} $t_0 = t^{\ast} - \tilde\upbeta f(x^{\ast} \!, 4 \rho)$, $s_2 = t^{\ast} + \tilde\upbeta f(x^{\ast} \!, 4 \rho)$,
$R = \tilde{r} = 5\rho$, $r = 4\rho$, $\varepsilon = 0$, $\tilde\theta = 0$ and
$\theta = \tilde\upbeta \, \frac{16}{25} \frac{h(x^{\ast}\!\!,4\rho)}{h(x^{\ast}\!\!,5\rho)}$
in \eqref{DGgamma+},
\begin{align*}
\int_{t^{\ast}}^{t^{\ast} + \tilde\upbeta \, f(x^{\ast}\!,4\rho)} & \!\!\! \int_{B_{4\rho}^+}
		\Big|D \Big(u - \frac{\eta h}{ 2^{m-1}}\Big)_-\Big|^2 (x,\tau) \, \la  \, dx d\tau  \leqslant												\\
\leqslant \gamma & \Bigg[ 
		\int_{I_{5\rho, \rho}^+}  \Big(u - \frac{\eta h}{ 2^{m-1}}\Big)_-^2 \big(x,t^{\ast} - \tilde\upbeta f(x^{\ast} \!, 4 \rho)\big) \mu_+(x) \, dx +		\\
& \hskip25pt + \sup_{t \in ( t^{\ast}\!,  \, t^{\ast}\! + \tilde\upbeta \, f(x^{\ast}\!, 4\rho))} 
						\int_{I^{5\rho, \rho}_+}  \Big(u - \frac{\eta h}{ 2^{m-1}}\Big)_-^2 (x,t) \mu_-(x) \, dx	 +							\\
& \hskip50pt + \, \frac{1}{\rho^2}
		\iint_{Q_{5\rho;5\rho, 0}^{+,\rho}} 
		 \Big(u - \frac{\eta h}{ 2^{m-1}}\Big)_-^2(x,\tau)\, \left( \frac{\mu_+}{h(x^{\ast}\!,5\rho)} + \la \right) \, dx d\tau  \Bigg] 	\leqslant			\\
\leqslant \gamma & \Bigg[  \left( \frac{\eta h}{ 2^{m-1}} \right)^2 \, \mu_+ (I_{5\rho, \rho}^+) +
		\left( \frac{\eta h}{ 2^{m-1}}\right)^2 \, \mu_- (I_+^{5\rho, \rho}) +															\\
& \hskip50pt + \left( \frac{\eta h}{ 2^{m-1}} \right)^2 \, \frac{1}{\rho^2} \, 
		\left( \frac{M_+}{h(x^{\ast}\!,5\rho)} + \Lambda \right) ( Q_{5\rho;5\rho, 0}^{+,\rho} ) \Bigg] \leqslant									\\
\leqslant \gamma  & \left( \frac{\eta h}{ 2^{m-1}}\right)^2 \frac{1}{\rho^2} \, \Bigg[ \rho^2 \, |\mu| (B_{5\rho}) + 
		2 \, \lambda (B_{5\rho}) \, 2 \, \tilde\upbeta \, f(x^{\ast}\!,4\rho) \Bigg]  .
\end{align*}
Then taking in \eqref{DGgamma-} $t_0 = t^{\ast} +2 \, \tilde\upbeta f(x^{\ast} \!, 4 \rho)$, $s_1 = t^{\ast}$,
$R = \tilde{r} = 5\rho$, $r = 4\rho$, $\varepsilon = 0$, $\tilde\theta = 0$ and 
$\theta = \tilde\upbeta \, \frac{16}{25} \frac{h(x^{\ast}\!\!,4\rho)}{h(x^{\ast}\!\!,5\rho)}$
in \eqref{DGgamma+},
\begin{align*}
\int_{t^{\ast}}^{t^{\ast} + \tilde\upbeta \, f(x^{\ast}\!,4\rho)} & \!\!\! \int_{B_{4\rho}^-}
		\Big|D \Big(u - \frac{\eta h}{ 2^{m-1}}\Big)_-\Big|^2 (x,\tau) \, \la  \, dx d\tau  \leqslant												\\
\leqslant \gamma & \Bigg[ 
		\int_{I_{5\rho, \rho}^-}  \Big(u - \frac{\eta h}{ 2^{m-1}}\Big)_-^2 \big(x,t^{\ast} + 2 \tilde\upbeta f^+(x^{\ast} \!, 4 \rho)\big) \mu_-(x) \, dx +	\\
& \hskip25pt + \sup_{t \in ( t^{\ast}\!,  \, t^{\ast}\! + \tilde\upbeta \, f(x^{\ast}\!, 4\rho))} 
						\int_{I^{5\rho, \rho}_-}  \Big(u - \frac{\eta h}{ 2^{m-1}}\Big)_-^2 (x,t) \mu_+(x) \, dx	 +							\\
& \hskip50pt + \, \frac{1}{\rho^2}
		\iint_{Q_{5\rho;5\rho, 0}^{-,\rho}} 
		 \Big(u - \frac{\eta h}{ 2^{m-1}}\Big)_-^2(x,\tau)\, \left( \frac{\mu_-}{h(x^{\ast}\!,5\rho)} + \la \right) \, dx d\tau  \Bigg] 	\leqslant			\\
\leqslant \gamma & \Bigg[  \left( \frac{\eta h}{ 2^{m-1}} \right)^2 \, \mu_- (I_{5\rho, \rho}^-) +
		\left( \frac{\eta h}{ 2^{m-1}}\right)^2 \, \mu_+ (I^{5\rho, \rho}_-) +															\\
& \hskip50pt + \left( \frac{\eta h}{ 2^{m-1}} \right)^2 \, \frac{1}{\rho^2} \, 
		\left( \frac{M_-}{h(x^{\ast}\!,5\rho)} + \Lambda \right) ( Q_{5\rho;5\rho, 0}^{-,\rho} ) \Bigg] \leqslant									\\
\leqslant \gamma  & \left( \frac{\eta h}{ 2^{m-1}}\right)^2 \frac{1}{\rho^2} \, \Bigg[ \rho^2 \, |\mu| (B_{5\rho}) + 
		4 \, \tilde\upbeta \, \lambda (B_{5\rho}) \, f(x^{\ast}\!,4\rho) \Bigg] .
\end{align*}
Finally taking in \eqref{DGgamma0} $s_1 = t^{\ast}$, $s_2 = t^{\ast} + \tilde\upbeta f(x^{\ast} \!, 4 \rho)$, 
$R = \tilde{r} = 5\rho$, $r = 4\rho$, $\varepsilon = 0$ in \eqref{DGgamma+},
\begin{align*}
\int_{t^{\ast}}^{t^{\ast} + \tilde\upbeta \, f(x^{\ast}\!,4\rho)} & \!\!\! \int_{B_{4\rho}^0}
		\Big|D \Big(u - \frac{\eta h}{ 2^{m-1}}\Big)_-\Big|^2 (x,\tau) \, \la  \, dx d\tau  \leqslant										\\
\leqslant \gamma & \Bigg[ \sup_{t \in ( t^{\ast}\!,  \, t^{\ast}\! + \tilde\upbeta \, f(x^{\ast}\!, 4\rho))} 
		\int_{I^{5\rho, \rho}_0}  \Big(u - \frac{\eta h}{ 2^{m-1}}\Big)_-^2 (x,t) \mu_-(x) \, dx +										\\
& \hskip25pt + \sup_{t \in ( t^{\ast}\!,  \, t^{\ast}\! + \tilde\upbeta \, f(x^{\ast}\!, 4\rho))} 
					\int_{I^{5\rho, \rho}_0}  \Big(u - \frac{\eta h}{ 2^{m-1}}\Big)_-^2 (x,t) \mu_+(x) \, dx	 +						\\
& \hskip50pt + \, \frac{1}{\rho^2}
		\int_{t^{\ast}}^{t^{\ast} + \tilde\upbeta \, f(x^{\ast}\!,4\rho)} \!\!\! \int_{(B_{4\rho}^0)^{\rho}} 
		 \Big(u - \frac{\eta h}{ 2^{m-1}}\Big)_-^2(x,\tau)\, \la \, dx d\tau  \Bigg] 	\leqslant										\\
\leqslant \gamma & \Bigg[
		\left( \frac{\eta h}{ 2^m}\right)^2 \, |\mu| (I^{5\rho, \rho}_0) +
			\left( \frac{\eta h}{ 2^{m-1}} \right)^2 \, \frac{1}{\rho^2} \, 
			\lambda \big( (B_{4\rho}^0)^{\rho} \big) \, \tilde\upbeta \, f(x^{\ast}\!,4\rho) \Bigg] \leqslant							\\
\leqslant \gamma  & \left( \frac{\eta h}{ 2^{m-1}}\right)^2 \frac{1}{\rho^2} \, \Bigg[ \rho^2 \, |\mu| (B_{5\rho}) + 
		\tilde\upbeta \, \lambda \big(B_{5\rho}\big) \, f(x^{\ast}\!,4\rho) \Bigg] .
\end{align*}
Summing up we get
\begin{equation}
\label{mo'lanumeriamo}
\begin{array}{l}
{\displaystyle
\int_{t^{\ast}}^{t^{\ast} + \tilde\upbeta \, f(x^{\ast}\!,4\rho)}
	\!\!\! \int_{B_{4\rho}} \Big|D \Big(u - \frac{\eta h}{ 2^{m-1}}\Big)_-\Big|^2 (x,\tau) \, \la  \, dx d\tau  \leqslant	}			\\	[1em]
\hskip50pt \leqslant \gamma 
	{\displaystyle \left( \frac{\eta h}{ 2^{m-1}}\right)^2 \frac{1}{\rho^2} \, \Bigg[ 3 \rho^2 \, |\mu| (B_{5\rho}) + 
		9 \, \tilde\upbeta \, f(x^{\ast}\!,4\rho)\, \lambda \big(B_{5\rho}\big) \Bigg]  }
\end{array}
\end{equation}
and so we can conclude, from \eqref{saraclara2}, that
\begin{align*}
a_{m+1}^+ \leqslant  64 \, \gamma_1 \gamma^{1/2} \q^2 \,  
		\frac{|\mu|_{\lambda} (B_{4\rho})}{(\lambda (B_{4\rho}))^{1/p}} & \, \big( \tilde\upbeta \, f(x^{\ast}\!,4\rho) \big)^{\frac{p-1}{p}}
	{\displaystyle	\cdot \left( d_{m-1}^{\texttt{\,>}} - d_m^{\texttt{\,>}} \right)^{\frac{2-p}{2p}} }	\cdot								\\
&	{\displaystyle \cdot \Bigg[ 3 \rho^2 \, |\mu| (B_{5\rho}) + 
		9 \, \tilde\upbeta \, f^+(x^{\ast}\!,4\rho)\, \lambda \big(B_{5\rho}\big) \Bigg]^{1/2} }.						
\end{align*}
Taking the power $\frac{2p}{2-p}$ and summing between $m = 1$ and $m = {m}^{\ast}$ we have
\begin{align*}
\sum_{m = 1}^{m^{\ast}} (a_{m+1}^+)^{\frac{2p}{2-p}}
	\leqslant  (64 \, \gamma_1 \gamma^{1/2} \q^2)^{\frac{2p}{2-p}} \, & 
	\frac{\left( |\mu|_{\lambda} (B_{4\rho}) \right)^{\frac{2p}{2-p}}}
		{(\lambda (B_{4\rho}))^{\frac{2}{2-p}}} \, \big( \tilde\upbeta \, f(x^{\ast}\!,4\rho) \big)^{\frac{2(p-1)}{2-p}} \cdot					\\
&	{\displaystyle \cdot \Bigg[ 3 \rho^2 \, |\mu| (B_{5\rho}) + 
		9 \, \tilde\upbeta \, f(x^{\ast}\!,4\rho)\, \lambda \big(B_{5\rho}\big) \Bigg]^{\frac{p}{2-p}} }
	{\displaystyle \left( d_{0}^{\texttt{\,>}} - d_{m^{\ast}}^{\texttt{\,>}} \right) } \, .
\end{align*}
Since the sequences $(a_{m}^+)_{m \in \N}$ and $(d_{m}^{\texttt{\,>}})_{m \in \N}$ are decreasing we can estimate 
$\sum_{m = 1}^{m^{\ast}} (a_{m+1}^+)^{\frac{2p}{2-p}}$ from below by $m^{\ast} (a_{m^{\ast}+1}^+)^{\frac{2p}{2-p}}$
and $d_{0}^{\texttt{\,>}} - d_{m^{\ast}}^{\texttt{\,>}}$ from above by $d_{0}^{\texttt{\,>}}$ and
$d_{0}^{\texttt{\,>}}$ by $\tilde\upbeta \, f(x^{\ast}\!,4\rho)\, \lambda ( B_{4\rho} )$ and get
\begin{align*}
(a_{m^{\ast}+1}^+)^{\frac{2p}{2-p}}
&	\leqslant  \frac{1}{m^{\ast}} \, (64 \, \gamma_1 \gamma^{1/2} \q^2)^{\frac{2p}{2-p}} \, 
	\frac{\left( |\mu|_{\lambda} (B_{4\rho}) \right)^{\frac{2p}{2-p}}}
		{(\lambda (B_{4\rho}))^{\frac{2}{2-p}}} \, \big( \tilde\upbeta \, f(x^{\ast}\!,4\rho) \big)^{\frac{2(p-1)}{2-p}} \cdot						\\
&	\qquad {\displaystyle \cdot \Bigg[ 3 \rho^2 \, |\mu| (B_{5\rho}) + 
	9 \, \tilde\upbeta \, f(x^{\ast}\!,4\rho)\, \lambda \big(B_{5\rho}\big) \Bigg]^{\frac{p}{2-p}} } 
	\tilde\upbeta \, f(x^{\ast}\!,4\rho)\, \lambda ( B_{4\rho} ) \leqslant																\\
&	\leqslant \frac{C^{\frac{2p}{2-p}}}{m^{\ast}} \, 
	\big( \tilde\upbeta   f(x^{\ast}\!,4\rho) \big)^{\frac{2p}{2-p}} \, \big( |\mu|_{\lambda} (B_{4\rho}) \big)^{\frac{2p}{2-p}} =					\\
&	= \frac{C^{\frac{2p}{2-p}}}{m^{\ast}} \, 
	\big( |M|_{\Lambda} \big( B_{4\rho} (x^{\ast}) \times (t^{\ast}\!, t^{\ast} + \tilde \upbeta f(x^{\ast}\!,4\rho) ) \big) \big)^{\frac{2p}{2-p}} \, ,
\end{align*}
where $C = 64 \, \gamma_1 \gamma^{1/2} \q^{5/2} (3 + 9 \, \tilde\upbeta)^{1/2} \tilde\upbeta^{-1/2}$,
by which finally
\begin{align*}
a_{m^{\ast}+1}^+
&	\leqslant  C \, \left( \frac{1}{m^{\ast}} \right)^{\frac{2-p}{2p}} \, 
		|M|_{\Lambda} \big( B_{4\rho} (x^{\ast}) \times (t^{\ast}\!, t^{\ast} + \tilde \upbeta f(x^{\ast}\!,4\rho) ) \big)  \, .
\end{align*}
Then for every $\epsilon > 0$ one can find $m^{\ast}$ such that $C / {m^{\ast}}^{\frac{2-p}{2p}} \leqslant \epsilon$.
Then we consider 
$$
m^{\ast} \geqslant \left( \frac{64 \, \gamma_1 \gamma^{1/2} \q^{5/2} (3 + 9 \, \tilde\upbeta)^{1/2}}{\tilde\upbeta^{1/2} \epsilon} \right)^{\frac{2p}{2-p}}
\qquad \text{and} \qquad
\eta_1 = \frac{\eta}{2^{m^{\ast}}}
$$
which depends on $\gamma_1 , \gamma , \q , \epsilon , \tilde\upbeta , \eta$. \\
Now by \eqref{carlettomio} we immediately get that
\begin{align*}
& \frac{\Lambda_+ \Big( \{u < \eta_1 h \} \, \cap \, \Big[B_{4\rho}^+ (x^{\ast}) \times (t^{\ast}, t^{\ast} + \tilde\upbeta \, \rho^2 h(x^{\ast}\!, 4\rho) ) \Big] \Big)}
	{\Lambda \Big( B_{4\rho} (x^{\ast}) \times (t^{\ast}, t^{\ast} + \tilde\upbeta \, \rho^2 h(x^{\ast} \!, 4\rho) ) \Big)}	\leqslant		\\
& \hskip30pt \leqslant \upkappa \left(
\frac{M_+ \Big( \{u < \eta_1 h \} \, \cap \, \Big[B_{4\rho}^+ (x^{\ast}) \times (t^{\ast}, t^{\ast} + \tilde\upbeta \, \rho^2 h(x^{\ast}\!, 4\rho) ) \Big] \Big)}
	{|M|_{\Lambda} \Big( B_{4\rho} (x^{\ast}) \times (t^{\ast}, t^{\ast} + \tilde\upbeta \, \rho^2 h(x^{\ast} \!, 4\rho) ) \Big)}
	\right)^{\uptau}
\leqslant \upkappa \, \epsilon^{\uptau} \, .
\end{align*}
\ \\
In a complete analogous way one can prove point $ii$). \\ \ \\
Point $iii$): the case where $\mu \equiv 0$, as usual, is a bit different. Let us see a sketch of the proof.
Integrating between
$t^{\ast} - \upbeta \, \rho^2 h(x^{\ast}\!, 4\rho)$ and $t^{\ast} + \upbeta \, \rho^2 h(x^{\ast}\!, 4\rho)$ and 
using Lemma \ref{lemma2.2} similarly as before but with $\nu = \lambda$ and $\bar{\nu} = \lambda_0|_{B_r (x^{\ast})}$, we get
\begin{align*}
\frac{1}{2 \q^2} \, \frac{\eta h}{2^{m+1}} a_{m+1}^0
& {\displaystyle \leqslant  8 \, \gamma_1 \, \rho
		\ (\lambda (B_{4\rho}))^{\frac{p-1}{p}} \cdot
		\left( \int_{t^{\ast} - \upbeta \, f(x^{\ast}\!, 4\rho)}^{t^{\ast} + \upbeta \, f(x^{\ast}\!,4\rho)} 
		\big[ \lambda (A_{m-1} (t)) - \lambda (A_m (t)) \big]\, d t	\right)^{\frac{2-p}{2p}} }	\cdot									\\
& {\displaystyle \hskip40pt \cdot \, \big( 2 \upbeta \, f(x^{\ast}\!,4\rho) \big)^{\frac{p-1}{p}}
		\left( \int_{t^{\ast} - \upbeta \,  f(x^{\ast}\!, 4\rho) }^{t^{\ast} + \upbeta \, f(x^{\ast}\!,4\rho)} \!\!\! \int_{B_{4\rho}} 
			\Big|D \Big(u - \frac{\eta h}{ 2^{m-1}}\Big)_-\Big|^2 (x,t) \, \la  \, dx dt \right)^{1/2} }.
\end{align*}		
Now estimating the part involving the gradient of $u - \frac{\eta h}{ 2^m}$ similarly as \eqref{mo'lanumeriamo} we get
\begin{align*}
\frac{1}{2 \q^2} \, \frac{\eta h}{2^{m+1}} a_{m+1}^0
& {\displaystyle \leqslant  8 \, \gamma_1 \, \gamma^{\frac{1}{2}}
		\ (\lambda (B_{4\rho}))^{\frac{p-1}{p}} \cdot
		\left( \int_{t^{\ast} - \upbeta \, f(x^{\ast}\!, 4\rho)}^{t^{\ast} + \upbeta \, f(x^{\ast}\!,4\rho)} 
		\big[ \lambda (A_{m-1} (t)) - \lambda (A_m (t)) \big]\, d t	\right)^{\frac{2-p}{2p}} }	\cdot									\\
& \hskip20pt \cdot \, \big( 2 \upbeta \, f(x^{\ast}\!,4\rho) \big)^{\frac{p-1}{p}} 
		\left( \frac{\eta h}{ 2^{m-1}}\right) \, \Bigg[ 3 \rho^2 \, |\mu| (B_{5\rho}) + 
		18 \, \upbeta \, f(x^{\ast}\!,4\rho)\, \lambda \big(B_{5\rho}\big) \Bigg]^{\frac{1}{2}}
\end{align*}		
and proceeding as before we reach
\begin{align*}
a_{m^{\ast}+1}^0
&	\leqslant  C' \, \left( \frac{1}{m^{\ast}} \right)^{\frac{2-p}{2p}} \, 
		\Lambda \big( B_{4\rho} (x^{\ast}) \times (t^{\ast}\! - \upbeta f(x^{\ast}\!,4\rho) , t^{\ast} + \upbeta f(x^{\ast}\!,4\rho) ) \big)
\end{align*}
with $C' = 64 \, \gamma_1 \gamma^{1/2} \q^{5/2} (3 + 18 \, \upbeta)^{1/2} (2 \upbeta)^{-1/2}$. The conclusion is as before. \\
\ \\
Finally let us see point $iv \, )$.
If $B_{4\rho} (x^{\ast}) \subset \Omega_0$ we have
\begin{align*}
\frac{1}{2 \q^2} \, \frac{\eta h}{2^{m+1}} & \lambda \big(A_{m+1} (t) \big) =
	\frac{1}{2 \q^2} \, \frac{\eta h}{2^{m+1}} \lambda \big(A_{m+1} (t) \big) \leqslant										\nonumber	\\
& {\displaystyle \leqslant  8  \, \gamma_1 \, \rho
		\ (\lambda (B_{4\rho}))^{\frac{p-1}{p}} \cdot \left(\int_{{A_{m-1} (t) \setminus A_m (t)}} 
								\!\!\!\!\!\!\!\!\!\!\!\!\!\!\!\!\! |D u|^p (x,t) \, \la  \, dx \right)^{1/p} } d\tau \leqslant		\nonumber	\\
& \leqslant \big( \lambda (A_{m-1} (t)) - \lambda (A_m (t)) \big)^{\frac{2-p}{2}}
	\left( \int_{B_{4\rho}} \Big|D \Big(u - \frac{\eta h}{ 2^{m-1}}\Big)_-\Big|^2 (x,t) \, \lambda (x)  \, dx  \right)^{1/2} \, .
\end{align*}
%
%
%
Since $B_{5\rho} (x^{\ast}) \subset \Omega_0$ taking $\tilde{r} = 5\rho$, $r = 4 \rho$ and $\varepsilon = 0$ in \eqref{tempofissato}, we get for almost
every $t$ that
\begin{align*}
\int_{B_{4\rho}} \Big|D \Big(u - \frac{\eta h}{ 2^{m-1}} & \Big)_-\Big|^2 (x,t) \, \lambda (x)  \, dx	 	\leqslant 									\\
&	\leqslant \gamma \, \frac{1}{\rho^2} 	\int_{B_{5 \rho}} \Big(u - \frac{\eta h}{ 2^{m-1}}\Big)_-^2 (x,t) \, \lambda (x) \, dx \leqslant
	\gamma \left( \frac{\eta h}{ 2^{m-1}}\right)^2 \frac{1}{\rho^2} \,  \lambda (B_{5\rho})
\end{align*}
and then $\lambda \big(A_{m+1} (t) \big) \leqslant	64 \, \q^2 \, \gamma_1 \, \gamma^{1/2} \, 
\big( \lambda (A_{m-1} (t)) - \lambda (A_m (t)) \big)^{\frac{2-p}{2}} \,  \big( \lambda (B_{5\rho}) \big)^{1/2}$.
By that we can conclude similarly as above.
\finedimo

\noindent
Now we state a result known as {\em expansion of positivity}. It will be a fundamental step
to prove the Harnack inequality.

\begin{lemma}
\label{esp_positivita}
Consider $(x^{\ast} \!, t^{\ast})$ such that
$B_{5\rho}(x^{\ast}) \times [t^{\ast} - 16 \, h(x^{\ast} \!, 4\rho) \, \rho^2, t^{\ast} + 16 \, h(x^{\ast} \!, 4\rho) \, \rho^2]
\subset \Omega \times (0,T)$. \\
Consider the value $\tilde\upbeta$
determined in Lemma \ref{lemma1} and used in in Lemma \ref{lemma2}.
Then for every $\hat\theta \in (0, 1)$ there is $\uplambda > 0$ depending only on
$\gamma_1 , \gamma , \q , \kappa, \tilde\upbeta , \hat\theta$
such that for every $h > 0$ and $u \geqslant 0$ in $DG(\Omega, T, \mu, \la, \gamma)$
points $i\, )$ and $ii\, )$ are true: \\ [0.3em]
$i\, )$ if $\mu_+ (B_{\rho}(x^{\ast})) > 0$ and
$$
u(\cdot , t^{\ast}) \geqslant h \qquad \text{a.e. in } B_{\rho}^+(x^{\ast})
$$
then
\begin{align*}
u \geqslant \uplambda h \qquad \text{a.e. in }
	& \, B_{2\rho}^+(x^{\ast}) \times 
		\big(t^{\ast} + \hat\theta \, \tilde\upbeta \, h(x^{\ast}, 4\rho) \rho^2, t^{\ast} + \tilde\upbeta \, h(x^{\ast}, 4\rho) \rho^2 \big) ;
\end{align*}
$ii\, )$ if $\mu_- (B_{\rho}^-(x^{\ast})) > 0$ and
$$
u(\cdot , t^{\ast}) \geqslant h \qquad \text{a.e. in } B_{\rho}^-(x^{\ast})
$$
then
\begin{align*}
u \geqslant \uplambda h \qquad \text{a.e. in }
	& \, B_{2\rho}^-(x^{\ast}) \times 
		\big(t^{\ast} + \hat\theta \, \tilde\upbeta \, h(x^{\ast}, 4\rho) \rho^2, t^{\ast} + \tilde\upbeta \, h(x^{\ast}, 4\rho) \rho^2 \big) .
\end{align*}
Moreover for every $\upbeta > 0$ for which
$B_{5\rho}(x^{\ast}) \times [t^{\ast} - \upbeta \, h(x^{\ast} \!, 4\rho) \, \rho^2, t^{\ast} + \upbeta \, h(x^{\ast} \!, 4\rho) \, \rho^2]
\subset \Omega \times (0,T)$
there is $\uplambda > 0$ depending only on $\gamma_1 , \gamma , \q , \kappa, \upbeta$
such that for every $h > 0$ and $u \geqslant 0$ in $DG(\Omega, T, \mu, \la, \gamma)$ point $iii\, )$ is true: \\ [0.3em]
$iii\, )$ if $\lambda_0 (B_{\rho}(x^{\ast})) > 0$ and
$$
u \geqslant h \qquad \text{a.e. in } B^0_{\rho}(x^{\ast}) \times 
	\big( t^{\ast} - \upbeta \, h(x^{\ast}, 4\rho) \rho^2, t^{\ast} + \upbeta \, h(x^{\ast}, 4\rho) \rho^2 \big)
$$
then
\begin{align*}
u \geqslant \uplambda h \qquad \text{a.e. in }
	B^0_{2\rho}(x^{\ast}) \times \big( t^{\ast} - \upbeta \, h(x^{\ast}, 4\rho) \rho^2, t^{\ast} + \upbeta \, h(x^{\ast}, 4\rho) \rho^2 \big) .
\end{align*}
If $B_{5 \rho} (x^{\ast}) \subset \Omega_0$ there is $\uplambda > 0$ depending only on
$\gamma_1 , \gamma , \q , \kappa$ such that
for every $h > 0$ and $u \geqslant 0$ in $DG(\Omega, T, \mu, \la, \gamma)$ point $iv\, )$ is true: \\ [0.3em]
$iv\, )$ for almost every $t \in (0,T)$ if
$$
u (\cdot , t) \geqslant h \qquad \text{a.e. in } B_{\rho} (x^{\ast})
$$
then
\begin{align*}
u (\cdot , t)  \geqslant \uplambda h \qquad \text{a.e. in } \, B_{2\rho}(x^{\ast}) .
\end{align*}
\end{lemma}
\noindent
\dimo
The proof is a consequence of Proposition \ref{prop-DeGiorgi2} and Lemma \ref{lemma2}.
We start from point $i\, )$: in Proposition \ref{prop-DeGiorgi2} we consider $\underline{m} = 0$, $R = 4 \rho$, $r = 2 \rho$, 
$\upbeta^{\diamond} = \tilde{\upbeta}$ (the value determined in Lemma \ref{lemma1} and used in Lemma \ref{lemma2}
and belonging to $(0,16]$),
$\theta^{\diamond}$ and $a \in (0,1)$ arbitrary;
from Proposition \ref{prop-DeGiorgi2} we derive the existence of $\underline{\nu}^{\diamond} \in (0,1)$ such that if, for $c > 0$ an arbitrary constant,
the following holds
\begin{align*}
& \frac{M_+ \Big( \big\{ u < c  \big\} \cap \big( B_{4\rho} (x^{\ast}) \times (t^{\ast}, t^{\ast} + \tilde{\upbeta} \rho^2 h(x^{\ast}, 4\rho) ) \big) \Big)}
{|M|_{\Lambda} \Big( B_{4\rho} (x^{\ast}) \times (t^{\ast}, t^{\ast} + \tilde\upbeta \, \rho^2 h(x^{\ast} \!, 4\rho) ) \Big)} +								\\
& \qquad \qquad \qquad
+ \frac{\Lambda_+ \Big( \big\{ u < c  \big\} \cap \big( B_{4\rho} (x^{\ast}) \times (t^{\ast}, t^{\ast} + \tilde{\upbeta} \rho^2 h(x^{\ast}, 4\rho) ) \big) \Big)}
{\Lambda \Big( B_{4\rho} (x^{\ast}) \times (t^{\ast}, t^{\ast} + \tilde\upbeta \, \rho^2 h(x^{\ast} \!, 4\rho) ) \Big)}
\leqslant \underline{\nu}^{\diamond}
\end{align*}
then
$$
u \geqslant a \, c \qquad \text{ in } B^+_{2\rho}(x^{\ast}) \times 
	\left(t^{\ast} + \theta^{\diamond} \tilde\upbeta \, h(x^{\ast}, 4\rho) \rho^2, t^{\ast} + \tilde\upbeta \, h(x^{\ast}, 4\rho) \rho^2 \right) \, .
$$
Now we use Lemma \ref{lemma2}: consider $\eta$ the value determined in Lemma \ref{lemma1} and used in in Lemma \ref{lemma2},
take $\upbeta = 16$ and $\epsilon$ such that $\epsilon + \upkappa \, \epsilon^{\uptau} = \underline{\nu}^{\diamond}$
and conclude that there is $\eta_1$
(depending on $\gamma_1 , \gamma , \q , \tilde\upbeta , \eta, \underline{\nu}^{\diamond}$ and then
on $\gamma_1 , \gamma , \q , \tilde\upbeta , \eta, \kappa, a , \theta^{\diamond}$, but $\eta$ depends only on $\gamma$ and $\q$) such that
\begin{align*}
& \frac{M_+ \Big( \big\{ u < \eta_1 h  \big\} \cap \big( B_{4\rho} (x^{\ast}) \times (t^{\ast}, t^{\ast} + \tilde{\upbeta} \rho^2 h(x^{\ast}, 4\rho) ) \big) \Big)}
{|M|_{\Lambda} \Big( B_{4\rho} (x^{\ast}) \times (t^{\ast}, t^{\ast} + \tilde\upbeta \, \rho^2 h(x^{\ast} \!, 4\rho) ) \Big)} +								\\
& \qquad \qquad \qquad
+ \frac{\Lambda_+ \Big( \big\{ u < \eta_1 h  \big\} \cap \big( B_{4\rho} (x^{\ast}) \times (t^{\ast}, t^{\ast} + \tilde{\upbeta} \rho^2 h(x^{\ast}, 4\rho) ) \big) \Big)}
{\Lambda \Big( B_{4\rho} (x^{\ast}) \times (t^{\ast}, t^{\ast} + \tilde\upbeta \, \rho^2 h(x^{\ast} \!, 4\rho) ) \Big)}
\leqslant \underline{\nu}^{\diamond} \, .
\end{align*}
Then
$$
u \geqslant a \, \eta_1 h \qquad \text{ in } B^+_{2\rho}(x^{\ast}) \times 
	\left(t^{\ast} + \theta^{\diamond} \tilde\upbeta \, h(x^{\ast}, 4\rho) \rho^2, t^{\ast} + \tilde\upbeta \, h(x^{\ast}, 4\rho) \rho^2 \right) \, .
$$
Taking $\hat{\theta} = \theta^{\diamond}$, $a = 1/2$ for simplicity and $\uplambda = \eta_1/2$ we conclude the proof of point $i \, )$.
In the same way one can prove point point $ii \, )$. \\ [0.3em]
Let us see point point $iii \, )$.
In Proposition \ref{prop-DeGiorgi2} we consider again $\underline{m} = 0$, $R = 4 \rho$, $r = 2 \rho$,
$\upbeta^{\star} = \upbeta \, h(x^{\ast}, 4\rho) / 8$ and $a \in (0,1)$ arbitrary.
We derive the existence of $\underline{\nu}^{\star} \in (0,1)$ such that if, for $c > 0$ an arbitrary constant,
the following holds
\begin{align*}
\Lambda_0 \bigg( \big\{ u < c  \big\} \, \cap \, 
\Big( & B_{4\rho} (x^{\ast}) \times 
	\big(t^{\ast} - \upbeta \, h(x^{\ast}, 4\rho) \, \rho^2, t^{\ast} + \upbeta \, h(x^{\ast}, 4\rho) \, \rho^2 \big) \Big) \bigg) \leqslant		\\
	 & \leqslant \, \underline{\nu}^{\star} \, 
	 \Lambda \Big( B_{4\rho} (x^{\ast}) \times 
	 	\big(t^{\ast} - \upbeta \, h(x^{\ast}, 4\rho) \, \rho^2, t^{\ast} + \upbeta \, h(x^{\ast}, 4\rho) \, \rho^2 \big) \Big)  \, ,
\end{align*}
then
$$
u \geqslant a \, c \qquad \text{ in } 
	B^0_{2\rho}(x^{\ast}) \times \big(t^{\ast} - \upbeta \, h(x^{\ast}, 4\rho) \, \rho^2, t^{\ast} + \upbeta \, h(x^{\ast}, 4\rho) \, \rho^2 \big) \, .
$$
Now in Lemma \ref{lemma2} take $\epsilon = \underline{\nu}^{\star}$ and conclude that there is $\eta_1$
(depending on $\gamma_1 , \gamma , \q , \kappa, a, \upbeta$) such that
$$
u \geqslant a \, \eta_1 h \qquad \text{ in } 
	B^0_{2\rho}(x^{\ast}) \times \big(t^{\ast} - \upbeta \, h(x^{\ast}, 4\rho) \, \rho^2, t^{\ast} + \upbeta \, h(x^{\ast}, 4\rho) \, \rho^2 \big) \, .
$$
Taking, e.g,  $a = 1/2$ we conclude. \\ [0.3em]
To prove point $iv$) we consider $\underline{m}$, $R$, $r$ and $a \in (0,1)$ as above and use point $iv$) of Proposition \ref{prop-DeGiorgi2}.
Then we get the existence of $\underline{\nu} \in (0,1)$ such that, for $c > 0$, if
\begin{align*}
\lambda \big(\{ x \in B_{4\rho} (x^{\star}) \, | \, u(x,t) < c \} \big) \leqslant \underline{\nu} \ \lambda (B_{4\rho} (x^{\star}) )
\end{align*}
then
$u(x,t) \geqslant a \, c$ for a.e. $x \in B_{2\rho} (x^{\star})$. Using Lemma \ref{lemma2} we conclude as above.
\finedimo

\section{The Harnack type inequality}
\label{secHarnack}

The following theorems (Theorem \ref{Harnack1} and Theorem \ref{Harnack2}) are the main results of the paper. \\

\begin{theorem}
\label{Harnack1}
Assume $u\in DG(\Omega, T, \mu, \la, \gamma)$, $u\geqslant 0$, $(x_o, t_o) \in \Omega \times (0,T)$ and fix $\rho > 0$.
\begin{itemize}
\item[$i\, $)]
Suppose $x_o \in \Omega_+ \cup I_+$.
For every $\vartheta_+ \in (0,1]$ for which
$B_{5\rho}(x_o) \times [t_o - h(x_o, \rho) \rho^2, t_o + 16 \, h(x_o, 4\rho) \rho^2 + \vartheta_+ h(x_o, \rho) \rho^2] \subset \Omega \times (0,T)$
there exists $c_+ > 0$ depending $($only$)$ on
$\gamma_1, \gamma, \q, \kappa, \alpha, \upkappa, \anna, K_1, K_2, K_3, q, \varsigma , \vartheta_+$ such that
$$u(x_o, t_o) \leqslant c_+ \, \inf_{B_{\rho}^+ (x_o)} u(x, t_o + \vartheta_+ \, \rho^2 h(x_o, \rho)) .$$
\item[$ii\, $)]
Suppose $x_o \in \Omega_- \cup I_-$.
For every $\vartheta_- \in (0,1]$ for which
$B_{5\rho}(x_o) \times [t_o - 16 \, h(x_o, 4\rho) \rho^2 + \vartheta_- h(x_o, \rho) \rho^2, t_o + h(x_o, \rho) \rho^2] \subset \Omega \times (0,T)$
there exists $c_- > 0$ depending $($only$)$ on
$\gamma_1, \gamma, \q, \kappa, \alpha, \upkappa, \anna, K_1, K_2, K_3, q, \varsigma , \vartheta_-$ such that
$$u(x_o, t_o) \leqslant c_- \, \inf_{B_{\rho}^- (x_o)} u(x, t_o - \vartheta_- \, \rho^2 h(x_o, \rho)) .$$
\item[$iii\, $)]
Suppose $x_o \in \Omega_0 \cup I_0$.
Suppose $B_{5\rho}(x_o) \times [t_o - 16 \, h(x_o, 4\rho) \rho^2, t_o + 16 \, h(x_o, 4\rho) \rho^2] \subset \Omega \times (0,T)$.
For every $s_1, s_2$ for which $s_2 - t_o = t_o - s_1 \leqslant 16 \, h(x_o, 4\rho) \rho^2$, suppose
$s_2 - t_o = t_o - s_1 = \upomega \, h(x_o, 4\rho) \rho^2$ for $\upomega \in (0,16]$,
there is $c_0$ depending $($only$)$ on 
$K_1, K_2, K_3, q, \varsigma, \kappa, \gamma_1, \gamma, \lela, h(x_o, 4\rho), \q$ such that
$$\sup_{B_{\rho}^+ (x_o) \times [s_1, s_2]}  u \leqslant c_0 \, \inf_{B_{\rho}^+ (x_o) \times [s_1, s_2]}  u.$$
\item[$iv\, $)]
Suppose $B_{5\rho}(x_o) \subset \Omega_0$. Then there is $c$ depending $($only$)$ on
$K_1, K_2, K_3, q, \varsigma, \kappa, \gamma_1, \gamma, \q$ such that
for almost every $t \in (0,T)$
$$\sup_{B_{\rho} (x_o)}  u (\cdot, t) \leqslant c \, \inf_{B_{\rho} (x_o)}  u (\cdot, t).$$
\end{itemize}
\end{theorem}
\noindent
\dimo
We start by proving the first of the three inequalities under the assumption that $B_{\rho}^+(x_o) \not= \emptyset$.
For some $r_1, r_2 > 0$ and  $(\bar{x}, \bar{t}) \in 
B_{5\rho}(x_o) \times [t_o - h(x_o, \rho) \rho^2, t_o + 16 \, h(x_o, 4\rho) \rho^2 + \vartheta_+ h(x_o, \rho) \rho^2] \subset \Omega \times (0,T)$
we define the sets 
\begin{gather*}
Q^{+, \texttt{\,<}}_{r_1, h(\bar{y}, r_2)} (\bar{x}, \bar{t}) := \Big( B_{r_1}^+ (\bar{x}) \times [\bar{t} - h(\bar{y}, r_2) r_1^2, \bar{t}] \Big)
\, , \quad	Q^{+, \texttt{\,>}}_{r, h(\bar{y}, r_2)} (\bar{x}, \bar{t}) := \Big( B_{r}^+ (\bar{x}) \times [\bar{t}, \bar{t} + h(\bar{y}, r_2) r_1^2] \Big)	\, ,		\\
Q^{\texttt{\,<}}_{r_1, h(\bar{y}, r_2)} (\bar{x}, \bar{t}) := \Big( B_{r_1} (\bar{x}) \times [\bar{t} - h(\bar{y}, r_2) r_1^2, \bar{t}] \Big)
\, , \quad	Q^{\texttt{\,>}}_{r,h(\bar{y}, r_2)} (\bar{x}, \bar{t}) := \Big( B_{r} (\bar{x}) \times [\bar{t}, \bar{t} + h(\bar{y}, r_2) r_1^2] \Big)	\, .
\end{gather*}
We may write $u(x_o, t_o) = b \, \rho^{-\xi}$ for some $b, \xi > 0$ to be fixed later.
Define the functions
$$
\emm (r) = \sup_{Q^{+, \texttt{\,<}}_{r, h (x_o, \rho)} (x_o, t_o)} u,		\qquad		\enn (r) = b (\rho - r)^{-\xi}, 	\qquad r \in [0,\rho) .
$$
Let us denote by $r_o \in [0,\rho)$ the largest solution of $\emm (r) = \enn (r)$.
Define
$$
\carlo := \enn (r_o) = b (\rho - r_o)^{-\xi} \, .
$$
We can find $(y_o, \tau_o) \in Q^{+, \texttt{\,<}}_{r_o, \, h(x_o,\rho)} (x_o, t_o)$ such that
\begin{equation}
\label{choicey0t0}
\frac{3\carlo}{4} < \sup_{Q^{+, \texttt{\,<}}_{\frac{\rho_o}{4}\!, \, h (y_o, \rho_o)} (y_o, \tau_o)} u \leqslant \carlo
\end{equation}
where $\rho_o \in (0, (\rho - r_o) / 2 ]$. If $\rho_o \leqslant (\rho - r_o) / 2$ then $B^+_{\rho_o} (y_o) \subset B_{\frac{\rho + r_o}{2}} (x_o)$.
We want the value of $\rho_o$ to be be chosen in such a way that
$$
Q^{+, \texttt{\,<}}_{\rho_o, \, h(y_o,\rho_o)} (y_o, \tau_o) \subset Q^{+, \texttt{\,<}}_{\frac{\rho + r_o}{2}, \, h(x_o,\rho)}(x_o,t_o)
$$
and the request $\rho_o \leqslant (\rho - r_o) / 2$ may be not sufficient. We also need
$\tau_o - h(y_o, \rho_o) \rho_o^2 \geqslant t_o - h(x_o,\rho) (\rho+r_o)^2/4$ and this is guaranteed if
\begin{equation}
\label{enumeriamopurequesta!}
h(y_o,\rho_o) \rho_o^2 \leqslant h(x_o, \rho) \left[ \frac{(\rho + r_o)^2}{4} - r_o^2 \right] ,
\end{equation}
which in turn is true, since $r_o^2 < \rho \, r_o$, if
$$
h(y_o,\rho_o) \rho_o^2 \leqslant h(x_o, \rho) \, \frac{(\rho - r_o)^2}{4} .
$$
so we will choose $\rho_o$ satisfying these two requests. 
Notice that this last request can be satisfied writing
$h(y_o, \rho_o) \, \rho_o^2 = h(y_o, \rho_o) \rho_o^{2\alpha} \rho_o^{2(1-\alpha)}$
because, thanks to Remark \ref{notaimportante}, point $\mathpzc{C}$, and (H.2)$'$ we have
\begin{align*}
h(y_o, \rho_o) \rho_o^{2\alpha} & \leqslant \tilde{K}_2^2 \, h(y_o, 2 \rho) (2\rho)^{2 \alpha}	\leqslant									\\
	& \leqslant \tilde{K}_2^2 \, \frac{|\mu|_{\lambda} (B_{4\rho} (x_o))}{\lambda (B_{2\rho} (y_o))} (2\rho)^{2 \alpha}	\leqslant				\\
	& \leqslant 4^{\alpha} \tilde{K}_2^2 \, \q^2 \, h (x_o, \rho) \rho^{2 \alpha}
\end{align*}
and then we have
\begin{align*}
h(y_o, \rho_o) \rho_o^2 \leqslant 4^{\alpha} \tilde{K}_2^2 \, \q^2 \, h (x_o, \rho) \rho^{2 \alpha} \rho_o^{2(1-\alpha)}.
\end{align*}
Then \eqref{enumeriamopurequesta!} holds if in particular
\begin{align*}
4^{\alpha} \tilde{K}_2^2 \, \q^2 \, h (x_o, \rho) \rho^{2 \alpha} \rho_o^{2(1-\alpha)} \leqslant
	h(x_o, \rho) \, \frac{(\rho - r_o)^2}{4}
\end{align*}
that is
\begin{align}
\label{rozero2}
\rho_o^{1-\alpha} \leqslant \frac{1}{2^{\alpha} \tilde{K}_2 \, \q} \, \frac{1}{\rho^{\alpha}} \, \frac{\rho - r_o}{2}
\end{align}
and it is always possible to choose $\rho_o$ small enough such that \eqref{rozero2} is satisfied. Therefore $\rho_o$ will be chosen satisfying
\begin{align}
\label{rozero}
\rho_o = \min \left\{ \frac{\rho - r_o}{2}, 
	\left[ \frac{1}{2^{\alpha} \tilde{K}_2 \, \q} \, \frac{1}{\rho^{\alpha}} \, \frac{\rho - r_o}{2} \right]^{\frac{1}{1-\alpha}} \right\} .
\end{align}
By this choice of $\rho_o$ and by the choice of $r_o$ we have
\begin{equation}
\label{oscillazione}
\sup_{Q^{+, \texttt{\,<}}_{\rho_o, h(y_o,\rho_o)} (y_o, \tau_o)} u \leqslant \sup_{Q^{+, \texttt{\,<}}_{\frac{\rho + r_o}{2}, h(x_o,\rho)}(x_o,t_o)} u < 
	\enn \left( \frac{\rho + r_o}{2}\right) = 2^{\xi} \carlo .
\end{equation}
We now proceed dividing the proof in six steps. \\ [0.3em]
\textsl{Step 1 - }
In this step we want to show that there is $\overline{\nu} \in (0,1)$,
depending only on $\kappa, \gamma_1, \gamma, \xi, \q$, such that
\begin{equation}
\label{mimancanogliamorimiei}
\begin{array}{c}
{\displaystyle
\frac{M_+ \left( \left\{ u > \frac{\carlo}{2} \right\} \cap 
	Q^{+, \texttt{\,<}}_{\rho_o/2,h (y_o, \rho_o)} (y_o, \tau_o) \right)}
	{|M|_{\Lambda} \left( Q^{\texttt{\,<}}_{\rho_o/2,h (y_o, \rho_o)} (y_o, \tau_o) \right)}	}			> \overline{\nu} \, ,			\\	[2em]
{\displaystyle
\frac{\Lambda_+ \left( \left\{ u > \frac{\carlo}{2} \right\} \cap 
	Q^{+, \texttt{\,<}}_{\rho_o/2,h (y_o, \rho_o)} (y_o, \tau_o) \right)}
	{\Lambda \left( Q^{\texttt{\,<}}_{\rho_o/2,h (y_o, \rho_o)} (y_o, \tau_o) \right)}
> 	\overline{\nu}	}
\end{array}
\end{equation}
and that
\begin{equation}
\label{mimancanogliamorimiei-1}
\iint_{Q^{+, \texttt{\,<}}_{\frac{\rho_o}{2},h(y_o, \rho_o)} (y_o, \tau_o)} |Du|^2 \, \lambda \, dx dt   \leqslant	
	9 \, \gamma \, (2^{\xi} N)^2 \,  h(y_o, \rho_o) \, \lambda \big( B_{\rho_o}(y_o) \big)	\, .
\end{equation}
To prove \eqref{mimancanogliamorimiei} first we show that there is $\nu \in (0,1)$ such that
\begin{equation}
\label{mimancanogliamorimiei-2}
{\displaystyle \frac{M_+ \left( \left\{ u > \frac{\carlo}{2} \right\} \cap 
	Q^{+, \texttt{\,<}}_{\rho_o/2,h(y_o, \rho_o)} (y_o, \tau_o) \right)}
	{|M|_{\Lambda} \left( Q^{\texttt{\,<}}_{\rho_o/2,h(y_o, \rho_o)} (y_o, \tau_o) \right)} \, + }	 \, 
{\displaystyle \frac{\Lambda_+ \left( \left\{ u > \frac{\carlo}{2} \right\} \cap 
	Q^{+, \texttt{\,<}}_{\rho_o/2,h(y_o, \rho_o)} (y_o, \tau_o) \right)}
	{\Lambda \left( Q^{\texttt{\,<}}_{\rho_o/2,h(y_o, \rho_o)} (y_o, \tau_o) \right)}
> 	\nu \,  . }
\end{equation}
Argue by contradiction and suppose that \eqref{mimancanogliamorimiei-2} is false. Since
\begin{align*}
& Q^{+, \texttt{\,<}}_{\frac{\rho_o}{2},h(y_o, \rho_o)} (y_o, \tau_o) = 
\left( B_{\frac{\rho_o}{2}}^+ (y_o) \times 
\left[\tau_o - h(y_o, {\textstyle \frac{\rho_o}{2}}) {\textstyle \frac{h (y_o, \rho_o)}{h(y_o, \rho_o/2)}{\textstyle \frac{\rho_o^2}{4}} }, \tau_o \right] \right) \, ,\\
& Q^{+, \texttt{\,<}}_{\frac{\rho_o}{4},h(y_o, \rho_o)} (y_o, \tau_o) = 
\left( B_{\frac{\rho_o}{4}}^+ (y_o) \times 
\left[\tau_o - h(y_o, {\textstyle \frac{\rho_o}{2}}) {\textstyle \frac{h (y_o, \rho_o)}{h(y_o, \rho_o/2)}{\textstyle \frac{\rho_o^2}{16}} }, \tau_o \right] \right) \, ,
\end{align*}
setting in Proposition \ref{prop-DeGiorgi1}
\begin{gather*}
\overline{m} = \omega = 2^\xi \carlo, \quad R = \frac{\rho_o}{2},		\quad	\rho = \frac{\rho_o}{4},		\quad
							\sigma = 1 - 2^{-\xi-1}, 	\quad 	a = \sigma^{-1}\biggl(1-\frac{3}{2^{\xi+2}}\biggr) \, , 			\\
x^{\diamond} = y_o \, , \qquad  t^{\diamond} = \tau_o - h(y_o, \rho_o) \frac{\rho_o^2}{4}	 \, , \qquad
\upbeta^{\diamond} = \frac{h (y_o, \rho_o)}{h(y_o, \rho_o/2)} \, ,	\qquad
{\theta}^{\diamond} = \frac{3}{4} \, ,
\end{gather*}
we obtain from Proposition \ref{prop-DeGiorgi1} that
$$
u\leqslant \frac{3\carlo}{4} \quad \textrm{in }\, Q^{+, \texttt{\,<}}_{\frac{\rho_o}{4}\!, \, h(y_o, \rho_o)} (y_o, \tau_o)
$$
which contradicts \eqref{choicey0t0}.  Notice that $\upbeta^{\diamond} \in [\q^{-1}, \q]$.
Now by \eqref{mimancanogliamorimiei-2} we derive that at least one of the two addends in \eqref{mimancanogliamorimiei-2}
is greater or equal to $\nu/2$.  Now we get \eqref{mimancanogliamorimiei} by \eqref{carlettomio} taking
$$
\overline{\nu} = \frac{1}{\upkappa} \, \left( \frac{\nu}{2} \right)^{\frac{1}{\alpha}} \, .
$$
To prove \eqref{mimancanogliamorimiei-1} we use \eqref{DGgamma+}. In $\eqref{DGgamma+}$ we choose
$x_0 = y_o$, $t_0 = \tau_o - h (y_o, \rho_o) \rho_o^2$,
$R = \rho_o$, $\tilde{r} = \rho_o$, $r = \rho_o / 2$, $\varepsilon = 0$, $\upbeta = 1$, $\theta = \frac{3}{4}$,
$\tilde\theta = \frac{1}{2}$, $k = 0$ and since $u \leqslant 2^{\xi} \carlo$ we get
\begin{align}
\label{est2.18}
& \iint_{Q^{+, \texttt{\,<}}_{\frac{\rho_o}{2}\!, \, h(y_o, \rho_o)} (y_o, \tau_o)} |Du|^2 \, \lambda \, dx dt   \leqslant				\nonumber	\\
& \qquad \leqslant \, \gamma \Bigg[ (2^{\xi} \carlo)^2 \, \mu_+ \left( I_{\frac{\rho_o}{2}, \frac{\rho_o}{2}}^+ (y_o) \right)
					+ (2^{\xi} \carlo)^2 \, \mu_- \left( I^{\frac{\rho_o}{2}, \frac{\rho_o}{2}}_+ (y_o) \right) +				\nonumber	\\
& \qquad \qquad  + \frac{4}{\rho_o^2}
		\iint_{\left(B_{\frac{\rho_o}{2}}^+(y_o)\right)^{\frac{\rho_o}{2}} \times [\tau_o - h(y_o, \rho_o) \frac{\rho_o^2}{2}, \tau_o] \cup 
		\left( I_{\frac{\rho_o}{2}}^+ (y_o) \right)^{\frac{\rho_o}{2}} \times [\tau_o - h(y_o, \rho_o) \rho_o^2, \tau_o]} 
		u^2\, \left( \frac{\mu_+}{h(y_o, \rho_o)} + \la \right) \, dx dt  \Bigg] 	\leqslant									\nonumber	\\
& \qquad \leqslant \, \gamma \Bigg[ (2^{\xi} \carlo)^2 \, \mu_+ \left( I_{\frac{\rho_o}{2}, \frac{\rho_o}{2}}^+ (y_o) \right)
					+ (2^{\xi} \carlo)^2 \, \mu_- \left( I^{\frac{\rho_o}{2}, \frac{\rho_o}{2}}_+ (y_o) \right)\Bigg] +			\nonumber	\\
& \qquad \qquad + \frac{4 \, \gamma}{\rho_o^2} (2^{\xi} \carlo)^2 \Bigg[ 
		\rho_o^2 \, \mu_+ \left( \left(B_{\frac{\rho_o}{2}}^+(y_o)\right)^{\frac{\rho_o}{2}} \right) +
		h(y_o, \rho_o) \, \rho_o^2 \, \lambda \left( \left(B_{\frac{\rho_o}{2}}^+(y_o)\right)^{\frac{\rho_o}{2}} \right)
		\Bigg] 	\leqslant																					\nonumber	\\
& \qquad \leqslant \, \frac{\gamma}{\rho_o^2} \Bigg[ (2^{\xi} \carlo)^2 \, \frac{h(y_o, \rho_o)}{h(y_o, \rho_o)} \, \rho_o^2 \, 
		|\mu| \left( \left(B_{\frac{\rho_o}{2}}^+(y_o)\right)^{\frac{\rho_o}{2}} \right)	\Bigg] +							\nonumber	\\
& \qquad \qquad + \frac{4 \, \gamma}{\rho_o^2} (2^{\xi} \carlo)^2 \Bigg[ 
		\frac{h(y_o, \rho_o)}{h(y_o, \rho_o)} \, \rho_o^2 \, \mu_+ \left( \left(B_{\frac{\rho_o}{2}}^+(y_o)\right)^{\frac{\rho_o}{2}} \right) +
		h(y_o, \rho_o) \, \rho_o^2 \, \lambda \left( \left(B_{\frac{\rho_o}{2}}^+(y_o)\right)^{\frac{\rho_o}{2}} \right)
		  \Bigg] 	\leqslant																					\nonumber	\\
& \qquad \leqslant \frac{9 \, \gamma}{\rho_o^2} (2^{\xi} \carlo)^2 \,  h(y_o, \rho_o) \, \rho_o^2 \, \lambda \big( B_{\rho_o}(y_o) \big)\, .	\nonumber
\end{align}
\textsl{Step 2 - }
The goal of this step is to show the existence of
$\bar{t} \in [\tau_o - h (y_o, \rho_o) {\textstyle \frac{\rho_o^2}{4}}, \tau_o]$ such that
\begin{equation}
\label{giochinipercarlettomio}
\begin{array}{c}
{\displaystyle
\frac{\mu_+ \left( \left\{ x \in B_{\rho_o / 2}^+ (y_o) \, \big| \, u(x,\bar{t}) > \frac{\carlo}{2} \right\} \right)}
	{{|\mu|}_{\lambda} \left( B_{\rho_o / 2} (y_o) \right)} > \frac{\overline{\nu}}{2} \, , }									\\	[1em]
{\displaystyle
\frac{\lambda_+ \left( \left\{ x \in B_{\rho_o / 2}^+ (y_o) \, \big| \, u(x,\bar{t}) > \frac{\carlo}{2} \right\} \right)}
	{\lambda \left( B_{\rho_o / 2} (y_o) \right)}
> 	\frac{\overline{\nu}}{2} \,  ,	}																			\\	[1em]
{\displaystyle
\int_{(B_{\frac{\rho_o}{2}}^+(y_o))^{\frac{\rho_o}{2}}}
	|Du(x,\bar t)|^2 \lambda (x) dx \leqslant  \frac{144 \gamma}{\overline{\nu}} \, (2^{\xi} \carlo)^2 \, \frac{\lambda (B_{\rho_o}(y_o))}{\rho_o^2} \, .
}
\end{array}
\end{equation}
%
%
To this aim we introduce the following sets ($b$ being a positive number to be fixed later)
\begin{gather*}
A^+(t) = \left\{ x\in B_{\rho_o/2}^+(y_o) \, \Big| \, u(x,t) > \frac{\carlo}{2} \right\}	 \, ,\qquad 
	t \in [\tau_o - h (y_o, \rho_o) {\textstyle \frac{\rho_o^2}{4}} , \tau_o]													\\
I^+_{\mu} = \left\{ t \in [\tau_o - h (y_o, \rho_o) {\textstyle \frac{\rho_o^2}{4}} , \tau_o] \, \Big| \, 
	\frac{\mu_+ (A^+(t))}{|\mu|_{\lambda}(B_{\rho_o/2}(y_o))} > \frac{\overline{\nu}}{2} \right\},								\\
J_b=\displaystyle \bigg\{t\in [\tau_o - h (y_o, \rho_o) {\textstyle \frac{\rho_o^2}{4}} , \tau_o] \, \Big| \,
\int_{(B_{\frac{\rho_o}{2}}^+(y_o))^{\frac{\rho_o}{2}}} |Du(x,t)|^2 \lambda(x) dx \leqslant 
	b \, (2^{\xi} \carlo)^2 \, \frac{\lambda (B_{\rho_o}(y_o))}{\rho_o^2}
\bigg\} \, .
\end{gather*}
Using \eqref{mimancanogliamorimiei} we can write
\begin{align*}
\overline{\nu} \, h(y_o, \rho_o) \frac{\rho_o^2}{4} & <  \int_{\tau_o - h (y_o, \rho_o) {\textstyle \frac{\rho_o^2}{4}}}^{\tau_o}
	\frac{\mu_+(A^+(t))}{{{|\mu|}_{\lambda} \left( B_{\rho_o / 2} (y_o) \right)}} \, dt	=										\\
& = \int_{I^+_{\mu}} \frac{\mu_+(A^+(t))}{{{|\mu|}_{\lambda} \left( B_{\rho_o / 2} (y_o) \right)}} \, dt	+
	\int_{[\tau_o - h (y_o, \rho_o) {\textstyle \frac{\rho_o^2}{4}} , \tau_o] \setminus I^+}
	\frac{\mu_+(A^+(t))}{{{|\mu|}_{\lambda} \left( B_{\rho_o / 2} (y_o) \right)}} \, dt	\leqslant									\\
& \leqslant | I^+_{\mu} | + \frac{\overline{\nu}}{2} \, h(y_o, \rho_o) \frac{\rho_o^2}{4}
\end{align*}
by which
\begin{align*}
| I^+_{\mu} | > \frac{\overline{\nu}}{2} \, h(y_o, \rho_o) \frac{\rho_o^2}{4} \, .
\end{align*}
Now from one hand we have \eqref{mimancanogliamorimiei-1}, on the other
\begin{align*}
\int_{[\tau_o - h (y_o, \rho_o) \frac{\rho_o^2}{4}, \tau_o ] \setminus J_{b}} 
	\int_{(B_{\frac{\rho_o}{2}}^+(y_o))^{\frac{\rho_o}{2}}} |Du|^2 \, \lambda \, dx dt   \geqslant
	b \, (2^{\xi} \carlo)^2 \, \frac{\lambda (B_{\rho_o}(y_o))}{\rho_o^2} \, 
	\Big|\Big[\tau_o - h (y_o, \rho_o) \frac{\rho_o^2}{4},\tau_o\Big]\setminus J_b\Big| .
\end{align*}
Then we get
$$
\big| J_b \big| \geqslant h(y_o, \rho_o) \frac{\rho_o^2}{4} \left( 1 - \frac{36 \gamma}{b} \right) .
$$
Choosing $b > 36 \gamma$ this inequality is not trivial. Choosing, e.g., $b = 144 \gamma / \overline{\nu}$ one gets
$$
| I^+_{\mu} \cap J_b | = | I^+_{\mu} | + | J_b | - | I^+_{\mu} \cup J_b | \geqslant  \frac{\overline{\nu}}{4} \, h(y_o, \rho_o) \frac{\rho_o^2}{4} \, .
$$
\ \\ [0.3em]
\textsl{Step 3 - }
Here we show that for every $\bar{\delta} \in (0,1)$ there are $\eta \in (0,1)$ and $y^{\ast} \in B_{\rho_o / 2}^+ (y_o)$,
$\eta = \eta (K_1, K_2, q, K_3, \varsigma, \bar{\delta})$,
$y^{\ast} = y^{\ast} (\gamma, 2^{\xi} N, \overline{\nu}, K_1, K_2, q, K_3, \varsigma, \bar{\delta}) = 
y^{\ast} (\gamma, 2^{\xi} N, \kappa, \gamma_1, \q, K_1, K_2, q, K_3, \varsigma, \bar{\delta})$, such that
$B_{\eta \frac{\rho_o}{2}} (y^{\ast}) \subset B_{\frac{\rho_o}{2}}^+(y_o)$ and such that
\begin{equation}
\label{estMis}
\mu_+ \left(\left\{u(\cdot,\bar t) \leqslant \frac{\carlo}{4}\right\} \cap B_{\eta \frac{\rho_o}{2}}(y^{\ast}) \right) \leqslant
	\bar{\delta} \, \mu_+ (B_{\eta \frac{\rho_o}{2}}(y^{\ast})).
\end{equation}
To see that it is sufficient to use the informations of the previous step and to apply Lemma \ref{lemmaMisVar} to the function $2u/\carlo$ 
with $\omega = \lambda$, $\nu = |\mu|_{\lambda}$, $\varepsilon = 1/2$, $\rho = \rho_o$, $x_0 = y_o$,
$\B = B_{\rho_o / 2}^+ (y_o)$, $\sigma = \rho_o / 2$, $\alpha = \frac{\overline{\nu}}{2}$,
$\beta = \frac{144 \gamma}{\overline{\nu}} \, (2^{\xi} \carlo)^2$ and we get
$$
\mu_+ \left(\left\{u(\cdot,\bar t) > \frac{\carlo}{4}\right\} \cap B_{\eta \frac{\rho_o}{2}}(y^{\ast}) \right) >
	(1 - \bar{\delta}) \, \mu_+ (B_{\eta \frac{\rho_o}{2}}(y^{\ast}))
$$
which is equivalent to \eqref{estMis}.  
Notice that $\eta$ depends on $K_1, K_2, q, K_3, \varsigma$, the constants of the weights, $\bar{\delta}$ and not on the value $N$. \\ [0.3em]
\textsl{Step 4 - }
Here we show that an estimate like that of the third step can be established also in a cylinder. Precisely we show that
for every $\delta \in (0,1)$ there is $\bar{x} \in B_{\eta \frac{\rho_o}{4}}(y^{\ast})$, $\varepsilon \in (0,1)$ which will depend only on $\delta$ and $\q$, and
$s^{\ast} = (\varepsilon \, \eta \, \rho_o/4)^2 \, h(\bar{x}, \varepsilon \eta \frac{\rho_o}{4})$ such that ($\bar{t}, \bar{\delta}, \eta, \rho_o$ as above)
\begin{align}
\label{cilindro_brutto_3}
M_+ \left( \left\{ u \leqslant \frac{\carlo}{8}\right \} \cap \big( B_{\varepsilon \eta \frac{\rho_o}{4}}(\bar{x}) \times [\bar{t}, \bar{t} + s^{\ast}] \big) \right)
	\leqslant \delta \, M_+ \big( B_{\varepsilon \eta \frac{\rho_o}{4}}(\bar{x}) \times [\bar{t}, \bar{t} + s^{\ast}] \big) \, .
\end{align}
Notice that $\bar{x}$ implicitely depends on $y^{\ast}$ and $\delta$ and then $\bar{x}$ depends on
$\gamma, 2^{\xi} N, \kappa, \gamma_1, \q, K_1, K_2, q, K_3, \varsigma, \bar{\delta}, \delta$. \\
To see this we consider $\varepsilon \in (0,1) $ and a disjoint family of balls
$\{B_{\varepsilon \eta \frac{\rho_o}{4}}(x_j)\}_{j=1}^m$ such that
\begin{align*}
& B_{\varepsilon \eta \frac{\rho_o}{4}}(x_j) \subset B_{\eta \frac{\rho_o}{4}}(y^{\ast})) \quad \text{for every } j=1,\ldots, m, \quad \text{ and }		\\
& B_{\eta \frac{\rho_o}{4}} (y^{\ast}) \subset \bigcup_{j=1}^m B_{\varepsilon \eta \frac{\rho_o}{2}}(x_j) \subset B_{\eta \frac{\rho_o}{2}} (y^{\ast})
\end{align*}
and define
$$
s^{\ast}_j := (\varepsilon \, \eta \, \rho_o/4)^2 \, h \left(x_j, \varepsilon \eta \frac{\rho_o}{4} \right) \, .
$$
If necessary one can choose $\varepsilon$ small enough so that $\bar{t} + s^{\ast}_j < T$.
We apply the energy estimate \eqref{DGgamma+_1} to the function $(u - \carlo/4)_-$ in each of the sets
$B_{\varepsilon \eta \frac{\rho_o}{4}}(x_j) \times [\bar{t}, \bar{t} + s^{\ast}_j]$.
Since $B_{\eta \frac{\rho_o}{2}}(y^{\ast}) \subset \Omega_+$ we get 
\begin{align*}
& \sup_{t \in [\bar{t}, \bar{t} + s^{\ast}_j]} \int_{B_{\varepsilon \eta \frac{\rho_o}{2}}(x_j)} \left(u - \frac{\carlo}{4}\right)^2_- (x,t) \mu_+ (x) dx	\leqslant	\\
& \hskip40pt \leqslant \int_{B_{\varepsilon \eta \frac{\rho_o}{2}}(x_j)} \left(u-\frac{\carlo}{4}\right)^2_- (x, \bar{t}) \mu_+ (x) dx +
	\frac{16 \gamma}{\eta^2 \rho_o^2}
	\int_{\bar{t}}^{\bar{t} + s^{\ast}_j} \!\!\!\! \int_{B_{\varepsilon \eta \frac{\rho_o}{2}}(x_j)} \left(u-\frac{\carlo}{4}\right)^2_-\, \la \, dx dt
\end{align*}
and summing over $j$ and using \eqref{estMis}
\begin{align*}
& \sum_{j = 1}^m 
	\sup_{t \in [\bar{t}, \bar{t} + s^{\ast}_j]} \int_{B_{\varepsilon \eta \frac{\rho_o}{2}}(x_j)} \left(u - \frac{\carlo}{4}\right)^2_- (x,t) \mu_+ (x) dx	\leqslant	\\
& \hskip30pt \leqslant \int_{B_{\eta \frac{\rho_o}{2}}(y^{\ast}))} \left(u-\frac{\carlo}{4}\right)^2_- (x, \bar{t}) \mu_+ (x) dx
	+ \sum_{j = 1}^m  \frac{16 \gamma}{\eta^2 \rho_o^2}
	\int_{\bar{t}}^{\bar{t} + s^{\ast}_j} \!\!\!\! \int_{B_{\varepsilon \eta \frac{\rho_o}{2}}(x_j)} \left(u-\frac{\carlo}{4}\right)^2_-\, \la \, dx dt	\leqslant		\\
& \hskip30pt \leqslant
	\frac{\carlo^2}{16} \, \mu_+ \left(\left\{u(\cdot,\bar t) \leqslant \frac{\carlo}{4}\right\} \cap B_{\eta \frac{\rho_o}{2}}(y^{\ast}) \right) +				\\
& \hskip70pt +	\sum_{j = 1}^m \frac{16 \gamma}{\eta^2 \rho_o^2} \, \frac{\carlo^2}{16} \, 
	\varepsilon^2 \, \eta^2 \, \frac{\rho_o^2}{16} \, h \left(x_j, \varepsilon \eta \frac{\rho_o}{4}\right)
	\, \lambda \big(B_{\varepsilon \eta \frac{\rho_o}{2}}(x_j) \big)	\leqslant																\\
& \hskip30pt \leqslant \frac{\carlo^2}{16} \, \bar{\delta} \, \mu_+ \big(B_{\eta \frac{\rho_o}{2}}(y^{\ast}) \big) +
	\sum_{j = 1}^m \frac{16 \, \gamma}{\eta^2 \rho_o^2} \, \frac{\carlo^2}{16} \, 
	\varepsilon^2 \eta^2 \, \frac{\rho_o^2}{16} \, h \left(x_j, \varepsilon \eta \frac{\rho_o}{4}\right)
	\, \q \, \lambda \big(B_{\varepsilon \eta \frac{\rho_o}{4}}(x_j) \big) \leqslant																\\
& \hskip30pt \leqslant \q \, \frac{\carlo^2}{16} \, \bar{\delta} \, |\mu|_{\lambda} \big(B_{\eta \frac{\rho_o}{4}}(y^{\ast})) \big) +
	\frac{\gamma \, \q \, \carlo^2 \varepsilon^2}{16} \, \sum_{j = 1}^m |\mu|_{\lambda} \big(B_{\varepsilon \eta \frac{\rho_o}{4}}(x_j) \big) \leqslant			\\
& \hskip30pt
	\leqslant \q \, \frac{\carlo^2}{16} ( \bar{\delta} + \gamma\, \varepsilon^2 ) \, |\mu|_{\lambda} \big(B_{\eta \frac{\rho_o}{4}}(y^{\ast}) \big) \, .
\end{align*}
On the other side, defining
\begin{gather*}
B_j(t) = \left\{ x \in B_{\varepsilon \eta \frac{\rho_o}{4}}(x_j) \, \bigg| \, u(x,t) \leqslant \frac{\carlo}{8}\right\} \, ,
\end{gather*}
we easily get (for $t \in [\bar{t}, \bar{t} + s^{\ast}_j]$)
$$
\int_{B_{\varepsilon \eta \frac{\rho_o}{4}}(x_j)} \Big( u - \frac{\carlo}{4} \Big)^2_- (x,t) \mu_+ (x) dx \geqslant
\int_{B_j(t)} \Big(u-\frac{\carlo}{2}\Big)^2_-(x,t) \mu_+(x) dx \geqslant \frac{\carlo^2}{64} |\mu|_{\lambda}( B_j(t)) \, .
$$
Now putting together these inequalities we get
\begin{align*}
\frac{\carlo^2}{64} \sum_{j = 1}^m  |M|_{\Lambda} &
\left( \left\{ u \leqslant \frac{\carlo}{8}\right \} \cap \big( B_{\varepsilon \eta \frac{\rho_o}{4}}(x_j) \times [\bar{t}, \bar{t} + s^{\ast}_j] \big) \right) \leqslant 	\\
& \hskip80pt \leqslant \q \, \frac{\carlo^2}{16} ( \bar{\delta} + \gamma\, \varepsilon^2 ) 
	\sum_{j = 1}^m |M|_{\lambda} \Big(B_{\eta \frac{\rho_o}{4}}(y^{\ast}) \times [\bar{t}, \bar{t} + s^{\ast}_j] \Big) \, .
\end{align*}
Once $\delta \in (0,1)$ is chosen we consider $\varepsilon$ and $\bar{\delta}$ in such a way that
$$
4 \, \q \, ( \bar{\delta} + \gamma\, \varepsilon^2 ) \leqslant \delta
$$
and then we get
\begin{align*}
\sum_{j = 1}^m  |M|_{\Lambda}
\left( \left\{ u \leqslant \frac{\carlo}{8}\right \} \cap \big( B_{\varepsilon \eta \frac{\rho_o}{4}}(x_j) \times [\bar{t}, \bar{t} + s^{\ast}_j] \big) \right) \leqslant
	\delta \sum_{j = 1}^m |M|_{\lambda} \Big(B_{\eta \frac{\rho_o}{4}}(y^{\ast}) \times [\bar{t}, \bar{t} + s^{\ast}_j] \Big) \, .
\end{align*}
Notice that $s^{\ast}_j$ depend on $\varepsilon$ and consequently on the choice of $\delta$.
To find a cylinder, independent of $\delta$, in which the estimate above holds true notice that, whatever the choice of $\delta$,
by the last inequality at least one among the $x_j$'s has to satisfy \eqref{cilindro_brutto_3}.
We call $\bar{x}$ that $x_j$ and $s^{\ast} := s^{\ast}_j$ . \\ [0.3em]
\textsl{Step 5 - }
Here we show that
\begin{equation}
\label{emanuela}
u \geqslant \frac{\carlo}{16} \qquad \qquad \text{a.e. in } B_{\varepsilon \eta \frac{\rho_o}{8}}(\bar{x}) \times
	\left[ \bar{t} + \frac{(\varepsilon \, \eta \, \rho_o)^2}{32} \, h \big( \bar{x}, \varepsilon \eta \frac{\rho_o}{4} \big) , 
		\bar{t} + \frac{(\varepsilon \, \eta \, \rho_o)^2}{16} \, h \big( \bar{x}, \varepsilon \eta \frac{\rho_o}{4} \big) \right] \, .
\end{equation}
First notice that $\varepsilon$ depens only on $\delta$ and $\q$.
By \eqref{cilindro_brutto_3} and \eqref{carlettomio} we also get
\begin{equation}
\label{cilindro_brutto_4}
\Lambda \left( \left\{ u \leqslant \frac{\carlo}{8}\right \} \cap \big( B_{\varepsilon \eta \frac{\rho_o}{4}}(\bar{x}) \times [\bar{t}, \bar{t} + s^{\ast}] \big) \right)
	\leqslant \upkappa \, \delta^{\anna} \, \Lambda \big( B_{\varepsilon \eta \frac{\rho_o}{4}}(y^{\ast}) \times [\bar{t}, \bar{t} + s^{\ast}] \big) \, .
\end{equation}
Now we want to apply Proposition \ref{prop-DeGiorgi2}, so first notice that,
by the choice of $\bar{x}$ and $\rho_o$, since $u \geqslant 0$ and by \eqref{oscillazione} we have, choosing $\varepsilon$ even smaller if necessary
so that $\bar{t} + s^{\ast} < \tau_o$, that
$$
\mathop{\rm osc}\limits_{B_{\varepsilon \eta \frac{\rho_o}{4}}(\bar{x}) \times [\bar{t}, \bar{t} + s^{\ast}]} \leqslant 2^{\xi} \carlo \, .
$$
Then taking in Proposition \ref{prop-DeGiorgi2}, point $i\, $), the following values
\begin{gather*}
\underline{m} = 0 , \qquad \omega = 2^{\xi} \carlo , \qquad r = \varepsilon \eta \frac{\rho_o}{8} , \qquad R = \varepsilon \eta \frac{\rho_o}{4} ,		\\
x^{\diamond} = \bar{x} , \qquad t^{\diamond} = \bar{t} , \qquad 
\upbeta^{\diamond} = 1 ,								\\
\sigma = \frac{1}{8} \frac{1}{2^{\xi}} ,  \qquad a = \frac{1}{2} , \qquad \theta^{\diamond} = \frac{1}{2}
\end{gather*}
we have the existence of $\underline{\nu}^{\diamond}$, which in this case depends only on $\kappa, \gamma_1, \gamma$, such that if
$$
\frac{M_+ \left( \left\{ u \leqslant \frac{\carlo}{8}\right \} \cap \big( B_{\varepsilon \eta \frac{\rho_o}{4}}(\bar{x}) \times [\bar{t}, \bar{t} + s^{\ast}] \big) \right)}
	{M_+ \big( B_{\varepsilon \eta \frac{\rho_o}{4}}(\bar{x}) \times [\bar{t}, \bar{t} + s^{\ast}] \big)} +
\frac
{\Lambda \left( \left\{ u \leqslant \frac{\carlo}{8}\right \} \cap \big( B_{\varepsilon \eta \frac{\rho_o}{4}}(\bar{x}) \times [\bar{t}, \bar{t} + s^{\ast}] \big) \right)}
	{\Lambda \big( B_{\varepsilon \eta \frac{\rho_o}{4}}(\bar{x}) \times [\bar{t}, \bar{t} + s^{\ast}] \big)}
\leqslant \underline{\nu}^{\diamond}
$$
then \eqref{emanuela} holds. Then, by \eqref{cilindro_brutto_3} and \eqref{cilindro_brutto_4}, it is sufficient to choose $\delta$ in the fourth step in such a way that
$$
\delta + \upkappa \, \delta^{\anna} = \underline{\nu}^{\diamond}
$$
to get that \eqref{emanuela} holds (so $\delta$ depends only on $\kappa, \gamma_1, \gamma, \upkappa, \anna)$. \\ [0.3em]
\textsl{Step 6 - }
Now, starting from \eqref{emanuela}, we are in the conditions to apply the expansion of positivity. \\
Before going on we recall the dependence of some parameters that are involved (and that we will need):
\begin{align*}
& \eta = \eta (K_1, K_2, q, K_3, \varsigma, \bar{\delta}) = \eta (K_1, K_2, q, K_3, \varsigma, \delta, \anna, \q)
	= \eta (K_1, K_2, q, K_3, \varsigma, \kappa, \gamma_1, \gamma, \upkappa , \q)	,								\\
& \varepsilon = \varepsilon (\bar{\delta}) = \varepsilon (\delta, \q) = \varepsilon (\kappa, \gamma_1, \gamma, \upkappa, \anna, \q) ,		\\
\end{align*}
We call just for simplicity
$$
r := \frac{\varepsilon \, \eta \, \rho_o}{8} \qquad \text{and} \qquad \bar{s} := \bar{t} + 4 \, h \Big( \bar{x}, \varepsilon \eta \frac{\rho_o}{4} \Big) r^2 \, .
$$
In Lemma \ref{esp_positivita} we consider
\begin{gather*}
x^{\ast} = \bar{x} , \qquad t^{\ast} = \bar{t} + \frac{(\varepsilon \, \eta \, \rho_o)^2}{16} \, h \Big(\bar{x}, \varepsilon \eta \frac{\rho_o}{4} \Big) = \bar{s},
\qquad \rho = r  , \qquad h = \frac{\carlo}{16} ,
\end{gather*}
and get that there is $\tilde\upbeta$ depending on $\gamma$ and for every $\hat\theta$ there is $\uplambda > 0$ depending on
$\gamma_1 , \gamma , \q , \kappa, \tilde\upbeta , \hat\theta$ such that
$$
u \geqslant \uplambda \, \frac{\carlo}{16} \qquad \qquad \text{a.e. in } B^+_{2 r}(\bar{x}) \times
	\left[ \bar{s} + \hat\theta \, \tilde\upbeta \, h (\bar{x}, 4 r) r^2, \bar{s} + \tilde\upbeta \, h (\bar{x}, 4 r) r^2 \right] \, .
$$
Since this holds for every $t \in [ \bar{s} + \hat\theta \, \tilde\upbeta \, h (\bar{x}, 4 r) r^2, \bar{s} + \tilde\upbeta \, h (\bar{x}, 4 r) r^2 ]$,
applying again this lemma we reach
$$
u \geqslant \uplambda^2 \, \frac{\carlo}{16} \qquad \text{a.e. in } B^+_{4 r}(\bar{x}) \times
	\left[ \bar{s} + \hat\theta \, \tilde\upbeta \, \big( h (\bar{x}, 4 r) + 4 h (\bar{x}, 8 r) \big) r^2, 
	\bar{s} + \tilde\upbeta \big( h (\bar{x}, 4 r) + 4 h (\bar{x}, 8 r) \big) r^2 \right] .
$$
Iterating this arguement $m$ times we get
$$
u \geqslant \uplambda^m \, \frac{\carlo}{16} \qquad \text{a.e. in }
B^+_{2^m r}(\bar{x}) \times \left[ \bar{s} + \hat\theta \, \tilde\upbeta \, r^2 \sum_{j=1}^m 4^{j-1} h(\bar{x}, 2^{j+2} r), 
	\bar{s} + \tilde\upbeta \, r^2 \sum_{j=1}^m 4^{j-1} h(\bar{x}, 2^{j+2} r) \right] .
$$
%
%
%
%
%
%
%
%
%
%
%
%
Now we define the quantities ($m \in \N$ is still to be fixed)
$$
\left\{
\arrst{1.5}
\begin{array}{l}
{\displaystyle {s}_m := \bar{s} + \hat\theta \, \tilde\upbeta \, r^2 \sum_{j=1}^m 4^{j-1} h(\bar{x}, 2^{j+2} r)	}	\, ,					\\
{\displaystyle {t}_m := \bar{s} + \tilde\upbeta \, r^2 \sum_{j=1}^m 4^{j-1} h(\bar{x}, 2^{j+2} r) }\, .
\end{array}
\right.
$$
Since $\bar{x} \in B_{\rho}(x_o)$ requiring that $2^m r \geqslant 2\rho$ provides that $B_{2^m r}(\bar{x}) \supset B_{\rho}(x_o)$
so we require that $m$ is such that
\begin{equation}
\label{emme}
2 \rho \leqslant 2^m r < 4 \rho \, , \qquad \text{ i.e.} \quad 1 + \log_2\frac{\rho}{r} \leqslant m < 2 + \log_2\frac{\rho}{r} \, .
\end{equation}
What we have still to fix in the times interval is the value of $\hat\theta$ and moreover the values of $b$ and $\xi$.
Now notice that for every $x,y \in \Omega$ and $\varrho > 0$ such that $B_{2\varrho} (x) \subset \Omega$ and $B_{2\varrho} (y) \subset \Omega$ 
and such that $| x - y | < \varrho$ we have
\begin{gather*}
|\mu|_{\lambda} (B_{\varrho}(x)) \leqslant |\mu|_{\lambda} (B_{2\varrho}(y))  \leqslant \q \, \mu|_{\lambda} (B_{\varrho}(y)) 				\\
\lambda (B_{\varrho}(y)) \leqslant \lambda (B_{2\varrho}(x)) \leqslant \q \, \lambda (B_{\varrho}(x))
\end{gather*}
by which we derive
$$
h (x, \varrho) \leqslant \q^2 h (y, \varrho) \, .
$$
Then, using this last estimate, (H.2)$'$ (see also Remark \ref{notaimportante}, point $\mathpzc{C}$) and \eqref{emme} we can estimate
\begin{align*}
\sum_{j=1}^m 4^{j-1} r^2 h(\bar{x}, 2^{j+2} r) & \, \leqslant \q^2 \sum_{j=0}^{m-1} 4^{j} r^2 h(x_o, 2^{j+3} r)	=							\\
& \, = \frac{\q^2}{4^3} \sum_{j=0}^{m-1} (4^{j+3} r^2)^{1 - \alpha} (4^{j+3} r^2)^{\alpha} h(x_o, 2^{j+3} r)  \leqslant							\\
& \, \leqslant \frac{\q^2}{4^3} \sum_{j=0}^{m-1} (4^{j+3} r^2)^{1 - \alpha} {\tilde{K}_2}^2 (4^{m+2} r^2)^{\alpha} h(x_o, 2^{m+2} r)  =				\\
& \, = \frac{\q^2}{4^3} {\tilde{K}_2}^2 (4^{m+2} r^2)^{\alpha} h(x_o, 2^{m+2} r) \sum_{j=0}^{m-1} (4^{3} r^2)^{1 - \alpha}(4^{1 - \alpha})^j  \leqslant	\\
& \, \leqslant \frac{\q^2 \, {\tilde{K}_2}^2}{4 - 4^{\alpha}} \, 4^{m} r^2 h(x_o, 2^{m+2} r) \leqslant												\\
& \, \leqslant \frac{4 \, \q^6 \, {\tilde{K}_2}^2}{4 - 4^{\alpha}} \, \rho^2 h(x_o, \rho)
\end{align*}
by which
\begin{gather*}
s_m \leqslant \bar{s} + \hat\theta \, \tilde\upbeta \, \frac{4 \, \q^6 \, {\tilde{K}_2}^2}{4 - 4^{\alpha}} \, \rho^2 h(x_o, \rho) \, .
\end{gather*}
Now for a fixed constant $\vartheta_+ \in (0,1]$ we can choose
$$
\hat\theta \leqslant \vartheta_+ \, \frac{4 - 4^{\alpha}}{4} \, \frac{1}{\tilde\upbeta \, \q^6 \, {\tilde{K}_2}^2} ,
$$
indipendent of $m$, and, since $\bar{s} < t_o$, we get
\begin{gather*}
s_m < t_o + \vartheta_+ \, \rho^2 h(x_o, \rho) \, .
\end{gather*}
Notice that once $\hat\theta$ is fixed $\uplambda$ depends only on $\gamma_1 , \gamma , \q , \kappa, \tilde\upbeta$. By the choice of $m$
and recalling the definition of $\carlo$ we have
\begin{align*}
u \geqslant \uplambda^m \, \frac{b (\rho - r_o)^{-\xi}}{16} \qquad \text{a.e. in } B^+_{\rho}(x_o) \times \left[ s_m, t_m \right] .
\end{align*}
By the choice we made of $\rho_o$ in \eqref{rozero} we have that
$$
\text{or }  \frac{1}{\rho - r_o} = \frac{1}{2 \, \rho_o} \qquad
	\text{either } \frac{1}{\rho - r_o} = \frac{1}{2^{1+\alpha} \, \q \, \tilde{K}_2} \, \frac{1}{\rho^{\alpha}} \, \frac{1}{\rho_o^{1 - \alpha}}\, .
$$
Then, by the definition of $r$ and since $u(x_o, t_o) = b \, \rho^{-\xi}$, in the first case we get
\begin{align*}
u (x,t) \geqslant \frac{(2^{\xi} \uplambda)^m}{16} \, \frac{b (\varepsilon \, \eta)^{\xi}}{(2^6 \rho)^{\xi}} =
	(2^{\xi} \uplambda)^m \, \frac{(\varepsilon \, \eta)^{\xi}}{2^{6\xi + 4}} \, u(x_o, t_o)
	\qquad \text{a.e. in } B^+_{\rho}(x_o) \times \left[ s_m, t_m \right] .
\end{align*}
In the second we get
\begin{align*}
u (x,t) \geqslant & \, \frac{\uplambda^m}{16} b \left( \frac{1}{2^{1+\alpha} \, \q \, \tilde{K}_2} \right)^{\xi} \, \frac{1}{\rho^{\alpha \xi}}
	\left( \frac{\varepsilon \eta}{8}\right)^{(1-\alpha)\xi} \left( \frac{2^m}{4 \rho}\right)^{(1-\alpha)\xi}							=			\\
= & \, (2^{(1-\alpha)\xi} \uplambda)^m \, b \, \frac{(\varepsilon \, \eta)^{(1-\alpha)\xi}}{(2^{6 - 4\alpha}\q \, \tilde{K}_2)^{\xi}} \, \frac{1}{\rho^{\xi}} =	\\
= & \, (2^{(1-\alpha)\xi} \uplambda)^m \, \frac{(\varepsilon \, \eta)^{(1-\alpha)\xi}}{(2^{6 - 4\alpha}\q \, \tilde{K}_2)^{\xi}} \, u(x_o, t_o)
	\qquad \text{a.e. in } B^+_{\rho}(x_o) \times \left[ s_m, t_m \right] .
\end{align*}
So we can get rid of the dependence of $m$ choosing now $\xi$ in such a way that
\begin{align*}
2^{\xi} \uplambda = 1 & \qquad \text{ in the first case} ,					\\
2^{(1-\alpha)\xi} \uplambda = 1 & \qquad \text{ in the second case} .
\end{align*}
Since $r$ depends on $\rho_o$, which depends on $r_o$, which depends on $\xi$, once we have fixed $\xi$ we have also chosen the value of $r$, and
consequently of $m$. Summing up, we have reached
\begin{align*}
u (x,t) \geqslant {c}_o \, u(x_o, t_o)
	\qquad \quad \text{a.e. in } B^+_{\rho}(x_o) \times \left[ s_m, t_m \right]
\end{align*}
with $s_m < t_o + \vartheta_+ \, \rho^2 h(x_o, \rho)$, where
$$
{c}_o = \frac{(\varepsilon \, \eta)^{\xi}}{2^{6\xi + 4}}
\qquad \text{ or } \qquad
{c}_o = \frac{(\varepsilon \, \eta)^{(1-\alpha)\xi}}{(2^{6 - 4\alpha}\q \, \tilde{K}_2)^{\xi}} .
$$
By the dependence of $\eta$, $\varepsilon$ and $\xi$ and since $\tilde{K}_2$ depends only on $K_2$ we have that
$$
c_o \quad \text{ depends on } \qquad \gamma_1, \gamma, \q, \kappa, \tilde\upbeta, \alpha, \upkappa, \anna, K_1, K_2, K_3, q, \varsigma \, .
$$
Now we are done if $t_m \geqslant t_o + \vartheta_+ \, \rho^2 h(x_o, \rho)$ and the constant $c_+$ is ${c}_o$. \\
If, otherwise, $t_m < t_o + \vartheta_+ \, \rho^2 h(x_o, \rho)$ we consider
$$
\hat{t} \in [s_m , t_m] \qquad \text{such that } \quad \hat{t} + \hat\theta \, \tilde\upbeta \, h(x_o, 4\rho) \rho^2 
	\leqslant t_o + \vartheta_+ h(x_o, \rho) \rho^2 \, .
$$
By \eqref{stimeacca} this is true, taking if necessary $\hat\theta$ smaller, if
$$
\hat\theta \leqslant \frac{\vartheta_+}{\q^2 \tilde\upbeta} \, .
$$
Applying again Lemma \ref{esp_positivita} and since $u(x,t) \geqslant c_o \, u(x_o, t_o)$ a.e. in $B^+_{\rho}(x_o) \times \left[ s_m, t_m \right]$
(and then also in $B^+_{\rho/4}(x_o) \times \left[ s_m, t_m \right]$)
we get, in particular, that both
$$
u(x,t) \geqslant \uplambda \, c_o \, u(x_o, t_o)	\qquad \text{a.e. in } B^+_{2\rho}(x_o) \times 
	\left[ \hat{t} + \hat\theta \, \tilde\upbeta \, h(x_o, 4\rho) \rho^2, \hat{t} + \tilde\upbeta \, h(x_o, 4\rho) \rho^2  \right]
$$
and
$$
u(x,t) \geqslant \uplambda \, c_o \, u(x_o, t_o)	\qquad \text{a.e. in } B^+_{\rho/2}(x_o) \times 
	\left[ \hat{t} + \hat\theta \, \tilde\upbeta \, h(x_o, \rho) \rho^2/ 16, \hat{t} + \tilde\upbeta \, h(x_o, \rho) \rho^2/16  \right] \, ;
$$
then in particular
$$
u(x,t) \geqslant \uplambda \, c_o \, u(x_o, t_o)	\qquad \text{a.e. in } B^+_{\rho}(x_o) \times 
	\left[ \hat{t} + \hat\theta \, \tilde\upbeta \, h(x_o, \rho) \rho^2/ 16, \hat{t} + \tilde\upbeta \, h(x_o, \rho) \rho^2/16  \right] \, .
$$
Repeating this argument for every $t$ in 
$\left[ \hat{t} + \hat\theta \, \tilde\upbeta \, h(x_o, \rho) \rho^2/ 16, \hat{t} + \tilde\upbeta \, h(x_o, \rho) \rho^2/16  \right]$
we get
$$
u(x,t) \geqslant \uplambda^2 \, c_o \, u(x_o, t_o)	\qquad \text{a.e. in } B^+_{\rho}(x_o) \times 
	\left[ \hat{t} + 2 \hat\theta \, \tilde\upbeta \, h(x_o, \rho) \rho^2 / 16, \hat{t} + 2 \tilde\upbeta \, h(x_o, \rho) \rho^2 / 16 \right] \, .
$$
If necessary, we add the requirement $2\hat\theta < 1$ so that
$\left[ \hat{t} + 2 \hat\theta \, \tilde\upbeta \, h(x_o, \rho) \rho^2 / 16, \hat{t} + 2 \tilde\upbeta \, h(x_o, \rho) \rho^2 / 16 \right] \cap
\left[ \hat{t} + \hat\theta \, \tilde\upbeta \, h(x_o, \rho) \rho^2/ 16, \hat{t} + \tilde\upbeta \, h(x_o, \rho) \rho^2/16  \right] \not= \emptyset$.
Going on, we get
$$
u(x,t) \geqslant \uplambda^3 \, c_o \, u(x_o, t_o)	\qquad \text{a.e. in } B^+_{\rho}(x_o) \times 
	\left[ \hat{t} + 3 \, \hat\theta \, \tilde\upbeta \, h(x_o, \rho) \rho^2 / 16, \hat{t} + 3 \tilde\upbeta \, h(x_o, \rho) \rho^2 / 16 \right]
$$
requiring $3 \hat \theta < 2$, which is free since we already imposed $2\hat\theta < 1$. 
We iterate $k$ times, without additional assumptions about $\hat\theta$,
till $\hat{t} + k \, \tilde\upbeta \, h(x_o, \rho) \rho^2 / 16 > t_o + \vartheta \, h(x_o, \rho) \rho^2$ and get
$$
u(x,t) \geqslant \uplambda^k \, c_o \, u(x_o, t_o)	\qquad \text{a.e. in } B^+_{\rho}(x_o) \times 
	\left[ \hat{t} + k \, \hat\theta \, \tilde\upbeta \, h(x_o, \rho) \rho^2 / 16, \hat{t} + k \, \tilde\upbeta \, h(x_o, \rho) \rho^2 /16  \right] \, .
$$
Since $t_o - \hat{t} > h(x_o, \rho) \rho^2$, the inequality
$$
\hat{t} + \frac{k \, \tilde\upbeta}{16} \, h(x_o, \rho) \rho^2 > t_o + \vartheta_+ h(x_o, \rho) \rho^2
$$
holds if we choose
$$
k > \frac{16}{\tilde\upbeta} \, (1 + \vartheta_+ ) \, .
$$
For instance we can choose $[\frac{16}{\tilde\upbeta} \, (1 + \vartheta_+)] + 1$, the minimum integer greater that $\frac{16}{\tilde\upbeta} \, (1 + \vartheta_+)$
and the constant $c_+$ is $\uplambda^k c_o$, where $k$ depends only on $\tilde\upbeta$ and $\vartheta_+$.
Since $\tilde\upbeta$ depends only on $\gamma$ we conclude that $c_+$ depends (only) on
$$
\gamma_1, \gamma, \q, \kappa, \alpha, \upkappa, \anna, K_1, K_2, K_3, q, \varsigma , \vartheta_+ \, .
$$
\ \\ [0.3em]
In a complete analogous way one can prove point $ii \, )$. \\ [0.3em]
We see now point $iii \, )$. Since $s_1$ and $s_2$ will remain fixed in the following we will use the simplified notations, for some $r > 0$ and $\bar{x} \in \Omega$,
$$
Q^{0}_{r} (\bar{x}) := B^0_r  (\bar{x}) \times [s_1, s_2] \, , \qquad
Q_{r} (\bar{x}) := B_r  (\bar{x}) \times [s_1, s_2] \, .
$$
Similarly as for point $i \, )$, we may write $u(x_o, t_o) = b \, \rho^{-\xi}$ for some $b, \xi > 0$ to be fixed later.
Define the functions
$$
\emm (r) = \sup_{Q^{0}_{r} (x_o)} u,		\qquad		\enn (r) = b (\rho - r)^{-\xi}, 	\qquad r \in [0,\rho) .
$$
Let us denote by $r_o \in [0,\rho)$ the largest solution of $\emm (r) = \enn (r)$.
Define
$$
\carlo := \enn (r_o) = b (\rho - r_o)^{-\xi} \, .
$$
We can find $y_o \in B^{0}_{r_o} (x_o)$ such that
\begin{equation}
\label{choicey_o}
\frac{3\carlo}{4} < \sup_{Q^{0}_{\frac{\rho_o}{4}} (y_o)} u \leqslant \carlo
\end{equation}
where $\rho_o = (\rho - r_o) / 2$, so $B^0_{\rho_o} (y_o) \subset B_{\frac{\rho + r_o}{2}} (x_o)$.
By this choice of $\rho_o$ and by the choice of $r_o$ we have
\begin{equation}
\label{oscillazione}
\sup_{Q^{0}_{\rho_o} (y_o)} u \leqslant \sup_{Q^{0}_{\frac{\rho + r_o}{2}}(x_o)} u < \enn \left( \frac{\rho + r_o}{2}\right) = 2^{\xi} \carlo .
\end{equation}
We now proceed dividing the proof in four steps. \\ [0.3em]
\textsl{Step 1 - }
In this step we want to show that there is $\overline{\nu} \in (0,1)$, depending on $\kappa, \gamma_1, \gamma$, such that
\begin{align*}
\Lambda_0 \left( \left\{ u > \frac{\carlo}{2} \right\} \cap  Q^{0}_{\rho_o/2}  (y_o) \right)
> 	\overline{\nu}	\, \Lambda \left( Q_{\rho_o/2}  (y_o) \right)
\end{align*}
and that
\begin{equation}
\label{chesonno!}
\iint_{Q^{0}_{\rho_o/2}  (y_o)} |Du|^2 \, \lambda \, dx dt   \leqslant
\gamma  \, (2^{\xi} \carlo)^2 \left( \frac{2 \, K_2^2 \, \q^2}{\lela} + 4 \right) (s_2 - s_1) \,  \frac{\lambda \big( B_{\rho_o} (y_o) \big)}{\rho_o^2} \, .
\end{equation}
Arguing by contradiction we immediatly get the first inequality: indeed if that were false,
setting in Proposition \ref{prop-DeGiorgi1}, point $iii \, )$,
\begin{gather*}
\overline{m} = \omega = 2^\xi \carlo, \quad R = \frac{\rho_o}{2},		\quad	\rho = \frac{\rho_o}{4},		\quad
							\sigma = 1 - 2^{-\xi-1}, 	\quad 	a = \sigma^{-1}\biggl(1-\frac{3}{2^{\xi+2}}\biggr) \, , 			\\
x^{\star} = y_o \, , \qquad  t^{\star} = t_o	 \, ,  \qquad \upbeta^{\star} = \frac{8 (s_2 - t_o)}{\rho_o} ,
				\qquad s_1^{\star} = s_1 \, , \qquad s_2^{\star} = s_2 \, ,
\end{gather*}
we would get that
$$
u \leqslant \frac{3\carlo}{4} \quad \textrm{in }\, B_{\rho_o/4}^0 (y_o) \times \left(s_1, s_2 \right)
$$
which contradicts \eqref{choicey_o}.
To prove \eqref{chesonno!} we choose in \eqref{DGgamma0}
$x_0 = y_o$, $R = \rho_o$, $\tilde{r} = \rho_o$, $r = \rho_o / 2$, $\varepsilon = 0$, 
$k = 0$ and since $u \leqslant 2^{\xi} \carlo$ we get
\begin{align*}
& \iint_{Q^{0}_{\rho_o/2}  (y_o)} |Du|^2 \, \lambda \, dx dt   \leqslant																\\
& \qquad \leqslant \, \gamma \Bigg[ 
		(2^{\xi} \carlo)^2 \, |\mu| \left( I^{\frac{\rho_o}{2}, \frac{\rho_o}{2}}_0 (y_o) \right) +
		\frac{4}{\rho_o^2} \iint_{\left(B_{\frac{\rho_o}{2}}^0(y_o)\right)^{\frac{\rho_o}{2}} \times [s_1, s_2]} 
		u^2\, \lambda \, dx dt  \Bigg] 	\leqslant																				\\
& \qquad \leqslant \, \gamma \Bigg[ 
		(2^{\xi} \carlo)^2 \, |\mu| \left( I^{\frac{\rho_o}{2}, \frac{\rho_o}{2}}_0 (y_o) \right) +
		(2^{\xi} \carlo)^2 \, \frac{4}{\rho_o^2} \, 2 \, \lela \, h(x_o, 4\rho) \rho^2 \, \lambda \big( B_{\rho_o} (y_o) \big) \Bigg] \leqslant		\\
& \qquad \leqslant \, \gamma  \, (2^{\xi} \carlo)^2
		\left[ |\mu|_{\lambda} \big( B_{\rho_o} (y_o) \big) + \frac{8 \, \lela}{\rho_o^2} \, h(x_o, 4 \rho) \rho^2 
		\, \lambda \big( B_{\rho_o} (y_o) \big) \right]	=																		\\
& \qquad = \, \gamma  \, (2^{\xi} \carlo)^2
		\left[ h(y_o, \rho_o) + \frac{8 \, \lela}{\rho_o^2} \, h(x_o, 4 \rho) \rho^2 \right] \lambda \big( B_{\rho_o} (y_o) \big)  \leqslant	\\
& \qquad \leqslant \, \gamma  \, (2^{\xi} \carlo)^2
		\left[ 4 \, K_2^2 \, \q^2 \, h(x_o, 4 \rho) \frac{\rho^2}{\rho_o^2} + 
			\frac{8 \, \lela}{\rho_o^2} \, h(x_o, 4 \rho) \rho^2 \right] \lambda \big( B_{\rho_o} (y_o) \big)  \, .
\end{align*}
\\ [0.3em]
\textsl{Step 2 - }
Here we show that for every $\delta \in (0,1)$ there are $\eta \in (0,1)$,
which will depend only on $K_1, K_2, K_3, q, \varsigma, \delta$,
and $y^{\ast} \in B_{\rho_o / 2}^0 (y_o)$,
which will depend only on $\gamma, 2^{\xi} N, \bar{\nu}, K_1, K_2, K_3, q, \varsigma, \lela, \delta$
($\delta$ will be chosen depending on $\kappa, \gamma_1, \gamma, \lela \, h(x_o, 4\rho))$), such that
$B_{\eta \rho_o/2}(y_o) \subset B_{\rho_o / 2}^0 (y_o)$ and
\begin{equation}
\label{cilindro_bello}
\Lambda \left( \left\{ u \leqslant \frac{\carlo}{4}\right \} \cap \big( Q_{\eta \frac{\rho_o}{2}} (y^{\ast}) \big) \right) 
												\leqslant \delta \, \Lambda \big( Q_{\eta \frac{\rho_o}{2}} (y^{\ast}) \big) \, .
\end{equation}
Indeed by \textsl{Step 1} and applying Corollary \ref{corollario3} to the function $2u/\carlo$ 
with $\omega = \nu = \lambda$, $\varepsilon = 1/2$, $\rho = \rho_o$, $x_0 = y_o$,
$\B = B_{\rho_o / 2}^0 (y_o)$, $\sigma = \rho_o / 2$, $a = s_1$, $b = s_2$,
$\alpha = \overline{\nu}$, $\beta = \gamma  \, (2^{\xi} \carlo)^2 \left( 2 \, K_2^2 \, \q^2 \, \lela^{-1} + 4 \right)$
we get the existence of $B_{\eta \frac{\rho_o}{2}} (y^{\ast}) \subset B_{\rho_o / 2}^0 (y_o)$ such that
$$
\Lambda \left( \left\{ u > \frac{\carlo}{4}\right \} \cap \big( Q_{\eta \frac{\rho_o}{2}} (y^{\ast}) \big) \right)
	> (1 - \delta) \, \Lambda \big( Q_{\eta \frac{\rho_o}{2}} (y^{\ast}) \big)
$$
which is equivalent to \eqref{cilindro_brutto_3}. \\ [0.3em]
\textsl{Step 3 - } Here we show that
\begin{equation}
\label{emanuela}
u \geqslant \frac{\carlo}{8} \qquad \qquad \text{a.e. in } Q_{\eta \frac{\rho_o}{4}} (y^{\ast}) \, .
\end{equation}
Now we want to apply Proposition \ref{prop-DeGiorgi2} so first notice that,
since $u \geqslant 0$ and by \eqref{oscillazione} we have that
$$
\mathop{\rm osc}\limits_{Q_{\eta \frac{\rho_o}{2}} (y^{\ast})} \leqslant 2^{\xi} \carlo \, .
$$
Then taking in Proposition \ref{prop-DeGiorgi2}, point $iii\, $), the following values
\begin{gather*}
\underline{m} = 0 , \qquad \omega = 2^{\xi} \carlo , \qquad r = \eta \frac{\rho_o}{4} , \qquad R = \eta \frac{\rho_o}{2} ,			\\
x^{\star} = y^{\ast} , \qquad t^{\star} = t_o , \qquad s_1^{\star} = s_1\qquad s_2^{\star} = s_2 ,							\\
\upbeta^{\star} = 8 \, \lela \, h(x_o, 4\rho) \, \frac{\rho^2}{\eta^2 \rho_o^2} ,
\qquad	\sigma = \frac{1}{4} \frac{1}{2^{\xi}} ,  \qquad a = \frac{1}{2} 
\end{gather*}
we have the existence of $\underline{\nu}^{\star} \in (0,1)$, which in this case depends only on $\kappa, \gamma_1, \gamma, \lela \, h(x_o, 4\rho)$, such that if
$$
\Lambda \left( \left\{ u \leqslant \frac{\carlo}{4}\right \} \cap \big( Q_{\eta \frac{\rho_o}{2}} (y^{\ast}) \big) \right) 
										\leqslant \underline{\nu}^{\star} \, \Lambda \big( Q_{\eta \frac{\rho_o}{2}} (y^{\ast}) \big) \, .
$$
then \eqref{emanuela} holds. Then it is sufficient to choose $\delta = \underline{\nu}^{\star}$ 
(so $\delta$ depends only on $\kappa, \gamma_1, \gamma, \lela \, h(x_o, 4\rho))$. \\ [0.3em]
\textsl{Step 4 - } Now we want to apply the expansion of positivity. We call, for simplicity
$$
r := \eta \frac{\rho_o}{4} \, .
$$
Taking in Lemma \ref{esp_positivita}, point $iii\, )$,
$$
\rho =  r \, , \qquad \upbeta = \upomega
$$
we get that
$$
u \geqslant \uplambda \, \frac{\carlo}{8} \qquad \qquad \text{a.e. in } B^0_{2 r}(y^{\ast}) \times \left[ s_1, s_2 \right]
$$
with $\uplambda$ depending on $\gamma_1, \gamma, \q, \kappa, \upomega$. Now taking in Lemma \ref{esp_positivita}, point $iii\, )$,
$$
\rho =  2r \, , \qquad \upbeta = \upomega
$$
we get that
$$
u \geqslant \uplambda^2 \, \frac{\carlo}{8} \qquad \qquad \text{a.e. in } B^0_{2 r}(y^{\ast}) \times
	\left[ s_1, s_2 \right] \, .
$$
We iterate this argument $m$ times getting
$$
u \geqslant \uplambda^m \, \frac{\carlo}{8} \qquad \qquad \text{a.e. in } B^0_{2^m r}(y^{\ast}) \times
	\left[ s_1, s_2 \right]
$$
till $B_{2^m r}(y^{\ast}) \supset B_{\rho}(x_o)$ and this is guaranteed if
$$
2 \rho \leqslant 2^m r < 4 \rho \, . 
$$
As done before, observe that
\begin{align*}
u (x,t) & \geqslant \uplambda^m \, \frac{\carlo}{8} = \uplambda^m \, \frac{b \, \eta^\xi}{8^{\xi+1}} \, \frac{2^{m\xi}}{(2^m r)^{\xi}}
	\geqslant \uplambda^m \, \frac{b \, \eta^\xi}{8^{\xi+1}} \, \frac{2^{m \xi}}{(4 \, \rho)^{\xi} } =								\\
   &  = (2^{\xi} \uplambda)^m \, \frac{\eta^\xi}{2^{5\xi + 3}} \, u(x_o, t_o) \, .
\end{align*}
Then, as before, choosing $\xi$ in such a way $2^{\xi} \uplambda = 1$ we get rid of the dependence of $m$ and then in particular we get
$$
u (x,t) \geqslant c_0 \, u(x_o, t_o) \qquad \qquad \text{a.e. in } B^0_{\rho}(x_o) \times
	\left[ s_1, s_2 \right]
$$
where $c_0 = \frac{\eta^\xi}{2^{5\xi + 3}}$ depends (only) on
$K_1, K_2, K_3, q, \varsigma, \kappa, \gamma_1, \gamma, \lela, h(x_o, 4\rho), \q$, the constants by which $\uplambda$ and $\eta$ depend. \\ [0.3em]
Finally the proof of point $iv \, )$ can be obtained similarly to that of point $iii \, )$, using in the order Proposition \ref{prop-DeGiorgi1}, point $iv \, )$,
Lemma \ref{lemmaMisVar}, Proposition \ref{prop-DeGiorgi2}, point $iv \, )$, Lemma \ref{esp_positivita}, point $iv \, )$.
\finedimo

\noindent
The previous theorem has an immediate consequence which we state here below.

\begin{theorem}
\label{Harnack2}
Assume $u\in DG(\Omega, T, \mu, \la, \gamma)$, $u\geqslant 0$. Fix $\rho > 0$ and $\vartheta \in (0,1]$ for which $B_{5\rho}(x_o) \times 
[t_o - 16 \, h(x_o, 4\rho) \rho^2 - \vartheta h(x_o, \rho) \rho^2, t_o + 16 \, h(x_o, 4\rho) \rho^2 + \vartheta h(x_o, \rho) \rho^2] \subset \Omega \times (0,T)$.
Suppose $x_o \in I$.
Then there exists $c > 0$ depending on
$\gamma_1, \gamma, \q, \kappa, \alpha, \upkappa, \anna, K_1, K_2, K_3, q, \varsigma , \vartheta$ such that
$$
u(x_o, t_o)\leqslant c \inf_{B_{\rho} (x_o)} \tilde{u} (x)
$$
where
\begin{align*}
\tilde{u} (x) =
	\left\{
	\begin{array}{ll}
	u(x, t_o + \vartheta \, h(x_o, \rho) \rho^2)			&	\text{ if } x \in B_{\rho}^+(x_o)		\\	[0.3em]
	u(x, t_o - \vartheta \, h (x_o, \rho) \rho^2)			&	\text{ if } x \in B_{\rho}^-(x_o)
	\end{array}
	\right.
& \qquad \text{if } \quad x_o \in I_+ \cap I_- , 																		\\
\tilde{u} (x) =
	\left\{
	\begin{array}{ll}
	u(x, t_o + \vartheta \, h(x_o, \rho) \rho^2)			&	\text{ if } x \in B_{\rho}^+(x_o)		\\	[0.3em]
	u(x, t_o)										&	\text{ if } x \in B_{\rho}^0(x_o)
	\end{array}
	\right.
& \qquad \text{if } \quad x_o \in I_+ \cap I_0 ,																		\\
\tilde{u} (x) =
	\left\{
	\begin{array}{ll}
	u(x, t_o - \vartheta \, h (x_o, \rho) \rho^2)			&	\text{ if } x \in B_{\rho}^-(x_o)		\\	[0.3em]
	u(x, t_o)										&	\text{ if } x \in B_{\rho}^0(x_o)
	\end{array}
	\right.
& \qquad \text{if } \quad x_o \in I_- \cap I_0 ,																		\\
\tilde{u} (x) =
	\left\{
	\begin{array}{ll}
	u(x, t_o + \vartheta \, h(x_o, \rho) \rho^2)			&	\text{ if } x \in B_{\rho}^+(x_o)		\\	[0.3em]
	u(x, t_o - \vartheta \, h (x_o, \rho) \rho^2)			&	\text{ if } x \in B_{\rho}^-(x_o)		\\	[0.3em]
	u(x, t_o)										&	\text{ if } x \in B_{\rho}^0(x_o) .
	\end{array}
	\right.
& \qquad \text{if } \quad x_o \in I_+ \cap I_- \cap I_0 \, .
\end{align*}
\end{theorem}
\noindent
\dimo
By Theorem \ref{Harnack1} we immediately get the result taking $\vartheta = \vartheta_+ = \vartheta_-$ and
$c = \max\{ c_+, c_-, c_0 \}$.
\finedimo

\noindent
One can give many different and equivalent formulations of the classical parabolic Harnack's inequality.
We conclude giving only one possible equivalent formulation, which can be proved by standard arguements, to the one given above.
Under the assumptions of Theorem \ref{Harnack2} one has for $u \in DG$, $u \geqslant 0$, and for instance
for $x_o \in \partial\Omega_+ \cap \partial\Omega_0 \cap \partial\Omega_-$ (and with obvious generalization in the other cases)
\begin{align}
\label{equivalent}
& \sup_{B_{\rho} (x_o)} \tilde{u} (x) \leqslant c \inf_{B_{\rho} (x_o)} u (x, t_o)							\nonumber	\\
& \text{where} \quad \tilde{u} (x) =
	\left\{
	\begin{array}{ll}
	u(x, t_o - \vartheta \, h(x_o, \rho) \rho^2)			&	\text{ if } x \in B_{\rho}^+(x_o)		\\	[0.3em]
	u(x, t_o + \vartheta \, h (x_o, \rho) \rho^2)			&	\text{ if } x \in B_{\rho}^-(x_o)		\\	[0.3em]
	u(x, t_o)										&	\text{ if } x \in B_{\rho}^0(x_o) .
	\end{array}
	\right.				
\end{align}
\ \\

\noindent
{\bf Some consequences of the Harnack inequality - }
An important and standard consequence for a function satisfying a Harnack's inequality is H\"older-continuity. By classical computations and
assuming (if necessary taking $\gamma$ bigger)
$$
\frac{\gamma}{\gamma - 1} < 2
$$
one can get that if $u \in DG(\Omega, T, \mu, \la, \gamma)$ then $u$ is locally $\alpha$-H\"older continuous with respect to $x$ and
$\alpha/2$-H\"older continuous with respect to $t$, where $\alpha = (\log_2 \frac{\gamma}{\gamma - 1})$,
in $\big( \Omega_+ \cup \Omega_- \cup I \big) \times (0,T)$. As regards $\Omega_0$ we can only get that for every $t \in (0,T)$
$u ( \cdot, t)$ is locally $\alpha$-H\"older continuous in $\Omega_0$.
Notice that in the interface $I$ separating $\Omega_0$ and $\Omega_+ \cup \Omega_-$ the function $u$ is regular also with respect to $t$. \\ [0.3em]
Another consequence is a strong maximum pronciple, which one can get, again by standard arguement, using \eqref{equivalent}.
One can derive a ``standard'' maximum principle from Theorem \ref{Harnack1}, which we do not state, and others from Theorem \ref{Harnack2}. \\
If, for instance, we suppose $x_o \in \partial\Omega_+ \cap \partial\Omega_0 \cap \partial\Omega_-$
(and again with obvious generalization in the other cases)
we could briefly state the maximum principles as follows: suppose $(x_o, t_o) \in \Omega \times (0,T)$ is a maximum point for $u$ in a set
\begin{align*}
& \Big( B_{\rho}^+ (x_o) \times (t_o - \vartheta \, h(x_o, \rho) \rho^2, t_o + \vartheta \, h(x_o, \rho) \rho^2) \Big)
	\cup \Big( B_{\rho}^0 (x_o) \times \{ t_o \} \Big) \, \cup															\\
& \quad \quad \Big( \cup B_{\rho}^- (x_o) \times (t_o - \vartheta \, h(x_o, \rho) \rho^2, t_o + \vartheta \, h(x_o, \rho) \rho^2) \Big)
\end{align*}
for some $\vartheta \in (0,1]$, then $u$ is constant in the set
\begin{align*}
\Big( B_{\rho}^+ (x_o) \times (t_o - \vartheta \, h(x_o, \rho) \rho^2, t_o] \Big)
	\cup \Big( B_{\rho}^0 (x_o) \times \{ t_o \} \Big) \, \cup
	\Big( \cup B_{\rho}^- (x_o) \times [t_o, t_o + \vartheta \, h(x_o, \rho) \rho^2) \Big) \, .
\end{align*}
\ \\

\newpage
\section{Examples}
\label{paragrafo9}

In this section we show some possible examples of $\mu$ (and consequently of $I$) and $\lambda$.
In all the examples, just for simplicity, we suppose $\Omega \subset \R^2$.

\begin{itemize}[itemsep=1ex, leftmargin=0.62cm]
\item[1.]
In the simplest situation when $\mu \equiv \lambda \equiv 1$ we get the classical case in which the De Giorgi class contains the solutions of
$$
\frac{\partial u}{\partial t} - \textrm{div} (a(x,t,u,Du)) = b(x,t,u,Du)
$$
with $a, b$ satisfying
\begin{align*}
\big(a (x,t,u,Du) , Du \big) \geqslant \lambda  |Du|^p			\, , 					\\
| a (x,t,u,Du) | \leqslant \Lambda |Du|^{p-1}				\, , 					\\
| b (x,t,u,Du) | \leqslant M |Du|^{p-1}						\, ,
\end{align*}
with $\lambda, \Lambda, M$ positive numbers.
Obviously if $\mu \equiv -1$ we have the analogous results for backward parabolic equations.
\item[2.]
If $\mu \equiv 0$ and $\lambda \equiv 1$ we have a family (in the parameter $t$) of elliptic equations for which
one cannot expect regularity in time, neither for ``solutions''. The same may happen if $\Omega_0$ is a proper subset of $\Omega$. \\
For example, in dimension $1$ consider the solutions of
$$
\frac{d}{dx} \left( a(x,t) \frac{du}{dx} \right) = 0 \, , \qquad u(0)=0 \, , \ u(2) = 1 \, ,
$$
with
$$
a(x,t) = \alpha (t) \text{ in } [0,1] \qquad \text{and} \qquad a(x,t) = \beta (t) \text{ in } [1,2]
$$
with $\alpha(t) \not= \beta(t)$ for every $t$ and $\alpha$ and $\beta$ discontinuous. The solutions are clearly discontinuous in time for $x \in (0,2)$.
\item[3.]
If $\mu > 0$ and $\lambda > 0$ we have the Harnack's inequality for doubly weighted equations, like for instance
\begin{equation}
\label{ultimavolta}
\mu \frac{\partial u}{\partial t} - \textrm{div} (\lambda Du) = 0  \, .
\end{equation}
In the particular case $\mu \equiv 1$ we rescue the result contained in \cite{chia-se3}
(and also contained in \cite{surnachev}),
while if $\mu \equiv \lambda$ we rescue the result contained in \cite{chia-se2}.
\item[4.]
Consider now an example where for simplicity $|\mu| \equiv \lambda \equiv 1$ in $\Omega$, but $\mu \not\equiv1$.
Suppose, for instance, that $\mu$ changes sign around an interface like that in the first of the two following pictures
where $I$ is a cross intersecting in a point $x_o$. This kind of interface clearly satisfies assumptions (H.4) and (H.5) and
then also in a neighbourhood of the points $(x_o, t)$, $t \in (0,T)$, the solution, e.g., of \eqref{ultimavolta} is H\"older-continuous. \\
Also an interface like that shown in the second of the two following pictures is admitted.
\ \\
\ \\
\ \\
\ \\
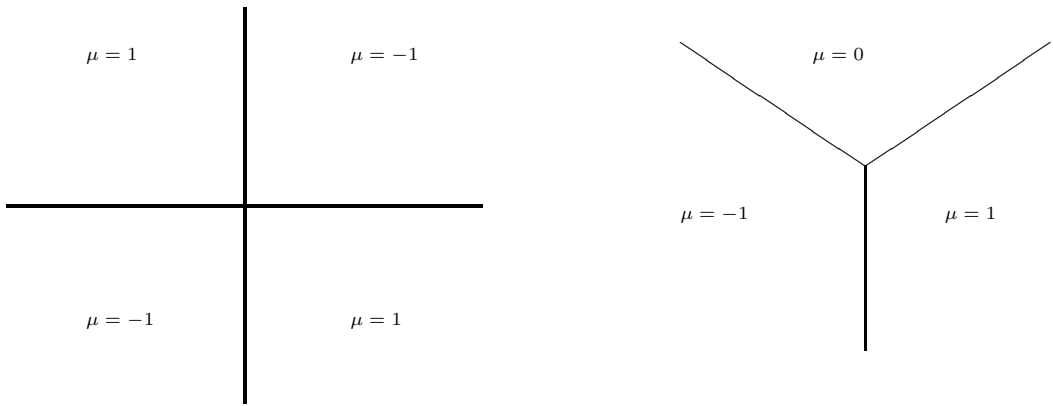
\begin{picture}(150,200)(-180,0)
\hspace{-3.5cm}
\put (-40,180){\tiny$\mu =1$}
\put (60,80){\tiny$\mu =1$}
\put (60,180){\tiny$\mu = -1$}
\put (-40,80){\tiny$\mu = -1$}

\put (20,50){\linethickness{1pt}\line(0,1){150}}

\put (-70,125){\linethickness{1pt}\line(1,0){180}}

\end{picture}
\begin{picture}(150,200)(-180,0)
\hspace{-1cm}
\put (10,180){\tiny$\mu =0$}
\put (60,120){\tiny$\mu =1$}
\put (-40,120){\tiny$\mu = -1$}

\put (30,70){\linethickness{1pt}\line(0,1){70}}
\put (30,140){\linethickness{1pt}\line(3,2){70}}
\put (30,140){\linethickness{1pt}\line(-3,2){70}}

\end{picture}

\item[5.]
Consider $\mu \geqslant 0$ and, for simplicity, suppose that $\mu$ takes only the values $1$ and $0$.
in the pictures below there are two simple examples: in the first one the interface is made by just a line, in the second one is made by two
intersecting lines. In both cases a function belonging to the De Giorgi class turns out to be H\"older-continuous in $(\Omega_+ \cup I ) \times (0,T)$. 
In particular it is continuous in the interface $I$ both in $x$ and $t$, even if it could not be continuous in $\Omega_0$ as shown in the second example.
\ \\
\ \\
\ \\
\ \\
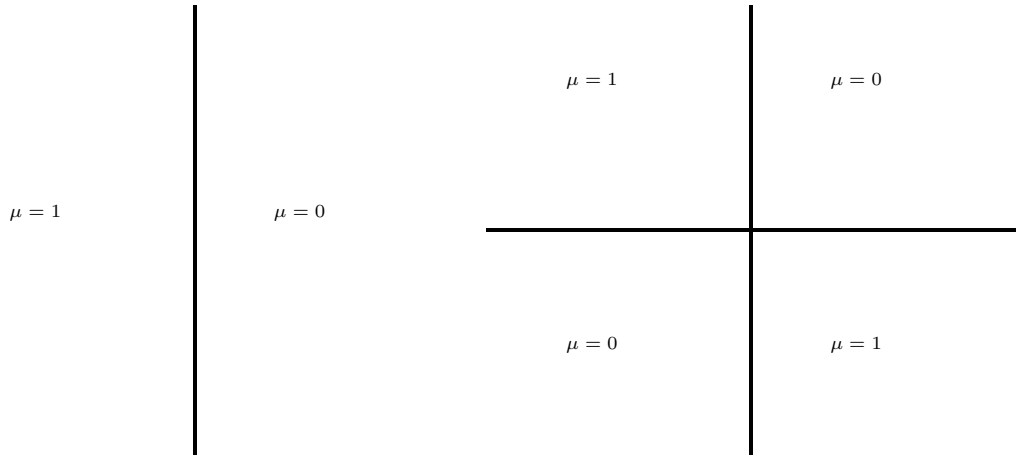
\begin{picture}(150,200)(-180,0)
\hspace{-4cm}
\put (-40,130){\tiny$\mu =1$}
\put (60,130){\tiny$\mu = 0$}

\put (30,40){\linethickness{1pt}\line(0,1){170}}

\end{picture}
\begin{picture}(150,200)(-180,0)
\hspace{-2cm}
\put (-40,180){\tiny$\mu =1$}
\put (60,80){\tiny$\mu =1$}
\put (60,180){\tiny$\mu = 0$}
\put (-40,80){\tiny$\mu = 0$}

\put (30,40){\linethickness{1pt}\line(0,1){170}}
\put (-70,125){\linethickness{1pt}\line(1,0){200}}

\end{picture}

\item[6.]
Also some cusps like the one in the picture below can be admitted, provided that assumption (H.4) is satisfied.
For example, suppose (part of) the interface is that in the picture below and the vertex is the point $(0,0)$
and suppose $\mu \not=0$.
If $\mu_+$ satisfies (H.4) then we are in the assumptions and the theorems of Section \ref{secHarnack} hold. \\
For instance, suppose $\lambda \equiv 1$ and consider $\mu \equiv -1$ on the left of the curve and $\mu \equiv 1$ on the other side of the curve
which is the union of the graphs of $f(x) = x^n$ and $g(x) = - x^n$
for $x \in [0, L]$, $L > 0$, and $n \in \N$, $n \geqslant 1$. We have that
$$
\mu_+ \big( {B_{2\rho} (0,0)} \big) \leqslant \q \, \mu_+ \big( {B_{\rho} ((0,0))} \big)
$$
for some $\q$ depending on $n$. \\
While if, for instance, we consider $f(x) = e^{-1/x}$ and $g(x) = - e^{-1/x}$ the above inequality does not hold any more.
\ \\ 
\ \\
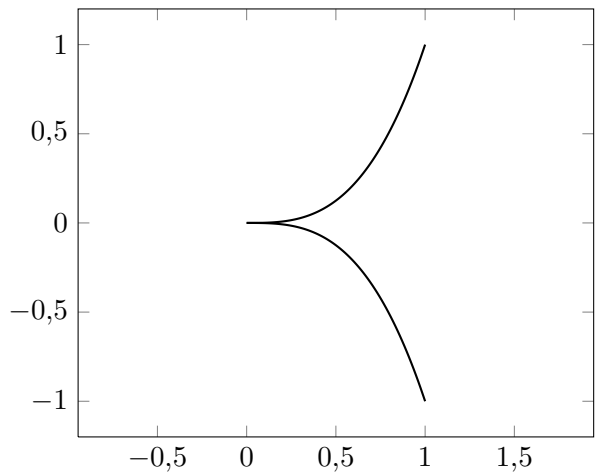
\begin{tikzpicture}
\hspace{3cm}
\begin{axis} [axis equal]
\addplot
[domain=0:1,variable=\t,
samples=40,smooth,thick,black]
({t},{t^3});
\addplot
[domain=0:1,variable=\t,
samples=40,smooth,thick,black]
({t},{-t*t*t});
\end{axis}
\end{tikzpicture}
\ \\
\ \\
If we consider different $\mu$, i.e. $\mu$ which can degenerate to zero,
the geometry of the interface can change depending also on how the weights $|\mu|$ and $\lambda$ degenerate near the interface.

\item[7.]
The final example is the following: again for simplicity suppose $|\mu| \equiv \lambda \equiv 1$ in $\R^2$ and suppose $\mu \equiv 1$ in the region
above the graphic of $f$, which we will call $\Omega_+$, and $\mu \equiv -1$ in the region below the graphic of $f$, which we will call $\Omega_-$, where
$$
f(y) = y \cos \frac{1}{y} \qquad (f(0) = 0) \, .
$$
In spite of the fact that the length of the graphic inside the ball $B := B_1(0,0)$ is infinite, the measure (the $2$-dimensional Lebesgue measure $\mathcal{L}^2$)
of the $\varepsilon$-neighbourhood of $I$ is of order $\varepsilon$ and then going to zero when $\varepsilon \to 0^+$.
Moreover, due to the simmetry of the graphic of $f$ we have that
$$
\mu_+ \big( B_{2\rho}(0,0) \big) = \frac{1}{2} \, \mathcal{L}^2 \big( B_{2\rho}(0,0) \big) \leqslant
	\frac{1}{2} \,  c \, \mathcal{L}^2 \big( B_{\rho}(0,0) \big) = \frac{1}{2} \,  c \, \mu_+ \big( B_{\rho}(0,0) \big)
$$
where $c$ denotes the doubling constant for $\mathcal{L}^2$. Therefore also in this case assumptions (H.4) and (H.5) are satisfied
and even if $I$ is not rectifiable can be an admissible interface.
\ \\ 
\ \\
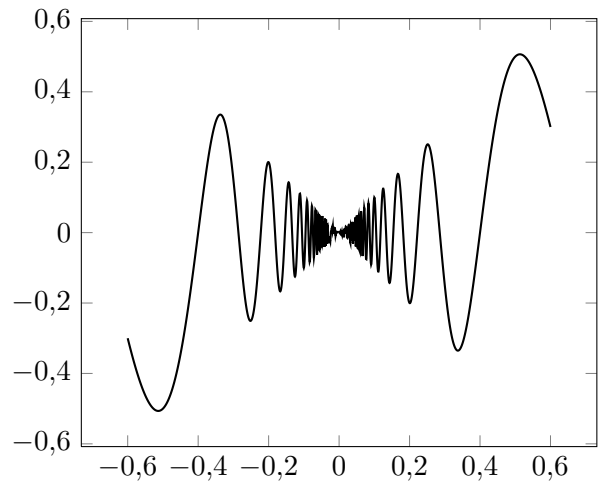
\begin{tikzpicture}
\hspace{3cm}
\begin{axis} [axis equal]
\addplot
[domain=-0.001:0.6,variable=\t,
samples=400,smooth,thick,black]
({t},{t*cos(deg(pi/t))});
\addplot
[domain=-0.6:-0.001,variable=\t,
samples=400,smooth,thick,black]
({t},{t*cos(deg(pi/t))});
\end{axis}
\end{tikzpicture}
\ \\ 
\ \\

\end{itemize}


\begin{thebibliography}{10}

\bibitem{alk-liske}
{\sc Y.~A. Alkhutov - V.~Liskevich}, {\em H\"older continuity of solutions to
  parabolic equations uniformly degenerating on a part of the domain}, Adv.
  Differential Equations, 17 (2012), pp.~747--766.

\bibitem{baogris}
{\sc M.~S. Baouendi - P.~Grisvard}, {\em Sur une \'equation d'evolution
  changeant de type}, J. Funct. Anal., 2 (1968), pp.~352--367.

\bibitem{beals}
{\sc R.~Beals}, {\em On an equation of mixed type from electron scattering
  theory}, J. Math. Anal. Appl., 58 (1977), pp.~32--45.

\bibitem{chanillo-wheeden}
{\sc S.~Chanillo - R.~L. Wheeden}, {\em {Weighted Poincar\'e and Sobolev
  inequalities and estimates for weighted Peano maximal functions}}, Amer.
  Jour. Math., 107 (5) (1985), pp.~1191--1226.

\bibitem{chia-se1}
{\sc F.~Chiarenza - R.~Serapioni}, {\em Degenerate parabolic equations and
  {Harnack} inequality}, Ann. Mat. Pura Appl., 137 (IV) (1984), pp.~139--162.

\bibitem{chia-se3}
{\sc F.~Chiarenza - R.~Serapioni}, {\em A {Harnack }inequality for degenerate
  parabolic equations}, {Comm. Partial Differential Equations}, 9 (8) (1984),
  pp.~719--749.

\bibitem{chia-se2}
{\sc F.~Chiarenza - R.~Serapioni}, {\em A remark on a {Harnack }inequality for
  degenerate parabolic equations}, Rend. Sem. Mat. Univ. Padova, 73 (1985),
  pp.~179--190.

\bibitem{dib1}
{\sc E.~DiBenedetto}, {\em Harnack estimates in certain function classes}, Atti
  Sem. Mat. Fis. Univ. Modena, 37 (1989), pp.~173--182.

\bibitem{articolo}
{\sc E.~DiBenedetto - U.~Gianazza - V.~Vespri}, {\em Harnack estimates for
  quasi-linear degenerate parabolic differential equations}, Acta Math., 200
  (2008), pp.~181--209.

\bibitem{dib-tru}
{\sc E.~DiBenedetto - N.~S. Trudinger}, {\em Harnack inequalities for
  quasiminima of variational integrals}, Ann. Inst. H. Poincar\'e Anal. Non
  Lin\'eaire, 1 (1984), pp.~295--308.

\bibitem{fa-ke-se}
{\sc E.~B. Fabes - C.~E. Kenig - R.~Serapioni}, {\em The local regularity of
  solutions of degenerate elliptic equations}, Comm. Partial Differential
  Equations, 7 (1) (1982), pp.~77--116.

\bibitem{fernandes}
{\sc J.~C. Fernandes}, {\em {Mean value and {Harnack} inequalities for a
  certain class of degenerate parabolic equation}}, Rev. Mat. Iberoamericana, 7
  (3) (1991), pp.~247--286.

\bibitem{gc-rdf}
{\sc J.~{Garcia Cuerva} - J.~L. {Rubio de Francia}}, {\em Weighted norm
  inequalities and related topics}, North-Holland, Amsterdam, 1985.

\bibitem{gianazza-vespri}
{\sc U.~Gianazza - V.~Vespri}, {\em Parabolic {D}e {G}iorgi classes of order
  {$p$} and the {H}arnack inequality}, Calc. Var. Partial Differential
  Equations, 26 (2006), pp.~379--399.

\bibitem{giaquinta}
{\sc M.~Giaquinta}, {\em Introduction to regularity theory for nonlinear
  elliptic systems}, Lectures in Mathematics ETH Z\"urich, Birkh\"auser Verlag,
  Basel, 1993.

\bibitem{giagiu}
{\sc M.~Giaquinta - E.~Giusti}, {\em Quasiminima}, Ann. Inst. H. Poincar\'e
  Anal. Non Lin\'eaire, 1 (1984), pp.~79--107.

\bibitem{giusti}
{\sc E.~Giusti}, {\em Direct methods in the calculus of variations}, World
  Scientific Publishing Co. Inc., River Edge, NJ, 2003.

\bibitem{gut-wheeden2}
{\sc C.~E. Guti\'errez - R.~L. Wheeden}, {\em Harnack's inequality for
  degenerate parabolic equations}, Comm. Partial Differential Equations, 16
  (483) (1991), pp.~745--770.

\bibitem{gut-wheeden3}
{\sc C.~E. Guti{\'e}rrez - R.~L. Wheeden}, {\em Sobolev interpolation
  inequalities with weights}, Trans. Amer. Math. Soc., 323 (1991),
  pp.~263--281.

\bibitem{ishige}
{\sc K.~Ishige}, {\em On the behavior of the solutions of degenerate parabolic
  equations}, Nagoya Math. J., 155 (1999), pp.~1--26.

\bibitem{mohammed}
{\sc A.~Mohammed}, {\em Harnack's inequality for solutions of some degenerate
  elliptic equations}, Rev. Mat. Iberoamericana, 18 (2002), pp.~325--354.

\bibitem{pag-tal}
{\sc C.~D. Pagani - G.~Talenti}, {\em On a forward-backward parabolic
  equation}, Ann. Mat. Pura Appl., 90 (1971), pp.~1--57.

\bibitem{fabio11}
{\sc F.~Paronetto}, {\em {A Harnack's inequality and H\"older continuity for
  solutions of mixed type evolution equations}}.
\newblock To appear in Atti Accad. Naz. Lincei Rend. Lincei Mat. Appl.

\bibitem{fabio9}
{\sc F.~Paronetto}, {\em Further existence results for evolution equations of
  mixed type and for a generalized {Tricomi} equation}.
\newblock Submitted.

\bibitem{fabio4}
{\sc F.~Paronetto}, {\em Existence results for a class of evolution equations
  of mixed type}, J. Funct. Anal., 212 (2) (2004), pp.~324--356.

\bibitem{fabio3}
{\sc F.~Paronetto}, {\em Homogenization of degenerate elliptic-parabolic
  equations}, Asymptotic Anal., 37 (2004), pp.~21--56.

\bibitem{fabio7}
{\sc F.~Paronetto}, {\em A time regularity result for forward-backward
  parabolic equations}, Boll. Unione Mat. Ital. (9), 4, no. 1 (2011),
  pp.~69--77.

\bibitem{show1}
{\sc R.~E. Showalter}, {\em Degenerate parabolic initial-boundary value
  problems}, J. Diff. Eq., 31 (1979), pp.~296--312.

\bibitem{surnachev}
{\sc M.~Surnachev}, {\em A {H}arnack inequality for weighted degenerate
  parabolic equations}, J. Differential Equations, 248 (2010), pp.~2092--2129.

\bibitem{tru1}
{\sc N.~Trudinger}, {\em On the regularity of generalized solutions of linear,
  non-uniformly elliptic equations}, Arch. Rat. Mech. Anal., 48 (1971),
  pp.~51--62.

\bibitem{tru2}
{\sc N.~Trudinger}, {\em Linear elliptic operators with measurable
  coefficients}, Ann. Scuola Norm. Sup. Pisa Cl. Sci., 87 (1973), pp.~265--308.

\bibitem{wieser}
{\sc W.~Wieser}, {\em Parabolic {$Q$}-minima and minimal solutions to
  variational flow}, Manuscripta Math., 59 (1987), pp.~63--107.

\end{thebibliography}

\end{document}